\documentclass[twoside,11pt]{article}

\usepackage{blindtext}

\usepackage[preprint]{jmlr2e}

\newcommand{\dataset}{{\cal D}}
\newcommand{\fracpartial}[2]{\frac{\partial #1}{\partial  #2}}

\usepackage[linesnumbered,ruled,vlined]{algorithm2e}
\usepackage{amsmath, dsfont}
\usepackage{graphicx}
\usepackage{enumerate}
\usepackage{natbib}
\usepackage{url} % not crucial - just used below for the URL 
\RequirePackage{enumitem}
\usepackage{booktabs}
\usepackage{subcaption}
\usepackage{multirow}
\usepackage{float}
\usepackage{placeins}
\usepackage{xr}
\newcommand{\lsup}{\theta}

\usepackage{mwe}
\newlength{\tempheight}
\newlength{\tempwidth}
\newcommand{\rowname}[1]% #1 = text
{\rotatebox{90}{\makebox[\tempheight][c]{\centering \textbf{#1}}}}
\newcommand{\rownameunrot}[1]% #1 = text
{{\makebox[][c]{\textbf{#1}}}}
\newcommand{\columnname}[1]% #1 = text
{\makebox[\tempwidth][c]{\textbf{#1}}}

\usepackage{color}
\usepackage{verbatim}
\usepackage{enumitem}
\newtheorem{assumption}[theorem]{Assumption}

\usepackage[varg]{txfonts}    % comment out for thesis
\renewcommand{\mathbb}{\varmathbb} %comment out for thesis

\renewcommand{\leq}{\leqslant}

\renewcommand{\geq}{\geqslant}

\usepackage{bm}

\usepackage{xspace}

\usepackage{boxedminipage}

\newcommand{\brac}[1]{[#1 ]}
\newcommand{\Brac}[1]{\left[#1\right]}

\newcommand{\norm}[1]{\left\lVert#1\right\rVert}

\newcommand{\N}{{\mathbb N}}
\newcommand{\R}{\mathbb R}

\newcommand{\Esymb}{\mathbb{E}}
\newcommand{\Psymb}{\mathbb{P}}

\DeclareMathOperator*{\E}{\Esymb}

\DeclareMathOperator*{\ProbOp}{\Psymb}

\newcommand{\Prob}[1]{\ProbOp\Brac{#1}}
\newcommand{\probz}[1]{\mathbb{P}_0\brac{#1}}
\newcommand{\Probz}[1]{\mathbb{P}_0\Brac{#1}}

\newcommand{\Probf}[1]{\mathbb{P}_f\Brac{#1}}

\newcommand{\Ex}[1]{\E\Brac{#1}}

\newcommand{\Exz}[1]{\mathbb{E}_0\Brac{#1}}

\newcommand{\Exf}[1]{\mathbb{E}_f\Brac{#1}}

\newcommand{\e}{\epsilon}

\let\e\epsilon

\newcommand{\eve}{\Tilde{\Omega}_T}

\PassOptionsToPackage{unicode}{hyperref}
\PassOptionsToPackage{naturalnames}{hyperref}

\usepackage{lastpage}
\jmlrheading{23}{2022}{1-\pageref{LastPage}}{1/21; Revised 5/22}{9/22}{21-0000}{Author One and Author Two}

\ShortHeadings{Scalable and adaptive variational Bayes methods for Hawkes processes}{D. Sulem, V. Rivoirard, and J. Rousseau}
\firstpageno{1}

\begin{document}

\title{Scalable and adaptive variational Bayes methods for Hawkes processes}

\author{\name D\'eborah Sulem  \email deborah.sulem@bse.eu \\
       \addr Department of Statistics \\
       University of Oxford \\
       \AND
       \name Vincent Rivoirard \email vincent.rivoirard@ceremade.dauphine.fr \\
       \addr Ceremade, CMRS, UMR 7534 \\
       Universit\'e Paris-Dauphine, PSL University\\
       \AND
              \name  Judith Rousseau \email judith.rousseau@stats.ox.ac.uk \\
       \addr Department of Statistics \\
       University of Oxford
}

\editor{}

\maketitle

\begin{abstract}%   <- trailing '%' for backward compatibility of .sty file
    Hawkes processes are often applied to model dependence and interaction phenomena in multivariate event data sets, such as neuronal spike trains, social interactions, and financial transactions. 
    In the nonparametric setting, learning the  temporal dependence structure of Hawkes processes is generally a computationally expensive task, all the more with Bayesian estimation methods. In particular, for generalised nonlinear Hawkes processes, Monte-Carlo Markov Chain methods applied to  compute the \emph{doubly intractable} posterior distribution are not scalable to high-dimensional processes in practice. Recently, efficient algorithms targeting a mean-field variational approximation of the posterior distribution have been proposed. In this work, we first unify existing variational Bayes approaches under a general nonparametric inference framework, and analyse the asymptotic properties of these methods under easily verifiable conditions on the prior, the variational class, and the nonlinear model.
    Secondly, we propose a novel sparsity-inducing procedure, and derive an adaptive mean-field variational algorithm %, based on thresholding appproach of the interaction ``strengths",
    for the popular sigmoid Hawkes processes. Our algorithm is parallelisable and therefore computationally efficient in high-dimensional setting. Through an extensive set of numerical simulations, we also demonstrate that our procedure is able to adapt to the dimensionality of the parameter of the Hawkes process, and is partially robust to some type of model mis-specification. 
\end{abstract}


\section{Introduction}\label{sec:intro}

Modelling point or event data with temporal dependence often implies to infer a local dependence structure between  events and to estimate interaction parameters. In this context, the multivariate Hawkes process is a widely used temporal point process (TPP) model, for instance, in seismology \citep{ogata1999seismicity}, criminology \citep{mohler11}, finance \citep{bacry2015second}, and social network analysis \citep{lemonnier2014nonparametric}. In particular, the generalised nonlinear multivariate Hawkes model, an extension of the classical self-exciting process \citep{hawkes1971point}, is able to account for different \emph{types} of temporal interactions, including \emph{excitation} and \emph{inhibition} effects, often found in event data \citep{hawkes18, bonnet2021maximum}. The \emph{excitation} phenomenon, sometimes named \emph{contagion} or \emph{bursting behaviour}, corresponds to empirical observation that the occurrence of an event, e.g., a post on a social media, increases the probability of observing similar events in the future, e.g., reaction comments. The \emph{inhibition} phenomenon refers to the opposite observation and is prominent in neuronal applications due to biological regulation mechanisms  \citep{bonnet2021maximum}, and in criminology due to the enforcement of policies \citep{olinde20}. 
%The generalised nonlinear Hawkes model can account for these two types of interactions  \citep{bremaud96}, which is of interest for modelling neuronal spike trains data \citep{bonnet2021maximum}, financial applications \citep{zhou2021nonlinear} and internet data \citep{lemonnier2014nonparametric}. 
%These processes have notably been used to model the times of occurrence of earthquakes  \citep{ogata1999seismicity}, terrorist attacks \citep{mohler11}, disease cases \citep{unwin2021using}, market orders \citep{Rambaldi_2016, hawkes18} and neuronal spike trains \citep{LAMBERT20189, bonnet2021maximum}. In the Hawkes model, the probability of an event depends on the past event occurrences through a functional parameter called interaction functions or triggering kernels.
Moreover, the Hawkes model has become popular for the interpretability of its parameter, in particular the \emph{connectivity} or \emph{dependence} graph parameter, which corresponds to a Granger-\emph{causal} graph for the multivariate point process \citep{eichler2016graphical}.
%This parameter can capture complex temporal dependence on previous events. 
%In the linear Hawkes model, this dependence is called the mutual-excitation phenomenon, while the nonlinear model, the dependence can also correspond to a mutual-inhibition phenomenon \citep{bremaud96}.

More precisely, in event data modelling, a multivariate TPP is often described as a counting process of events (or points), $N = (N_t)_{t \in [0,T]} = (N_t^1, \dots, N_t^K)_{t \in [0,T]}$, where $K \geq 1$ is the number of components (or dimensions) of the process, observed over a period $[0,T]$ of length $T > 0$. Each component of a TPP can represent a specific type of event (e.g., an earthquake), or a particular location where events are recorded (e.g., a country). For each  $k=1,\dots, K$ and time $t \in [0,T]$, $N_t^k \in \mathbb{N}$ counts the number of events that have occurred until $t$ at component $k$, therefore, $(N_t^k)_{t \in [0,T]}$ is an integer-valued, non-decreasing, process. %We note that a multivariate TPP is also equivalent to a \emph{marked TPP} where the marks belong to the set $\{1,2,\dots,K\}$ \citep{daley2007introduction}
In particular, multivariate TPP models are of interest for jointly modelling the occurrences of events separated into distinct types, or recorded at multiple places, by specifying a multivariate conditional intensity function (or, more concisely, intensity). The latter,  denoted $(\lambda_t)_t = (\lambda_t^1, \dots, \lambda_t^K)_{t\in \R}$, characterises the probability distribution of events, for each component. It is informally defined as the infinitesimal probability rate of event, conditionally on the history of the process, i.e,
\begin{equation*}
\lambda_t^k dt = \Prob{\text{event at dimension $k$ in } \: [t,t+dt] \Big|\mathcal{G}_{t}}, \quad k=1, \dots, K, \quad t \in [0,T],
\end{equation*}
where $\mathcal{G}_{t} = \sigma(N_s, 0 \leq s<t)$ denotes the history of the process until time $t$. %We note that this intensity function is in general a stochastic process.
In the generalised nonlinear Hawkes model, the intensity is defined as
\begin{equation}\label{def:NLintensity}
    \lambda_t^k = \phi_k \left( \nu_k + \sum_{l=1}^{K} \int_{-\infty}^{t^-} h_{lk}(t-s)dN_s^l \right), \quad k=1,\dots K,
\end{equation}
where for each $k$, $\phi_k: \mathbb{R}\to \mathbb{R}^+$ is a \textit{link} or \emph{activation} function, $\nu_k > 0$  is a \textit{background} or \textit{spontaneous} rate of events, and for each $l=1, \dots, K$, $h_{lk}: \R^+ \to \R$ is an \emph{interaction function} or \emph{triggering kernel}, modelling the influence of $N^l$ onto $N^k$. We note that in this model, the parameter $\nu = (\nu_k)_k$ %is called the \textit{background} or \textit{spontaneous} rate of events. It
characterises the external influence of the environment on the process, here,  assumed constant over time, while the functions $h = (h_{lk})_{l,k=1, \dots, K}$ %are called the \textit{interaction functions} or \emph{triggering kernels} and
parametrise the \emph{causal} influence of past events, that depends on each ordered pair of dimensions. In particular, % the function $h_{lk}$ parametrises the influence of the component $N^l$ onto the component $N^k$, and implies
for any $(l,k)$, there exists a \emph{Granger-causal} relationship from $N^l$ to $N^k$, or in other words, $N^k$ is \emph{locally-dependent} on $N^l$, if and only if $h_{lk} \neq 0$ \citep{eichler2016graphical}. Moreover, defining for each $(l,k)$, $\delta_{lk} := \mathds{1}_{h_{lk}\neq 0}$, the parameter $\delta := (\delta_{lk})_{l,k} \in \{0,1\}^{K \times K}$ defines a Granger-causal graph, called the \emph{connectivity} graph.

Finally, the link functions $\phi = (\phi_k)_k$'s are in general nonlinear and monotone non-decreasing, so that a value $h_{lk}(x) > 0$ can be interpreted as an excitation effect, and $h_{lk}(x) < 0$ as an inhibition effect, for some $x \in \R^+$. Link functions are an essential part of the model chosen by the practitioner, and frequently set as ReLU functions $\phi_k(x) = \max(x,0) = (x)_+$ \citep{Hansen:Reynaud:Rivoirard, chen17b, costa18, lu2018high, bonnet2021maximum, deutsch2022}, sigmoid-type functions, e.g., $\phi_k(x) = \lsup_k (1 + e^{x})^{-1}$ with a scale parameter $ \lsup_k > 0$ \citep{zhou2021efficient, zhou2021nonlinear, malemshinitski2021nonlinear}, softplus functions $\phi_k(x) = \log (1+e^x)$ \citep{mei2017neural}, or clipped exponential functions, i.e.,  $\phi_k(x) = \min(e^x, \Lambda_k)$ with a clip parameter $\Lambda_k > 0$ \citep{gerhard17,  carstensen2010multivariate}. When all the interaction functions are non-negative and $\phi_k(x) = x$ for every $k$, the intensity \eqref{def:NLintensity} corresponds to the linear Hawkes model. Defining the \emph{underlying} or \emph{linear} intensity as
\begin{align}\label{def:lin_intensity}
        \Tilde \lambda_t^k = \nu_k + \sum_{l=1}^{K} \int_{-\infty}^{t^-} h_{lk}(t-s)dN_s^l, \quad k=1,\dots, K,
\end{align}
for any $t \in \R$, the nonlinear intensity \eqref{def:NLintensity} can be re-written as $\lambda_t^k = \phi_k(\Tilde \lambda_t^k)$.

Estimating the parameter of the Hawkes model, denoted $f = (\nu, h)$, and the graph parameter $\delta$, %and also possibly some additional parameter of the activation functions $(\phi_k)_k$, say $\theta = (\theta_k)_k$ (e.g., the scale in the sigmoid model)
can be done via Bayesian nonparametric methods, by leveraging standard prior distributions such as random histograms, B-splines, mixtures of Beta densities \citep{donnet18, sulem2021bayesian}, or Gaussian processes \citep{malemshinitski2021nonlinear}, which enjoy asymptotic guarantees under mild conditions on the model. However, Monte-Carlo Markov Chain (MCMC) methods to compute the posterior distribution are too computationally expensive in practice, even in linear Hawkes models with a moderately large number of dimensions \citep{donnet18}. In contrast, frequentist methods such as maximum likelihood estimates \citep{bonnet2021maximum} and penalised projection estimators \citep{Hansen:Reynaud:Rivoirard, bacry2020sparse, cai2021latent} are more computationally efficient but do not provide uncertainty quantification on the parameter estimates. Yet, in practice, most methods rely on estimating a parametric exponential form of the interaction functions, i.e., $h_{lk}(x) = \alpha_{lk}e^{-\beta_{lk}x}$ \citep{bonnet2021maximum, wang16, deutsch2022}. % in the ReLU model.

%Moreover, general nonlinear Hawkes models have been theoretically studied in \cite{ sulem2021bayesian} under the Bayesian nonparametric framework, nonetheless, it is not yet understood how these methods can be efficiently implemented. 

%For modelling high-dimensional processes, i.e., $K >> 1$, it is often assumed that the model is sparse \cite{Hansen:Reynaud:Rivoirard, bacry2015second}, which means that many interaction functions are null.  In the nonparametric setting, \cite{Hansen:Reynaud:Rivoirard} and \cite{bacry2020sparse} use penalised projection estimators with group Lasso penalty. Penalised maximum likelihood methods have also been considered in the linear model.
The implementation of Bayesian methods using MCMC algorithms is computational intensive for two reasons: the high dimensionality of the parameter space ($K^2$ functions and $K$ parameters to estimate) and the non linearity induced by the link function. 
Recently, data augmentation strategies have been used to answer the second difficulty, jointly with variational Bayes algorithms in sigmoid Hawkes processes \citep{malemshinitski2021nonlinear, zhou2021jmlr}. These novel methods leverage the conjugacy of an augmented  mean-field  variational posterior distribution with certain families of Gaussian priors.
In particular,  \cite{ zhou2021nonlinear} propose an efficient iterative mean-field variational inference (MF-VI) algorithm in a semi-parametric multivariate model. A similar type of algorithm is introduced by \cite{malemshinitski2021nonlinear}, based on a nonparametric Gaussian process prior construction. Nonetheless, these methods do not consider the high-dimensional nonparametric setting. They do not address either the problem of estimating the connectivity graph $\delta$, which is of interest in many applications and which also allows to reduce the computational complexity. In fact, the connectivity graph also determines the dimensionality and the sparsity of the estimation problem, similarly to the structure parameter in high-dimensional regression \citep{Ray_2021}. Moreover, variational Bayes approaches have not been yet theoretically analysed. 

In this work, we make the following contributions to the variational Bayes estimation of multivariate Hawkes processes.
\begin{itemize}
    \item First, we provide a general nonparametric variational Bayes estimation framework for multivariate Hawkes processes %that encompasses existing approaches, such as \cite{zhou2021jmlr} and \cite{malemshinitski2021nonlinear},
    and analyse the asymptotic properties of variational methods in this context. We notably establish concentration rates for variational posterior distributions, leveraging the general methodology of \cite{Zhang2017ConvergenceRO}, based on verifying a prior mass, a testing, and a variational class condition. Moreover, we apply our general results to variational classes of interest in the Hawkes model, namely mean-field and model-selection variational families.
    \item Secondly, we propose a novel adaptive and sparsity-inducing variational Bayes procedure, based on a estimate of the connectivity graph using thresholding of the $\ell_1$-norms of the interaction functions, and relying on \emph{model selection} variational families \citep{Zhang2017ConvergenceRO, ohn2021adaptive}. For sigmoid Hawkes processes, we additionally leverage a mean-field  approximation to derive an efficient adaptive variational inference algorithm. In addition to being theoretically valid in the asymptotic regime, we show that this approach performs very well in practice. % that is able to estimate the dimensionality of  the interaction functions.
    \item In addition to the previous theoretical guarantees and proposed methodology, we empirically demonstrate the effectiveness of our algorithm in an extensive set of simulations. We notably show that, in low-dimensional settings, our adaptive variational algorithm is more computationally efficient than MCMC methods, while enjoying comparable estimation performance. Moreover, our approach is scalable to high-dimensional and sparse processes, and provides good estimates.  In particular, %of the connectivity graph parameter,  Finally
    our algorithm is able to uncover the causality structure of the true generating process given by the graph parameter,
    even in some type of model mis-specification. 
\end{itemize}
%We also leverage the theory of \cite{dennis21}, in particular to analyse sparse variational approximations of Gaussian process distributions. 
% In comparison to the equivalent conditions for the exact posterior \cite{ghosal2000, ghosal:vdv:07}, the neighborhood coverage in the prior mass condition is expressed in Renyi $\rho$-divergence with $\rho > 1$, which is stronger than the standard Kullback-Leibler condition. 
%, which notably obtain slightly weaker conditions in the context of sparse Gaussian process regression.
%A related type of works focus on the \emph{tempered} posterior distributions and derive similar set of conditions \cite{alquier2019concentration, pati2017statistical, zhang2006}.

%Additionally, we note that in the context of sigmoid Hawkes models with link function $\phi_k(x) = \lsup_k (1 + e^{-x})^{-1}$ with $ \lsup_k > 0$, $k \in [K]$, existing algorithms  also  aim at estimating the scale parameter $\theta = (\theta_k)_k$ \citep{apostolopoulou2019mutually, zhou2021efficient, malemshinitski2021nonlinear, zhou2021nonlinear}. However, the latter estimation problem has not been thoroughly analysed yet, neither in the Bayesian nor the frequentist frameworks. 
%Therefore,  % of \cite{zhou2021nonlinear} and \cite{malemshinitski2021nonlinear}.
%we also extend the posterior concentration results of \cite{sulem2021bayesian} to the latter model with unknown scale parameter $\theta$, and validate the use of Bayesian methods in this setup. \textcolor{magenta}{Ca y est ?}

\paragraph{Outline } In the remaining part of this section, we introduce some useful notation. Then, in Section \ref{sec:setup}, we describe our general model and inference setup, and present our novel adaptive and sparsity-inducing variational algorithm in Section \ref{sec:DA_sigmoid}. Moreover, Section \ref{sec:main_results} contains our general results, and their applications to prior and variational families of interest in the Hawkes model. Finally, we report in Section \ref{sec:numerical} the results of an in-depth simulation study. Besides, the proofs of our main results are reported in Appendix \ref{sec:proofs}.

\textbf{Notations. } For a function $h$, we denote $\norm{h}_1 = \int_{\mathbb R} |h(x)|dx$ the $L_1$-norm, $\norm{h}_2 = \sqrt{\int_{\mathbb R}h^2(x)dx}$ the $L_2$-norm, $\norm{h}_\infty = \sup \limits_{x\in \mathbb R} |h(x)|$ the supremum norm, and $h^+ = max(h,0), \: h^- = max(-h,0)$ its positive and negative parts. For a $K \times K$ matrix $M$, we denote $r(M)$ its spectral radius, $\norm{M}$ its spectral norm, and $tr(M)$ its trace. %and $\norm{A}_\infty = \max_i \sum_{j=1}^K |A_{ij}|$  its $L_\infty$-norm.
For a vector $u \in \R^K, \norm{u}_1 = \sum_{k=1}^K |u_k|$. The notation $k \in [K]$ is used for $k \in \{ 1, \ldots, K\}$. For a set $B$ and $k \in [K]$, we denote $N^k(B)$ the number of events of $N^k$ in $B$ and $N^k|_B$ the point process measure restricted to the set $B$. For random processes, the notation $ \overset{\mathcal{L}}{=}$ corresponds to equality in distribution. 
We also denote $\mathcal{N}(u, \mathcal{H}_0,d)$ the covering number of a set $\mathcal{H}_0$ by balls of radius $u$ w.r.t. a metric $d$. For any $k \in [K]$, let $\mu_k^0 = \mathbb{E}_0[\lambda_t^k(f_0)]$ be the mean of $\lambda_t^k(f_0)$ under the stationary distribution $\mathbb{P}_0$. For a set $\Omega$, its complement is denoted $\Omega^c$. We also use the notations $u_T \lesssim v_T$ if $|u_T/v_T|$ is bounded when $T \to \infty$, $u_T \gtrsim v_T$ if $|v_T/u_T|$ is bounded and $u_T \asymp v_T$ if $|u_T/v_T|$ and $|v_T/u_T|$ are bounded. We recall that a function $ \phi$ is $L$-Lipschitz, if for any $(x,x')\in\R^2$,
   $|\phi(x) - \phi(x')| \leq L |x - x'|$. We denote $\mathds{1}_n$ and  $\mathbf{0}_n$ the all-ones and all-zeros vectors of size $n$.  Finally, we denote $\mathcal{H}(\beta, L_0)$ the H\"{o}lder class of $\beta$-smooth functions with radius $L_0$.

\section{Bayesian nonparametric inference of multivariate Hawkes processes}\label{sec:setup}

%\subsection{Multivariate Hawkes processes} \label{sec:hawkes}
%\paragraph{Nonlinear Hawkes processes}

%An alternative definition of Hawkes processes can be formulated via a system of stochastic equations driven by a marked Poisson point process (see for instance \cite{bremaud96}).

%\subsection{Bayesian and variational Bayesian  inference framework} %\label{sec:bayes}

\subsection{The Hawkes model and Bayesian framework}\label{sec:bayes}
Formally a $K$-dimensional temporal point process $N = (N_t)_{t \in \R} = (N_t^1, \dots, N_t^K)_{t \in \R}$, defined as a process on the real line $\R$ and on a probability space  $(\mathcal{X}, \mathcal{G}, \mathbb{P})$, is a Hawkes process if it satisfies the following properties.
 \begin{enumerate}[label={\roman*)}]
    
    \item Almost surely, $\forall k,l \in [K]$, $(N^k_t)_t$ and $(N^l_t)_t$ never jump simultaneously.
    
    \item For all $k \in [K]$, the $\mathcal{G}_t$-predictable conditional intensity function of $N^k$ at $t \in \R$ is given by \eqref{def:NLintensity}, where $\mathcal{G}_t = \sigma(N_s, s < t) \subset \mathcal{G}$.

\end{enumerate}
From now on, we assume that $N$ is a stationary, finite-memory, $K$-dimensional Hawkes process $N$ with parameter $f_0 = (\nu_0, h_0)$, link functions $(\phi_k)_k$, and  memory parameter $A>0$, defined as $A = \sup \{x \in \R^+; \max_{l,k}  |h_{lk}^0(x)| > 0 \}$. We note that $A$ characterises the temporal length of interaction of the point process and that this inference setting is commonly used in previous work on Hawkes processes \citep{Hansen:Reynaud:Rivoirard, donnet18, sulem2021bayesian,cai2021latent}. We assume that $f_0$ is the unknown parameter, and that $(\phi_k)_k$ and $A$ are known to the statistician.
%and $\theta_0$
%such that $((\phi_k)_k, h_0)$ verifies condition \textbf{(C1)}, i.e., $\norm{S_0^+} < 1$ with $S_0^+ = ( L \norm{h_{lk}^{0+}}_1)_{l,k} $. %and $(\psi, f_0)$ verify Assumption \ref{ass-psi}. We denote $\mathbb{P}_0(.|\mathcal{G}_0)$ the conditional distribution of $N$ with initial condition $\mathcal{G}_0$.

Similarly to \cite{donnet18}, we consider that our data is an observation of $N$ over a time window $[-A,T]$, with $T > 0$, but our inference procedure is based on  the log-likelihood function corresponding to the observation of $N$ over $[0,T]$. For a parameter $f = (\nu, h)$, this log-likelihood  is given by
\begin{equation}\label{eq:loglik}
        L_T(f)  := \sum_{k=1}^K L_T^k(f), \quad L_T^k(f) = \left[\int_0^T \log (\lambda_t^k(f)) dN_t^k - \int_0^T \lambda_t^k(f) dt\right].
    \end{equation}
We denote by $\mathbb{P}_0(.|\mathcal{G}_0)$ the true conditional distribution of $N$, given the initial condition $\mathcal{G}_0$, and by $\mathbb{P}_f(.|\mathcal{G}_0)$ the distribution defined as
$
    d\mathbb{P}_f(.|\mathcal{G}_0) = e^{L_T(f) - L_T(f_0)}\mathbb{P}_0(.|\mathcal{G}_0).
$
We also denote $\mathbb{E}_0$ and $\mathbb{E}_f$ the expectations associated to $\mathbb{P}_0(.|\mathcal{G}_0)$ and $\mathbb{P}_f(.|\mathcal{G}_0)$. With a slight abuse of notation, we drop the notation $\mathcal{G}_0$ in the subsequent expressions.

We consider a nonparametric setting for estimating the parameter $f$, within a parameter space $\mathcal{F}$. Given a prior distribution $\Pi$ on $\mathcal{F}$, the posterior distribution, for any subset $B \subset \mathcal{F}$, is defined as
\begin{equation}\label{def:pposterior_dist}
    \Pi(B|N) = \frac{\int_{B} \exp(L_T(f)) d\Pi(f)}{\int_{\mathcal{F}} \exp(L_T(f)) d\Pi(f)} =: \frac{N_T(B)}{D_T}, \quad D_T := \int_{\mathcal{F}} \exp(L_T(f)) d\Pi(f).
\end{equation}
This posterior distribution \eqref{def:pposterior_dist} is often said to be \emph{doubly intractable}, because of the integrals in the log-likelihood function \eqref{eq:loglik} and in the denominator $D_T$. Before studying the problem of computing the posterior distribution, we explicit our construction of the prior distribution. % which is important to understand the variational classes of approximation considered in the next section.

%\subsection{Hierarchical prior distribution}

Firstly, our prior distribution $\Pi$ is built so that it puts mass 1 to finite-memory processes, i.e., to parameter $f$ such that the interaction functions $(h_{lk})_{l,k}$ have a bounded support included in $[0,A]$. % with $A > 0$ a known constant. We call $A$ a \emph{memory} parameter, since it characterises the temporal length of interaction of the point process.
Moreover, we use a hierarchical spike-and-slab prior  based on the connectivity graph parameter $\delta$ similar to \cite{donnet18, sulem2021bayesian}. For each $(l,k) \in [K]^2$, we consider the following parametrisation
\begin{align*}
    h_{lk} = \delta_{lk} \bar h_{lk}, \quad \delta_{lk} \in \{0,1\}, \quad \text{ with } \quad \bar h_{lk}=0 \quad \iff \quad \delta_{lk}=0\end{align*}
so that $\delta = (\delta_{lk})_{lk} \in \{0,1\}^{K^2}$ is the connectivity graph associated to $f$. We therefore consider $\delta \sim \Pi_\delta$, where $\Pi_\delta$ is a prior distribution on the space  $\{0,1\}^{K^2}$, and, for each $(l,k)$ such that $\delta_{lk} = 1$, $\bar h_{lk} \sim \tilde \Pi_h$  where $\tilde \Pi_h$ is a prior distribution on functions with support included in $[0,A]$. In this paper we will mostly consider the case where the functions $\bar h_{lk}$, when non null, are developed on a dictionary of functions $(e_j)_{j\geq 1}$, such that $e_j: [0,A] \to \R, \: \forall j$, and
\begin{equation} \label{dictionary}
\bar h_{lk} = \sum_{j=1}^{J_{k}} h_{lk}^j e_j, \quad h_{lk}^j \in \mathbb R, \: \quad \forall j \in [J_{k}], \quad J_{k}\geq 1, \quad (l,k) \in [K]^2.
\end{equation}
Then, choosing a prior distribution $\Pi_J$ on $J = (J_k)_{k \in [K]}$, our hierarchical prior on $f$ finally writes as
\begin{align}\label{eq:prior-distribution}
    d\Pi(f) = d\Pi_\nu(\nu) d\Pi_\delta(\delta) d\Pi_J(J) d \Pi_{h|\delta, J}(h),
\end{align}
where $\Pi_\nu$ is a prior distribution on $\R_+^K$, suitable to the nonlinear model (see \cite{sulem2021bayesian} for some examples), and 
\begin{align*}
     d \Pi_{h|\delta, J}(h) = \prod_{l,k} (1- \delta_{lk})\delta_{(0)}(\bar h) +  \delta_{lk}d\Tilde  \Pi_{h|\delta, J}(\bar h),
\end{align*}
where $\delta_{(0)}$ denotes the Dirac measure at 0 and $\Tilde  \Pi_{h|\delta, J}(\bar h)$ is a prior distribution on non-null functions decomposed over $J_k$ functions from the dictionary. From the previous construction, one can see that the graph parameter  $\delta \in \{0,1\}^{K^2}$ defines the sparsity structure of $h = (h_{lk})_{l,k}$. This parameter plays a crucial role when performing inference on high dimensional Hawkes processes, either in settings when sparsity is a reasonable assumption, or as the only parameter of interest \citep{bacry2020sparse, Chen2017}.

As previously noted, it is generally expensive to compute the posterior distribution \eqref{def:pposterior_dist}, %since the parameter $f$  includes (at most) $K^2$ functions, and because of the nonlinearities $(\phi_k)_k$ in the definition of the intensity function \eqref{def:NLintensity},
which does not have an analytical expressions. However, we note that when the prior on $f$ is a product of probability distributions on the dimension-restricted parameters
$f_k = (\nu_k, (h_{lk})_{l=1, \dots, K}) \in \mathcal{F}_k$, for $k \in [K]$, so that $f = (f_k)_k$, $\mathcal{F} = \mathcal{F}_1 \times \dots \times \mathcal{F}_K$ and $d\Pi(f) = \prod_k d\Pi_k(f_k)$, then, given the expressions of the log-likelihood function \eqref{eq:loglik} and the intensity function \eqref{def:NLintensity}, we have that each term $L_T^k(f)$ in \eqref{eq:loglik} only depends on $f_k$, i.e., $L_T^k(f) = L_T^k(f_k)$. Furthermore, the posterior distribution can be written as 
\begin{align} \label{rem:factorisation}
    d\Pi(f|N) = \prod_k d\Pi_k(f_k|N), \quad d\Pi_k(f_k|N) = \frac{\exp(L_T^k(f_k)) d\Pi_k(f_k)}{\int_{\mathcal{F}_k} \exp(L_T^k(f_k)) d\Pi_k(f_k)}.
\end{align}

In particular, the latter factorisation implies that each factor $\Pi_k(.|N)$ of the posterior distribution  can be  computed in parallel, nonetheless,  given the whole data $N$. Despite this possible parallelisation, implementation of MCMC methods for computing the posterior distribution in the context of multivariate nonlinear Hawkes processes remains very challenging  \citep{donnet18, zhou2021nonlinear, malemshinitski2021nonlinear}. %\textcolor{red}{ see Fig blabla, even in small dimensions.}  
To alleviate this computational bottleneck, we consider in the next section a family of variational algorithms, together with a two-step procedure to handle high-dimensional processes.

% , using the simpler construction 
% $
%     d\Pi_h( h)  \propto \left(\prod_{l,k} d\Tilde \Pi_h(h_{lk}) \right) \mathds{1}_{ \norm{S^+} <1}(h), %\quad \text{ or } \quad d\Pi_h( h | \delta)  \propto d\pi_h^{\otimes |I(\delta)|}(h) \mathds{1}_{ \|S^+\|_1<1}(h),
% $
% or  $d\Pi_h( h)  = \prod_{l,k} d\Tilde \Pi_h(h_{lk})$ in bounded models, 

\subsection{Variational Bayes inference}\label{sec:var_bayes}

To scale up Bayesian nonparametric methods to high-dimensional processes, we consider a variational Bayes approach. The latter consists of approximating the posterior distribution within a variational class of distributions on $\mathcal{F}$, denoted $\mathcal{V}$. Then, the variational Bayes (VB) posterior distribution, denoted $\hat Q$, is defined as the best approximation of the posterior distribution within $\mathcal{V}$, with respect to the Kullback-Leibler divergence, i.e.,
\begin{align}\label{eq:var_posterior}
    \hat Q := \arg \min_{Q \in \mathcal{V}} KL \left(Q|| \Pi(.|N)\right),
\end{align}
where the Kullabck-Leibler divergence between $Q$ and $Q'$ is defined as
\begin{align*}
    KL(Q||Q') := \begin{cases}
    \int \log \Big(\frac{dQ}{dQ'}\Big) dQ, & \text{if } Q \ll Q' \\
    +\infty, & \text{otherwise}
    \end{cases}.
\end{align*}
For a more in-depth introduction to this framework in the context of Hawkes processes, we refer to the  works of \cite{Zhang_2020, zhou2021jmlr, malemshinitski2021nonlinear}.

In the variational Bayes approach, there are many possible families for  $\mathcal{V}$. Interestingly, we note that under a product posterior \eqref{rem:factorisation},  the variational distribution also factorises in $K$ factors, $\hat Q = \prod_k \hat Q_k$ where each factor $\hat Q_k$ approximates $\Pi_k(.|N)$. Therefore, one can choose a variational class $\mathcal{V}'$ of distributions on $\mathcal{F}_1$, and define $\mathcal{V} := \mathcal{V}'^{\otimes K}$. In the case of multivariate Hawkes processes, we combine mean-field variational approaches \citep{zhou2021jmlr, malemshinitski2021nonlinear} with different versions of model selection variational methods \citep{Zhang2017ConvergenceRO, ohn2021adaptive}. Some important notions related to the two latter inference strategies are recalled in Appendix \ref{app:mean-field}. Before presenting our method, we introduce additional concepts and notation.

%With a slight abuse of notation, we use the notation $\mathcal{V}$ also for $\mathcal{V}'$.

%\textcolor{red}{ Onn'a pas besoin de faire l'hypothese a priori car avec KL si $Q$ n'a pas de densite par rapport a une mesure qui domine la posterior KL est infinie . Donc en minimisant on va forcement aller chercher une densite, il me semble. J'ai donc enleve la remarque}
% From now on, we assume that the mean-field variational posterior distribution has a density with respect to a dominating measure $\mu = \prod_d \mu_d$, and with a slight abuse of notation, we denote $\hat Q$ both the distribution and density with respect to $\mu$. Moreover, 
We consider a general model where the log-likelihood function of the nonlinear Hawkes process can be augmented with some latent variable $z \in \mathcal{Z}$, with $\mathcal{Z}$ the latent parameter space. This  approach is notably used by \cite{malemshinitski2021nonlinear, zhou2021nonlinear} in  the sigmoid Hawkes model, for which $\phi_k(x) \propto (1 + e^{x})^{-1}, \forall k \in [K]$. Denoting $L_T^A(f,z)$ the augmented log-likelihood, we define the \emph{augmented} posterior distribution  as
\begin{align*}%\label{def:aug_posterior_dist}
    \Pi_A(B |N) = \frac{\int_{B} \exp(L_T^A(f,z)) d(\Pi(f) \times \mathbb{P}_{A}(z))}{\int_{\mathcal{F} \times \mathcal{Z}} \exp(L_T^A(f,z)) d(\Pi(f) \times \mathbb{P}_{A})(z)}, \quad B \subset \mathcal{F} \times  \mathcal{Z},
\end{align*}
where $\mathbb{P}_{A}$ is a prior density on $\mathcal{Z}$ with respect to a dominating measure $\mu_z$. One can then define an approximating mean-field family of $\Pi_A(. |N)$ as
\begin{align}\label{eq:aug-mean-field}
        \mathcal{V}_{AMF} = \left \{Q: \mathcal{F} \times \mathcal{Z} \to [0,1] ; \: Q(f, z) =  Q_1(f)Q_2(z) \right \},
\end{align}
by only ``breaking'' correlations between parameters and latent variables. The corresponding mean-field variational posterior distribution is then 
\begin{align}\label{eq:var-mf}
     \hat Q_{AMF} = \arg \min_{Q \in \mathcal{V}_{AMF}} KL(Q|| \Pi_A(. |N)).
\end{align}

%Moreover, by construction, our prior distribution on $f$ has the form $d\Pi(f) = d\Pi(\delta)d\Pi_{f|\delta}(\nu_k, \bar h_{lk}, l,k\in [K])$ where $\bar h_{lk} = \sum_{j=1}^{J_{lk}} h_{lk}^j e_j(x)$ whenever $\delta_{lk}=1$. 
Moreover, our hierarchical prior construction \eqref{eq:prior-distribution} implies that a parameter $f$ is indexed by a set of hyperparameters in the form $m=(\delta, J_{lk}; (l,k) \in \mathcal I(\delta))$, where $\mathcal I(\delta) := \{ (l,k) \in [K]^2;\, \delta_{lk}=1\}$ is the set of ``edges", i.e., pair indices corresponding to non-null interaction functions in $f$. Moreover, $J_{lk}$ is the number of functions in the dictionary used to decompose $h_{lk}$. We note that $m$ characterises the dimensionality of the parameter $f$, and we call it a \emph{model}. We can then re-write our parameter space as
\begin{align}\label{eq:space-decomposition}
    \mathcal{F} = \bigcup_{m \in \mathcal{M}} \mathcal{F}_m, \quad \mathcal{F}_m = \left\{f' \in  \mathcal{F}; \: \delta' = \delta, \: J' = J \right\}, \quad m = (\delta, J), \: \delta = (\delta_{lk})_{l,k}, \: J = (J_{lk})_{l,k},
\end{align}
where $\mathcal{M}$ is the set of models
\begin{align*}
    \mathcal{M} = \left \{ m = (\delta, J); \delta \in \{0,1\}^{K\times K}, \: J \in \mathbb N^{K \times K} \right \}.
\end{align*}
From now on, we assume that for each $k$, $J_{lk} = J_k, \: forall l$ and re-define $J = (J_1, \dots, J_K) \in \mathbb N^{K}$.

The decomposition \eqref{eq:space-decomposition} of the parameter space is key to compute a variational distribution that has support on the whole space $\mathcal{F}$, and that in particular, provides a distribution on the space of graph parameter. Next, we can construct an \emph{adaptive} variational posterior distribution by considering an approximating family of variational distributions within each subspace $\mathcal{F}_m$, denoted $\mathcal{V}^m$. 
%We denote $\mathcal{F}_m$ the latter subspace so that, if $\mathcal M$ contains all possible models in $\mathcal{F}$, then we have $\mathcal{F} = \cup_{m \in \mathcal{M}} \mathcal{F}_m$.
We leverage two types of adaptive variational posterior distributions,  $\hat Q_{A1}$ and $\hat Q_{A2}$, considered respectively by \cite{Zhang2017ConvergenceRO} and \cite{ohn2021adaptive}, and defined as 
\begin{align}
    &\hat Q_{A1} := \hat Q_{\hat m}, \quad  \hat m :=  \arg \max_{m \in \mathcal{M}} ELBO(\hat Q^m),\label{eq:ms_var_post}  \\ 
    &\hat Q_{A2} := \sum_{m \in \mathcal{M}} \hat \gamma_{m} \hat Q_m \label{eq:ms_var_post_avg},
\end{align}
 where  $\hat Q_m = \arg \min_{Q \in \mathcal{V}^m} KL(Q|||\Pi(.|N))$ is the variational posterior distribution in model $m$ (defined on $\mathcal{F}_m$), $ELBO(\cdot)$ is the \textit{evidence lower bound (ELBO)} (defined in our context in \eqref{eq:elbo_1} in Appendix \ref{app:comp_elbo}),  and $ \{\hat \gamma_{m}\}_{m \in \mathcal{M}} $ are the model marginal probabilities defined as 
\begin{align*}
    \hat \gamma_{m} = \frac{ \Pi_m(m) \exp \left \{ ELBO(\hat Q_m)  \right \}}{\sum_{m \in \mathcal{M}}\Pi_m(m) \exp \left \{ ELBO(\hat Q_m)  \right \}}, \quad m \in \mathcal{M}.%\label{eq:ms_var_post}
\end{align*}

\begin{remark}
    We note that in practice one might prefer using the adaptive VB posterior \eqref{eq:ms_var_post} rather than \eqref{eq:ms_var_post_avg}, to avoid manipulating a distribution mixture. In our simulations in Section \ref{sec:numerical}, we often find that one or two models only have significant marginal probabilities $\hat \gamma_k^m$, and therefore the two  adaptive variational posteriors \eqref{eq:ms_var_post}  and \eqref{eq:ms_var_post_avg} are often close. 
\end{remark}

To leverage the computational benefits of the augmented mean-field variational class \eqref{eq:aug-mean-field}, we can set the variational family $\mathcal{V}^m$ as 
\begin{align}\label{eq:mean-field-m}
     \mathcal{V}^m_{AMF} = \left \{Q: \mathcal{F}_m \times \mathcal{Z} \to [0,1] ; \: Q(f, z) =  Q_1(f)Q_2(z) \right \}, \quad \forall m \in \mathcal{M}.
\end{align}

Nonetheless, in the case of moderately large to large values of $K$, it is not computationally feasible to explore all possible models in  $\mathcal{M}$, which number is greater than $2^{K^2}$, the cardinality of the graph space $\{0,1\}^{K^2}$. Even with parallel inference on each dimension, the number of models per dimension is greater than $2^K$ and remains too large. Therefore, for this dimensionality regime, we propose an efficient two-step procedure in the next section. This procedure consists first in estimating $\delta$ using a thresholding procedure, then computes the adaptive mean-field variational Bayes posterior in a restricted set of models with $\delta$ fixed at this estimator.

\subsection{Adaptive two-step procedure} \label{two-step}

In this section, we propose an adaptive and sparsity-inducing variational Bayes procedure for estimating the parameter of Hawkes processes with a moderately large or large number of dimensions $K$. 

Firstly, we note that in Section \ref{sec:main_results}, we will provide theoretical guarantees for the above types of variational approaches in nonlinear multivariate Hawkes processes. In particular, we show that under easy to verify assumptions on the prior and on the parameters, the variational posterior concentrates, in $L_1$-norms at some rate $\epsilon_T$, which typically depends on the smoothness of the interaction functions. Moreover, this concentration rate is the same as for the true  posterior distribution. For instance, using \cite{sulem2021bayesian}, for Lipshitz link functions and well-behaved priors, such as hierarchical Gaussian processes, histogram priors, or Bayesian splines, if the interaction functions belong to a H\"older or Sobolev class with smoothness parameter $\beta$, we obtain that $\epsilon_T\asymp T^{-\beta/(2\beta+1)}$, up to $\log T$ terms. 

A consequence of this result is that for each $(l,k) \in [K]^2$, the (variational) posterior distribution of 
$S_{lk} := \|h_{lk}\|_1$ concentrates around the true value $S_{lk}^0 := \|h_{lk}^0\|_1$ at the same rate $\epsilon_T$. Hence, if for all $(l,k)$ such that $\delta_{lk}^0=1$, $S_{lk}^0$ is large compared to $\epsilon_T$, then the following thresholding estimator of $\delta$ is consistent
\begin{align}\label{eq:graph-estimator}
    \hat \delta_{lk} = 1 \quad \Leftrightarrow \quad \hat S_{lk}> \eta_0,
\end{align}
where $\hat S_{lk}$ is the variational posterior mean or median on $S_{lk}$ and $\epsilon_T<<\eta_0< \min_{lk}S_{lk}^0$.

\begin{figure}[hbt!]
\centering
    \includegraphics[width=0.8\textwidth, trim=0.cm 0.cm 0cm  0.cm,clip]{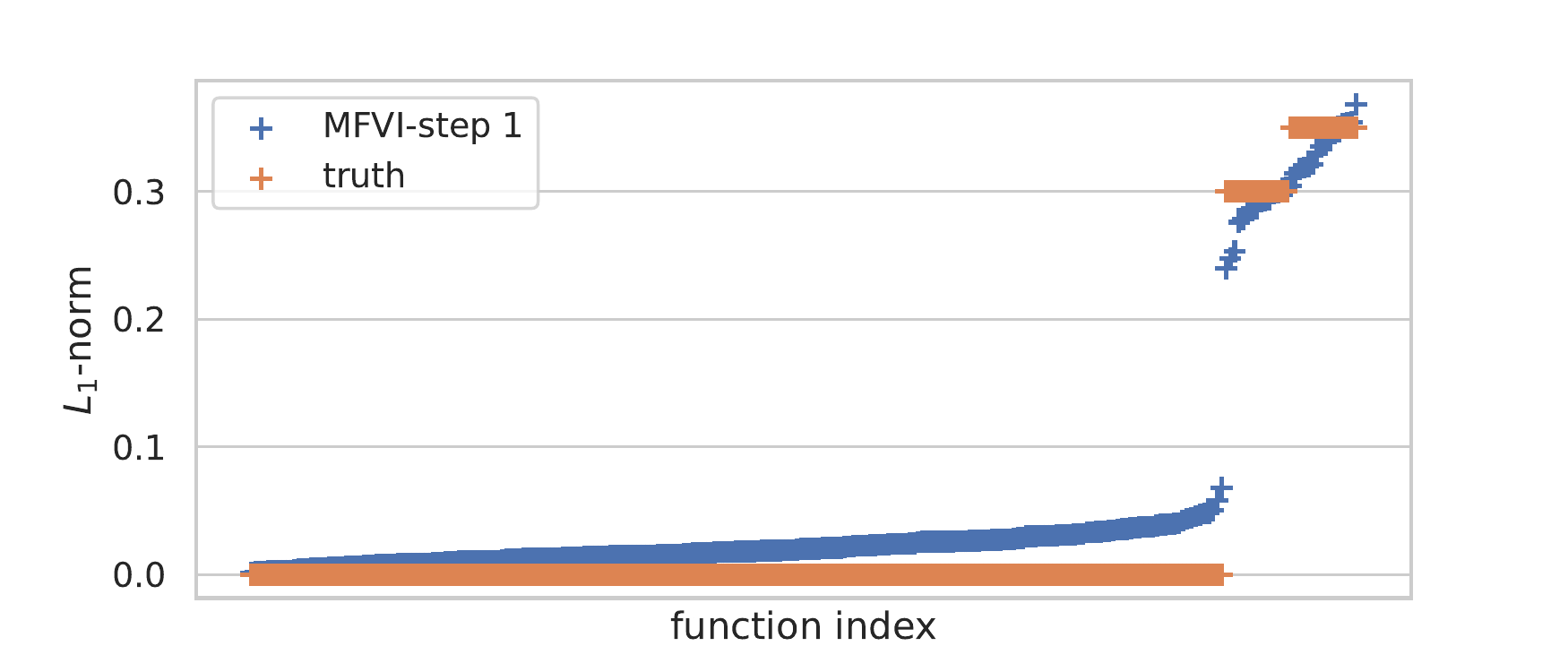}
\caption{Estimated $L_1$-norms $(\hat S_{(i)})_{i \in [K^2]}$ (blue dots), based on the mean-field adaptive variational posterior mean and the set of models $\mathcal{M}_C$ containing models with complete graph $\delta_C = \mathds{1}\mathds{1}^T$, plotted in increasing order. The orange dots correspond to the true values $S_{lk}^0 = \|h_{lk}^0\|_1$. These results correspond to one realisation of the \emph{Excitation} scenario of Simulation 4, for the Hawkes processes with $K=16$ dimensions. }
\label{fig:threshold_example}
\end{figure}

In particular, the above results hold for the adaptive variational Bayes posterior  with the set $\mathcal{M}_C$ of candidate models  with the complete graph $\delta_C$, defined as
\begin{align}\label{eq:set_models_complete}
    \mathcal{M}_C := \{m = (\delta_C = \mathds{1}\mathds{1}^T, J = (J_{k})_{k}); \: J_{k} \geq 1, \: \forall k \in [K] \}.   
\end{align}
In this case, to choose the threshold $\eta_0$ in a data-driven way, we order the estimators $\hat S_{lk}, \: (l,k) \in [K]^2$, say $\hat S_{(1)} \leq  \hat S_{(2)} \leq \cdots \leq \hat S_{(K^2)}$, and set $\eta_0 \in (S_{(i_0)},  S_{(i_0+1)})$ where $i_0$ is the index of the first \textit{significant} gap in $(\hat S_{(i)})_i$, i.e., the first significant values of $ S_{(i+1)}- S_{(i)}$. In Figure \ref{fig:threshold_example}, we plot the estimates $(\tilde S_{(i)})_i$ (blue dots) in one of the simulation settings of Section \ref{sec:simu6}. In this case, the true graph $\delta_0$ is sparse and many $S_{lk}^0$ (orange dots) are equal to 0. From this picture, we can see that by choosing $\eta_0$ anywhere between $0.1$ and $0.2$, we can correctly estimate the true graph  $\delta_0$. More details on these results and their interpretation are provided  in Section \ref{sec:simu6}.

Therefore, once $\hat \delta$ is obtained, we compute an adaptive  variational Bayes posterior, conditional on $\delta=\hat \delta$, by considering the set of models
\begin{align}\label{eq:set_models_restricted}
    \mathcal{M}_E := \{m = (\hat \delta, J = (J_{k})_{k}); \: J_{k} \geq 1, \: \forall k \in [K]\}.   
\end{align}
In summary, our adaptive two-step algorithm writes as:

\begin{enumerate}
    \item \textbf{Complete graph VB}: 
    
    (a) compute the VB posterior associated to the set of models $\mathcal{M}_C$, i.e., to the complete graph $ \delta_C = \mathds{1}\mathds{1}^T$, and compute the posterior mean of $S_{lk}=\|h_{lk}\|_1$,  denoted $\hat S_{lk}$, $\forall (l,k) \in[K]^2$.
    
    (b) order the values $\hat S_{lk}$ in increasing order, say $\hat S_{(1)} \leq  \hat S_{(2)} \leq \cdots \leq \tilde S_{(K^2)}$, and define 
    $\hat \delta_{lk}= 1$ iff $\hat S_{lk} > \eta_0$, where $\eta_0$ is a threshold defined by the first significant value of $\hat S_{(i+1)} -  \hat S_{(i)}$.
    
    \item \textbf{Graph-restricted VB}: compute the VB posterior associated to the set of models $\mathcal{M}_E$, i.e., to models with $\delta = \hat \delta$. 
\end{enumerate}
Theoretical validation of our procedure is provided in Section \ref{sec:conc_rate}. We also note that different variants of our two-step strategy are possible. In particular one can choose a different threshold for each dimension $k \in [K]$, since different convergence rates could be obtained in the different dimensions. Moreover, one can potentially remove the model selection procedure to choose the $J_{k}$'s, $k\in [K]$ in the first step 1(a), and compute a variational posterior in only one model $m \in \mathcal{M}_C$.

In the next section we consider the case of the sigmoid Hawkes processes, for which a data augmentation scheme allows to efficiently compute a mean-field approximation of the posterior distribution within a model $m$. 

\begin{remark}
    In recent work, \cite{bonnet2021maximum} also propose a thresholding approach for estimating the connectivity graph $\delta$ in the context of parametric maximum likelihood estimation. In fact, an alternative strategy to our procedure derived from their work would consist in defining the graph estimator as $\hat \delta_{lk} = 1  \iff \tilde S_{lk}> \varepsilon \sum_{l,k} \tilde S_{lk}$, where $\varepsilon \in (0,1)$ is a pre-defined or data-driven threshold.
\end{remark}

\section{Adaptive variational Bayes algorithms in the sigmoid model}\label{sec:DA_sigmoid}

In this section, we focus on the \emph{sigmoid} Hawkes model, for which the link functions in \eqref{def:NLintensity} are sigmoid-type functions.  We consider the following parametrisation of this model: for each $ k \in [K]$,
\begin{align}\label{eq:sigmoid_tilde}
    \phi_k(x) =  \theta_k\Tilde \sigma(x), \quad \Tilde \sigma(x) = \sigma \left(\alpha(x-\eta)\right), \quad \sigma(x) :=  (1 + e^{-x})^{-1}, \quad \alpha > 0, \: \eta > 0, \: \theta_k > 0.
\end{align} 
Here, we assume that the hyperparameters $\alpha, \eta$ and $\theta = (\theta_k)_k$ are known; however, our methodology can be directly extended to estimate an unknown $\theta$, similarly to \cite{zhou2021jmlr} and \cite{malemshinitski2021nonlinear}. We first note that for $\alpha = 0.1, \eta = 10$ and $\theta_k = 20$, the nonlinearity $\phi_k$ is similar to the ReLU  % $\max(x,0)$
and softplus functions on $[-\infty, 20]$ (see Figure \ref{fig:links} in Section \ref{sec:numerical}). This is helpful to compare the impact of the link functions on the inference in our numerical experiments in Section \ref{sec:numerical}.

For sigmoid-type of link functions,  efficient mean-field variational inference methods based on data augmentation and Gaussian priors have been previously proposed, notably by \cite{malemshinitski2021nonlinear, zhou2021nonlinear, zhou2021jmlr}. We first recall this latent variable augmentation scheme, which allows to obtain a conjugate form for the variational posterior distribution in a fixed model $m = (\delta, J)$ (see Section \ref{sec:bayes}). %, leading  to  conjugate forms of the augmented likelihood function with Gaussian priors and an iterative mean-field variational algorithm.
Then, building on this prior work, we provide two explicit algorithms based on the adaptive and sparsity-inducing methodology presented in Section \ref{two-step}.

%In particular, our approach consists in reducing the dimensionality of the problem by inferring the connectivity graph parameter using the graph selection approach described in Section \ref{conc-sas}.

% This technique has been used to derive Gibbs samplers, a variational EM and efficient mean-field variational inference algorithms, in semi-parametric \cite{zhou2021efficient, zhou2021nonlinear} and nonparametric \cite{malemshinitski2021nonlinear} Bayesian frameworks. Building on this prior work, we first prove new concentration guarantees on the posterior distribution when the scale parameter $\theta = (\theta_k)_k$  is unknown. Then we recall the data augmentation strategy leading to the augmented likelihood function in order to  present our novel adaptive semi-parametric and sparsity-induce approach. In particular, we aim at providing an efficient variational Bayes methodology for selecting the connectivity graph parameter and the dimensionality of the problem.

\subsection{Augmented mean-field variational inference 
in a fixed model}\label{sec:aug-mf-vi}

In our method, we leverage existing latent variable augmentation strategy and Gaussian prior construction, which allows to efficiently compute a mean-field variational posterior distribution on $\mathcal{F}_m \subset \mathcal{F}$, the parameter subspace within a model $m=(\delta, J_{lk}; (l,k) \in \mathcal I(\delta))$. The details of this construction are provided in Appendix \ref{app:data-aug} and we recall that in this context, the set of latent variables $\omega, \bar Z$ correspond respectively to marks at each point of the point process $N$ and to a marked Poisson point process on $[0,T] \times \R^+$.

Then, the augmented mean-field variational family \eqref{eq:aug-mean-field} approximating the augmented posterior distribution corresponds to
$$
\mathcal{V}_{AMF} = \left \{ Q : \mathcal{F} \times \mathcal{O} \times \mathcal{Z} \to \R^+
    ; \: dQ(f, \omega, \bar Z) = dQ_1(f) dQ_2(\omega, \bar Z) \right \},
$$
where $\mathcal{O}$ and  $\mathcal{Z}$ denote the latent variable spaces.  More precisely, in our method, we use the mean-field approach within a fixed model $m$, and therefore define the \emph{model-restricted mean-field variational class} as
$$
\mathcal{V}_{AMF}^m = \left \{  Q: \mathcal{F}_m \times \mathcal{O} \times \mathcal{Z} 
    ; \: dQ(f, \omega, \bar Z) = dQ_1(f) dQ_2(\omega, \bar Z) \right \},
$$
leading to the model-restricted variational posterior
    $\hat Q_{AMF}^m(f, \omega , \bar Z) = \hat Q_1^m(f) \hat Q_2^m(\omega, \bar Z)$.
%\begin{align}\label{eq:var_factor_1} 
 %   \hat Q_{1}(f) &\propto \exp \left \{ \mathbb{E}_{\hat Q_{2}} [\log p(f, \omega, \bar Z, N ) ] \right \}, 
 %   \quad \hat Q_{2}(\omega, \bar Z) \propto \exp \left \{ \mathbb{E}_{\hat Q_{1}} [\log p(f, \omega, \bar Z, N ) ] \right \}.
%\end{align}

Then, we introduce a family of Gaussian prior distributions $\Pi_{h|\delta, J}(h)$ on $\mathcal{F}_m$  such that the factors of $\hat Q_{AMF}^m$, $\hat Q_1^m$ and $\hat Q_2^m$, are conjugate. This conjugacy leads to an iterative variational inference algorithms with closed-forms updates, using \eqref{eq:var_factor_1}. %We recall from Section \ref{sec:bayes} that our hierarchical prior construction on $h$ is based on the graph parameter $\delta = (\delta_{lk})_{l,k} \in \{0,1\}^{K \times K}$ and the representation \eqref{dictionary} on a dictionary with $J_{lk}$ functions, for all non null interaction functions $h_{lk}$.  We now assume that for each $k$, $J_{lk} = J_k, \: \forall l$ and we fix both $\delta$ and the sizes $J := (J_1, \ldots, J_K)$.  We  call $J$ the dimensionality of $h$, and
Let $|J| = \sum_k J_k$. %. With this notation, we have that $m = (\delta, J)$ corresponds to a model in our model-selection approach. Furthermore, we define
We define %our space of functions
\begin{align*}
    \mathcal{H}_{e}^J = \left \{ h = (h_{lk})_{l,k} \in \mathcal{H}; \: h_{lk}(x) = \sum_{j=1}^{J_k} h^j_{lk} e_j(x), \: x \in [0,A], \: \underline{h}_{lk}^{J_k} = (h_{lk}^1, \dots, h_{lk}^{J_k}) \in \R^{J_k}, \: \forall (l,k) \in [K]^2 \right \}.
\end{align*}

Now, for each $(l,k)$, if $\delta_{lk} = 1$, we consider a normal prior distribution %$\Tilde \Pi_{h|\delta} (h)$
on $\underline{h}_{lk}^{J_k}$, with mean vector $\mu_{J_k} \in \R^{J_k}$ and covariance matrix $\Sigma_{J_k} \in \R^{J_k \times J_k}$, i.e.,
$\underline{h}_{lk}^{J_k} \sim  \mathcal{N}(\mu_{J_k}, \Sigma_{J_k})$, and if $\delta_{lk} = 0$, we set  $\underline{h}_{lk}^{J_k} = \mathbf{0}_{J_k}$. We then denote $\mu_m = (\mu_k^m)_k$ with $\mu_k^m = (\delta_{lk} \mu_{J_k})_{l} \in \R^{K J_k}$ and $\Sigma_m = Diag((\Sigma_k^m)_k)$ with $\Sigma_k^m = Diag((\delta_{lk} \Sigma_{J_k})_{l}) \in \R^{KJ_k \times KJ_k}$. We also consider a normal prior on the background rates, i.e., $\nu_k \overset{i.i.d}{\sim}  \mathcal{N}( \mu_\nu,  \sigma_\nu^2)$ with hyperparameters $\mu_\nu,  \sigma_\nu > 0$. We finally denote  by $f_m := (f_k^m)_k \in \mathcal{F}_m$ where for each $k$,
$ f_k^m = (\nu_k, \underline{h}_{1k}^{J_k}, \dots,  \underline{h}_{Kk}^{J_k}) \in \R^{KJ_k+1}$, and define  $H(t) = (H^0(t), H^1(t), \dots, H^K(t)) \in \R^{|J| + 1}$, where $H_0(t) = 1$ and for $k \in [K]$, $H^k(t) = (H_j^k(t))_{j=1, \dots, J_k}$ with
\begin{align}\label{eq:notation_Hjk}
    %&H_j(t) := \int_{t-A}^t e_j(t-s)dN^{k_j}_s, \quad k_j \in [K], \quad j = 1, \dots, K J, \\
    &H_j^k(t) := \int_{t-A}^t e_j(t-s)dN^{k}_s, \quad j \in [J_k].
\end{align}

Using similar computations as in \cite{donner2019efficient, zhou2021nonlinear, malemshinitski2021nonlinear}, we can derive analytic forms for $\hat Q_{1}^m$ and $\hat Q_{2}^m$. In particular, we have that $\hat Q_{1}^m(f_m) = \prod_{k} \hat Q_{1}^{m,k}(f_k^m)$, and for each $k$,  $\hat Q_{1}^{m,k}(f_k^m)$ is a normal distribution  with mean vector $ \tilde \mu_k^m \in \R^{KJ_k+1}$ and covariance matrix  $\tilde{\Sigma}_k^m \in \R^{(KJ_k +1) \times (K J_k +1) }$ given by %$Q_k^s(f_k^s) = \mathcal{N}(f_k^s; \tilde \mu_k^s, \tilde \Sigma_k^s)$, and
\begin{align}
            \tilde{\Sigma}_k^m &= \left[ \alpha^2  \sum_{i \in [N_k]} \mathbb{E}_{\hat Q_{2}^{m,k}}[\omega_i^k]  H(T_i^k)  H(T_i^k)^T +  \alpha^2  \int_0^T \int_0^{+\infty} \bar \omega_t^k H(t) H(t)^T  \Lambda^k(t,\bar \omega) d\bar \omega dt + (\Sigma_k^m)^{-1 } \right]^{-1}, \label{eq:update-sigma}\\
            %&=  \left[ \alpha^2  \sum_{i \in [N_k]}\mathbb{E}[\omega_i^k]  H_k(t_i^k)  H_k(t_i^k)^T +  \alpha^2 \lsup \sum_q v_q \Ex{\bar \omega_q^k}  \frac{\exp(-\frac{1}{2}\mathbb{E}[h(p_q,f_k)])}{2\cosh \frac{\tilde{h}(p_q,f_k)}{2}} H(p_q) H(p_q)^T   + \Sigma^{-1 } \right]^{-1}
    \tilde{\mu}_k^m &= \frac{1}{2 } \tilde{\Sigma}_k^m \left[ \alpha  \sum_{i \in [N_k]} (2 \mathbb{E}_{\hat Q_{2s}^k}[\omega_i^k] \alpha \eta  + 1 ) H(T_i^k)  +  \alpha \int_{0}^T \int_0^{+\infty}    \left(2 \bar \omega^k \alpha \eta  - 1 \right) H(t) \Lambda^k(t,\bar \omega) d\bar \omega dt + 2(\Sigma_k^m)^{-1 } \mu_k^m  \right], \label{eq:update-mu}
    %&= \frac{1}{2 } \tilde{\Sigma}_k \left[ \alpha  \sum_{i \in [N_k]} (2 \mathbb{E}[\omega_i^k] \alpha \eta  + 1 ) H_k(t_i^k)^T  +  \alpha \lsup \sum_q v_q    (2 \Ex{\bar \omega^k_q} \alpha \eta  - 1 ) \frac{\exp(-\frac{1}{2}\mathbb{E}[h(p_q,f_k)])}{2\cosh \frac{\tilde{h}(p_q,f_k)}{2}} H(p_q)^T  + 2\Sigma^{-1 } \mu  \right]
\end{align}
where $N_k := N^k[0,T]$ and
\begin{align*}
     &\Lambda^k(t,\bar \omega) := \lsup_k \frac{\exp \left \{-\frac{1}{2}\mathbb{E}_{Q_{1}^{m,k}}[\Tilde{\lambda}^k_t(f_k^s)] \right \}}{2\cosh \frac{ c^k_{t}}{2}}  p_{PG}(\bar \omega;1,  c^k_{t}), \quad c^k_{t} := \sqrt{ \mathbb{E}_{Q_{1}^{m,k}}[\Tilde{\lambda}^k_t(f)^2]}.
\end{align*}
Besides, we also have that $\hat Q_{2}^m(\omega, \bar Z) = \hat Q_{21}^m(\omega) \hat Q_{22}^m (\bar Z)$ with 
$   \hat Q_{21}^m(\omega) = \prod_k \prod_{i \in [N_k]} p_{PG}(\omega_i^k;1, c^k_{T_i^k})$
and $\hat Q_{22}^m = \prod_k \hat Q^{m,k}_{22}$ where for each $k$, $ \hat Q^{m,k}_{22}$ is the probability distribution of a marked Poisson point process on $[0,T] \times \R^+$ with intensity measure $ \Lambda^k(t,\bar \omega) $. %(see Appendices B and D in \cite{donner2019efficient})
 The full derivation of these formulas can be found in Appendix \ref{app:fixed_updates}.
%where $\mu_D = (\mu_k^D)_k$ is the prior mean vector. 

From the previous expression, we can compute $\hat Q_{2}^m$ given an estimate of $\hat Q_{1}^m$, and conversely. Therefore, to compute the model-restricted mean-field variational posterior $\hat Q^m$, we use an iterative algorithm that updates each factor alternatively, a procedure summarised in Algorithm \ref{alg:cavi}. We note that the updates of the mean vectors and covariance matrices require to compute an integral, which we perform using the  Gaussian quadrature method \citep{golub1969calculation}, where the number of points, denoted $n_{GQ}$, is a hyperparameter of our method. We finally recall that in this algorithm, each variational factor $\hat Q_k^m$ can be computed independently %, and therefore, these computations can be parallelised to further accelerate posterior inference. Besides, each factor 
and only depends on a subset of the parameter $f_k$, and hence, of the sub-model, $m_k := (\delta_k, J_k)$.

\begin{remark}\label{rem:nb-iterations}
The number of iterations $n_{iter}$ in Algorithm \ref{alg:cavi} is another hyperparameter of our method. In practice, we implement an early-stopping procedure, where we set a maximum number of iterations, such as 100, and stop the algorithm whenever the increase of the ELBO is small, e.g., lower than $10^{-3}$, indicating that the algorithm has converged.
    
\end{remark}

\begin{remark}\label{rem:gibbs}
    Similarly to \cite{zhou2021nonlinear, malemshinitski2021nonlinear}, we can also derive analytic forms of the conditional distributions of the augmented posterior \eqref{eq:augm_posterior}. Therefore, the latter could be computed  via a Gibbs sampler, which is provided in  Algorithm \ref{alg:gibbs} in Appendix \ref{app:gibbs_sampler}. However, in this Gibbs sampler, one needs to sample the latent variables - in particular a $K$-dimensional inhomogeneous Poisson point process. This is therefore computationally much slower than the variational inference counterpart, which only implies to compute \emph{expectation} wrt to the latent variables distribution.   % and is also tested in our numerical experiments in Section \ref{sec:numerical}.
\end{remark}

\begin{small}
\begin{algorithm}
\caption{Mean-field variational inference algorithm in a fixed model}\label{alg:cavi}
%\begin{algorithmic}
\KwIn{ $N = (N^1, \dots, N^K)$, $m = (\delta, J), \: J = (J_1, \dots, J_K)$, $\mu_m = (\mu_k^m)_k, \Sigma_m = (\Sigma_k^m)_k$, $n_{iter}$, $n_{GQ}$.}
\KwOut{$\Tilde \mu_m = (\Tilde \mu_k^m)_k, \Tilde \Sigma_m = (\Tilde \Sigma_k^m)_k$.} % Variational posterior mean and covariance $\tilde{\mu} = (\tilde{\mu}_k)_k, \tilde{\Sigma} = (\tilde{\Sigma}_k)_k$.
Precompute $(H(T_i^k))_{i, k}$. \\
Precompute $(p_q, v_q)_{q \in [n_{GQ}]}$ (points and weights for Gaussian quadrature) and $(H(p_q))_{q \in [n_{GQ}]}$ . \\
\DontPrintSemicolon
\SetKwBlock{DoParallel}{do in parallel for each $k = 1, \dots, K$}{end}
\DoParallel{
 Initialisation: $\tilde{\mu}_k^m \gets \mu_k^m$, $\tilde{\Sigma}_k^m \gets \Sigma_k^m$. \\
\For{$t \gets 1$ to $n_{iter}$}{
        \For{$i \gets 1$ to $N_k$}{
        $\mathbb{E}_{\hat Q_{1}^{m,k}}[\tilde{\lambda}^k_{T_i^k}(f_k^m)^2] = \alpha  \left( H(T_i^k)^T \Tilde \Sigma_k^s H(T_i^k) + (H(T_i^k)^T \Tilde \mu_k^s)^2 - 2 \eta H(T_i^k)^T \Tilde \mu_k^m + \eta^2 \right)$ \\
        $\mathbb{E}_{\hat Q_{2}^{m,k}}[\omega_i^k] =  \tanh\left(\sqrt{\mathbb{E}_{\hat Q_{1}^{m,k}}[\tilde{\lambda}^k_{T_i^k}(f_k^s)^2]}\right) / \left(2 \sqrt{\mathbb{E}_{\hat Q_{1}^{m,k}}[ \tilde{\lambda}^k_{T_i^k}(f_k^m)^2]} \right) $
        }
        \For{$q \gets 1$ to $n_{GQ}$}{
        $\mathbb{E}_{\hat Q_{1}^{m,k}}[\tilde{\lambda}_{p_q}^k ( f_k^m)^2] = \alpha \left( H(p_q)^T \Tilde \Sigma_k^m H(p_q) + (H(p_q)^T \Tilde \mu_k^m)^2 - 2 \eta H(p_q)^T \Tilde \mu_k^m + \eta^2 \right)$ \\
        $\mathbb{E}_{\hat Q_{2}^{m,k}}[\omega_q^k] =  \tanh\left(\sqrt{\mathbb{E}_{\hat Q_{1}^{m,k}}[\tilde{\lambda}^k_{p_q}(f_k^s)^2]}\right) / \left(2 \sqrt{\mathbb{E}_{\hat Q_{1}^{m,k}}[ \tilde{\lambda}^k_{p_q}(f_k^m)^2]} \right) $ \\
        $\mathbb{E}_{\hat Q_{1}^{m,k}}[\tilde{\lambda}_{p_q}^k ( f_k^m)  ] = \alpha \left((\tilde{\mu}_k^m)^T H(p_q) -\eta\right)$ \\
        }
        Compute  $\tilde{\Sigma}_k^m$ and $\tilde{\mu}_k^m$ using \eqref{eq:update-sigma} and \eqref{eq:update-mu}
        %$\tilde{\Sigma}_k^s %&=  \alpha^2  \sum_i \mathbb{E}[\omega_i^k]  H_k(t_i^k)  H_k(t_i^k)^T +  \alpha^2  \int_0^T \bar \omega_t^k H(t) H(t)^T  \Lambda^k(t,\bar \omega) d\bar \omega dt + \Sigma^{-1 } \\
        %     = \left[ \alpha^2  \sum \limits_{i \in [N_k]} \mathbb{E}_{\hat Q_{2s}}[\omega_i^k]  H(T_i^k)  H(T_i^k)^T +  \alpha^2 \lsup_k \sum \limits_{q \in [n_{GQ}]} v_q \mathbb{E}_{\hat Q_{2s}}[\bar \omega_q^k]  \frac{\exp(-\frac{1}{2}\mathbb{E}_{\hat Q_{1s}}[\tilde{\lambda}_{p_q}^k ( f_k^s)])}{2\cosh \frac{1}{2}\mathbb{E}_{\hat Q_{1s}}[\tilde{\lambda}^k_{p_q}(f_k^s)^2]} H(p_q) H(p_q)^T   + \Sigma_s^{-1 } \right]^{-1}.
        % $
        % \State $
        %     \tilde{\mu}_k^s %&= \frac{1}{2 } \tilde{\Sigma}_k \left[ \alpha  \sum_i (2 \mathbb{E}[\omega_i^k] \alpha \eta  + 1 ) H_k(t_i^k)^T  +  \alpha \int_{0}^T \int   (2 \bar \omega^k_t \alpha \eta  - 1 ) H(t)^T \Lambda^k(t,\bar \omega) d\bar \omega dt + 2\Sigma^{-1 } \mu  \right] \\
        %     = \frac{1}{2 } \tilde{\Sigma}_k^s \left[ \alpha  \sum \limits_{i \in [N_k]} (2 \mathbb{E}_{\hat Q_{2s}}[\omega_i^k] \alpha \eta  + 1 ) H(T_i^k)^T  +  \alpha \lsup_k \sum \limits_{q \in [n_{GQ}]} v_q    (2 \mathbb{E}_{\hat Q_{2s}}[\bar \omega^k_q] \alpha \eta  - 1 ) \frac{\exp(-\frac{1}{2}\mathbb{E}[\tilde{\lambda}_{p_q}^k ( f_k^s)])}{2\cosh \frac{1}{2}\mathbb{E}_{\hat Q_{1s}}[\tilde{\lambda}^k_{p_q}(f_k^s)^2]} H(p_q)^T  + 2\Sigma_s]^{-1 } \mu_s  \right].
        % $
}
}
%\end{algorithmic}
\end{algorithm}
\end{small}

\subsection{Adaptive variational algorithms}\label{sec:adapt-mf-vi}

Using Algorithm \ref{alg:cavi} for computing a model-restricted mean-field variational posterior, we now leverage the model-selection and two-step approach from Section \ref{sec:bayes} to design two adaptive variational Bayes algorithms. The first one, denoted \emph{fully-adaptive}, is only based on the model-selection strategy from Section \ref{sec:var_bayes} and is suitable for low-dimensional settings.  The second one, denoted \emph{two-step adaptive}, relies on a partial model-selection strategy and the two-step approach from Section \ref{two-step}, and is more efficient for moderately large to large dimensions of the point process.

\subsubsection{Fully-adaptive variational algorithm}\label{sec:fully-adaptive}

%Using the notations of \cite{ohn2021adaptive}, in our setting, a ``model" is a pair of graph parameter and depth, i.e., $m = = \{\delta, D\}$, using that 
%For simplicity, we assume that the depth is kept the same for all non-null $h_{lk}$.

From now on, we assume that the number of functions $(e_j)$ in the dictionary is bounded by $J_T \in \mathbb{N}$. We then define the set of models
\begin{align}\label{eq:set-models-bounded}
    \mathcal{M}_T = \big\{m = (\delta, J=(J_k)_k); \: \delta \in \{0,1\}^{K \times K}, \: 1 \leq J_k \leq J_T, \: k\in [K] \big\}.
\end{align}
We can easily see that in this case $|\mathcal{M}_T | \sim 2^{K^2} J_T$, and that %the number of models to explore grows exponentially with $K$. Moreover, let $s = (s_k)_{k=1,\dots K}$ with $s_k = (\delta_{\cdot k}, D)$, so that $|s_k| = (D+1) \sum_{l=1}^K \delta_{lk} + 1$ and
for any $m = (\delta, J) \in \mathcal{M}_T$, the number of parameters in $m$ is equal to $\sum_{l,k} \delta_{lk}(J_k+1) + 1$. Therefore, we recall that exploring all models in $ \mathcal{M}_T $ is only computationally  feasible for low-dimensional settings, e.g., $K \leq 3$. %Nonetheless, since the variational posterior can be computed dimension per dimension, the cardinality of the set of models \emph{per dimension} is of order $2^{K} J_T$.
We also recall our notation $m = (m_k)_k$ with $m_k = (\delta_{\cdot k}, J_k), \forall k$.

Let $\Pi_{m}$ be a prior distribution on $\mathcal{M}_T$ of the form 
\begin{align*}
    \Pi_{m}(m) = \prod_k \Pi_{m}(m_k) = \prod_k \Pi_{k,\delta}(\delta_{\cdot k}) \Pi_{k,J}(J_k).
\end{align*}
%$\Pi_{m}(m) = \Pi_\delta(\delta) \Pi_J(J)$. 
For instance, one can choose $\Pi_{k,\delta}$ as a product of Bernoulli distribution with parameter $p \in (0,1)$ and $\Pi_{k,J}$ as the uniform distribution over $[J_T]$. 
Using Algorithm \ref{alg:cavi}, for each $m = (m_k)_k$, we  compute  $\hat Q_k^m $ together with the corresponding $ELBO(\hat Q_k^m )$ for each $k$. % (see Appendix \ref{app:comp_elbo} for  the latter derivation).
We note that the computations for each model can be computed independently, and therefore be parallelised to further accelerate posterior inference.  %we can leverage parallelisation to increase the efficient of our algorithm.

Then, we recall that the model-selection adaptive variational approach consists in either selecting $\hat m$ which maximises the ELBO over $m \in  \mathcal{M}_T$ (see \eqref{eq:ms_var_post}) or in averaging over the different models $m$ (see \eqref{eq:ms_var_post_avg}). In the first case, %the maximizer accross $m_k$ of $ELBO(\hat Q_k^m)$ is then
with $\hat m_k = \arg \max_{m_k} ELBO(\hat Q_k^m)$, the VB posterior is $\hat Q_{MS}=\otimes_{k=1}^{K}\hat Q_k^{\hat m_k} $.
%where for each $k$, $\hat Q_{MV}^k $  maximizes $ELBO(\hat Q_{s_k})$, since  as explained in Section \ref{sec:aug-mf-vi}, we have that  $ELBO(Q_s) = \sum_{k=1}^K ELBO(Q_{s_k}^k)$. 
In the second case, the \emph{model-averaging} adaptive variational posterior is given by
    \begin{align}
      &\hat Q_{AV} = \otimes_{k=1}^K \hat Q_k^{AV},     \quad \hat Q^{AV}_k = \sum_{ m_k } \hat \gamma_k^{m} \hat Q_k^{m},\quad \hat \gamma_k^{m} = \frac{\Tilde \gamma_k^{m}}{\sum_m \Tilde \gamma_k^{m}} \nonumber \\
      &\Tilde \gamma_k^{m} = \Pi_{k,\delta}(\delta_{\cdot k}) \Pi_{k,J}(J_k) \exp \left \{ ELBO(\hat Q_k^{m})  \right \}.  \label{eq:vpost_hat}
    \end{align}
% $
%     \hat Q_{AV} =\otimes_{k=1}^{K}\hat Q_k^{\hat m_k}  \sum_{ s  \in \mathcal{S}_T} \hat \gamma_{s} \hat Q_s,
% $
% where $\hat \gamma_{s} $ is the marginal probability on $s$ defined as
% \begin{align}
% \hat \gamma_{s} =   \frac{\Tilde \gamma_{s}}{\sum_{s \in \mathcal{S}_T}   \Tilde \gamma_{s}},    \quad&\Tilde \gamma_{s} = \Pi_{\delta}(\delta) \Pi_J(J) \exp \left \{ ELBO(\hat Q_s)  \right \}.
% \label{eq:vpost_hat}
% \end{align}
% Similarly if  $\Pi_J(J) = \otimes \Pi_{k,J}(J_k)$ and $\pi_\delta (\delta)  = \otimes \pi_{k,\delta}(\delta_{\cdot k})$, then the averages can be computed along each dimension:

We call this procedure (exploring all models in $\mathcal{M}_T$) the \emph{fully-adaptive mean-field variational inference} algorithm, and  summarise its steps in Algorithm \ref{alg:adapt_cavi}. In the next section, we propose a faster algorithm that avoids the exploration of all models in $\mathcal{M}_T$.

\begin{algorithm}
\caption{Fully-adaptive mean-field variational inference}\label{alg:adapt_cavi}
%\begin{algorithmic}
\KwIn{ $N = (N^1, \dots, N^K)$, $\mathcal{M}_T$,  $\mu = (\mu_m)_{m \in \mathcal{M}_T}, \Sigma = (\Sigma_m)_{m \in \mathcal{M}_T}$, $n_{iter}$, $n_{GQ}$.} % $\{(\mu_D, \Sigma_D)\}_{D=1,\dots, D_M}$,
\KwOut{ $\hat Q_{AV}$ or $\hat Q_{MV}$.}
\DontPrintSemicolon
\SetKwBlock{DoParallel}{do in parallel for each $m = (\delta, D) \in \mathcal{M}_T$}{end}
\DoParallel{
 Compute the variational posterior $\hat Q_m$ using Algorithm \ref{alg:cavi} with $\mu_m, \Sigma_m$, $n_{iter}$ and $n_{GQ}$ as hyperparameters. \\
 Compute $(ELBO(\hat Q_k^m)_k)$ and $(\Tilde{\gamma}_{k}^m)_k$ using \eqref{eq:vpost_hat}.
}
 Compute $\{\hat{\gamma}_{m}\}_{m \in \mathcal{M}_T}$ and $\hat Q_{AV}$ or $\hat Q_{MS}$.
%\end{algorithmic}
\end{algorithm}

\subsubsection{Two-step adaptive mean-field algorithm}\label{sec:two_step_mfvi}

As discussed in the above section, for moderately large values of $K$, the model-averaging or model-selection procedures  in Algorithm \ref{alg:adapt_cavi} become prohibitive. In this case, we instead use the two-step approach introduced in Section \ref{two-step}.

We recall that this strategy corresponds to starting with a maximal graph $\delta_C$, typically the  complete graph $\delta_C = \mathds{1}\mathds{1}^T$, and considering the set of models 
$
    \mathcal{M}_C = \big \{m = (\delta_C, J=(J_k)_k); \: 1\leq J_k \leq J_T, \: k\in [K] \big  \},
$
where here as well we assume that the number of functions in the dictionary is bounded by $J_T$. Then, after computing a graph estimator $\hat \delta$, we consider the second set of models $
    \mathcal{M}_E = \big \{m = (\hat \delta, J=(J_k)_k); \: 1\leq J_k \leq J_T, \: k\in [K] \big  \}
$. We note that  both $ \mathcal{M}_C$ and  $ \mathcal{M}_E$  have cardinality of order $ K J_T$, and the cardinality of models \emph{per dimension} is $J_T$. Therefore, as soon as the computation for each model is fast and $J_T$ is not too large, optimisation procedures over these two sets are feasible, even for large values of $K$.

In the first step of our fast algorithm, we compute the model-selection adaptive VB posterior $\hat Q_{MS}^{C}$ using Algorithm \ref{alg:adapt_cavi}, replacing $\mathcal{M}_T$ by $\mathcal{M}_C$. Then, we use $\hat Q_{MS}^{C}$ to estimate the norms  $(\norm{h_{lk}}_1)_{l,k}$ and the graph parameter, with the thresholding method described in Section \ref{two-step}: \\
% $\{\hat \gamma_{\delta_C, D};  \: 1\leq D \leq D_T \}$ using Algorithm \ref{alg:adapt_cavi}, replacing $\mathcal{M}_T$ by $\mathcal{M}_T^{\delta_C}$. Now, let
% \begin{align*}
%     \hat D_C = \arg \max_{D} \hat \gamma_{\delta_C, D}, \quad J_C = 2^{D_C - 1},
% \end{align*}
% and for any $l,k$, let $\tilde{\mu}_{lk}^{D_C}, \tilde \Sigma_{lk}^{D_C}$ be the mean vector and covariance matrix of the weights associated to $h_{lk}$. 
(a) denoting $\hat J_C = (J_{k,C})_k$ the selected dimensionality in $\hat Q_{MS}^{C}$, we compute our estimates of the norm $\hat S_{lk} = \mathbb{E}_{\hat Q_{MS}^{C}}[\norm{h_{lk}}_1], \: \forall (l,k)$, and define $\hat S= (\hat S_{lk})_{l,k} \in \R_+^{K\times K}$; \\
(b) we order our estimates $\hat S_{(1)} < \dots < \hat S_{(K^2)}$ and choose a threshold $\eta_0$ in the first significant gap between $\hat S_{(i)}$ and $\hat S_{(i+1)}$, $i \in [K^2]$; \\
(c) we compute the graph estimator $ \hat \delta = (\hat \delta_{lk})_{l,k}$ defined for any $k$ and $l$ by
$
   \hat \delta_{lk} = \mathds{1}_{\{\hat S_{lk} > \eta_0\}}.
$

% \textcolor{red}{ To be put in the simu section when we specialize in the histogram prior:
% \begin{align*}
%    \Tilde S_{lk} = \mathbb{E}_{\hat Q_1^{\delta_C}}[\norm{h_{lk}}_1] &=  \sum_{j=1}^{J_{k,C}}  \mathbb{E}_{\hat Q_1^{\delta_C}}[|h_{lk}^j|] \\
%     &= \sum_{j=1}^{J_{k,C}} \sqrt{\frac{2}{\pi} [\Sigma_{lk}^{J_{k,C}}]_{jj}}  \exp \left \{ - \frac{[\tilde{\mu}_{lk}^{J_{k,C}}]_{j}^2}{[\Sigma_{lk}^{D_{k,C}}]_{jj} } \right \} - [\tilde{\mu}_{lk}^{J_{k,C}}]_{j} \left[ 1 - 2 \Phi\left(- \frac{[\tilde{\mu}_{lk}^{J_{k,C}}]_{j} }{\sqrt{[\Sigma_{lk}^{J_{k,C}}]_{jj} }}\right) \right].
% \end{align*}}
% with $\Phi$ is the normal cumulative distribution. We then define $\Tilde S_{\hat D_C} = ([\Tilde S_{\hat D_C}]_{lk})_{lk} \in \R_+^{K\times K}$,
% the associated matrix of $\ell_1$-norms of the mean estimator $\Tilde{\mu}_{lk} = (\Tilde{\mu}_{lk}^0, \dots, \Tilde{\mu}_{lk}^{J-1})$ of the interaction functions $h_{lk}$, i.e., for any $l,k \in [K]$, 
% \begin{align*}
%     [\Tilde S_{\hat D_C}]_{lk} = \frac{A}{J} \sum_{j=0}^{J-1} \tilde{\mu}_{lk}^j.
% \end{align*}

In the second step,  we compute the adaptive model-selection  VB posterior $\hat Q_{MS}$  or model-averaging VB posterior $\hat Q_{AV}$  using Algorithm \ref{alg:adapt_cavi}, replacing $\mathcal{M}_T$ by $\mathcal{M}_E$.

This procedure is summarised in Algorithm \ref{alg:2step_adapt_cavi}. In the next section, we provide theoretical guarantees for general variational Bayes approaches, and apply them to our adaptive and mean-field algorithms.

 %We note that the latter is a heuristic approach to deal with the computational cost of exploring all the possible graph parameter, which is not covered by our theoretical results. Nonetheless, we show in our numerical experiments in Section \ref{sec:numerical} that this is an effective approach which allows to scale up to high-dimensional processes.  

\begin{algorithm}
\caption{Two-step adaptive mean-field variational inference}\label{alg:2step_adapt_cavi}
%\begin{algorithmic}
\KwIn{$N = (N^1, \dots, N^K)$, $\mathcal{M}_T$,  $\mu = (\mu_m)_m, \Sigma = (\Sigma_m)_m$, $n_{iter}$, $n_{GQ}$.}
\KwOut{$\hat Q_{MS}$ or $\hat Q_{AV}$}
Compute $\hat Q_{MS}$ using Algorithm \ref{alg:adapt_cavi} with input set $ \mathcal{M}_C$ and hyperparameters  $\mu = (\mu_m)_m, \Sigma = (\Sigma_m)_m$, $n_{iter}$, $n_{GQ}$.
Compute $\hat \delta$ using the thresholding of the estimate $\Tilde S$.
Compute $\hat Q_{MS}^C$ or $\hat Q_{AV}$ using Algorithm \ref{alg:adapt_cavi} with input set  $\mathcal{M}_E$  and hyperparameters  $\mu = (\mu_m)_m, \Sigma = (\Sigma_m)_m$, $n_{iter}$, $n_{GQ}$.
%\end{algorithmic}
\end{algorithm}

%In practice, we only need to compute the values $\{ELBO(\hat Q^m)\}_{m \in \mathcal{M}_T}$ up to a constant in order to compute the model probabilities $\{\hat \gamma_m\}_{m \in \mathcal{M}_T}$.

\section{Theoretical properties of the variational posteriors}\label{sec:main_results}

This section contains general results on variational Bayes methods for estimating the parameter of Hawkes processes, and theoretical guarantees for our adaptive and mean-field approaches proposed in Section \ref{sec:setup} and Section \ref{sec:DA_sigmoid}. In particular, we derive the concentration rates of  variational Bayes posterior distributions, under general conditions on the model, the prior distribution, and the variational family. Then, we apply our general result to variational methods of practical interest, in particular our model-selection adaptive and mean-field methods. % proposed in Section \ref{sec:setup} and Section \ref{sec:DA_sigmoid}.

We recall that in our problem setting, the link functions $\phi := (\phi_k)_k$ in the nonlinear intensity \eqref{def:NLintensity} are fixed by the statistician and therefore known  \textit{a-priori}. Throughout the section  we assume that these functions %$(\phi_k)_k$
are monotone non-decreasing, $L$-Lipschitz, $L > 0$, and that one of the two following  conditions is satisfied:
\begin{itemize}
     \item [\textbf{(C1)}]  For a parameter $f = (\nu, h)$, the matrix defined by $S^+ = (S^+_{lk})_{l,k} \in \R_+^{K \times K}$ with $S_{lk}^+ = L \norm{h_{lk}^+}_1, \forall l,k$, satisfies $\norm{S^+} < 1$;
       \item [\textbf{(C2)}] For any $k \in [K]$, the link function $\phi_k$ is bounded, i.e., $\exists \Lambda_k > 0, \forall x \in \R$, $0 \leq \phi_k(x) \leq \Lambda_k$. 
\end{itemize}
These conditions are sufficient to prove that the Hawkes process is stationary (see for instance \cite{bremaud96}, \cite{deutsch2022}, or \cite{sulem2021bayesian}).

\subsection{Variational posterior concentration rates }\label{sec:conc_rate}

%the setting where the parameter $\theta$ is  known, so that we only need to estimate $f$. In other words, the activation functions $(\phi_k)_k$ are supposed to be fixed  \textit{a priori}. %With a slight abuse of notation, we willl still denote $f = (\nu, h)$ the parameter of interest.
To establish our general concentration result on the VB posterior distribution, we need to introduce the following assumption, also used to prove the concentration of the posterior distribution \eqref{def:pposterior_dist} in the nonlinear Hawkes model %\eqref{def:pposterior_dist} concentrates in our setup
in \cite{sulem2021bayesian}.
\begin{assumption}\label{ass-psi}
For a parameter $f$, we assume that %the activation functions $\phi_k$'s verifies that
there exists $\varepsilon > 0$ such that for each $k \in [K]$, the link function $\phi_k$ restricted to  $I_k = (\nu_k - \max \limits_{l \in [K]} \norm{h_{lk}^{-}}_\infty - \varepsilon, \nu_k + \max \limits_{l \in [K]} \norm{h_{lk}^{+}}_\infty + \varepsilon)$ is bijective from $I_k$ to $J_k = \phi_k(I_k)$ and its inverse is $ L'$- Lipschitz on $J_k$, with $L' > 0$. We also assume that at least one of the two following conditions is satisfied.
\begin{enumerate}[label={\roman*)}]
	\item[(i)] For any $k \in [K]$, $ \inf \limits_{x \in \R}  \phi_k(x) >0$.
	\item[(ii)]  For any $k \in [K]$, $\phi_k > 0$, and $\sqrt{\phi_k}$  and  $\log \phi_k$ are $L_1$-Lipschitz  with $L_1 > 0$ . %on $\mathbb{R}^-$
\end{enumerate}
\end{assumption}
In \cite{sulem2021bayesian}, Assumption \ref{ass-psi} is used to obtain general posterior concentration rates, and is verified for commonly used link functions (see Example 1 in \cite{sulem2021bayesian}). In particular, it holds for sigmoid-type link functions, such as the ones considered in Section \ref{sec:DA_sigmoid}, when the parameter space is bounded (see below). %We also need this assumption to obtain the concentration of the variational posterior distribution since our proofs leverage the existing theory on posterior concentration rates. 

We  now define our parameter space $\mathcal{F}$ %and the functional space $\mathcal{H}$
as follows
\begin{align*}
&\mathcal{H}' = \left \{ h: [0,A] \to \R; \: \|h \|_\infty < \infty \right \}, \quad \mathcal{H} = \left \{ h =  (h_{lk})_{l,k=1}^K \in \mathcal{H}'^{K^2}; \:  (h, \phi) \text{ satisfy  \textbf{(C1)} or  \textbf{(C2)} }  \right \}, \\
&\mathcal{F} = \left \{  f = (\nu, h) \in (\R_+\backslash\{0\})^K \times \mathcal{H}; \: (f, \phi) \text{ satisfies Assumption ~\ref{ass-psi} } \right \}.
\end{align*}
We also define the $L_1$-distance for any $f,f' \in \mathcal{F}$ as
\begin{align*}
    &\|f-f'\|_1 := \norm{\nu - \nu'}_1 + \norm{h -h'}_1, \quad \norm{h -h'}_1 := \sum_{l,k=1}^K \norm{h_{lk} -h_{lk}'}_1, \quad \norm{\nu - \nu'}_1 := \sum_{k} |\nu_{k} -\nu_{k}'|.
\end{align*}
In particular, for the sigmoid function $\phi_k(x) = \theta_k \sigma( \alpha(x - \eta))$, we can choose
$
    \mathcal{F} = \left \{  f = (\nu, h) \in [0,B]^K \times \mathcal{H}  \right \}
$, with $B > 0$.
Moreover, we introduce %and $L_2$-norms respectively, for $B > 0$,
\begin{align*}
&B_\infty(\epsilon) = \left \{f \in \mathcal{F}; \:  \nu_k^0 \leq \nu_k \leq  \nu_k^0 +  \epsilon, \, h_{lk}^0 \leq h_{lk} \leq  h_{lk}^0 + \epsilon, \:  (l,k) \in [K]^2 \right \}, \quad  \epsilon > 0,
%&B_2(\epsilon_T, B) = \{f \in \mathcal{F}; \: \max_k |\nu_k - \nu_k^0| \leq \epsilon_T,  \: \max_{l,k} \|h_{lk} - h_{lk}^0\|_2 \leq \epsilon_T,  \: \max_{l} \nu_l + \max_k \|h_{kl}\|_\infty < B \},
\end{align*}
a neighbourhood around $f_0$ in  supremum norm, and a sequence $(\kappa_T)_T$ defined as
\begin{align} \label{kappaT}
\kappa_T := 10 (\log T)^r, 
\end{align}
with $r=0$ if $(\phi_k)_k$ satisfies Assumption~\ref{ass-psi} (i), and $r=1$ if $(\phi_k)_k$ satisfies Assumption~\ref{ass-psi} (ii).
%We can now present our variational posterior concentration theorem, which is weaker than the general result of \cite{Zhang2017ConvergenceRO} (Theorem 2.1), but is also based on a milder condition on the variational class.
We can now state our general theorem.

\begin{theorem}\label{thm:cv_rate_vi}
Let $N$ be a Hawkes process with link functions $\phi = (\phi_k)_k$ and parameter $f_0 = (\nu_0, h_0)$ such that $(\phi, f_0)$  satisfy Assumption~\ref{ass-psi} and \textbf{(C1)} or \textbf{(C2)}. Let $\epsilon_T = o(1/\sqrt{\kappa_T})$ be a positive sequence verifying $\log^3 T=O(T \epsilon_T^2)$, $\Pi$ a prior distribution on $\mathcal{F}$ and $\mathcal{V}$ a variational family of distributions on $\mathcal{F}$. We assume that the following conditions are satisfied for $T$ large enough.

\textbf{(A0)} There exists $c_1 > 0$ such that $\Pi(B_\infty(\epsilon_T)) \geq e^{-c_1T\epsilon_T^2}.$

\textbf{(A1)} There exist $\mathcal{H}_T \subset \mathcal{H}$, $\zeta_0 > 0$, and $x_0 > 0$ such that %, with $\Theta_T = \{\theta \in \Theta, \: 0 < \theta_k \leq e^{c_2T\epsilon_T^2}, \forall\, k\}$,
$$\Pi(\mathcal{H}_T^{c})  = o(e^{- (\kappa_T +c_1) T\epsilon_T^2})\quad\mbox{and}\quad\log \mathcal{N}\left(\zeta_0 \epsilon_T, \mathcal{H}_T, ||\cdot||_1\right) \leq x_0 T \epsilon_T^2.$$
%    $$ \frac{\Pi(\mathcal{H}_T^{c})}{\Pi(B(\epsilon_T, B))} = o(e^{- \kappa T\epsilon_T^2)}, $$

\textbf{(A2)} There exists $Q \in \mathcal{V}$ such that $supp(Q) \subset B_\infty(\epsilon_T)$ and $KL(Q||\Pi) = O( \kappa_T T \epsilon_T^2).$
% \begin{align}\label{eq:kl_prior}
%     KL(Q||\Pi) \leq c_v \kappa_T T \epsilon_T^2,
% \end{align}
% with $c_v > 0$.

Then, for any $M_T \to \infty$ and $\hat Q$ defined in \eqref{eq:var_posterior}, we have that
 \begin{align*}
       \Exz{ \hat Q \left( \norm{f - f_0}_1 > M_T \sqrt{\kappa_T} \epsilon_T \right) } \xrightarrow[T\to \infty]{} 0.
    \end{align*}
    % or equivalently,
    % $
    %     \Exz{\hat Q( \norm{f - f_0}_1 > M_T \sqrt{\kappa_T}\epsilon_T )} \xrightarrow[T\to \infty]{} 0.
    % $
\end{theorem}
The proof of Theorem \ref{thm:cv_rate_vi} is reported in Appendix \ref{sec:proof_main_thm} and leverage existing theory on posterior concentration rates. We now make a few remarks related to the previous results.

%\begin{remark}\label{rem:l2_norm}
Firstly, similarly to \cite{donnet18} \cite{sulem2021bayesian}, Theorem \ref{thm:cv_rate_vi}  also holds when the neighborhoods $B_\infty(\epsilon_T)$ around $f_0$ in supremum norm, considered in Assumptions \textbf{(A0)}  and \textbf{(A2)}, are replaced neighborhoods in $L_2$-norm, defined as
\begin{align*}
    B_2(\epsilon_T, B) = \left \{f \in \mathcal{F}; \: \max_k |\nu_k - \nu_k^0| \leq \epsilon_T,  \: \max_{l,k} \|h_{lk} - h_{lk}^0\|_2 \leq \epsilon_T,  \: \max_{l} \nu_l + \max_k \|h_{kl}\|_\infty < B \right \},
\end{align*}
with $B>0$, and when $\kappa_T$ replaced by $\kappa_T' = 10 (\log \log T) (\log T)^r$.
%\end{remark}

%\begin{remark}\label{rem:ass-a2}
Secondly, Theorem \ref{thm:cv_rate_vi} also holds under the more general condition on the variational family: \\
$\quad$ \textbf{(A2')}  The variational family $\mathcal{V}$ verifies $ \min_{Q \in \mathcal{V}} KL(Q||\Pi(.|N)) = O(\kappa_T T \epsilon_T^2).$
% The variational family is such that
% \begin{align*}
%     \min_{Q \in \mathcal{V}} KL(Q||\Pi(.|N)) = o(\kappa_T T \epsilon_T^2).
% \end{align*}
However, in practice, one often verifies \textbf{(A2)} and deduces \textbf{(A2')}  using the following steps from \cite{Zhang2017ConvergenceRO}. For any $Q \in \mathcal{V}$, we have that
\begin{align*}
    KL(Q||\Pi(.|N)) \leq %KL(Q||\Pi) + Q(KL(d\mathbb{P}_{0},d\mathbb{P}_f)) = 
    KL(Q||\Pi) + Q(KL(\mathbb{P}_{T,f_0}, \mathbb{P}_{T,f})), %Q(\Exz{L_T(f_0) - L_T(f)}),
\end{align*}
where we denote $\mathbb{P}_{T,f_0} = e^{L_T(f_0)}$ and $\mathbb{P}_{T,f} = e^{L_T(f)}$.
Using Lemma S6.1 from \cite{sulem2021bayesian}, for any $f \in B_\infty(\epsilon_T)$, % (or $B_2(\epsilon_T,B)$), 
we also have that
\begin{align*}
    \Exz{L_T(f_0) - L_T(f)} \leq \kappa_T T \epsilon_T^2.
\end{align*}
Therefore,  under \textbf{(A2)},  there exists $Q \in \mathcal{V}$ such that %and $c_v > 0$ such that $supp(Q) \subset B_\infty(\epsilon_T)$ and $KL(Q||\Pi) \leq c_v \kappa_T T \epsilon_T^2$, then for this $Q$,
$
    KL(Q||\Pi(.|N)) =  O(\kappa_T T \epsilon_T^2),
$
which implies \textbf{(A2')} . Besides, \textbf{(A2)} (or \textbf{(A2')}), is the only condition on the variational class, and informally states that this family of distributions can approximate the true posterior conveniently. Nonetheless, under  \textbf{(A2)}, we may still have $\min_{Q \in \mathcal{V}} KL(Q||\Pi(.|N)) \xrightarrow[T\to \infty]{} \infty$, as has been observed by \cite{dennis21}.

Finally, Assumptions \textbf{(A0)} and \textbf{(A1)} are similar to the ones of Theorem 3.2 in \cite{sulem2021bayesian}. They are sufficient conditions for proving that the posterior concentration rate is at least as fast as $\sqrt{\kappa_T} \epsilon_T$. %in this setting, and to obtain a concentration inequality on the posterior distribution of the type

\subsection{Applications to variational classes and prior families of interest}

In this section, we apply the previous result to variational inference methods of interest in nonlinear Hawkes models, in particular, the mean-field and model-selection variational families, introduced in Section \ref{sec:var_bayes} and used in  our algorithms. We also verify our general conditions on the prior distribution on two common examples of nonparametric prior families, namely random histograms and Gaussian processes, see for instance in \cite{donnet18, malemshinitski2021nonlinear}. %for which we propose adaptive algorithms in the sigmoid model in Section \ref{sec:DA_sigmoid}. We consider two prior families, namely the Gaussian process prior and the random histogram prior
We then obtain explicit concentration rates for the variational posterior distribution and for  H\"{o}lder classes of functions.

First, we re-write our hierarchical spike-and-slab prior distribution  from Section~\ref{sec:bayes} as 
\begin{align}\label{eq:prior-distribution-rewritten}
    d\Pi(f) = d\Pi_\nu(\nu) d\Pi_\delta(\delta)  d \Pi_{h|\delta}(h), \quad d \Pi_{h|\delta}(h) = \prod_{l,k} d \Tilde \Pi_{h|\delta}(h_{lk})
\end{align}
and recall that from \cite{sulem2021bayesian}, we know that Assumption  \textbf{(A0)} of Theorem \ref{thm:cv_rate_vi} can be replaced by

\medskip
\textbf{(A0')} There exists $c_1 > 0$ such that $\Pi(  B_\infty(\epsilon_T) | \delta= \delta_0)   \geq e^{ -c_1 T\epsilon_T^2/2} $ and $ \Pi_\delta ( \delta = \delta_0) \geq e^{ -c_1 T\epsilon_T^2/2}$.
\medskip

Furthermore, one can choose for instance $\Pi_\delta = \mathcal B(p)^{K^2}$ with $p \in  (0,1)$, implying that the $\delta_{lk} $'s are i.i.d. Bernoulli random variables. Then, for any fixed $p$, one only needs to verify $ \Pi_{h|\delta}( B_\infty(\epsilon_T)|\delta = \delta_0)  \geq e^{ -c_1 T\epsilon_T^2/2} $.

\subsubsection{Mean-field variational family}\label{conc-mf-vi}

Here, we consider the mean-field variational inference method with general latent variable augmentation as described in Section \ref{sec:bayes}.
%can be used in the Hawkes model within a conveniently chosen latent variable augmentation scheme, and using approximating distributions that remove the correlation between the parameter $f$ and the latent variable $z$.
We recall that for some latent variable $z \in \mathcal{Z}$, the mean-field family  $\mathcal{V}_{AMF}$ for approximating the augmented posterior $\Pi_A(. |N)$ is defined as
\begin{align*}
        \mathcal{V}_{AMF} = \left \{Q: \mathcal{F} \times \mathcal{Z} \to [0,1] ; \: Q(f, z) =  Q_1(f)Q_2(z) \right \},
\end{align*}
and the corresponding mean-field variational posterior is $\hat Q_{AMF} = \arg \min_{Q \in \mathcal{V}_{AMF}} KL(Q|| \Pi_A(. |N))$. We also recall our notation $ \mathbb{P}_{A}$, for the prior distribution on the latent variable.  %and two nonparametric  prior distributions of interest. 
% Using an augmentation scheme is particularly useful in the sigmoid Hawkes model (see Section \ref{sec:DA_sigmoid}), %one can define $w := (\omega, \bar N)$, 
% however, more generally, one could define a mean-field family in the Hawkes without the latent variables. We define the \emph{augmented} mean-field variational posterior as
% \begin{align}\label{eq:mean-field-vp}
%     \hat Q_A = \arg \min_{Q \in \mathcal{V}_{MF, A}} KL \left(Q ||  \Pi_A(. |N) \right) =: \hat Q_1 \hat Q_2. 
% \end{align}
% In the mean-field variational scheme, one considers factorising each distribution on $f_k$ into $p\geq 2$ variational factors. We define this variational family in a general Hawkes model, possibly with added latent variables (see for instance the data augmentation scheme in the sigmoid model in Section \ref{sec:DA_sigmoid}).
% \begin{definition}[Mean-field variational class in the Hawkes model]
% In the Hawkes model with parameter $f = (f_k)_k \in \mathcal{F}$ with $f_k = (\theta_k, \nu_k,h_k)$ and latent variables $\omega = (\omega_k)_k \in \mathcal{O}$, a mean-field variational family is defined as
% \begin{align*}
%      \mathcal{V}_{MF} = \left \{ Q: \mathcal{F} \times  \mathcal{O} \to \R ; \: dQ(f) = \prod_{k=1}^K \prod_{j=1}^p dQ_j(x_k^j), \: x_k^j \subset (\theta_k, \nu_k, h_{1k}, \dots, h_{Kk}, \omega_k), \: p\geq 2 \right \}.
% \end{align*}
% \end{definition}
We  note that here the augmented prior distribution is $\Pi \times \mathbb{P}_{A} \in \mathcal{V}_{AMF}$, therefore, assumption \textbf{(A2)} is equivalent to the prior mass condition (see for instance \cite{Zhang2017ConvergenceRO}). Therefore, we only need to verify the assumptions \textbf{(A0')} and \textbf{(A1)}. Besides, these assumptions are the same as in \cite{sulem2021bayesian} and therefore can be applied to any prior family discussed there. In particular, priors on the $h_{lk}$'s based on decompositions on dictionaries like in \eqref{dictionary} have been studied in \cite{arbeletal:13} or \cite{shen:ghosal14} and their results can be applied to prove  assumptions \textbf{(A0')} and \textbf{(A1)}. Below, we apply Theorem \ref{thm:cv_rate_vi} in two examples,  random histogram priors and hierarchical Gaussian process priors. 

%We  now apply Theorem \ref{thm:cv_rate_vi} to $\hat Q_{AMF}$ and two family of prior distributions.

%the notation $\Tilde \Pi_h$ from Section \ref{sec:bayes} for the prior on the slabs of the interaction functions $h_{lk}$.
%we recall that $\Pi(f) = \Pi_h(h) \prod_k \Pi_\nu(\nu_k)$, with $\Pi_h(h)=  \otimes_{l,k} \Tilde \Pi_h(h_{lk}) $ the nonparametric prior on the interaction functions.

\paragraph{Random histogram prior}
 We consider a random histogram prior for $\Pi_{h|\delta}(h)$, using a similar construction as in Section \ref{sec:aug-mf-vi}. %the prior distribution on the "slabs", i.e., the non-null interaction functions, from Section \ref{sec:bayes}.
 This prior family is notably  used in \cite{donnet18, sulem2021bayesian}, and is similar to the basis decomposition prior in \cite{zhou2021efficient, zhou2021nonlinear}. For simplicity, we assume here that $J=J_1=\dots=J_k$ and consider a  regular partition of $(0,A]$ based on
%\subsection{Random histogram basis}
% Prior distribution based on basis decomposition have been widely used in the context of linear and nonlinear Hawkes processes \cite{donnet18, sulem2021bayesian, zhou2021efficient, zhou2021nonlinear}. In this section, we consider the prior construction \eqref{eq:parametrisation2} and a random histogram prior for $\tilde \Pi_{h}$. For simplicity, we consider regular partitions
$(t_j)_{j=0,\dots,J}$ with $t_j = j A/J, \: j=0, \dots, J$, $J \geq 1$, and define piecewise-constant interaction functions as
$$ h_{lk}^{w}(x) =  \sum_{j=1}^{J} w_{lk}^j e_j(x), \quad  e_j(x) =\frac{ J }{A}\mathds{1}_{(t_{j-1}, t_j]}(x), \quad w_{lk}^j \in \R \quad \forall j \in [J],  \forall l,k \in [K].  % \quad J \sim \mathcal P(\lambda), \quad \lambda > 0 .
$$
Note that $\|e_j\|_2=\sqrt{J/A}$ but $\|e_j\|_1 = 1, \: \forall j \in [J]$, therefore, the functions of the dictionary, $(e_j)_j$ are orthonormal in terms of the $L_1$-norm.
In this general construction, we also consider a prior on the number of pieces $J$ with exponential tails, for instance we can choose $J \sim \mathcal P(\lambda)$ with $ \lambda > 0 $, or $J= 2^D$ where $2^D \leq J_D < 2^{D+1}$ and $J_D\sim \mathcal P(\lambda)$. Finally, given $J$, we consider a normal prior distribution on each weight $w_{lk}^j$, i.e.,
$$w_{lk}^j |J \overset{\mathrm{i.i.d.}}{\sim} \mathcal{N}(0_J,K_J), \quad K_J = \sigma_0^2 I_J, \quad \sigma_0 > 0.$$
With this prior construction, assumptions  \textbf{(A0')} and  \textbf{(A1)} are easily checked. For instance, this Gaussian random histogram prior is a particular case of the spline prior family in \cite{sulem2021bayesian}, with a spline basis of order $q=0$. We note that these conditions are also verified easily for other prior distributions on the weights, for instance, the shrinkage prior of  \cite{zhou2021efficient} based on the  Laplace distribution $p_{Lap}(w_{lk}^j ;0,b) = (2b)^{-1} \exp \{ - |w_{lk}^j |/b \}$ with $ b> 0$, and a ``local" spike-and-slab prior inspired by the construction in \cite{donnet18, sulem2021bayesian}:
\begin{align*}
    w_{lk}^j | J \overset{\mathrm{i.i.d.}}{\sim} p \delta_{(0)} + (1 - p) p_{Lap}(. ;0,b), \quad p \in (0,1), \quad b>0,
\end{align*}
where $\delta_{(0)} $ is the Dirac measure at 0.

In the following proposition, we further assume that the true functions in $h_0$ belong to a Holder-smooth class of functions $\mathcal{H}(\beta, L_0)$ with $\beta \in  (0,1)$, so that explicit variational posterior concentration rates $\epsilon_T$ for the mean-field family and the  random histogram prior can be derived. 
% \begin{align*}
%     w_{lk}^j = z_{lk}^j \Tilde w_{lk}^j, \quad z_{lk}^j \in \{-1,0,1\}, \quad  z_{lk}^j  \overset{\mathrm{i.i.d.}}{\sim} \Pi_z, \quad \Tilde w_{lk}^j \overset{\mathrm{i.i.d.}}{\sim} \mathcal{D}(\alpha_{1}, \dots, \alpha_{J}),
% \end{align*}
% with $\Pi_z$ is a prior distribution on $\{-1,0,1\}$, $\alpha_{1}, \dots, \alpha_{J} > 0$ and $\mathcal{D}$ is the Dirichlet distribution. 
% \begin{align*}
%     &w_{lk}^j \sim Lap(0,b). % \quad b>0, \\
%     %&p_{Lap}(w_{lk}^j;0,b) = \frac{1}{2b} \exp \{ - \frac{|w_{lk}^j |}{b}\}.
% \end{align*}
 
%and assuming that the true functions $h_0$ are H\"{o}lder-smooth, we can thus obtain an explicit variational concentration rate. 

% Similarly to \cite{sulem2021bayesian}, we can verify assumptions (A0)-(A2) with these priors and obtain an explicit posterior concentration rate for $\beta$-smooth functions with $\beta \in (0,1]$. We consider the mean-field variational family
% \begin{align*}
%     \mathcal{V}_{MF} = \left \{  Q ; \: dQ(f_k) = dQ_1(\theta_k) dQ_2(\nu_k, h_k), \: h_{k} = (h_{lk})_{l=1}^K \right \},
% \end{align*}
% where $Q_1$ and $Q_2$ are distributions respectively on $(0,+\infty)$ and $\R \cup \mathcal{H}'^{K}$. In particular, $ \mathcal{V}_{MF} $ contains distributions with the same form as the prior. 

\begin{proposition}\label{prop:histo} 
Let $N$ be a Hawkes process with link functions $\phi = (\phi_k)_k$ and parameter $f_0 = (\nu_0, h_0) $ such that $(\phi,f_0)$ verify Assumption \ref{ass-psi}. Assume that for any $l,k \in [K]$, $h_{lk}^0 \in \mathcal{H}(\beta, L_0)$ with $\beta \in  (0,1)$ and $L_0 > 0$. Then, under the above Gaussian random histogram prior, the mean-field variational distribution $\hat Q_1$ defined in \eqref{eq:mean-field-vp} satisfies, for any $M_T \to +\infty$,
 \begin{align*}
        \Exz{\hat Q_1 \left( \norm{f - f_0}_1 > M_T(\log T)^q (T / \log T)^{-\beta / (2\beta + 1)} \right)} \xrightarrow[T\to \infty]{} 0,
\end{align*}
with $q = 0$ if $\phi$ verifies Assumption \ref{ass-psi}(i) and $q = 1/2$ if $\phi$ verifies Assumption \ref{ass-psi}(ii). 
\end{proposition}

The proof of Proposition \ref{prop:histo} is omitted since it is a direct application of Theorem \ref{thm:cv_rate_vi} to mean-field variational families in the context of a latent variable augmentation scheme. We note that the variational concentration rates also match the true posterior concentration rates (see \cite{sulem2021bayesian}).

\paragraph{Gaussian process prior}  We now consider a prior family $\Pi_{h|\delta}$ based on Gaussian processes which is commonly used for nonparametric estimation of Hawkes processes (see for instance \cite{Zhang_2020, zhou2020, malemshinitski2021nonlinear}).
%\subsubsection{Gaussian processes priors}
We define a centered Gaussian process distribution  with %mean function $m_{GP}$ and
covariance function $k_{GP}$ as the prior distribution $\Tilde \Pi_{h|\delta}$ on each $h_{lk}$ such that $\delta_{lk} = 1$, $l,k \in [K]$, i.e., for any $n \geq 1$ and $x_1,\dots, x_n \in [0,A]$, we have
\begin{align*}
    (h_{lk}(x_i))_{i=1,\dots, n} \sim \mathcal{N}\left(0_{n},  (k_{GP}(x_i,x_j))_{i,j=1,\dots, n}\right).
\end{align*}
% \begin{align*}
%     \begin{pmatrix}
% h_{lk}(x)\\
% h_{lk}(x')
% \end{pmatrix} \sim  \mathcal{N} \left(\begin{pmatrix}
% 0\\
% 0
% \end{pmatrix},
% \begin{pmatrix}
% k_{GP}(x,x) & k_{GP}(x,x')\\
% k_{GP}(x,x') & k_{GP}(x',x')
% \end{pmatrix} \right).
% \end{align*}
% \begin{align*}
%     h_{lk} \sim GP(m_{GP}, k_{GP}).
% \end{align*}
We then verify assumptions \textbf{(A0')}  and \textbf{(A1)}  based on the $L_2$-neighborhoods (see comment after Theorem \ref{thm:cv_rate_vi}), i.e.,  we check that there exist $ \mathcal{H}_T \subset \mathcal{H}$ and $c_1, x_0, \zeta_0 > 0$, such that 
\begin{align*}
    &\Pi(\mathcal{H}_T^c) \leq e^{-(\kappa_T + c_1) T\epsilon_T^2}, \quad \log \mathcal{N}(\zeta_0 \epsilon_T, \mathcal{H}_T, \|.\|_1) \leq x_0 T \epsilon_T^2,\quad \Pi(B_2(\epsilon_T,B)) \geq e^{-c_1 T\epsilon_T^2}.
\end{align*}
It is therefore enough to find $ \mathcal{B}_T \subset L_2([0,A])$ such that
\begin{align*}
    &\Tilde \Pi_h(\mathcal{B}_T^c) \leq e^{-(\kappa_T + c_1) T\e_T^2}, \quad \log \mathcal{N}(\zeta_0 \epsilon_T, \mathcal{B}_T, \|.\|_1) \leq \frac{x_0  T \epsilon_T^2 }{K^2}, \quad \Tilde \Pi_h(\norm{h_{lk} - h_{lk}^0}_2 < \epsilon_T) \geq   e^{-c_2 T\e_T^2} / K^2,
\end{align*}
and define $\mathcal{H}_T = \mathcal{B}_T^{\otimes K^2}$, since for all $\zeta>0$, there exists $\zeta_2 >0$   (independent of $T$) such that $\Pi(\mathcal{H}_T^c)  \leq \Tilde \Pi(\mathcal{B}_T^c),$ and 
\begin{align*}
  \log \mathcal{N}(\zeta \epsilon_T, \mathcal{H}_T, \|.\|_1) \leq  K^2\log  \mathcal{N}(\zeta_2 \epsilon_T, \mathcal{B}_T, \|.\|_1), \quad \Pi(B_2(\epsilon_T, B)) \geq \prod_{l,k} \Tilde \Pi_h\left(\norm{h_{lk} - h_{lk}^0}_2 < \epsilon_T\right).
\end{align*}
%Our Banach space is the space of functions $\mathbb{B} = \{f: [0,A] \to \R, \: f \in \ell^2([0,A])\}$ equipped with the $\ell^2$-norm. Let $K(.,.)$ be the covariance kernel of the Gaussian process and
These conditions are easily deduced from Theorem 2.1  in \cite{van_der_Vaart_2009} that we recall here. Let $\mathbb{H}$ be the Reproducing Kernel Hilbert Space of $k_{GP}$ 
% We recall that $\mathbb{H}$ is the completion of the following set of functions
% \begin{align*}
%     \left \{ f: x \in [0,A] \to \sum_{i=1}^n \alpha_i k(s_i,x) = \Ex{W_t H}, \quad H = \sum_{i=1}^n \alpha_i W_{s_i}, \: \alpha_i \in \R, \: s_i \in [0,A], \: \forall i \in [n], \: n \in \N, \right \},
% \end{align*}
% where $W$ is the squared exponential centered Gaussian process on $\R$. Moreover, the inner product on $\mathbb{H}$ is defined, for any $f,g$, as
% \begin{align*}
%     <f,g>_{\mathbb{H}} = <\Ex{WH_f}, \Ex{WH_g}>_{\mathbb{H}}  = \Ex{H_f H_g},
% \end{align*}
% with $H_f, H_g$ the random elements associated to $f$ and $g$. Now let
 and $\phi_{h_0}(\e)$ be the concentration function associated to $\Tilde \Pi_{h|\delta}$ defined as
\begin{align*}
    \phi_{h_0}(\e) = \inf_{h \in \mathbb{H}, \norm{h_{lk}-h_{lk}^0}_2 \leq \e} \norm{h_{lk}-h_{lk}^0}_{\mathbb{H}} - \log \Tilde \Pi(\norm{h_{lk}}_2 \leq \e), \quad \e > 0.
\end{align*}
%From Theorem 2.1  in \cite{van_der_Vaart_2009}
For any $\epsilon_T > 0$ such that $\phi_{h_0}(\epsilon_T) \leq T \epsilon_T^2$, there exists $ \mathcal{B}_T \subset L_2([0,A])$ satisfying
\begin{align*}
    &\Tilde \Pi_h(\mathcal{B}_T^c) \leq e^{-C T\epsilon_T^2}, \quad \log \mathcal{N}(3 \epsilon_T, \mathcal{B}_T, \|.\|_2) \leq 6C  T \e_T^2, \quad \Tilde \Pi_h(\norm{h_{lk} - h_{lk}^0}_\infty < 2\epsilon_T) \geq   e^{- T\epsilon_T^2},
\end{align*}
for any $C > 1$ such that $e^{-CT \epsilon_T^2}<1/2$. Since $\norm{h_{lk}}_1 \leq \sqrt{A} \norm{h_{lk}}_2$, we then obtain that
\begin{align*}
    &\log \mathcal{N}(3 \sqrt{A} \epsilon_T, \mathcal{B}_T, \|.\|_1) \leq  \log \mathcal{N}(3 \epsilon_T, \mathcal{B}_T, \|.\|_2) \leq 6C  T \epsilon_T^2,
\end{align*}
and finally, that $\log \mathcal{N}(\zeta_0 \epsilon_T, \mathcal{H}_T, \|.\|_1) \leq 6C  K^2  T \epsilon_T^2 \leq x_0 T \epsilon_T^2$ with $\zeta_0 = 3 \sqrt{A}, x_0 = 12 C K^2$.

Although more general kernel functions $k_{GP}$ could be considered, we focus on the hierarchical squared exponential kernels for which 
\begin{align*}
   \forall x,y \in \R, \quad k_{GP}(x,y; \ell) = \exp \left \{- (x-y)^2 / \ell^2 \right \}, \quad  \ell \sim IG(\ell; a_0, a_1 ), \quad a_0,a_1 > 0,
\end{align*}
where $IG(.; a_0, a_1 )$ with $ a_0,a_1 > 0$ is the Inverse Gamma distribution. The hierarchical squared exponential kernel is notably chosen in the variational method of \cite{malemshinitski2021nonlinear}, and its adaptivity and near-optimality has been proved by \cite{vzanten:vdv:09}.

\begin{proposition}\label{prop:gp}
Let $N$ be a Hawkes process with link functions $\phi = (\phi_k)_k$ and parameter $f_0 = (\nu_0, h_0) $ such that $(\phi,f_0)$ verify Assumption \ref{ass-psi}.  Assume that for any $l,k \in [K]$, $h_{lk}^0 \in \mathcal{H}(\beta, L_0)$ with $\beta > 0$ and $L_0 > 0$. Let $\Tilde \Pi_{h|\delta}$ be the above Gaussian Process prior with hierarchical squared exponential kernel $k_{GP}$. Then, under our hierarchical prior, the mean-field variational distribution $\hat Q_1$ defined in \eqref{eq:mean-field-vp} satisfies, for any $M_T \to +\infty$,
 \begin{align*}
        \Exz{\hat Q_{1} \left( \norm{f - f_0}_1 > M_T(\log \log T)^{1/2}(\log T)^q (T/ \log T)^{-\beta / (2\beta + 1)} \right)} \xrightarrow[T\to \infty]{} 0,
\end{align*}
with $q = 1$ if $\phi$ verifies Assumption \ref{ass-psi}(i) and $q = 3/2$ if $\phi$ verifies Assumption \ref{ass-psi}(ii).
\end{proposition}

 Given Theorem \ref{thm:cv_rate_vi}, Proposition \ref{prop:gp} is then a direct consequence of Theorem \ref{thm:cv_rate_vi} and  \cite{vzanten:vdv:09}, therefore its proof is omitted.

\begin{remark}
The Gaussian process prior has been used in variational methods for Hawkes processes when there exists a conjugate form of the mean-field variational posterior distribution, i.e., $\hat Q_1$ is itself a Gaussian process with mean function $m_{VP}$ and kernel function $k_{VP}$. This is notably the case in the sigmoid Hawkes  model under the latent variable augmentation scheme described in Section \ref{sec:aug-mf-vi} and used for instance by\cite{malemshinitski2021nonlinear}. Since the computation of the  Gaussian process variational distribution is often expensive for large data set, the latter is often further approximated using the sparse Gaussian process approximation via inducing variables \citep{titsias_NIPS2011}. %This leads to an ``inducing variable" mean-field variational posterior, which is itself a variational approximation of the original mean-field variational posterior \citep{dennis21}.
Using results of \cite{dennis21}, we conjecture that our result in Proposition \ref{prop:gp} would also hold for the mean-field variational posterior with inducing variables.   
\end{remark}

\subsubsection{Model-selection variational family}\label{conc-sas}

In this section, we consider the model-selection adaptive variational posterior distributions \eqref{eq:ms_var_post} and \eqref{eq:ms_var_post_avg}, and similarly obtain their concentration rates. We recall that these two types of adaptive variational posterior correspond to the following variational families (see also Appendix \ref{app:model-selection})
\begin{align*}
    &\mathcal{V}_{A1} = \cup_{m \in \mathcal M} \{ \{m\}\times \mathcal{V}^m \},
    &\mathcal{V}_{A2} = \left \{ \sum_{m \in \mathcal{M}} \alpha_{m} Q_m; \sum_{m} \alpha_{m}  = 1, \: \alpha_m \geq 0, \: Q_m \in \mathcal{V}^m, \: \forall m \in \mathcal{M} \right\},
\end{align*}
where here, $\mathcal{M}$ is the set of all possible models, i.e.,
\begin{align*}
    \mathcal{M} = \left \{ m = (\delta, J=(J_1, \dots, J_K)); \delta \in \{0,1\}^{K\times K}, \: J_k \in \mathbb N, \: \forall k \in [K]\right \},
\end{align*}
and for a model $m \in \mathcal{M}$, the variational family $\mathcal{V}^m$ corresponds to a set of distributions on the subspace $\mathcal{F}_m \subset \mathcal{F}$ and $\bigcup_{m \in \mathcal{M}} \mathcal{F}_m = \mathcal{F}$. In the data augmentation context and with the mean-field approximation, $\mathcal{V}^m$ is the set of distributions $Q: \mathcal{F}_m  \times \mathcal{Z} \to [0,1]$ such that $Q(f,z) = Q(f)Q(z)$. We further recall that for each $k$, $J_k$ corresponds to the number of functions in the dictionary used to construct $(h_{lk})_{l \in [K]}$.  

In this context, the general results from \cite{Zhang2017ConvergenceRO} can be applied, and here, it is enough to replace the prior assumption \textbf{(A0)} by
\begin{align}\label{eq:A0''}
    \textbf{(A0'')} \quad \exists c_1 > 0, \: &\Pi \left(  B_\infty(\epsilon_T) \left| \delta= \delta_0, J= (J_k^0 )_k J_T \right. \right)   \geq e^{ -c_1 T\epsilon_T^2/3}, \nonumber \\
    &\Pi_\delta ( \delta = \delta_0) \geq e^{ -c_1 T\epsilon_T^2/3}, \quad  \Pi_J \left( J=\left(J_k^0 \right)_k J_T\right) \geq e^{ -c_1 T\epsilon_T^2/3},
\end{align}
where $J_T = \left(\frac{T}{\log T}\right)^{\beta / (2 \beta +1)}$,
assuming that,  for any $l,k \in [K]$, $h_{lk}^0 \in \mathcal{H}(\beta, L_0)$. Indeed,  \textbf{(A0'')} implies that
\begin{align*}
    - \log \Pi(m = m_0) - \log \Pi \left(  B_\infty(\epsilon_T) | m=m_0 \right) \leq c_1 T \epsilon_T^2, \quad m_0 = \left(\delta_0, \left(J_k^0 \right)_k J_T\right),
\end{align*}
which also implies \textbf{(A0)}. For example, under the random histogram prior of Section \ref{conc-mf-vi}, it is enough to choose $\Pi_J$ such that, for some sequence $(x_n)_{n\geq 1}$ such that $x_n \xrightarrow[n \to \infty]{}\infty$,
\begin{align*}
    \Pi_J(J_l > x_n) \lesssim e^{-cx_n}, \quad \Pi_J(J_l = x_n) \gtrsim e^{-cx_n}, \quad \forall n\geq 1, \quad c > 0,
\end{align*}
which is the case for instance when $\Pi_J$ is a Geometric distribution. In the next proposition, we state our result on the model-selection variational family, when using the random histogram prior distribution; however, this result also holds for other prior distributions based on decomposition over dictionaries such as the ones in \cite{arbeletal:13, shen:ghosal14}.

\begin{proposition}\label{prop:ss-conc}
Let $N$ be a Hawkes process with link functions $\phi = (\phi_k)_k$, parameter $f_0 = (\nu_0, h_0) $ such that $(\phi,f_0)$ verify Assumption \ref{ass-psi}. Assume that for any $l,k \in [K]$, $h_{lk}^0 \in \mathcal{H}(\beta, L_0)$ with $\beta \in  (0,1)$ and $L_0 > 0$. %Let $\epsilon_T = o(1/\sqrt{\kappa_T})$ be a positive sequence verifying $\log^3 T=O(T \epsilon_T^2)$ and $\Pi$ be a prior distribution on $\mathcal{F}$ satisfying \textbf{(A0')}  and \textbf{(A1)}.
Then, under the random histogram prior distribution, for the model selection variational posterior \eqref{eq:ms_var_post} , we have that, for any $M_T \to +\infty$,
 \begin{align*}
       \Exz{\hat Q_{A1} \left( \norm{f - f_0}_1 >  M_T(\log T)^q (T / \log T)^{-\beta / (2\beta + 1)} \right)} \xrightarrow[T\to \infty]{} 0,
    \end{align*}
with $q = 0$ if $\phi$ verifies Assumption \ref{ass-psi}(i) and $q = 1/2$ if $\phi$ verifies Assumption \ref{ass-psi}(ii). 
\end{proposition}
%We note that explicit concentration rates for H\"{o}lder-smooth functions can then be derived when using the prior families of Section \ref{conc-mf-vi}. %for the slab distribution $\Tilde \Pi_{h|\delta}$.

Since Proposition \ref{prop:ss-conc} is a direct consequence of Theorem \ref{thm:cv_rate_vi} and Theorem 4.1 in \cite{Zhang2017ConvergenceRO}, its proof is omitted. Finally, we note that we can obtain similar guarantees for the model-averaging adaptive variational posterior \eqref{eq:ms_var_post_avg}, by adapting Theorem 3.6 from \cite{ohn2021adaptive}, which directly holds under the same assumptions as Proposition \ref{prop:ss-conc}. 

\subsection{Convergence rate associated to the two-step algorithm} 
 As discussed in Section \ref{sec:var_bayes}, when the number of dimensions $K$ is moderately large,  both the $\hat Q_{A1}$ and $\hat Q_{A2}$ are intractable, due to the necessity of exploring all models in $\mathcal{M}_T$, defined in \eqref{eq:set-models-bounded}. For this setting, we have proposed a two-step procedure (Algorithm \ref{alg:2step_adapt_cavi}) that first constructs the estimator of the graph with \eqref{eq:graph-estimator}, then constructs a restricted set of models $\mathcal{M}_E$ and computes the corresponding variational distribution $\hat Q^{\hat \delta}$. We now show that this two-step procedure is theoretically justified. We recall our notation $S_{lk}^0=\norm{h_{lk}^0}_1, \: \forall l,k \in [K]$.

Firstly, since the complete graph $\delta_C = \mathds{1}\mathds{1}^T$ is larger than the true graph $\delta_0$, the subspace $\bigcup_{m \in \mathcal{M}_C} \mathcal{F}_m$ contains the true parameter $f_0$. %\textcolor{magenta}{Qu'est ce que ca veut dire ? : the model conditional on $\delta_C$ is correctly specified.}
Hence Theorem \ref{thm:cv_rate_vi} remains valid with $\mathcal{V}_C = \cup_{m \in \mathcal M_C} \{ \{m\}\times \mathcal{V}^m \}$. %$\mathcal V = \mathcal V^{(\delta_C)}$. 
In particular the rates $\epsilon_T = (\log T)^q T^{-\beta/(2\beta+1)}$ obtained in Propositions \ref{prop:histo} and \ref{prop:gp} apply to the corresponding variational posterior $\hat Q_{MS}^{C}$, under the assumption that $K$ is large but fixed. In particular, for each $(l,k)$ and $\hat S_{lk} = \int\norm{h_{lk}}_1dQ_{MS}^{C}(h_{lk}) $, Theorem \ref{thm:cv_rate_vi} implies that 
  $$\mathbb P_0 \big(|\hat S_{lk} - S_{lk}^0| > \epsilon_T \big) =o(1).$$
In our two-step procedure, we consider the two following thresholding strategies:
  \begin{enumerate}
      \item[(i)] given a threshold $\eta_0 > 0$ defined \emph{a-priori}, we compute $\hat \delta = (\hat \delta_{lk})_{l,k}, \: \delta_{lk} = \mathds{1}_{\{\hat S_{lk} > \eta_0\}}, \forall l,k$. 
      \item[(ii)] we choose a data-dependent threshold $\eta_0 \in (\hat S_{(i_0)}, \hat S_{(i_0+1)})$, where $(\hat S_{(i)})_{i \in [K^2]}$ corresponds to the values $(\hat S_{lk})_{l,k}$ in increasing order and $i_0$ is the first index such that $\hat S_{(i+1)} - \hat S_{(i)}$ is \emph{large}. We then compute $\hat \delta$ as in (i).
  \end{enumerate}
Let $i^* := K^2 - \sum_{l,k}\delta_{lk}^0 = \min \left \{i \in [K^2]; \: S^0_{(i+1)} -S^0_{(i)} \neq 0 \right \}$ be  the first index of non-zero such that $S^0_{(i+1)} > 0$, %where $ S^0_{(i+1)} -S^0_{(i)} \neq 0$,
where $(S^0_{(i)})_{i \in [K^2]}$ corresponds to the values of $(S_{lk}^0)_{l,k}$ in increasing order. %\textcolor{magenta}{on a besoin de $i^*$ ?}
We recall our notation $\mathcal I(\delta_0)$ for the set of index pairs $(l,k)$ such that $S_{lk}^0> 0$. We now assume that $f_0$ is such that  
  \begin{equation}\label{signal_detection}
 S_{lk}^0 \geq u_T, \quad  \forall l,k \in \mathcal I(\delta_0),
  \end{equation} 
 where $u_T>> \epsilon_T$. We note that \eqref{signal_detection} is a mild requirement on $f_0$ since we allow $u_T$ to go to 0 almost as fast as $\epsilon_T$. 
Now, for the thresholding strategy (i), for any $\eta_0$ (possibly depending on $T$) such that $u_T \leq \eta_0 <\min_{(l,k) \in \mathcal I(\delta_0)}\norm{h_{lk}^0}_1/2  $, we obtain that 
\begin{align}\label{eq:graph-consistency}
    \mathbb P_0( \hat \delta \neq \delta_0 ) = o(1).
\end{align}
Moreover, for the data-dependent thresholding strategy (ii), as soon as the gap $\hat S_{(i+1)} - \hat S_{(i)}$ is larger than $u_T$ but smaller than $\min_{lk\in \mathcal I(\delta_0)}\norm{h_{lk}^0}_1/2$, then \eqref{eq:graph-consistency} also holds.
%$$\mathbb P_0( \hat \delta \neq \delta_0 ) = o(1) $$
This is verified since
   \begin{align*}
      \mathbb P_0( \hat \delta \neq \delta_0)
       \leq \sum_{l,k \in [K]} \mathbb P_0( |\hat S_{lk}-S_{lk}^0|> u_T/2 )=o(1).
   \end{align*}

\section{Numerical results}\label{sec:numerical}

In this section, we perform a simulation study to evaluate our variational Bayesian method in the context of nonlinear  Hawkes processes, and demonstrate its efficiency, scalability, and robustnessin various estimation setups. In low-dimensional settings ($K = 1$ and $K = 2$), we can compare our variational posterior to the posterior distribution obtained from an MCMC method. As a preliminary experiment, we additionally analyse the performance of a Metropolis-Hastings sampler  in commonly used nonlinear Hawkes processes, namely with ReLU, sigmoid and softplus link functions (Simulation 1). In the subsequent simulations, we focus on the sigmoid model and test our adaptive variational  algorithms, in well-specified (Simulations 2-5) and mis-specified settings (Simulation 6), high-dimensional data sets, and for different connectivity graphs (Simulation 4).

In each setting, we sample one observation of a Hawkes process with dimension $K$, link functions $(\phi_k)_k$ and parameter $f_0 = (\nu_0, h_0)$ on $[0,T]$, using the thinning algorithm of \cite{adams09}.
In most simulated settings, the true interaction functions $(h_{lk}^0)_{l,k}$ will be piecewise-constant, and we use the random histogram prior described in Section \ref{sec:aug-mf-vi} in our variational Bayes method. For  $D\geq 1$, we  introduce the notation
\begin{align*}
    \mathcal{H}_{histo}^D = \left \{ h_k = (h_{lk})_{l} ; \: h_{lk}(x) = \sum_{j=1}^{2^D} w^j_{lk} e_j(x), \: x \in [0,A], \: l \in [K], \: e_j(x) = \frac{2^D}{A}\mathds{1}_{[\frac{j A}{2^D},\frac{(j+1)A}{2^D}}))(x) \right \},
\end{align*}
and for the remaining of this section, we index functions $h_{lk}$ by the histogram depth $D$.

In the next sections, we report the results of the following set of simulations. %\textcolor{magenta}{Peut etre s'assurer que pour chaque Simulation le valeur de $T$ est bien specifiee} 
\begin{itemize}
    \item \textbf{Simulation 1: Posterior distribution in parametric, univariate,  nonlinear Hawkes models.}
    We analyse the posterior distribution computed from a Metropolis-Hasting sampler (MH) in several nonlinear univariate Hawkes processes ($K=1$), with ReLU, sigmoid, and softplus link functions. For this sampler, we consider that the dimensionality $D_0$ such that $h_0 \in  \mathcal{H}_{histo}^{D_0}$ is known, and therefore, the posterior inference is non-adaptive. %We note that the MCMC iterations are computationally expensive, which prevents from scaling this algorithm to large dimensions. % where $h_0 \in \mathcal{H}_{hist}^{D_0}$ and $D_0 \geq 1$ is known.

    \item \textbf{Simulation 2: Variational and true posterior distribution in parametric, univariate sigmoid Hawkes models.} In a univariate setting with  $h_0 \in  \mathcal{H}_{histo}^{D_0}$ and the dimensionality $D_0$ is known (non-adaptive), we compare the variational posterior obtained from Algorithm \ref{alg:cavi} to the posterior distribution obtained from two MCMC samplers, i.e., the MH sampler of Simulation 1, and a Gibbs sampler available in the sigmoid model (Algorithm \ref{alg:gibbs}). % (see Section \ref{sec:DA_sigmoid} and Appendix \ref{app:gibbs_sampler}), and the fixed-dimension mean-field  variational inference algorithm (Algorithm \ref{alg:cavi}).
    % therefore our variational  algorithm (Algorithm \ref{alg:cavi}) reduces to the semi-parametric approach of \cite{zhou2021nonlinear}. Therefore, our numerical results are related to those of \cite{zhou2021nonlinear} but our setting involves a different function basis. 
    
    \item \textbf{Simulation 3: Fully-adaptive variational algorithm in univariate and bivariate sigmoid models.} This experiment evaluates our first adaptive variational algorithm (Algorithm \ref{alg:adapt_cavi}) in sigmoid Hawkes processes with $K=1$ and $K=2$, in nonparametric settings where the true interaction functions are either piecewise-constant functions with unknown dimensionality or continuous. %We also compare to the MH algorithm. We note that in the bivariate setting $K=2$, there are $2^{K^2} = 16$ possible graph parameters, and therefore $16 \times D_T$ ``models" to include in the adaptive variational posterior.
    %In this case, we  use the adaptive variational posterior on the set of histogram sizes $\{D=1,\dots, D_T\}$.
    % We will simulate from two settings of the true parameter $f_0$, where $D_0 \geq 1$ is an unknown parameter:
    % \begin{enumerate}
    %     \item $h_0 \in \mathcal{H}_{hist}^{D_0}$;
    %     \item $h_0$ is a continuous functions with non-zero $D_0$ first components in the Fourier basis.
    %      \end{enumerate}
         
    \item \textbf{Simulation 4: Two-step adaptive variational algorithm in high-dimensional sigmoid models.} This experiment evaluates the performance and scalability of our fast adaptive variational algorithm (Algorithm \ref{alg:2step_adapt_cavi}), for sigmoid Hawkes processes with $K \in \{2,4,8,10,16, 32, 64\}$, in sparse and less sparse settings of the true parameter $h_0 \in  \mathcal{H}_{histo}^{D_0}$ with unknown dimensionality $D_0$.

    \item  \textbf{Simulation 5:  Convergence of the two-step adaptive variational posterior for varying data set sizes.} In this experiment, we evaluate the asymptotic performance of our two-step  variational procedure (Algorithm \ref{alg:2step_adapt_cavi}), with respect to the number of observations, i.e., the length of the observation horizon $T$, for sigmoid Hawkes processes with $K =10$.

    \item  \textbf{Simulation 6: Robustness of the variational posterior to some types of mis-specification of the Hawkes model.} This experiment aims at evaluating the performance our variational algorithm for the sigmoid Hakwes model (Algorithm \ref{alg:2step_adapt_cavi}) on data sets generated from Hawkes processes with mis-specified nonlinear link functions and memory parameter of the interaction functions.
   
\end{itemize}

%In all simulations, we set $A=0.1$.

% \ds{We note that in our implementation of the Gibbs sampler, we use the package polyagamma. We also compare to
% \begin{itemize}
%     \item the Nonlinear Hawkes Process with Gaussian Process Self-Effects algorithm (NH-GPS)  \cite{malemshinitski2021nonlinear}
%     \item the Cross-Covariance estimator for estimating the graph \cite{Chen2017, chen17b, cai2021latent}? (ask code to Shizhe Chen or Biao Cai?)
% \end{itemize}}
%We will simulate a Hawkes process with parameter $f_0$ and estimate the posterior distribution on $f$ using a Metropolis-Hasting and a HMC algorithms, and in the sigmoid model, also with a mean-field variational inference (MF-VI) algorithm and Gibbs sampler.

% We first introduce some notation. For an observation $N$, we define $H(t) = (H_b(t))_{b=0, \dots, K B}$ with
% \begin{align*}
%     &H_b(t) = \int_{t-A}^t \alpha_b(t-s)dN^{k_b}_s, \quad k_b \in [K], \quad b = 1, \dots, K B, \\
%     &H_0(t) = 1.
% \end{align*}
% In all models \eqref{eq:nonlinearity}, we consider estimating the parameter $f = (f_k)_k, \: f_k = (\nu_k, h_{1k}^1, \dots, h_{Kk}^B) \in \R^{KB+1}, \: k \in [K]$. We also denote $ h(t,f_k) = f_k^T H(t)$, therefore, $\Tilde{\lambda}_t^k = \alpha ( h(t,f_k) - \eta)$.

In all simulations, we set the memory parameter as $A=0.1$, and we evaluate the performance visually, in low-dimensional settings, or with the $L_1$-risk on the continuous parameter and $\ell_0$-error on the graph parameter (defined below), in moderately large to large-dimensional settings. %Additional details on these experiments are reported in Appendix \ref{app:details_exp}. 

\begin{remark}\label{rem:excursions}
    One important quantity in these synthetic experiments is the number of \emph{excursions} in the generated data, formally defined in \cite{costa18} and Lemma \ref{lem:excursions} in Appendix \ref{app:main_lemmas}. Intuitively, the observation window of the data $[0,T]$ can be partitioned into contiguous  intervals $\{[\tau_{i-1},\tau_i)\}_{i=1,\dots, I}$, $\tau_0=0, \tau_{I}=T$, $I \in \mathbb N$, called \emph{excursions}, where the point process measures are i.i.d. The main properties of these intervals are that $N[\tau_{i-1}, \tau_i) \geq 1$ and $N[\tau_{i}-A, \tau_i)=0$. For our multivariate contexts, we additionally introduce a new concept of excursions, that we call \emph{local} excursions, defined for each dimension $k$ as a partition of $[0,T] = \bigcup_{i=1}^{I_k} [\tau^k_{i-1},\tau^k_i)$ such that $N^k[\tau^k_{i-1}, \tau^k_i) \geq 1$ and $N^k[\tau^k_{i}-A, \tau^k_i)=0$. To the best of our knowledge, this quantity has not yet been introduced for Hawkes processes, although we observe in our experiments that it is an important statistical property, as will be shown below.
\end{remark}

\FloatBarrier

\subsection{Simulation 1: Posterior distribution in univariate nonlinear Hawkes models }\label{sec:exp_1D_fixed}

\begin{figure}[hbt!]
    \centering
    \includegraphics[width=0.5\textwidth, trim=0.cm 0.cm 0cm  0.cm,clip]{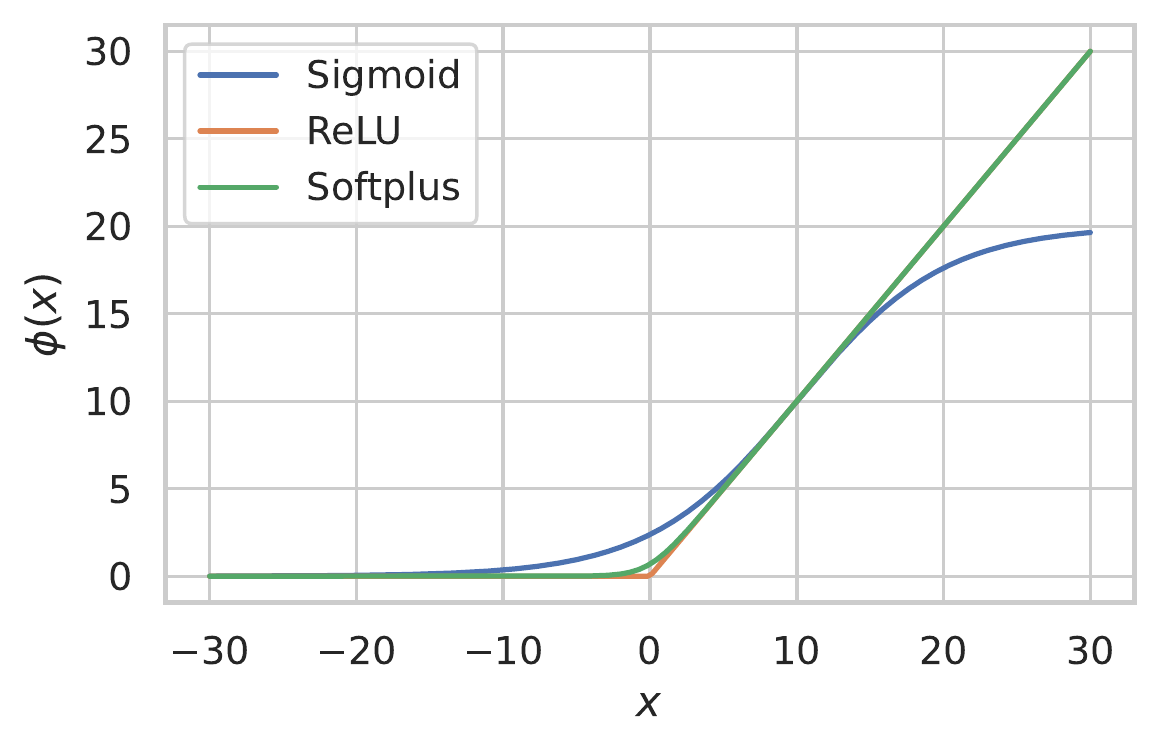}
    \caption{Link functions $\phi$ of the Hawkes model considered in Simulation 1, namely the sigmoid (blue), ReLU (red), and softplus (green) functions.}
    \label{fig:links}
\end{figure}

\begin{figure}[hbt!]
\setlength{\tempwidth}{.3\linewidth}\centering
\settoheight{\tempheight}{\includegraphics[width=\tempwidth, trim=0.cm 0.cm 0cm  0.65cm,clip]{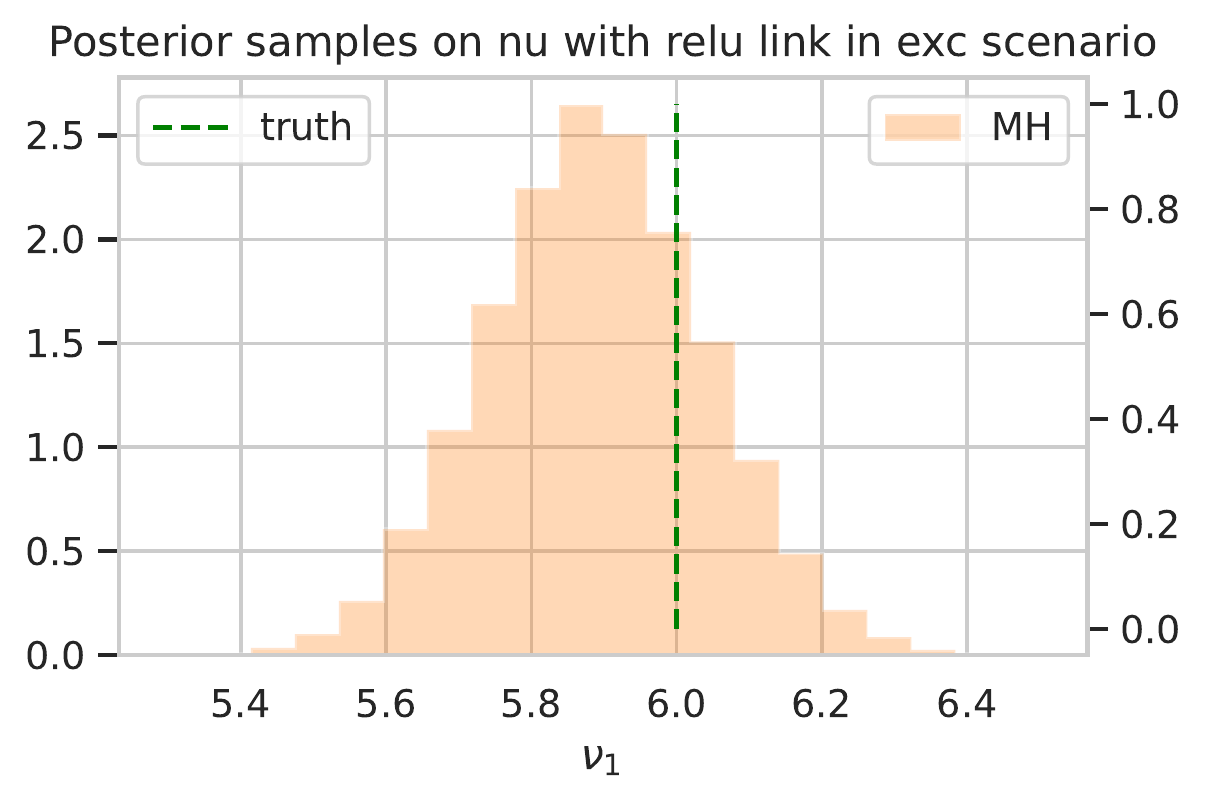}}%
\hspace{-5mm}
\fbox{\begin{minipage}{\dimexpr 15mm} \begin{center} \itshape \large  \textbf{ReLU} \end{center} \end{minipage}}
\hspace{-5mm}
\columnname{\emph{Excitation} only}\hfil
\columnname{Mixed effect }\hfil
\columnname{Inhibition only}\\
\rowname{Background}
    \includegraphics[width=\tempwidth, trim=0.cm 0.cm 0cm  0.65cm,clip]{figs/relu_d4_exc_nu_mh_True_hmc_False.pdf}\hfil
    \includegraphics[width=\tempwidth, trim=0.cm 0.cm 0cm  0.65cm,clip]{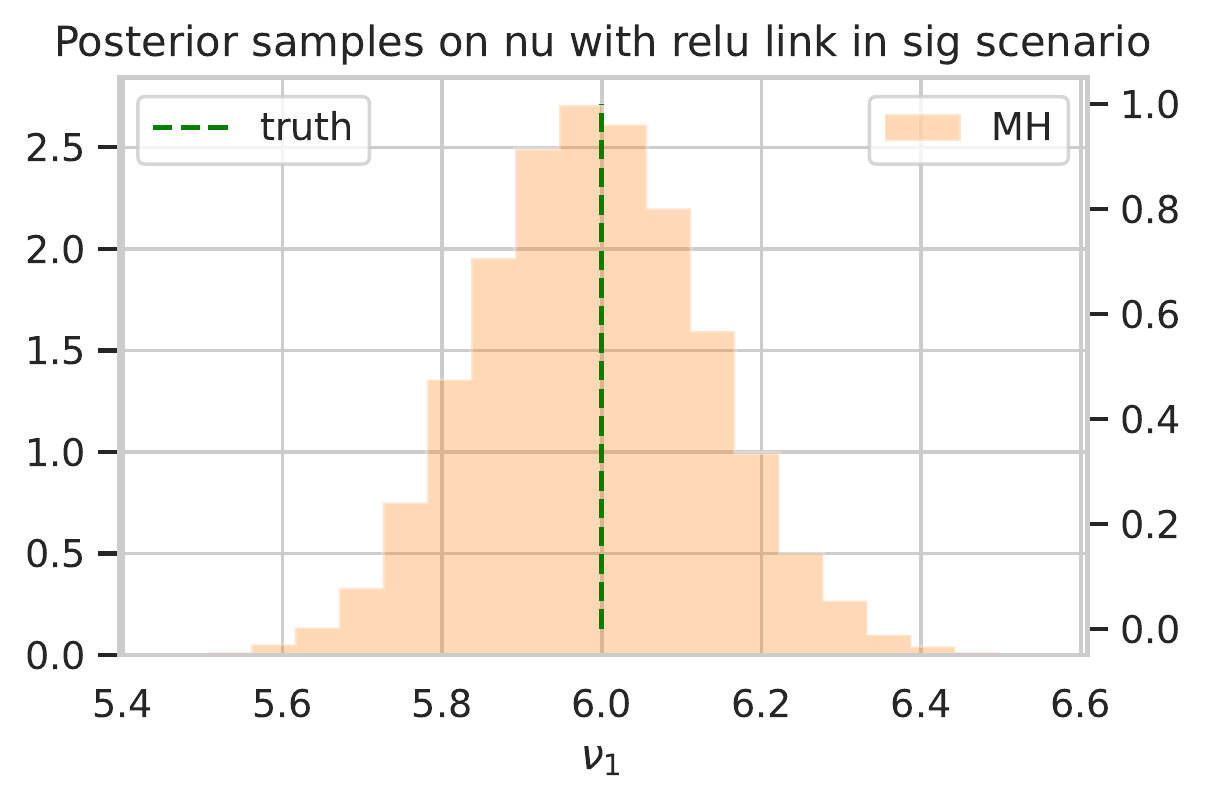}\hfil
    \includegraphics[width=\tempwidth, trim=0.cm 0.cm 0cm  0.65cm,clip]{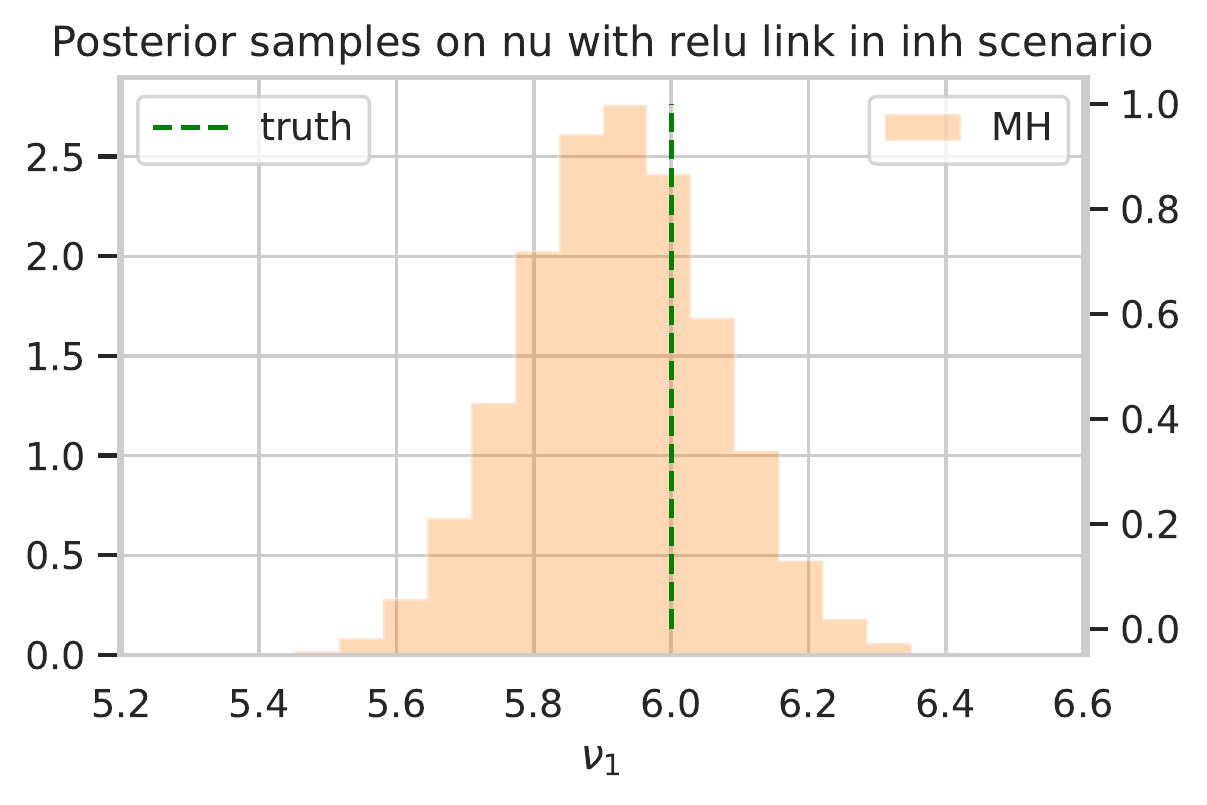}\\
\rowname{Interaction}
    \includegraphics[width=\tempwidth, trim=0.cm 0.cm 0cm  0.65cm,clip]{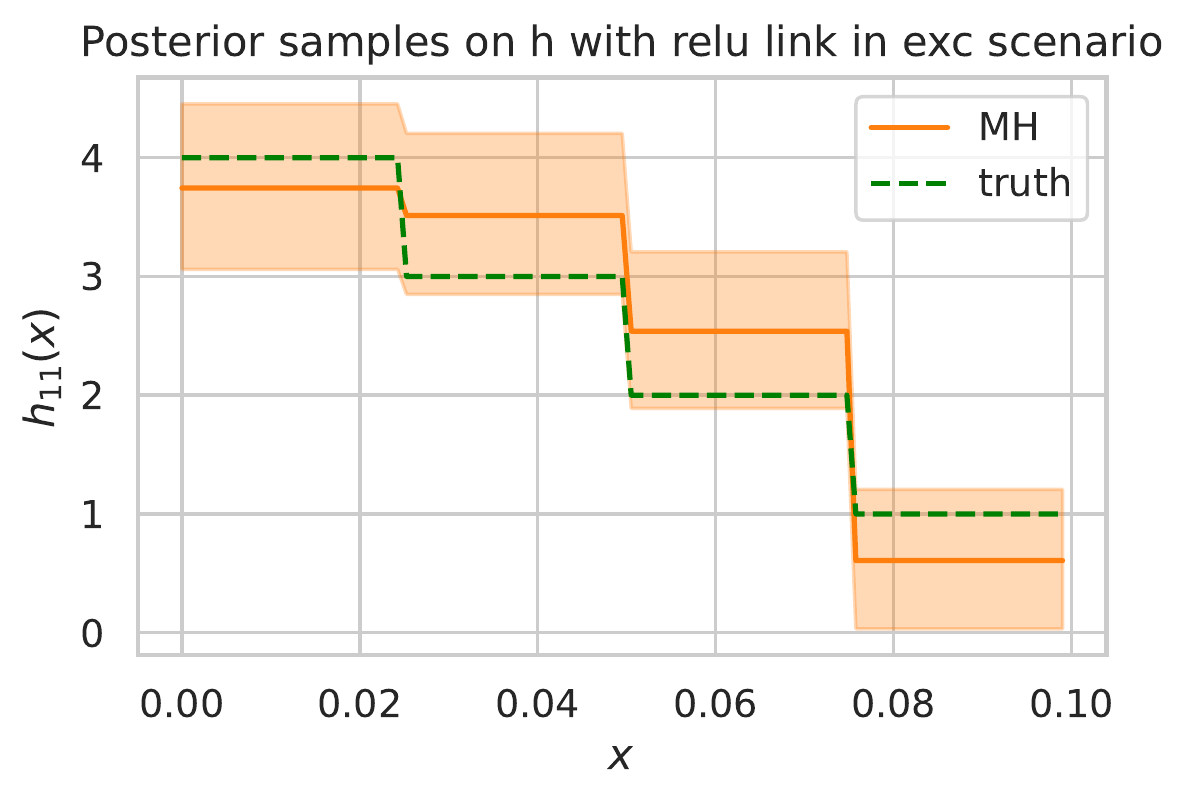}\hfil
    \includegraphics[width=\tempwidth, trim=0.cm 0.cm 0cm  0.65cm,clip]{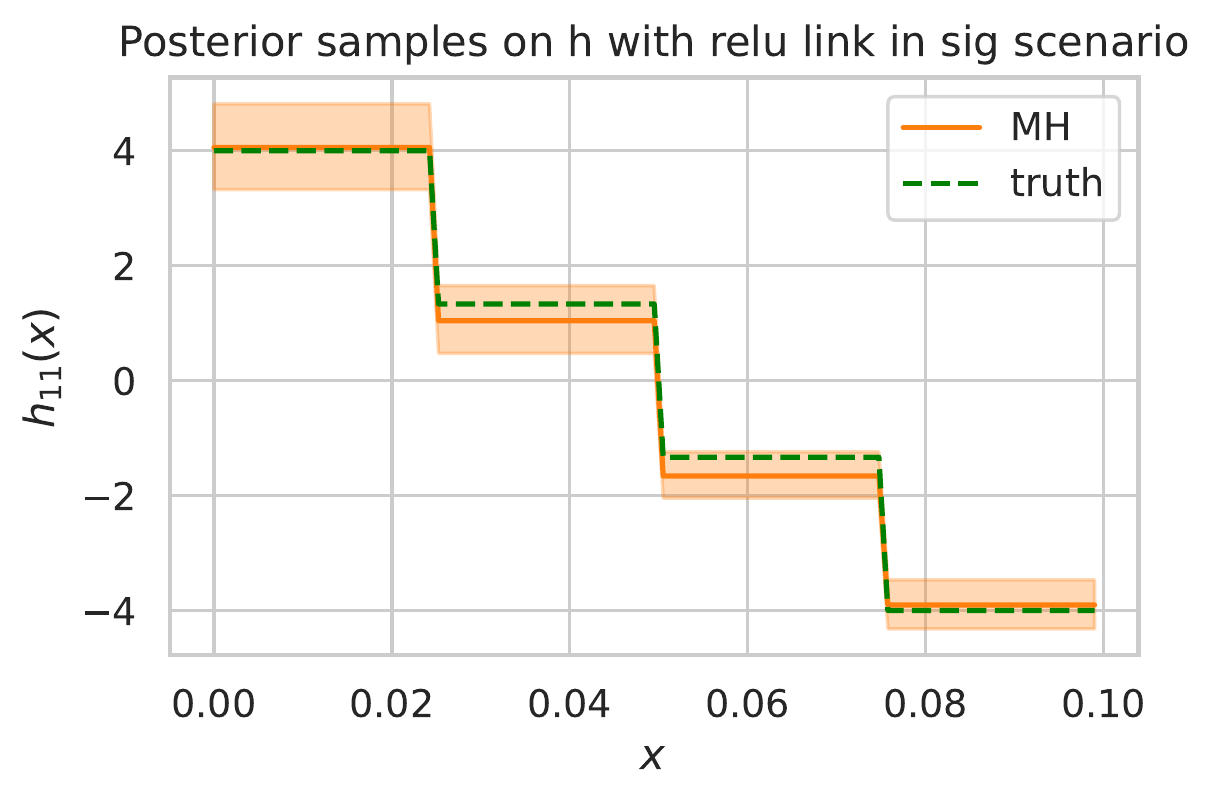}\hfil
    \includegraphics[width=\tempwidth, trim=0.cm 0.cm 0cm  0.65cm,clip]{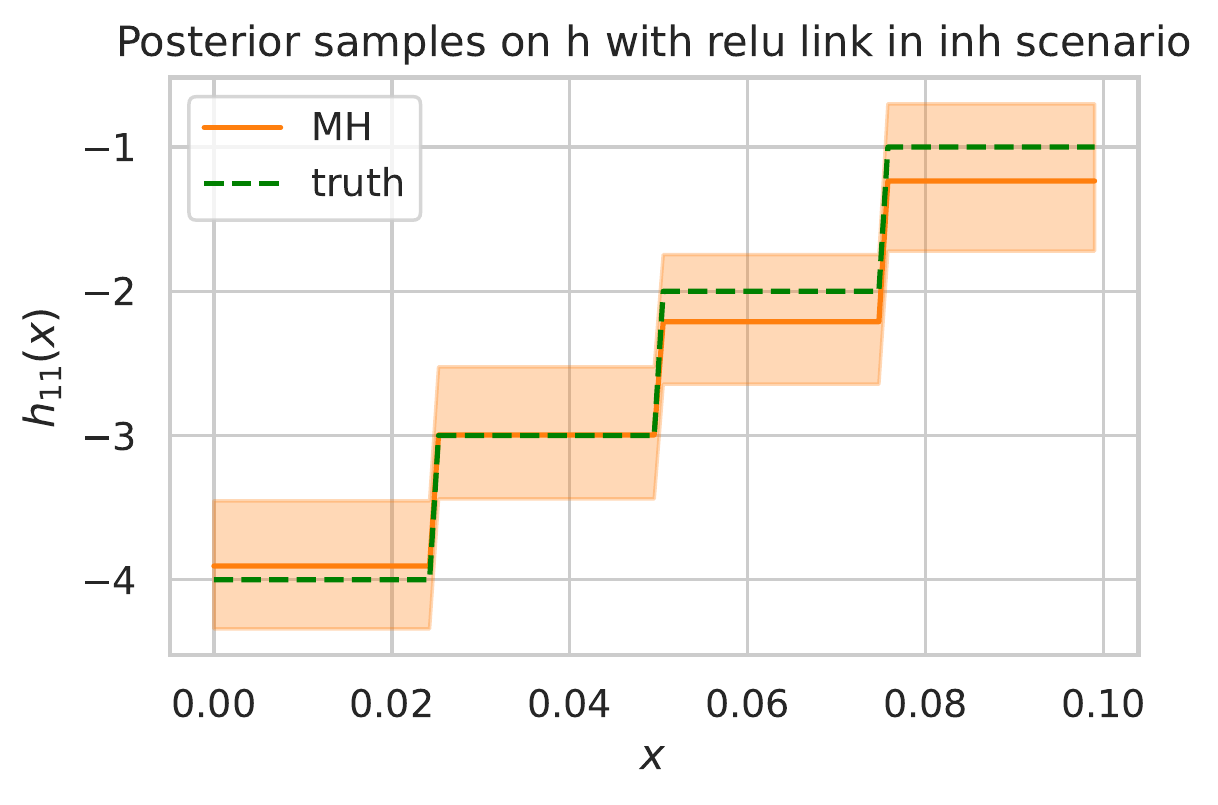}\\
\caption{Posterior distribution on $f = (\nu_1, h_{11})$ obtained with the Metropolis-Hastings sampler (MH), in the univariate ReLU models of Simulation 1. The three columns correspond to the \emph{Excitation only} (left), \emph{Mixed effect} (center), and \emph{Inhibition only} (right) scenarios. On the first row, we plot the marginal posterior distribution on the background rate $\nu_1$, and on the second row, the posterior mean (solid orange line) and  95\% credible sets (orange areas) on the interaction function $h_{11}$, here piecewise-constant with dimensionality $2^{D_0} = 4$. The true parameter $f_0 = (\nu_1^0, h_{11}^0)$ is plotted in dotted green line.}
\label{fig:relu_mcmc_D4}
\end{figure}

In this simulation, we consider univariate Hawkes processes ($K=1$) with link function $\phi = \phi_1$ of the form 
\begin{align}\label{eq:nonlinearity}
    &\phi(x) = \theta +   \Lambda \psi(\alpha (x - \eta)), 
\end{align}
where $\xi = (\theta, \Lambda, \alpha, \eta) $ and $\psi:\R \to \R^+$ are known and chosen as:
\begin{itemize}
    \item Sigmoid: $\psi(x) = (1 + e^{-x})^{-1}$ and $\xi = (0.0, 20.0, 0.2, 10.0)$;
    \item ReLU: $\psi(x) =  \max(x,0)$ and $\xi = (0.001, 1.0, 1.0, 0.0)$;
    \item Softplus: $\psi(x) = \log(1 + e^x)$ and $\xi = (0.0, 40.0, 0.1, 20.0)$.
\end{itemize}
Note that the corresponding link functions $\phi$ have similar shapes on a range of values between -20 and 20 (see Figure \ref{fig:links}). In all models, we consider a Hawkes process with %$\nu_0 = 6$ and
$h_0 = h_{11}^0  \in \mathcal{H}_{histo}^{D_0}$ with $D_0=2$, and three scenarios, called \emph{Excitation only}, \emph{Mixed effect}, and \emph{Inhibition only}, where $h_0$ is respectively non-negative, signed, and non-positive (see Figure \ref{fig:relu_mcmc_D4} for instance). In each of the nine settings, we set $T=500$ and in  Table \ref{tab:simu_data}, we report the corresponding number of events and excursions observed in each scenario and model. Note that, as we may expect, more events and less excursions are observed in the data generated in \emph{Excitation only} scenario than in the \emph{Mixed effect} and \emph{Inhibition only} scenarios.

Here, we assume that $D_0$ is known and we consider a normal prior on $ \mathcal{H}_{histo}^{D_0}$ such that %mon $w_{11}$ such that %$\nu_1$ and $h_{11}$ 
$
    w_{11}  \sim \mathcal{N}(0, \sigma^2 I),
$
and  for $\nu_1$, 
$
     \nu_1 \sim \mathcal{N}(0, \sigma^2),
$
with $\sigma=5.0$. 
%$\Pi(f) = \mathcal{N}(f; \mu, \Sigma)$, with  $\mu = 0_{KB+1}, \Sigma =  \sigma^2*I_{KB+1}$,  $\sigma=3.0$. % $\mu \in \R^{B+1}, \Sigma \in \R^{(B+1) \times (B+1)} $ with $B = 2^{D_0} = 4$. In fact,  we choose a centered, isotropic Gaussian prior with $\mu = 0_{KB+1}, \Sigma =  \sigma^2*I_{KB+1}$ with  $\sigma=5.0$.
%The number of unknown parameters is thus $N_{param} = 2^{D_0} + 1 = 5$.
%with $(h_0^b)_{b=1,\dots, 2^{D_0}}$ such that $h_0 = \sum_b h_0^b \alpha_b$,
%We simulate one realisation of the process until time $T=500$.
 %two MCMC samplers in the package PyMC3 \footnote{\url{https://docs.pymc.io/en/v3/index.html}},
To compute the (true) posterior distribution, we run a  Metropolis-Hasting (MH) sampler implemented via the Python package PyMC4\footnote{\url{https://www.pymc.io/welcome.html}} with  4 chains, 40 000 iterations, and a burn-in time of 4000 iterations. We also use the Gaussian quadrature method \citep{golub1969calculation}  for evaluating the log-likelihood function, except in the ReLU model and \emph{Excitation only} scenario, where the integral term is computed exactly. We note that we also tested a Hamiltonian Monte-Carlo sampler in this simulation, and obtained similar posterior distributions, but within a much larger computational time, therefore these results are excluded from this experiment.

The posterior distribution on $f = (\nu_1, h_{11})$ in the ReLU model and our three scenarios are plotted in Figure \ref{fig:relu_mcmc_D4}. For conciseness purpose in this section, our results for the sigmoid and softplus models are reported in Appendix~\ref{app:simu1}. We note that in almost all settings, the ground-truth parameter $f_0$ is included in the 95\% credible sets of the posterior distribution. %, except in the \emph{Excitation only} scenario in the softplus model.
Nonetheless, the posterior mean is sometimes biased, %in particular in the Excitation scenario, 
possibly due to the numerical integration errors in the log-likelihood computation. Moreover, we conjecture that the estimation quality depends on the number of events and the number of excursions, which could explain the  differences between the \emph{Excitation only}, \emph{Mixed effect}, and \emph{Inhibition only} scenarios. In particular, the credible sets seem consistently smaller for the second scenario, which realisations have more excursions than the other ones.  %In fact, this bias seems to disappear when the integral in the log-likelihood is computed exactly (i.e., in the ReLU model and \emph{Excitation only} scenario). 

This simulation therefore shows that the posterior distribution in commonly used nonlinear univariate Hawkes models behaves well and can be sampled from using a simple  MH sampler. Nonetheless, we note that the MH iterations are computationally expensive, which prevents from scaling this algorithm to large dimensions. Therefore, we will only  use the MH sampler to compute the posterior distribution in the low-dimensional settings, i.e., Simulations 2 and 3, with respectively  $K=1$ and $K=2$. 
%Besides, these results can be seen as an illustration of the theoretical results of \cite{sulem2021bayesian} on the asymptotic behaviour of the posterior distribution for general nonlinear Hawkes models. 
%We note that the posterior distributions obtained with HMC and MH are very similar in all settings and they both reasonably well estimate the ground-truth parameter $f_0$. However, the MH sampler is about 40 times faster than HMC (see the computational times in Table \ref{tab:comptime_mcmc}), therefore, in the subsequent simulations (Simulations 2,3, and 4), we will only consider the MH sampler.

%therefore this is the sampler that we compare to the Gibbs and Mean-Field VI algorithm in the sigmoid model below. 

\begin{table}[hbt!]
    \centering
\begin{tabular}{c|c|c|c|c}
\toprule
  Scenario & &  Sigmoid & ReLU &  Softplus \\
 \midrule
 \multirow{2}{*}{\emph{Excitation} only} & \# events & 5250 &   5352 & 4953 \\
 & \# excursions &  1558 &    1436 &  1373 \\
\midrule
 \multirow{2}{*}{Mixed effect} & \# events &  3876 &   3684 & 3418 \\
 & \# excursions &  1775  &     1795  &    1650 \\
 \midrule
  \multirow{2}{*}{Inhibition only} & \# events & 3047 &   2724  &  2596  \\
 & \# excursions &    1817 &    1693  &   1588 \\
\bottomrule
\end{tabular}
    \caption{Number of events and excursions in the  simulated data  of Simulation 1 with $T=500$. We refer to Remark \ref{rem:excursions} and Lemma \ref{lem:excursions} in Appendix \ref{app:main_lemmas} for the definition of an excursion in Hawkes processes.}
    \label{tab:simu_data}
\end{table}

% \begin{table}[hbt!]
%     \centering
% \begin{tabular}{c|c|c|c|c}
% \toprule
%   Scenario & Sampler & ReLU &  Softplus &  Sigmoid  \\
%  \midrule
%  \multirow{2}{*}{Excitation only} & MH & 37 &   105 &  63  \\
%  & HMC & 482 &  2,064     & 1,963  \\
% \midrule
%  \multirow{2}{*}{Mixed effect} & MH & 39 &   101 & 100 \\
%  & HMC &   931 &     2,073  &   1441 \\
%  \midrule
%   \multirow{2}{*}{Inhibition only} & MH &  33.40  &  101   &  128  \\
%  & HMC &  1,386 &    2,115 & 2298 \\
% \bottomrule
% \end{tabular}
%     \caption{Computational times (in seconds) of the MH and HMC samplers in the three models and scenarios of Simulation 1.}
%     \label{tab:comptime_mcmc}
% \end{table}

\FloatBarrier

\subsection{Simulation 2: Parametric variational posterior and posterior distribution in the univariate sigmoid  model.}

%\textcolor{magenta}{$T=500$  non ?}

In this simulation,  we consider the same univariate scenarios as Simulation 1, but only for the sigmoid Hawkes model and compare the variational and true posterior distributions. Here, the dimensionality $D_0$ of the true function $h_0$ is assumed to be known, therefore, the samplers are non-adaptive. Specifically, we compare the performance of the previous MH sampler, % (4 chains with 50 000 iterations),
the Gibbs sampler (introduced in Remark \ref{rem:gibbs} and described in Algorithm \ref{alg:gibbs} in Appendix ~\ref{app:gibbs_sampler}), % with 3000 iterations
and our mean-field variational algorithm in a fixed model (Algorithm \ref{alg:cavi}) - here, we fix the dimensionality of $h_{11}$ to $J=2^{D_0}=4$.
We run 4 chains for 40 000 iterations for the MH sampler, 3000 iterations of the Gibbs sampler, and use our early-stopping procedure for the mean-field variational algorithm. 

\begin{figure}[hbt!]
\setlength{\tempwidth}{.3\linewidth}\centering
\settoheight{\tempheight}{\includegraphics[width=\tempwidth, trim=0.cm 0.cm 0cm  0.65cm,clip]{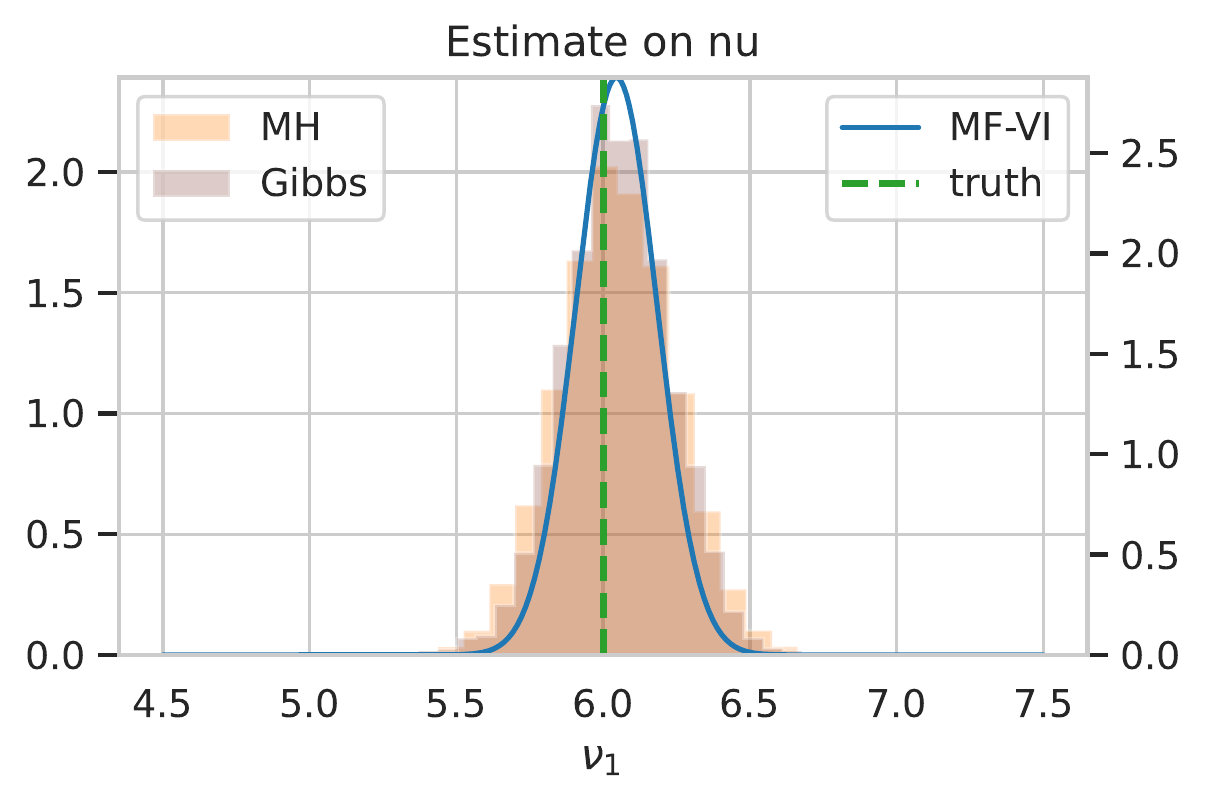}}%
\hspace{-5mm}
\fbox{\begin{minipage}{\dimexpr 20mm} \begin{center} \itshape \large  \textbf{Sigmoid} \end{center} \end{minipage}}
\hspace{-5mm}
\columnname{\emph{Excitation} only}\hfil
\columnname{Mixed effect}\hfil
\columnname{\emph{Inhibition} only}\\
\rowname{Background}
    \includegraphics[width=\tempwidth, trim=0.cm 0.cm 0cm  0.6cm,clip]{figs/sigmoid_D1_B4_exc_nu.pdf}\hfil
    \includegraphics[width=\tempwidth, trim=0.cm 0.cm 0cm  0.6cm,clip]{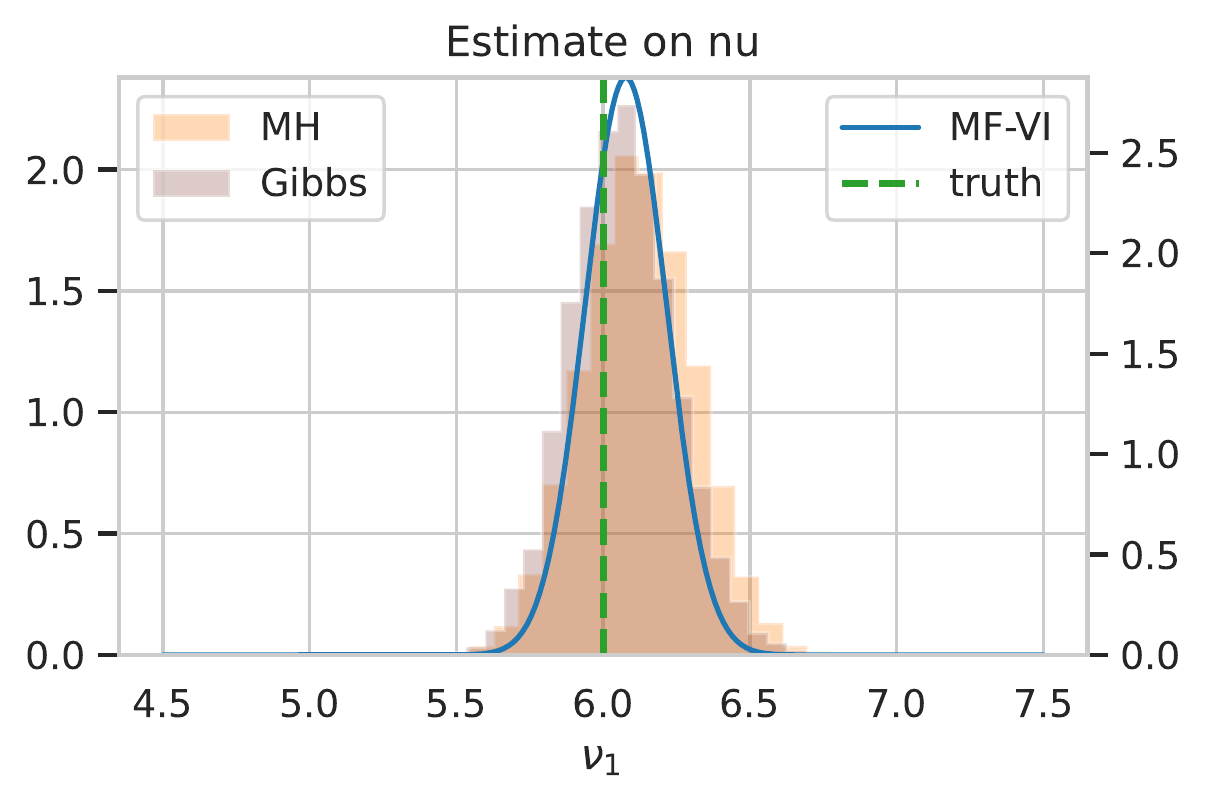}\hfil
    \includegraphics[width=\tempwidth, trim=0.cm 0.cm 0cm  0.6cm,clip]{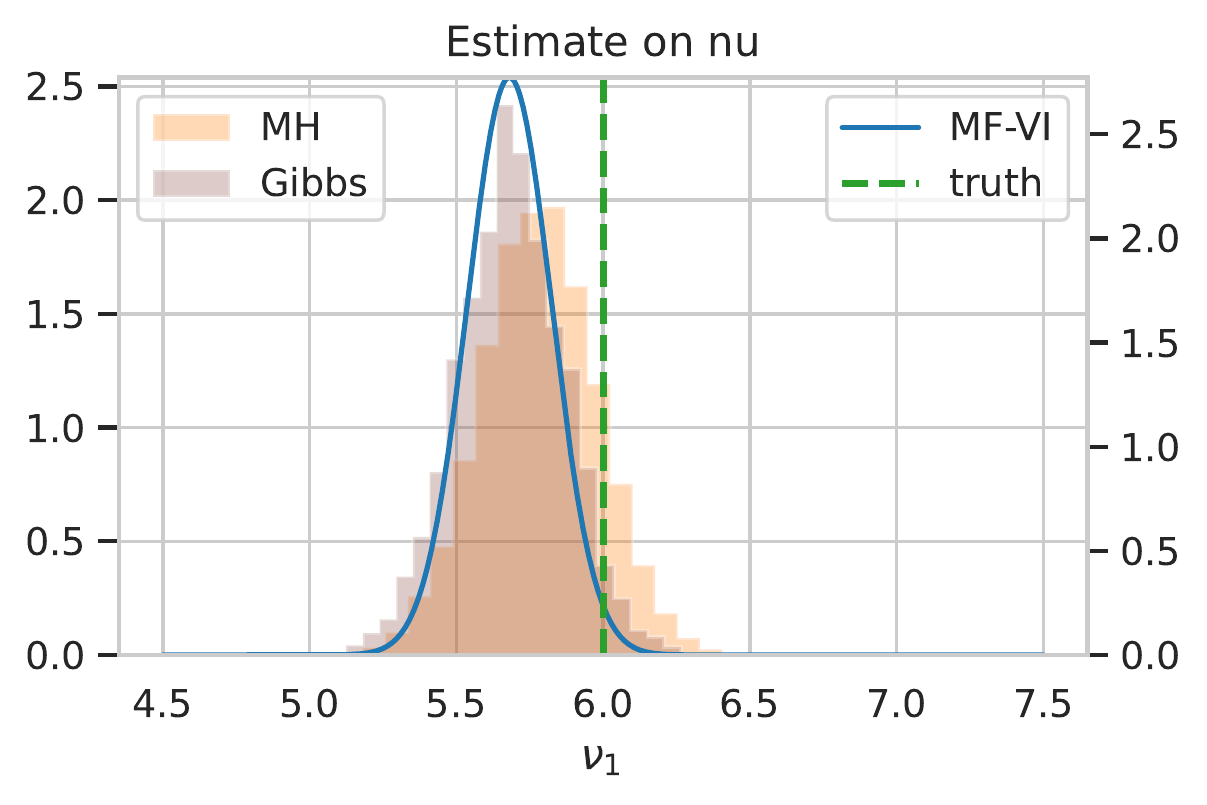}\\
\rowname{Interaction}
    \includegraphics[width=\tempwidth, trim=0.cm 0.cm 0cm  0.6cm,clip]{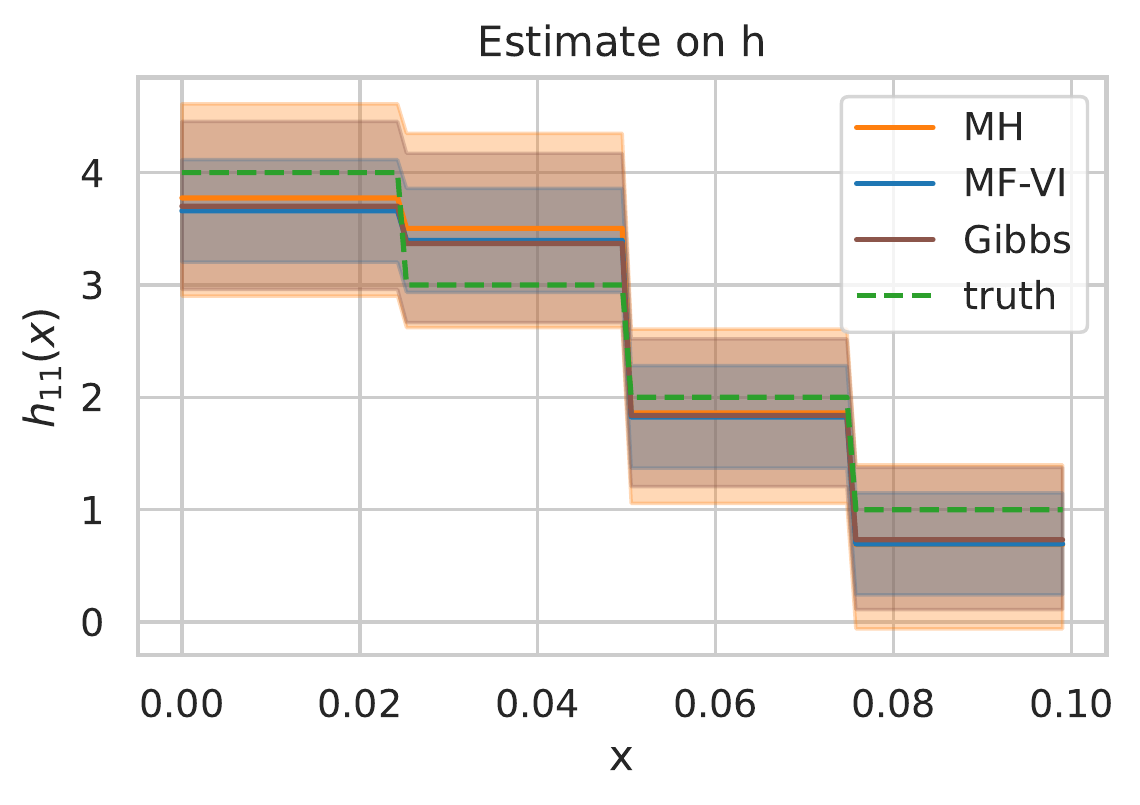}\hfil
    \includegraphics[width=\tempwidth, trim=0.cm 0.cm 0cm  0.6cm,clip]{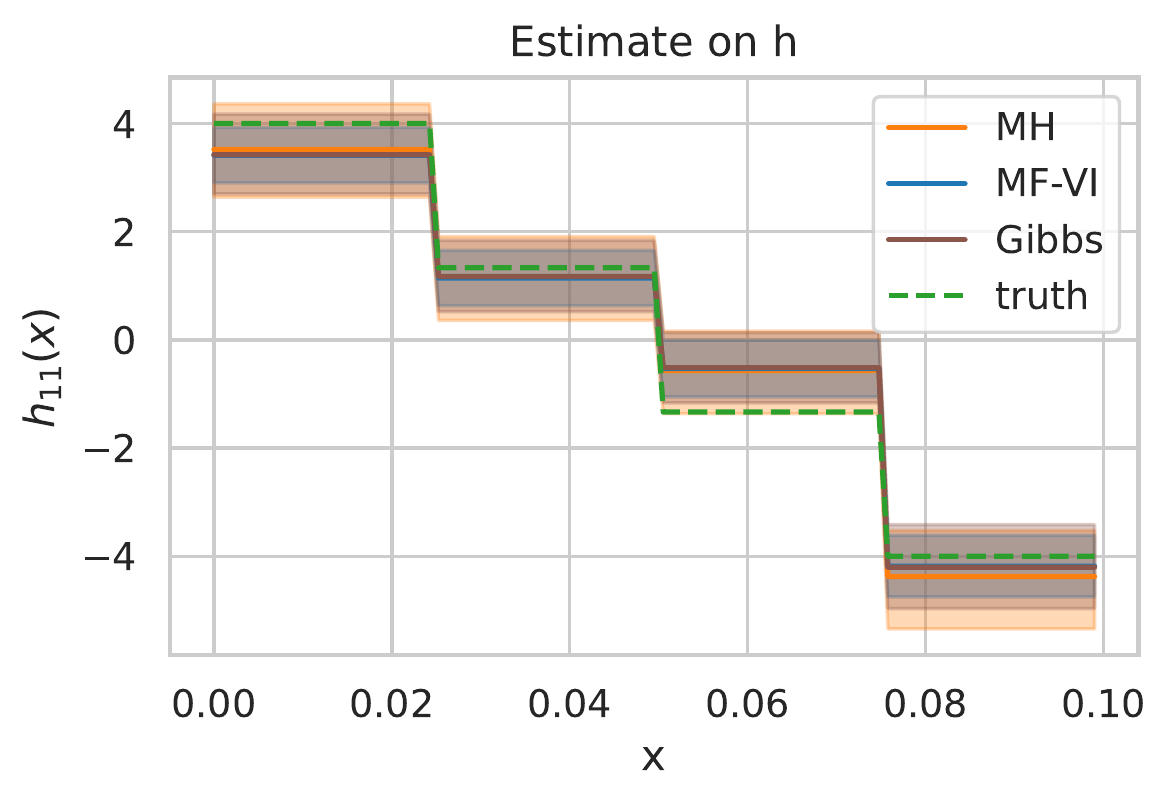}\hfil
    \includegraphics[width=\tempwidth, trim=0.cm 0.cm 0cm  0.6cm,clip]{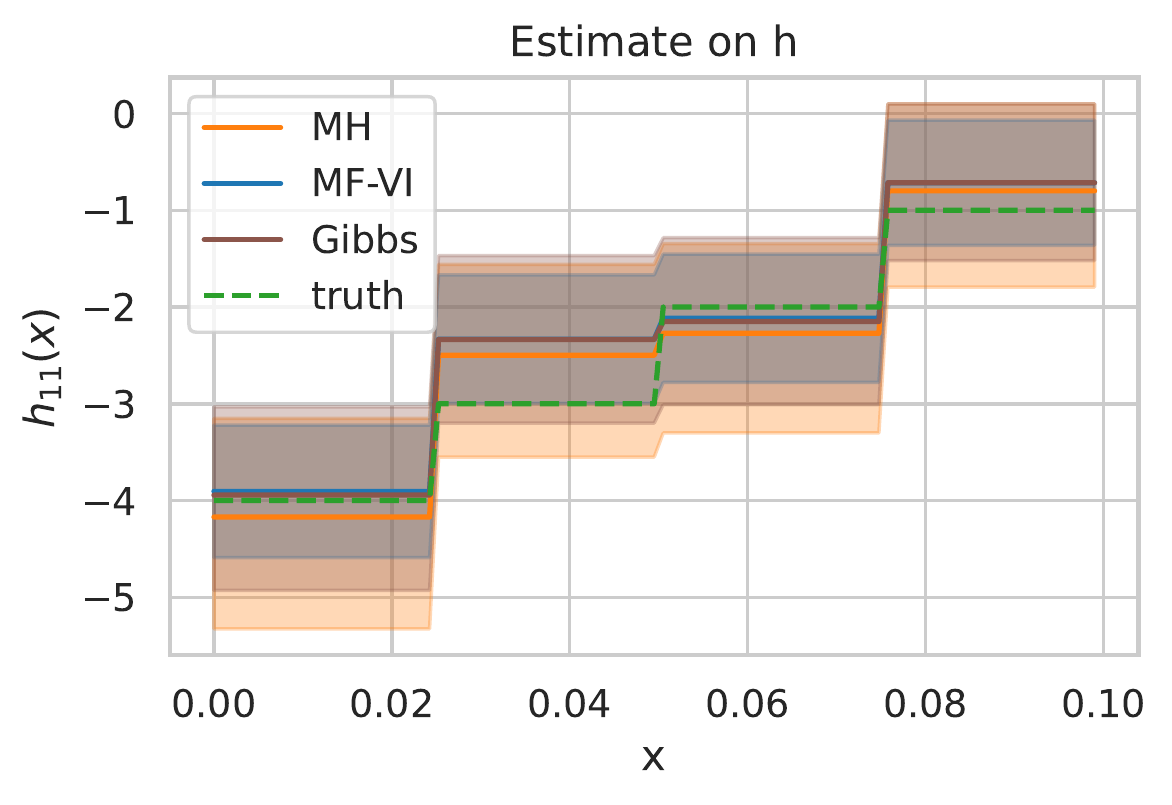}
\caption{Posterior and variational posterior distributions on $f = (\nu_1, h_{11})$ in the univariate sigmoid model of Simulation 2, evaluated by the MH sampler, the mean-field variational  (MF-VI) algorithm in a fixed model (Algorithm \ref{alg:cavi}) and the Gibbs sampler (Algorithm \ref{alg:gibbs}).  The three columns correspond to the \emph{Excitation only} (left), \emph{Mixed effect} (center), and \emph{Inhibition only} (right) scenarios. The true parameter $f_0$ is plotted in dotted green line. The first row contains the marginal distributions (VB, MH and Gibbs) on the background rate $\nu_1$, and the second row represents the posterior means (solid lines) and  95\% credible sets (colored areas) on the (self) interaction function $h_{11}$. We note that the variational posterior is close to the Gibbs posterior distribution, nonetheless, has smaller credible bands.}% \textcolor{red}{On ne voit pas trop pour MF-VI}  }
\label{fig:sigmoid_D4}
\end{figure}

In Figure~\ref{fig:sigmoid_D4}, we can compare the variational posterior on  $f = (\nu_1, h_{11})$ to the posterior distributions, computed either with the Gibbs or MH samplers, in the three estimation scenarios. We note that variational posterior mean is always close to the posterior mean, in particular when computed with the Gibbs sampler. Nonetheless, its credible sets are generally smaller, %than the ones of the posterior distribution, 
which is a common empirical observation of mean-field variational approximations. %, and one potential solution is to post-process the variational credible bands to obtain better coverage of the truth. 

\begin{figure}[hbt!]
    \centering
     \begin{subfigure}[b]{0.49\textwidth}
    \includegraphics[width=\textwidth, trim=0.cm 0.cm 0cm  0.cm,clip]{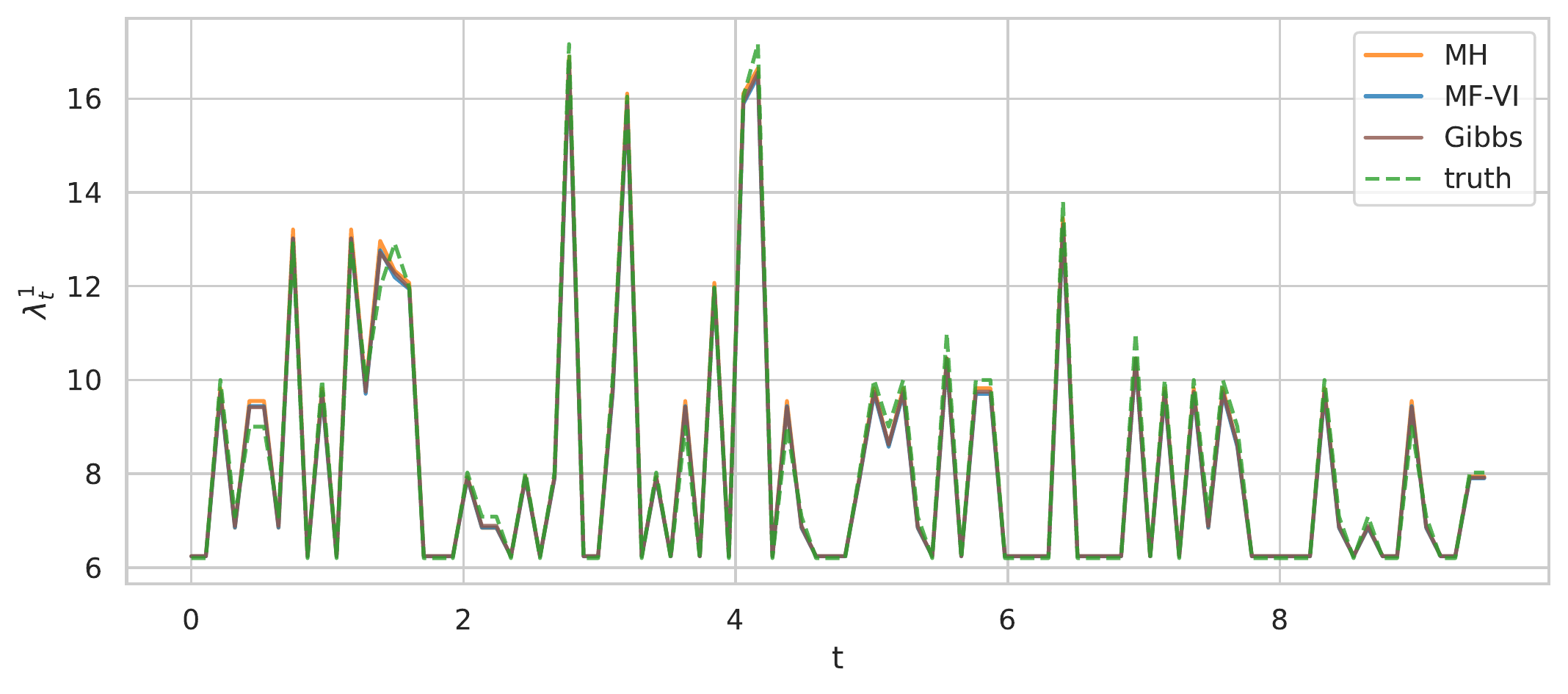}
    \caption{\emph{Excitation only} }
    \end{subfigure}
    \hfill
         \begin{subfigure}[b]{0.49\textwidth}
    \includegraphics[width=\textwidth, trim=0.cm 0.cm 0cm  0.cm,clip]{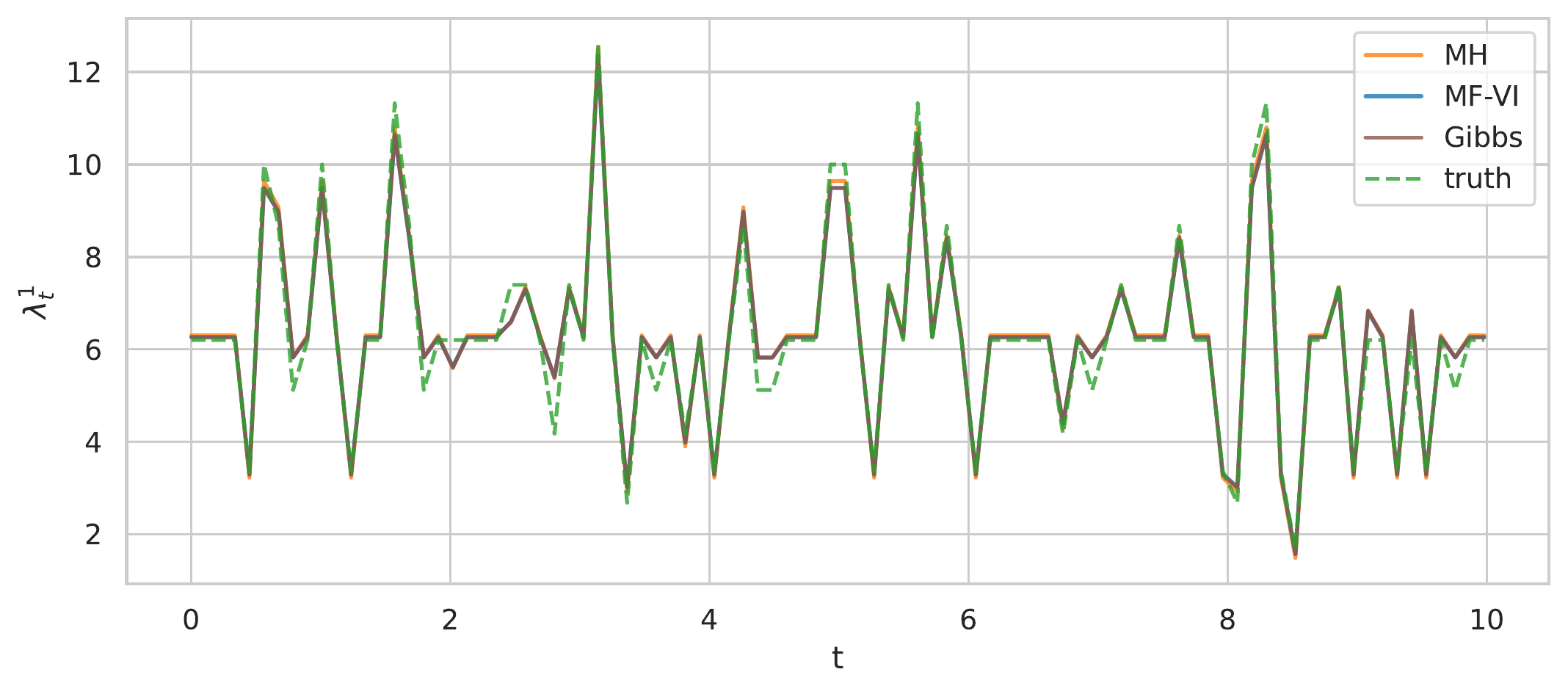}
    \caption{\emph{Mixed effect}}
    \end{subfigure}
    \hfill
         \begin{subfigure}[b]{0.49\textwidth}
    \includegraphics[width=\textwidth, trim=0.cm 0.cm 0cm  0.cm,clip]{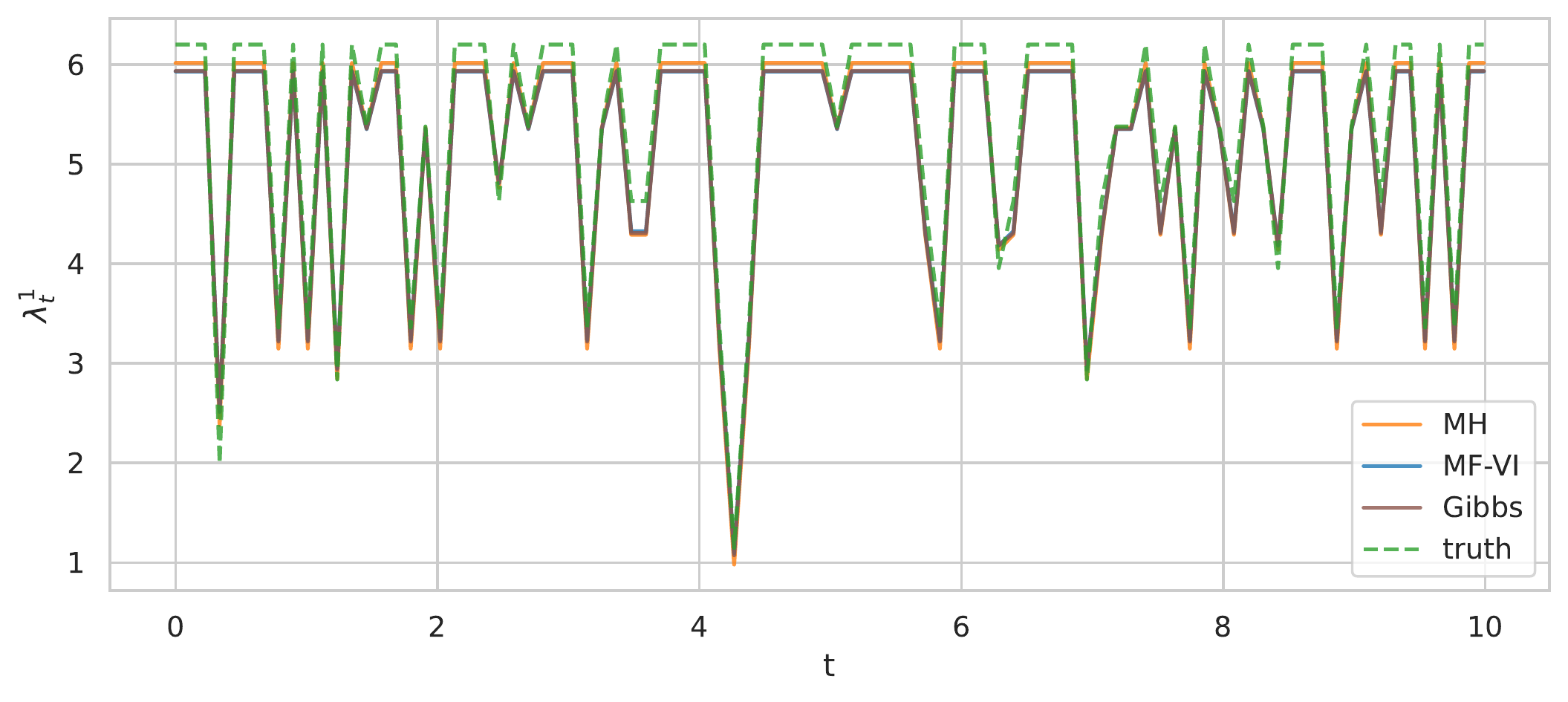}
    \caption{\emph{Inhibition only} }
    \end{subfigure}
    \hfill
\caption{Intensity function on a sub-window of the observation window estimated via the variational posterior mean (blue) or via the posterior mean, computed with the MH sampler (orange) or the Gibbs sampler (purple), in each scenario of Simulation 2. The true intensity $\lambda_t^1(f_0)$ is plotted in dotted green line. We note that all estimates are close in this simulation.} %, while the MH, Gibbs, and mean-field variational (MF-VI) estimates are plotted in solid orange (resp., red and blue) lines.}
\label{fig:intensity_D4}
\end{figure}

Besides, the variational posterior seems to be similarly biased as the posterior distribution, as can be seen for the background rate $\nu_1$ in the \emph{Inhibition} scenario. One could therefore test if this bias decreases with more data observations, i.e., larger $T$; however, the Gibbs sampler has a large computational time (between 3 and 5 hours), %seems slightly more precise than the other two algorithms; in spite of the small number of Gibbs iterations, it 
which is about 6 (resp. 40) times longer than the MH sampler (resp. our mean-field algorithm), due to the expensive latent variable sampling scheme (see Table \ref{tab:comptime_sigmoid}).
%quick to run and the variational posterior mean is a reasonably good estimator. However, the credible regions of the variational posterior are sometimes too narrow, in particular for estimating the background parameter $\nu$. 
Finally, we also compare the estimated intensity function using the (variational) posterior means, on a sub-window of the observations in Figure \ref{fig:intensity_D4}. The latter plot shows that all three methods provide fairly equivalent estimates on the nonlinear intensity function.

From this simulation, we conclude that, in the univariate and parametric sigmoid Hawkes model, the mean-field variational algorithm in a fixed model provides a good approximation of the posterior distribution. Moreover, we note that  although the Gibbs sampler is slightly better than MH, it is much slower than the latter and therefore cannot be applied to multivariate Hawkes processes in practice. Therefore, in the bivariate simulation in the next section, we only compare to the posterior distribution computed with the MH sampler, which can still be computed within reasonable time for $K=2$ .   

% is a good algorithm that provides well calibrated credible regions within a quite short amount of time. Nonetheless, the MF-Vi algorithm is 10 times faster to run, therefore seems more able to scale to multidimensional and high-dimensional Hawkes processes (Simulation 4). Besides, it can also be used in a nonparametric setting using the adaptive versions (Algorithms \ref{alg:adapt_cavi} and \ref{alg:2step_adapt_cavi}), that we first test in the univariate setting in Simulation 3. 

%\textcolor{red}{ Ju: What is $T$ in Table 2, is it still 500?}

\begin{table}[hbt!]
    \centering
\begin{tabular}{c|c|c|c}
\toprule
  Scenario  &  MH & Gibbs & MF-VI\\
 \midrule
\emph{Excitation} only  &  2169 & 16 092 &   416 \\
Mixed effect &  2181  & 13 097 &   338 \\
\emph{Inhibition} only &  2222 &  9 318 &   400 \\
\bottomrule
\end{tabular}
    \caption{Computational times (in seconds) of the Gibbs sampler (Algorithm \ref{alg:gibbs}), our mean-field variational (MF-VI) algorithm (Algorithm \ref{alg:cavi}), and the Metropolis-Hastings (MH) sampler in each parametric univariate scenario of Simulation 2 with $T = 500$. We note that the Gibbs sampler is much slower than the MH sampler, which is also slower than the mean-field variational algorithm. }
    \label{tab:comptime_sigmoid}
\end{table}

\FloatBarrier

\subsection{Simulation 3: Fully-adaptive variational method in the univariate and bivariate sigmoid models.}

%\textcolor{magenta}{Le setting du Scenario 3 semble similaire a celui du scenario 2, mais du coup, on ne comprend pas pourquoi les temps de calculs pour MH different. Est-ce parce que $T$ est plus grand ?}
%\ds{Oui les T et les parametres etaient differents, il faut que je corrige}

\begin{table}[hbt!]
    \centering
\begin{tabular}{c|c|c|c|c}
\toprule
 \# dimensions & Scenario & T & FA-MF-VI &  MH  \\
 \midrule
 \multirow{2}{*}{$K=1$}  & Excitation & 2000 & 32  &   417  \\
 & Inhibition & 3000 & 33 & 445    \\
\midrule
 \multirow{2}{*}{$K=2$} &  Excitation  & 2000 & 189 & 2605 \\
 & Inhibition &  3000 &  197 & 2791 \\
\bottomrule
\end{tabular}
    \caption{Computing times (in seconds) of our fully-adaptive mean-field variational method  (FA-MF-VI) (Algorithm \ref{alg:adapt_cavi}) and the Metropolis-Hastings (MH) sampler in the univariate and bivariate sigmoid models and the scenarios of Simulation 3.} 
    \label{tab:comptime_simu3}
\end{table}

\begin{figure}[hbt!]
    \centering
    %      \begin{subfigure}[b]{0.15\textwidth}
    % \includegraphics[width=\textwidth, trim=0.cm 0.cm 0cm  0.cm,clip]{figures/graph_2D_1.pdf}
    % \caption{$\delta_0$}
    % \label{fig:graph}
    % \end{subfigure}%
     \begin{subfigure}[b]{0.49\textwidth}
    \includegraphics[width=\textwidth, trim=0.cm 0.cm 0.cm  0.7cm,clip]{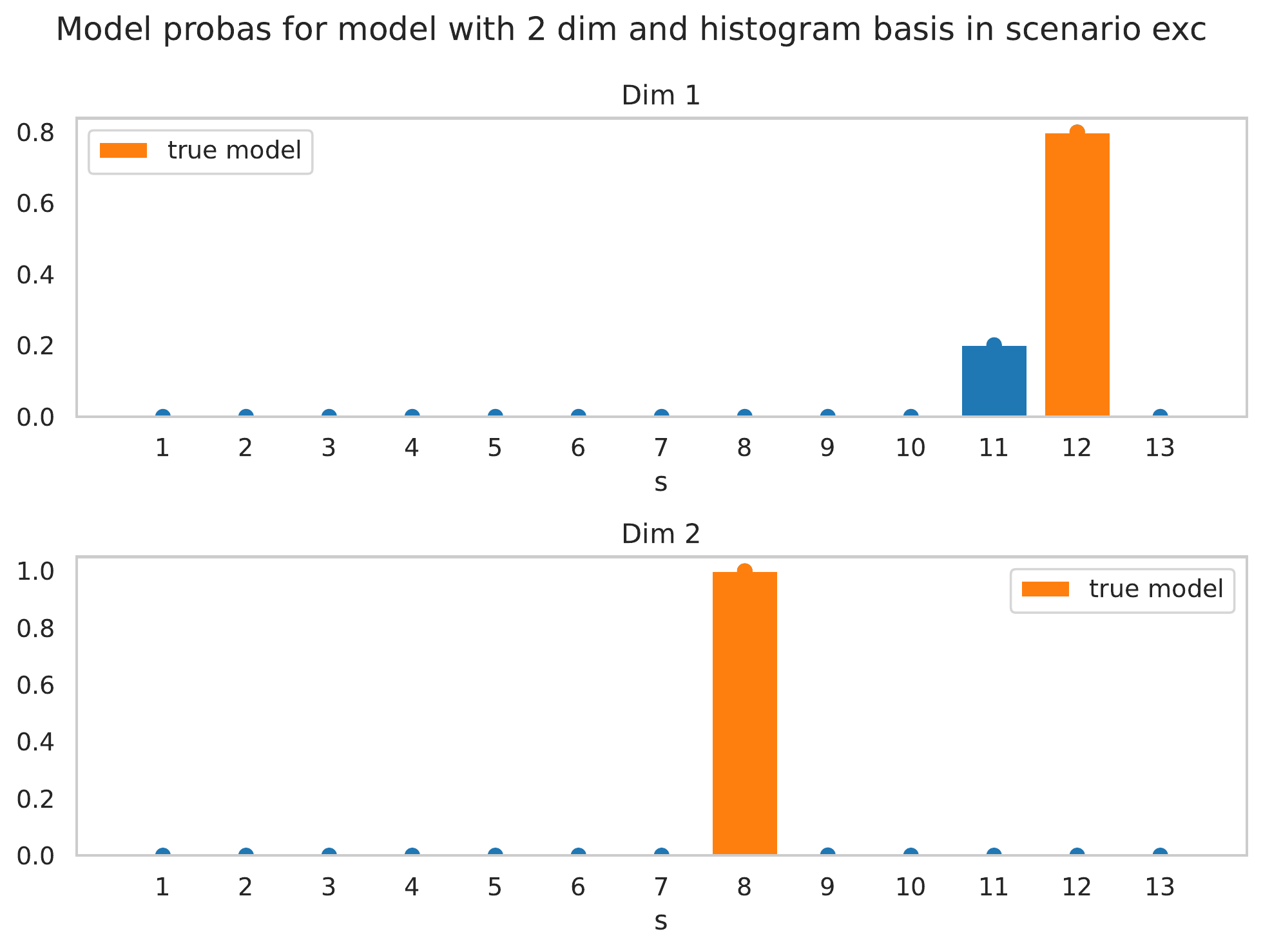}
    \caption{\emph{Excitation}}
    \label{exc}
    \end{subfigure}%
         \begin{subfigure}[b]{0.49\textwidth}
    \includegraphics[width=\textwidth, trim=0.cm 0.cm 0.cm  0.7cm,clip]{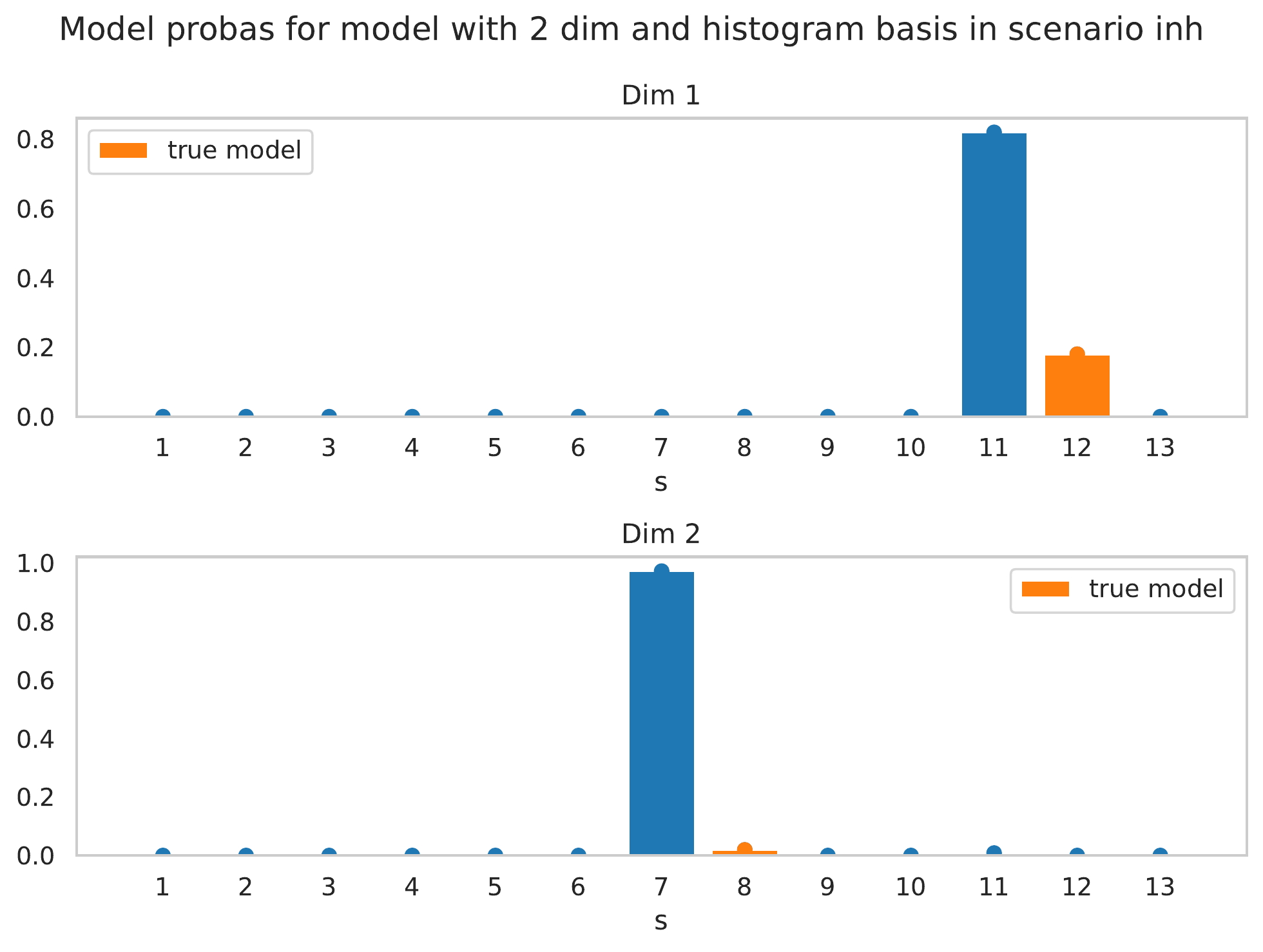}
    \caption{Inhibition}
    \label{inh}
    \end{subfigure}
\caption{ %True graph parameter $\delta_0$ (black=0, white=1) (\subref{fig:graph}) and 
Marginal probabilities on the graph and dimensionality parameter  $s_k = (\delta_{\cdot, k}, D_k)$ at each dimension, i.e., $(\hat \gamma_{s_k}^k)_{s_k \in \mathcal{S}_2}$ in the fully-adaptive averaged mean-field variational posterior, in the well-specified setting of Simulation 3 with $K=2$. % (i.e., with $h_0$ piecewise-constant, $\delta_0 = 1$ and $D_0=4$) of Simulation 3.
The  \emph{Excitation} scenario (\subref{exc}) corresponds to $h_0 \geq 0$, while in the  \emph{Inhibition} scenario (\subref{inh}) , $h_{11}^0, h_{22}^0  \leq 0$. The elements in $\mathcal{S}_2$ are indexed from 1 to 13 and the true model in this set is indicated in orange.}
\label{fig:graph_2D}
\end{figure}

In this simulation, we %consider that the dimensionality of the true interaction functions is unknown and 
test our fully-adaptive variational inference algorithm (Algorithm \ref{alg:adapt_cavi}), in the one-dimensional  ($K=1$) and two-dimensional ($K=2$) sigmoid models, and in two estimation settings:
%estimation scenarios, denoted \emph{Excitation} and \emph{Inhibition}, and for each scenario, we consider the two following settings for the true parameter $f_0$:
\begin{enumerate}
    \item \emph{Well-specified:}  $h_0 \in \mathcal{H}_{hist}^{D_0}$ (with $D_0 = 2$); % In this setting, we compare our variational posterior to the posterior distribution obtained with a non-adaptive MH sampler run with the true $s_0 = (\delta_0, D_0)$;
    \item \emph{Mis-specified:}  $h_0 \notin \mathcal{H}_{hist}^{D_0}$, and $h_{lk}^0$ is a continuous function, for all $(l,k) \in [K^2]$. % with non-zero $D_0 = 5$ first components in the Fourier basis.
     \end{enumerate}
Note that in the well-specified case, $m_0 := (\delta_0, 2^{D_0})$ is unknown for the variational method, nonetheless, we also compute the posterior distribution with the non-adaptive MH sampler using the true $m_0$.  %We set $T = 1500$ and
In the  bivariate model, we choose a true graph parameter $\delta_0$ with one zero entry (see Figure \ref{fig:graph_2D}a). %, i.e., three of the four interaction functions are non-null. 
We also consider an \emph{Excitation} scenario where all the true interaction functions $(h_{lk}^0)_{l,k}$ are non-negative and with $T = 2000$, and an \emph{Inhibition} scenario where the self-interaction functions $(h_{kk}^0)_{k=1,2}$ are non-positive with $T = 3000$. The latter setting aims at imitating the so-called self-inhibition phenomenon in neuronal spiking data, due to the refractory period of neurons \citep{bonnet2021maximum}.
%To estimate the dimensionality of the problem, we consider a family of nested histogram models with $2^{d}$ bins, $d=0,\dots,D$, or in other words a dyadic regular partition with maximum depth $D$.
In our adaptive variational algorithm, we set a maximum histogram depth $D_1 = 5$ for $K=1$, and  $D_1 = 4$ for $K=2$, so that the number of models per dimension is respectively 7 and 76. 
%In the univariate setting, we consider two estimation scenarios, denoted \emph{Excitation} and \emph{Inhibition}. In both scenarios, the true self-interaction function is continuous (see Figure \ref{fig:adaptive_VI_1D}) and we fix $D = 5$.
 %We will also compare to the parametric Metropolis-Hastings sampler, which is also available in the low-dimensional and well-specified settings. 

\begin{figure}[hbt!]
\setlength{\tempwidth}{.4\linewidth}\centering
\settoheight{\tempheight}{\includegraphics[width=\tempwidth, trim=0.cm 0.cm 0cm  1.cm,clip]{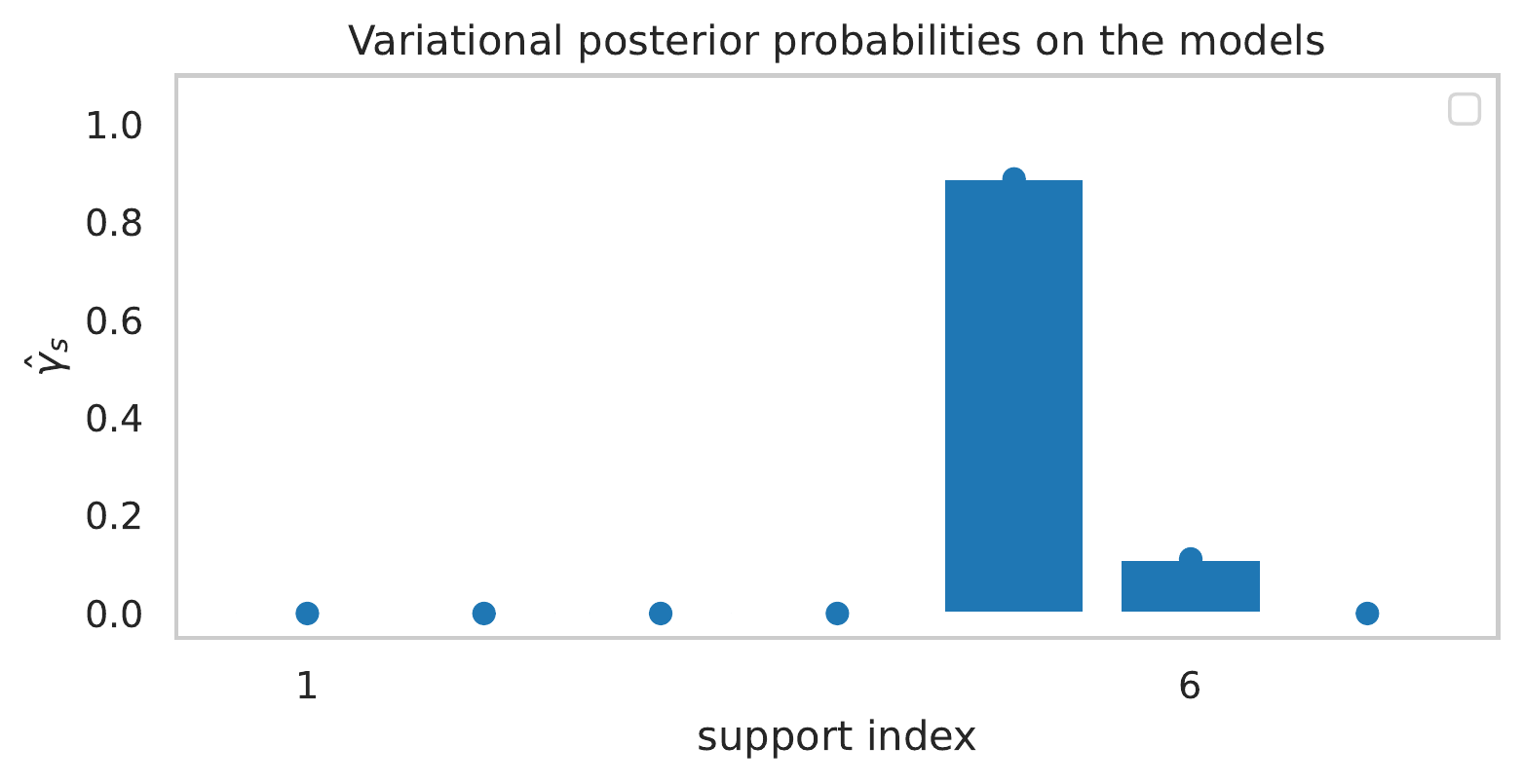}}%
\hspace{-5mm}
\fbox{\begin{minipage}{\dimexpr 15mm} \begin{center} \itshape \large  \textbf{$K = 1$} \end{center} \end{minipage}}
\hspace{-5mm}
\columnname{Excitation}\hfil
\columnname{Inhibition}\\
\rowname{Well-specified}
    \centering
    \includegraphics[width=\tempwidth, trim=0.cm 0.cm 0cm  0.7cm,clip]{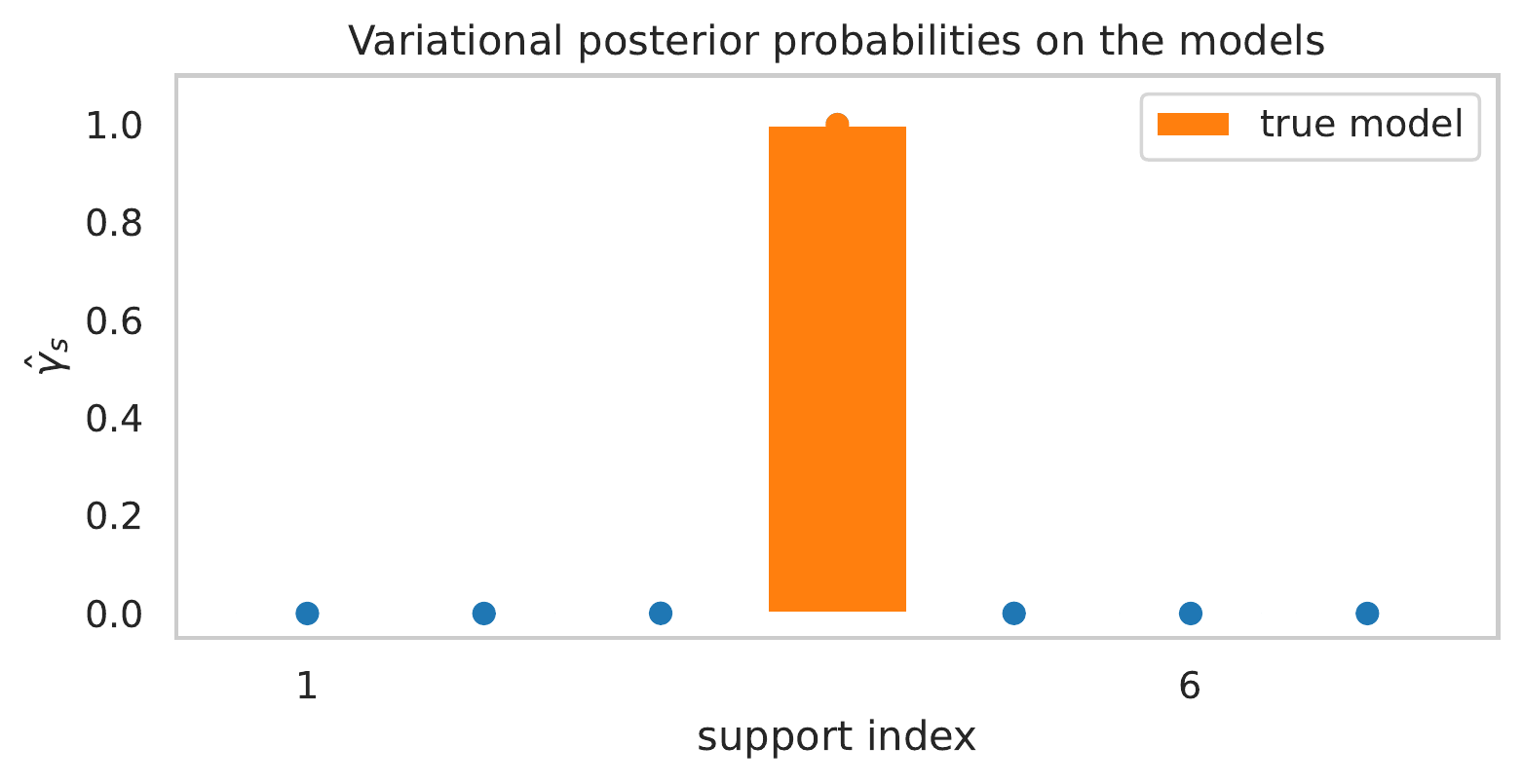}\hfil
     %\caption{\emph{Excitation} scenario }
    \includegraphics[width=\tempwidth, trim=0.cm 0.cm 0.cm  0.7cm,clip]{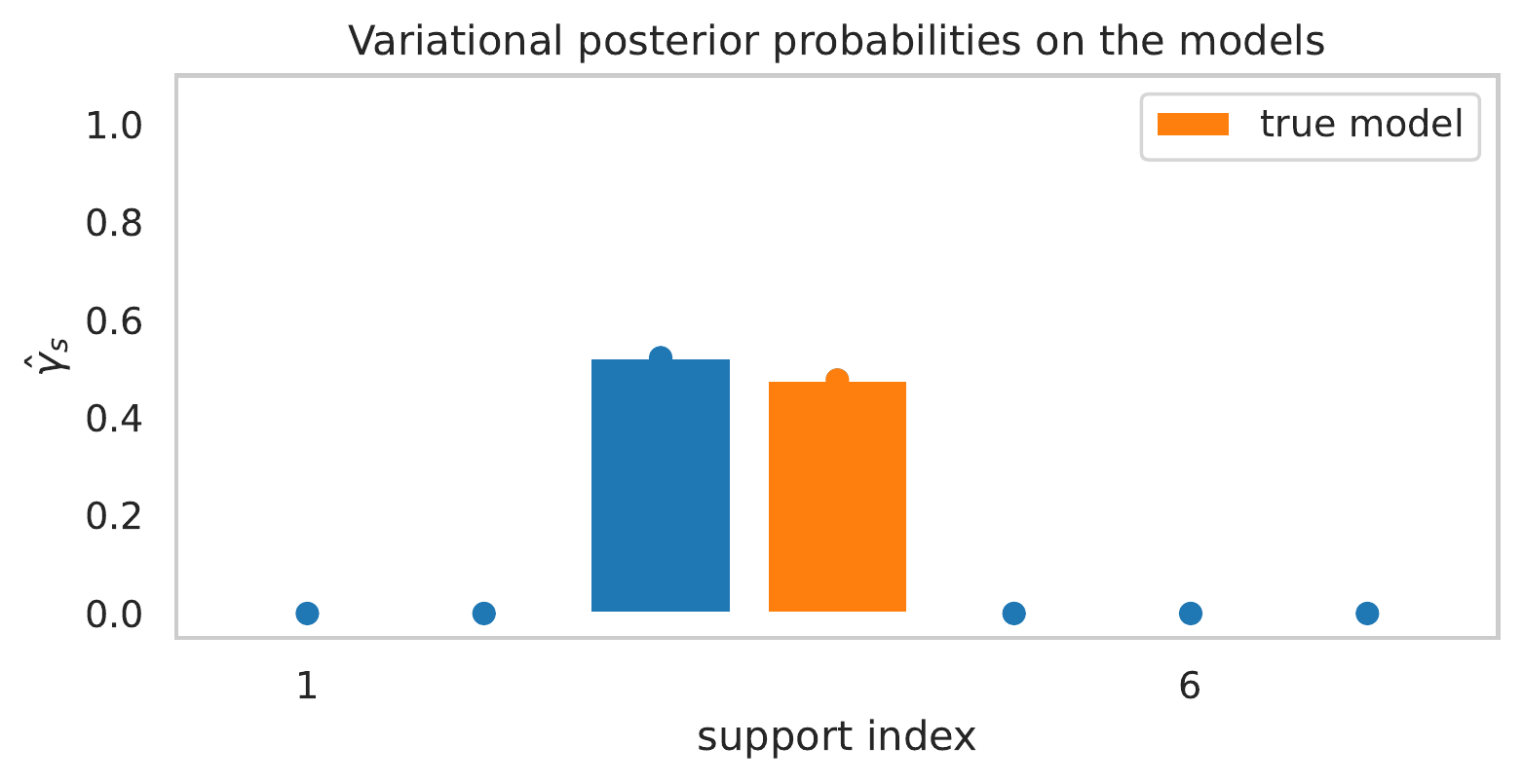}\\
     %\caption{\emph{Inhibition} scenario}
\rowname{Mis-specified}
    \centering
    \includegraphics[width=\tempwidth, trim=0.cm 0.cm 0cm  0.7cm,clip]{figs/sigmoid_exc_adaptive_1D_fourier_probas.pdf}\hfil
     %\caption{\emph{Excitation} scenario }
    \includegraphics[width=\tempwidth, trim=0.cm 0.cm 0.cm  0.7cm,clip]{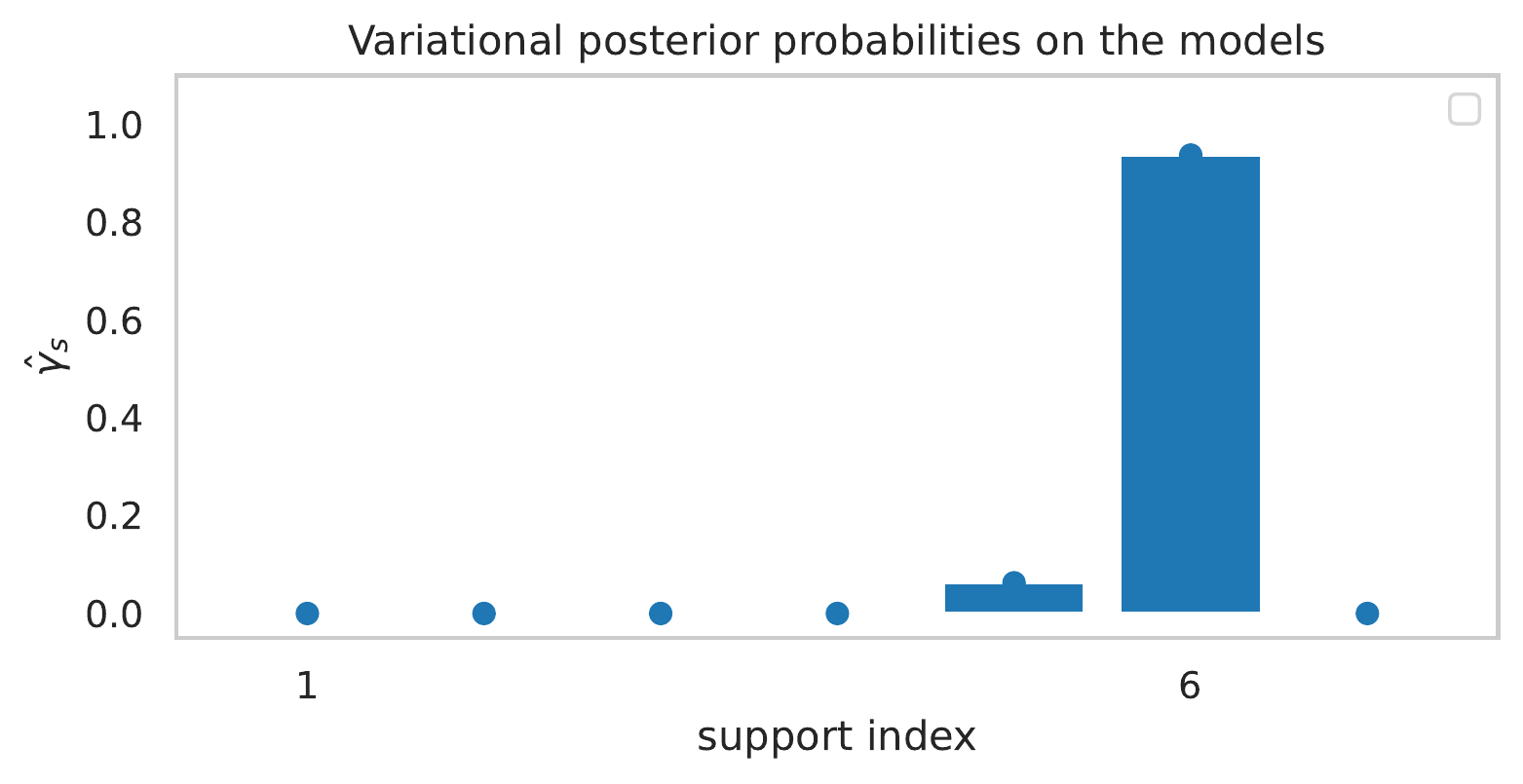}\\
     %\caption{\emph{Inhibition} scenario}
\caption{Model marginal probabilities $(\hat \gamma_{m})_{m}$ in the adaptive mean-field variational posterior, in the well-specified and mis-specified settings of Simulation 3 with $K=1$. % (i.e., with $h_0$ piecewise-constant, $\delta_0 = 1$ and $D_0=4$) of Simulation 3.
The left and  right panels correspond to the \emph{Excitation} (resp. \emph{Inhibition}) setting. The elements in $\mathcal{S}_1$ are indexed from 1 to 7, and correspond respectively to  $m = (\delta=0,2^{D}=1)$, and $m = (\delta=1,2^D)$ with $D=0, \dots, 5$.}
\label{fig:elbo_1D_fourier}
\end{figure}

In the well-specified setting, we first analyse the ability of Algorithm \ref{alg:adapt_cavi} to recover the true connectivity graph and dimensionality of $h_0$. In Figure \ref{fig:elbo_1D_fourier}, we plot the model marginal probabilities $(\hat \gamma_m)_{m}$ in our adaptive variational posterior and in the univariate setting. % and well-specified settings.
In the \emph{Excitation} scenario, the largest marginal probability $\hat \gamma_{\hat s}$ is on the true model, i.e., $\hat m = m_0 = (\delta_0 = 1,2^{D_0}=2)$, and  all the other marginal probabilities are negligible. Therefore, in this case, the model-averaging VB posterior  \eqref{eq:ms_var_post_avg} is essentially equivalent to the model-selection VB posterior \eqref{eq:ms_var_post}. In the \emph{Inhibition} scenario, the dimensionality $\hat D$ is not well inferred in the model selection variational posterior (maximising the ELBO), which is over-regularizing in this case, since $\hat m = (\hat \delta = \delta_0= 1, \hat D = 1)$. However, as seen in Figure \ref{fig:elbo_1D_fourier}, the ELBO for $D=1$ and for $D=2=D_0$ are very close, %$\hat \gamma_{s_0}$ being close to $\hat \gamma_{\hat s}$, 
therefore, the model-averaging variational posterior better captures the model since it% the marginal probability on $s_0$ is the second largest. %In the Self-inhibition scenario, $\hat s = (1,2)$, therefore corresponds the true dimensionality  $D_0=2$.
is essentially a mixture of two components, one corresponding to $\hat D = 1$, and the second one corresponding to the true model $D_0 + 2$.

Nonetheless, comparing the estimated nonlinear intensity based on the model-selection variational posterior mean and the posterior mean in Figure  \ref{fig:intensity_D1_adaptive} in Appendix \ref{app:details_exp}, we note that the model selection variational estimate is very close to the true intensity and the non-adaptive MH estimate, despite the error of dimensionality in the \emph{Inhibition} scenario.

\begin{figure}[hbt!]
\setlength{\tempwidth}{.32\linewidth}\centering
\settoheight{\tempheight}{\includegraphics[width=\tempwidth, trim=0.cm 0.cm 0cm  1.cm,clip]{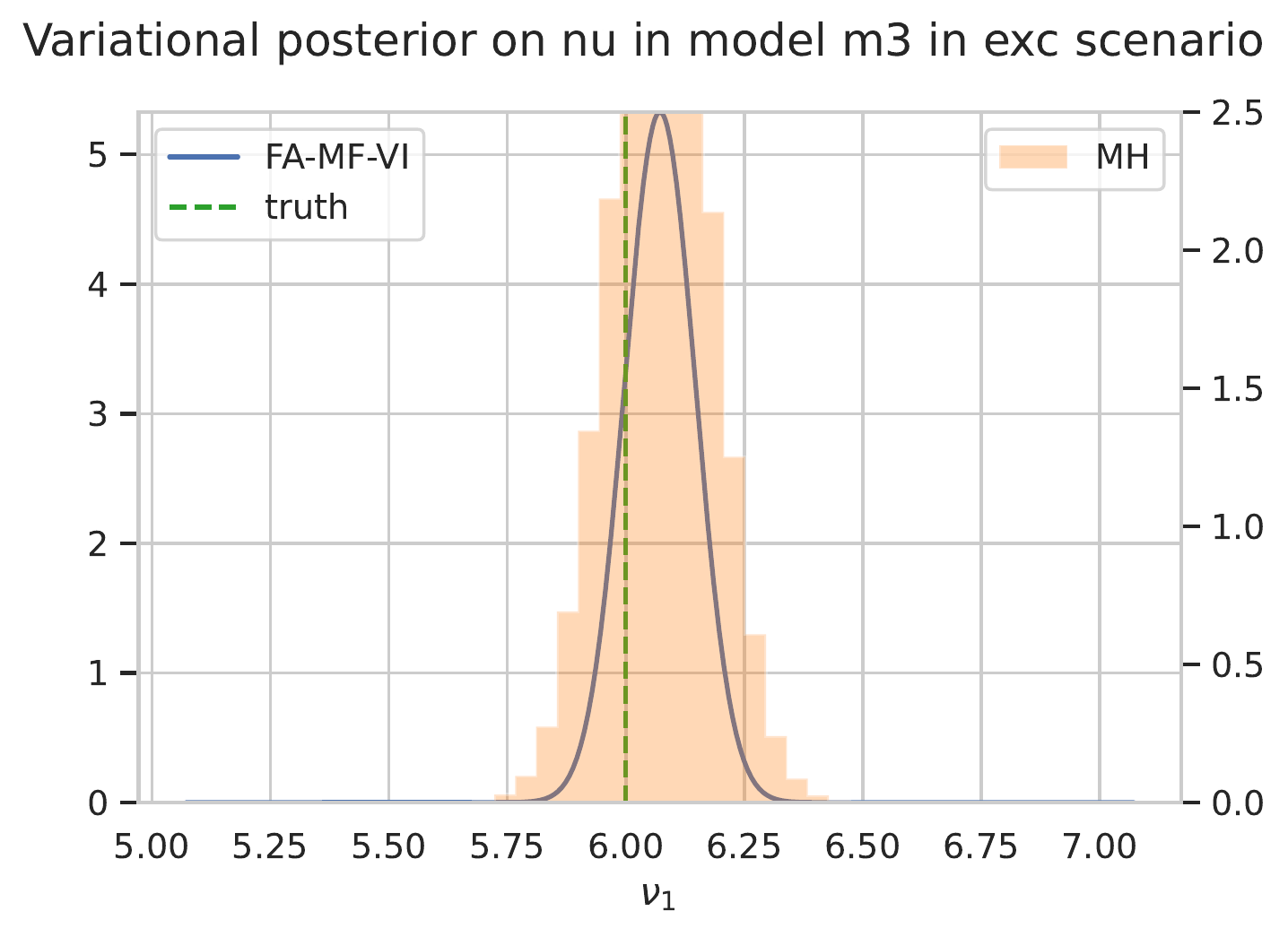}}%
\hspace{-15mm}
\fbox{\begin{minipage}{\dimexpr 15mm} \begin{center} \itshape \large  \textbf{$K = 1$} \end{center} \end{minipage}}
\hspace{-5mm}
\columnname{Well-specified-Exc}
\columnname{Well-specified-Inh}
\columnname{Mis-specified-Exc}\\
\hspace{-3mm}
\rowname{Background}
    \includegraphics[width=0.96\tempwidth, trim=0.cm 0.cm 0cm  1.cm,clip]{figs/adaptive_vi_1D_histogram_exc_nu_smod.pdf}
    \includegraphics[width=0.96\tempwidth, trim=0.cm 0.cm 0cm  1.cm,clip]{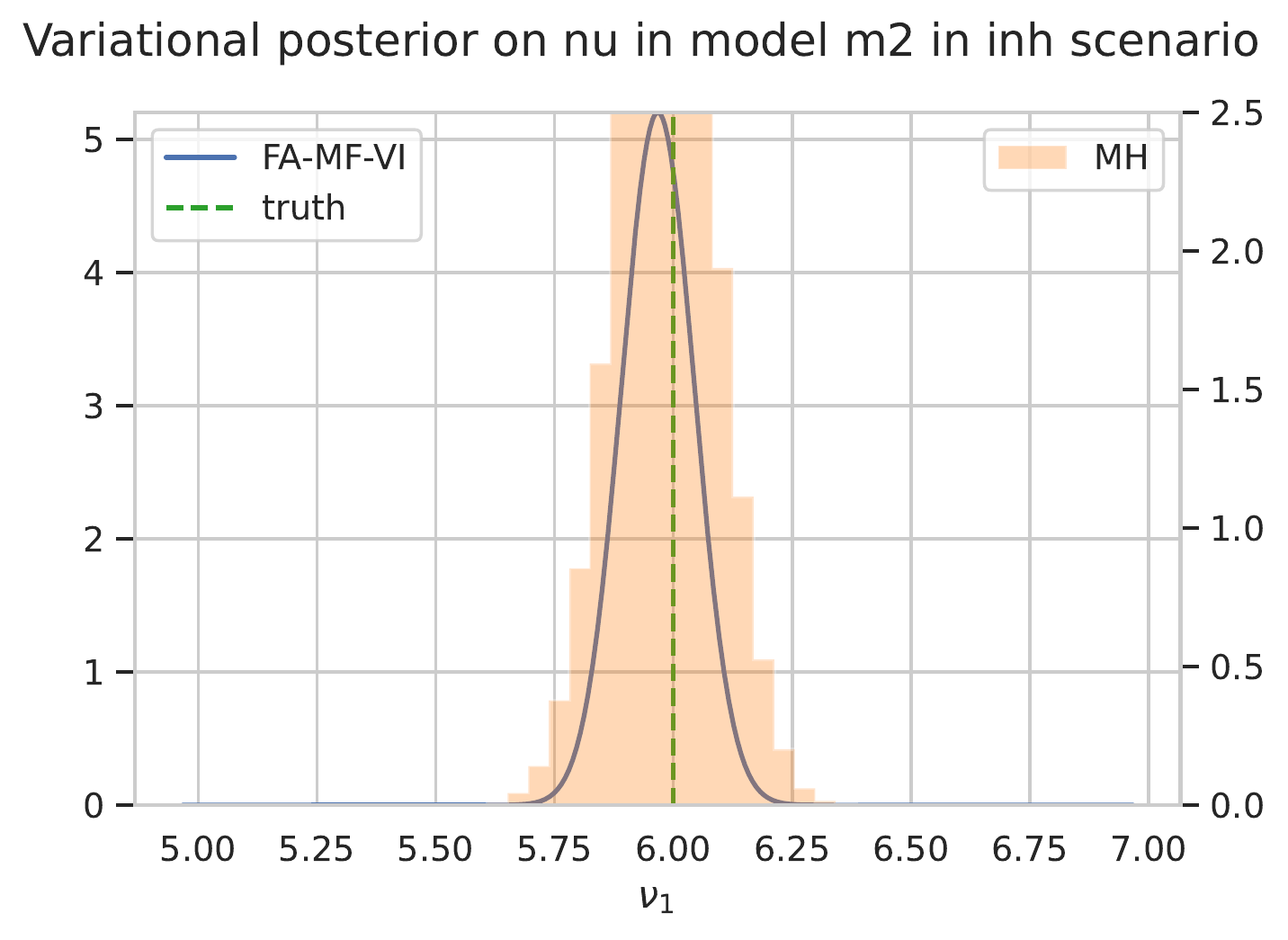}
        \includegraphics[width=0.98\tempwidth, trim=0.cm 0.cm 0cm  1.cm,clip]{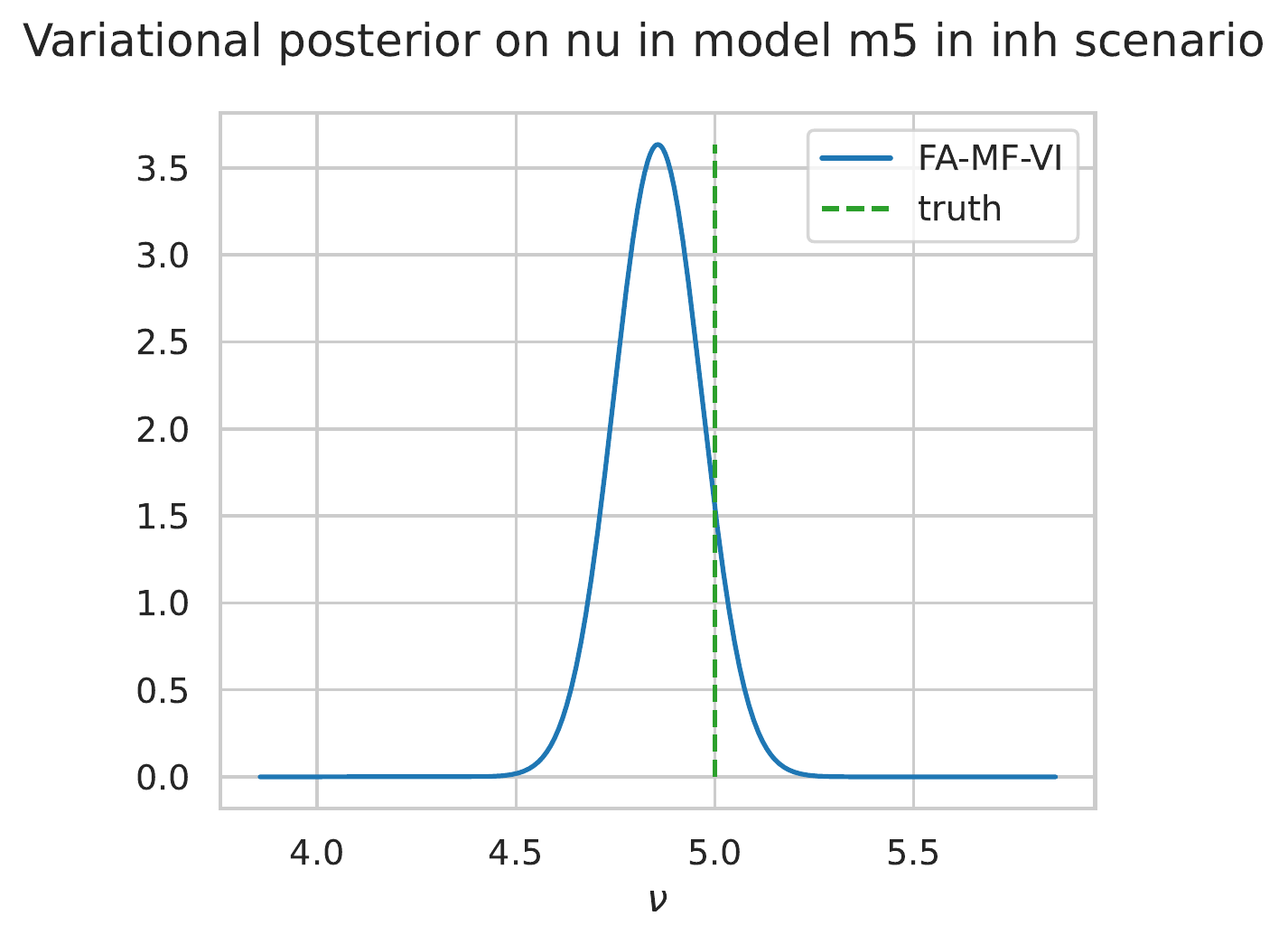}\\
\hspace{-7mm}
\rowname{Interaction}
    \includegraphics[width=0.95\tempwidth, trim=0.cm 0.cm 1.5cm  1.cm,clip]{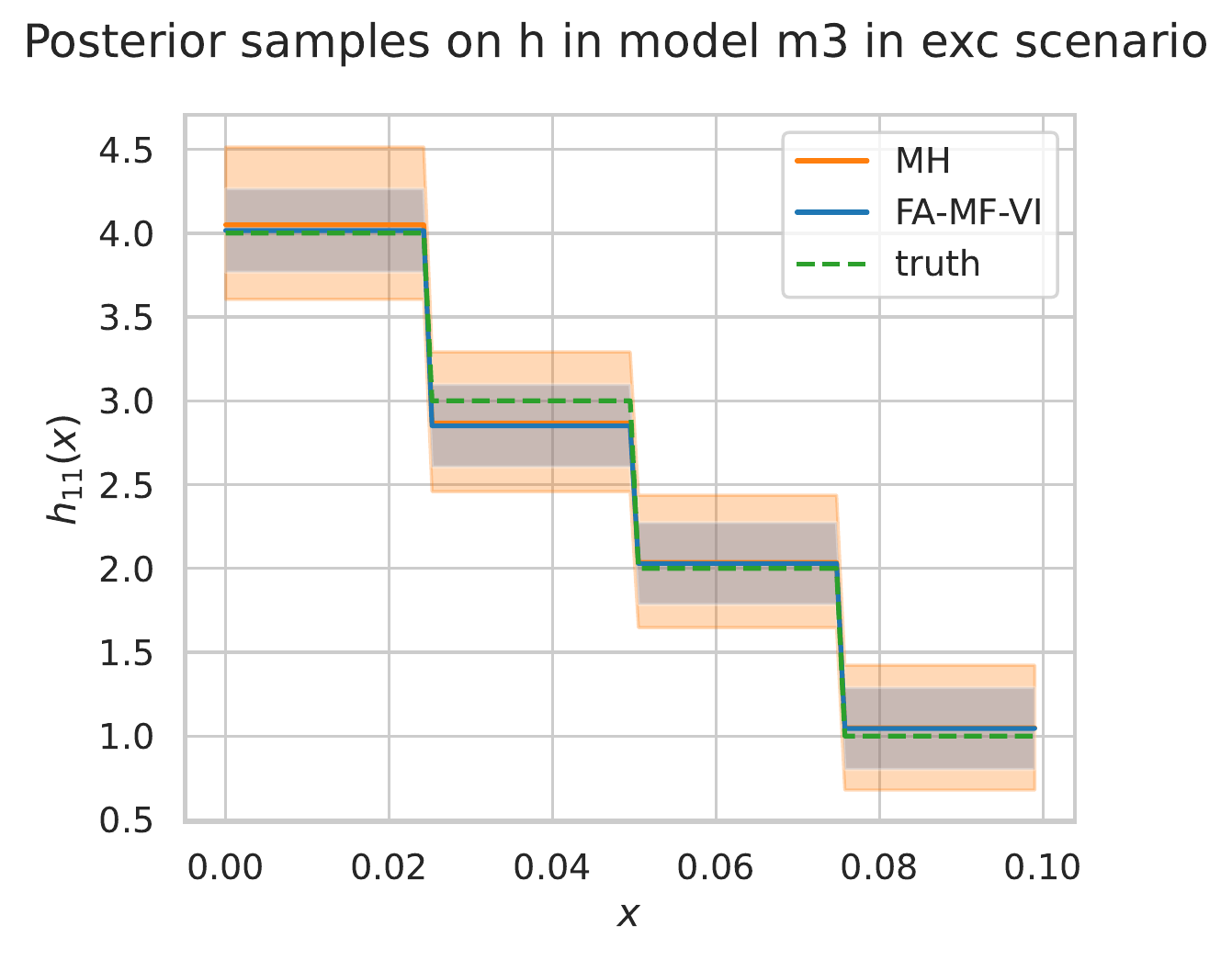}
    \includegraphics[width=0.95\tempwidth, trim=0.cm 0.cm 1.5cm  1.cm,clip]{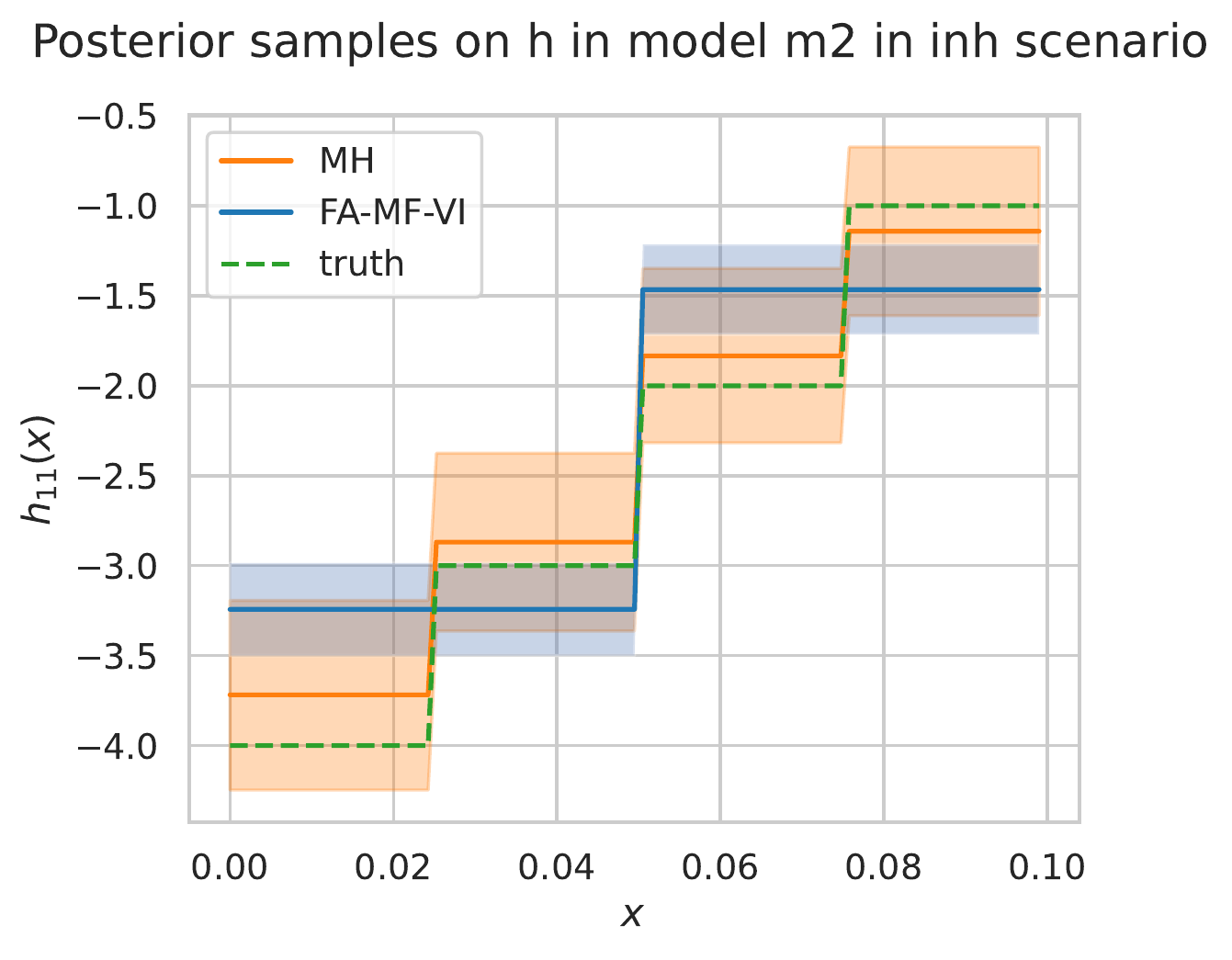}
        \includegraphics[width=0.95\tempwidth, trim=0.cm 0.cm 1.5cm  1.cm,clip]{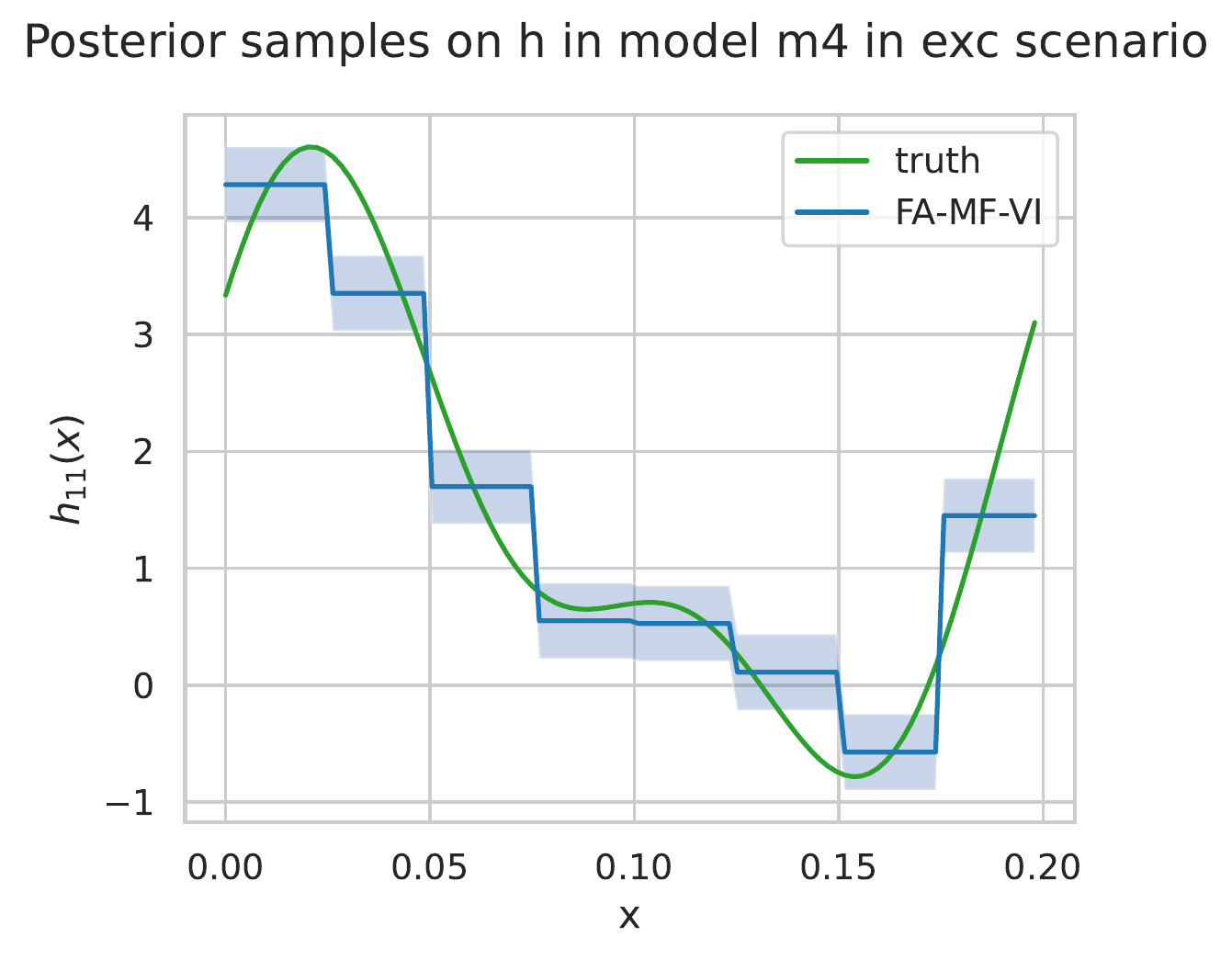}
\caption{Posterior and model-selection variational posterior distributions on $f = (\nu_1, h_{11})$ in the univariate sigmoid model and settings of Simulation 3, evaluated by the MH sampler and the fully-adaptive mean-field variational (FA-MF-VI) algorithm (Algorithm \ref{alg:adapt_cavi}).  The three columns correspond respectively to the two well-specified settings, i.e., the \emph{Excitation} (Well-specified-Exc)  and \emph{Inhibition} (Well-specified-Exc) scenarios, and one mis-specified setting (Mis-specified-Exc). The first row contains the marginal distribution on the background rate $\nu_1$, and the second row represents the (variational) posterior mean (solid line) and  95\% credible sets (colored areas) on the (self) interaction function $h_{11}$.  The true parameter $f_0$ is plotted in dotted green line.}
\label{fig:adaptive_VI_1D_histo}
\end{figure}

We then compare the model selection adaptive variational posterior distribution on the parameter with the true posterior distribution computed with the non-adaptive MH sampler in Figure  \ref{fig:adaptive_VI_1D_histo}.  We note that in the \emph{Excitation} scenario, the variational posterior mean is very close to the posterior mean, however, its 95\%  credible bands are significantly smaller. Note also that, in the \emph{Inhibition} scenario, in spite of the wrongly selected histogram depth, the estimated interaction function is still not too far from the truth.

In the mis-specified setting, all the marginal probabilities are negligible but one, in both the \emph{Excitation} and \emph{Inhibition} scenarios (see Figure \ref{fig:elbo_1D_fourier}), although there is no true $m_0$ in this case. In Figure \ref{fig:adaptive_VI_1D_fourier} in Appendix \ref{app:details_exp}, we note that the model selection adaptive variational posterior mean approximates quite well the true parameter. Moreover, its 95\% credible bands often cover the truth but are once again slightly too narrow.

\begin{figure}[hbt!]
\setlength{\tempwidth}{.49\linewidth}\centering
\settoheight{\tempheight}{\includegraphics[width=0.5\tempwidth, trim=0.cm 0.cm 0cm  1.cm,clip]{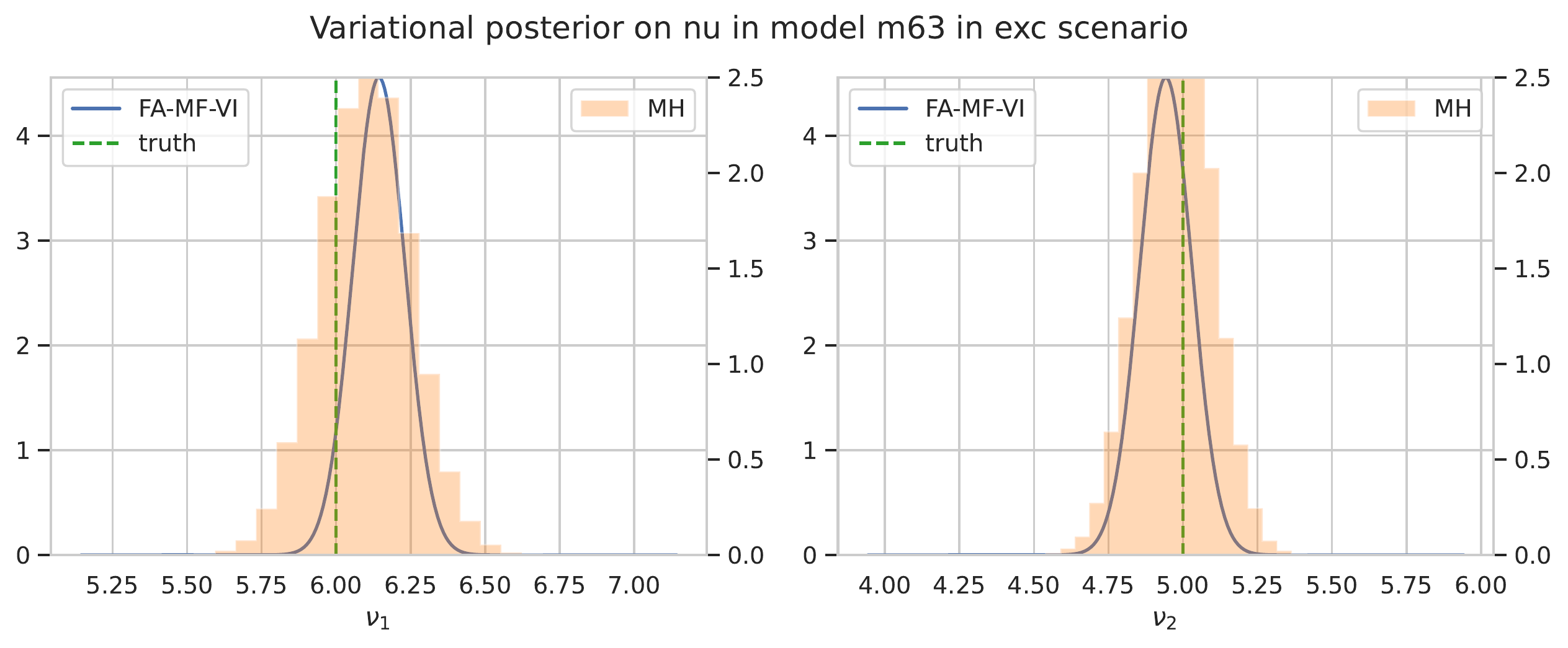}}%
\hspace{-20mm}
\fbox{\begin{minipage}{\dimexpr 13mm} \begin{center} \itshape \large  \textbf{$K = 2$} \end{center} \end{minipage}}
\columnname{Well-specified}
\columnname{Mis-specified}\\
\hspace{-10mm}
\rowname{\hspace{20mm} Background}
    \includegraphics[width=\tempwidth, trim=0.cm 0.cm 0cm  0.8cm,clip]{figs/adaptive_vi_2D_histogram_exc_nu_smod.pdf}
    \includegraphics[width=0.85\tempwidth, trim=0.cm 0.cm 0cm  0.8cm,clip]{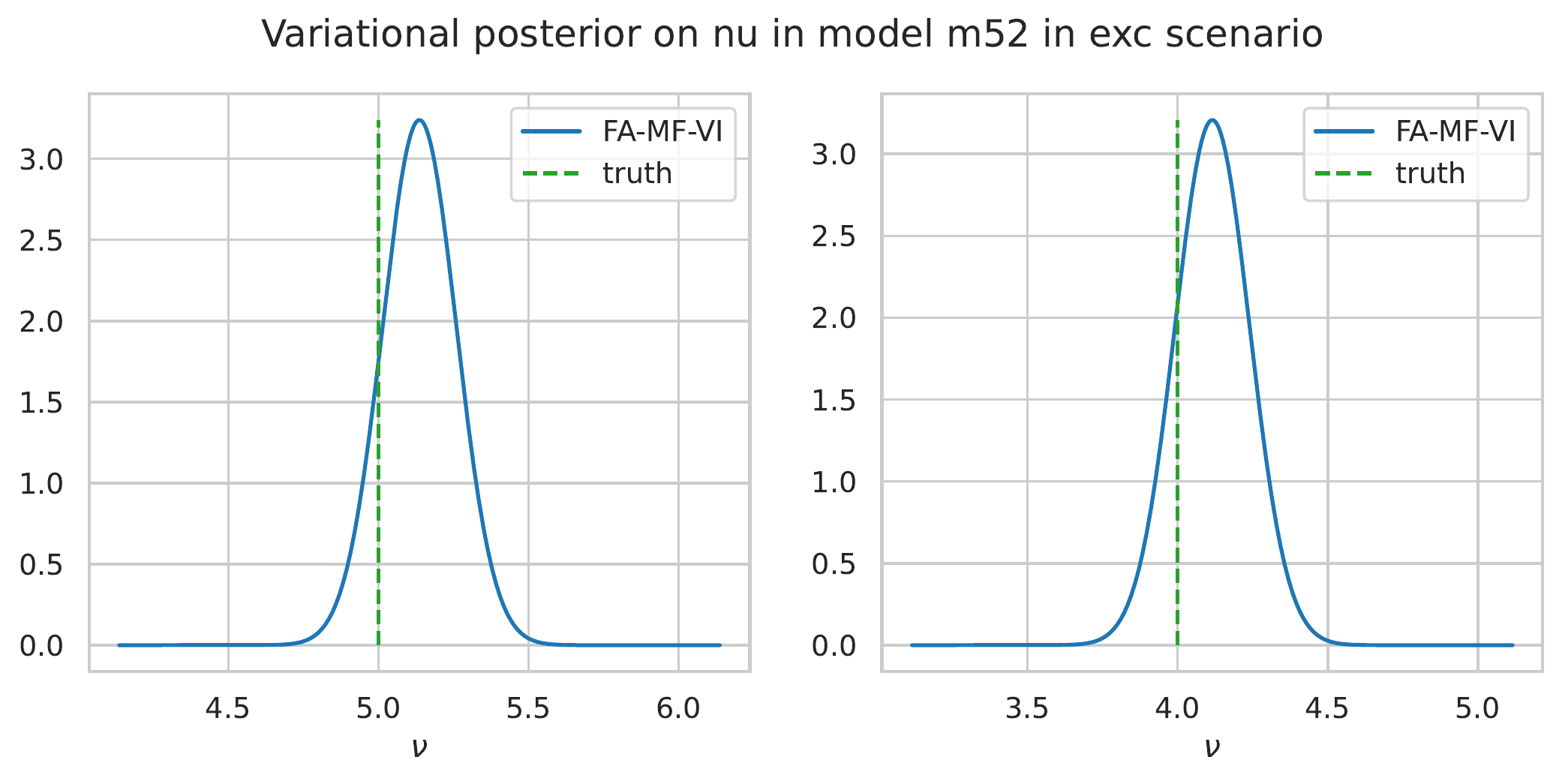}\\
\hspace{-10mm}
\rowname{\hspace{40mm} Interaction functions}
 \includegraphics[width=\tempwidth, trim=0.cm 0.cm 0cm  1.5cm,clip]{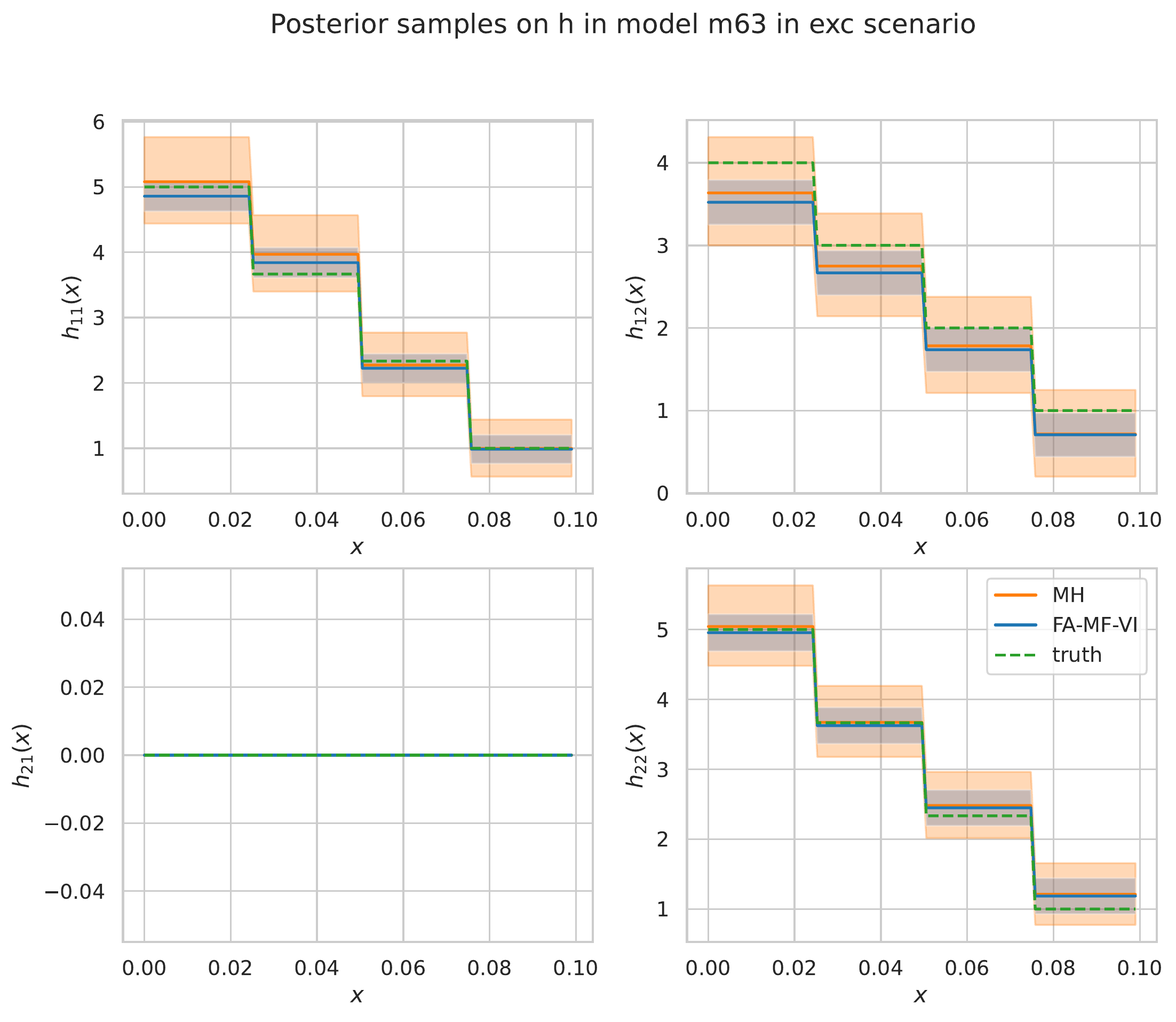}
    \includegraphics[width=\tempwidth, trim=0.cm 0.cm 0cm  1.5cm,clip]{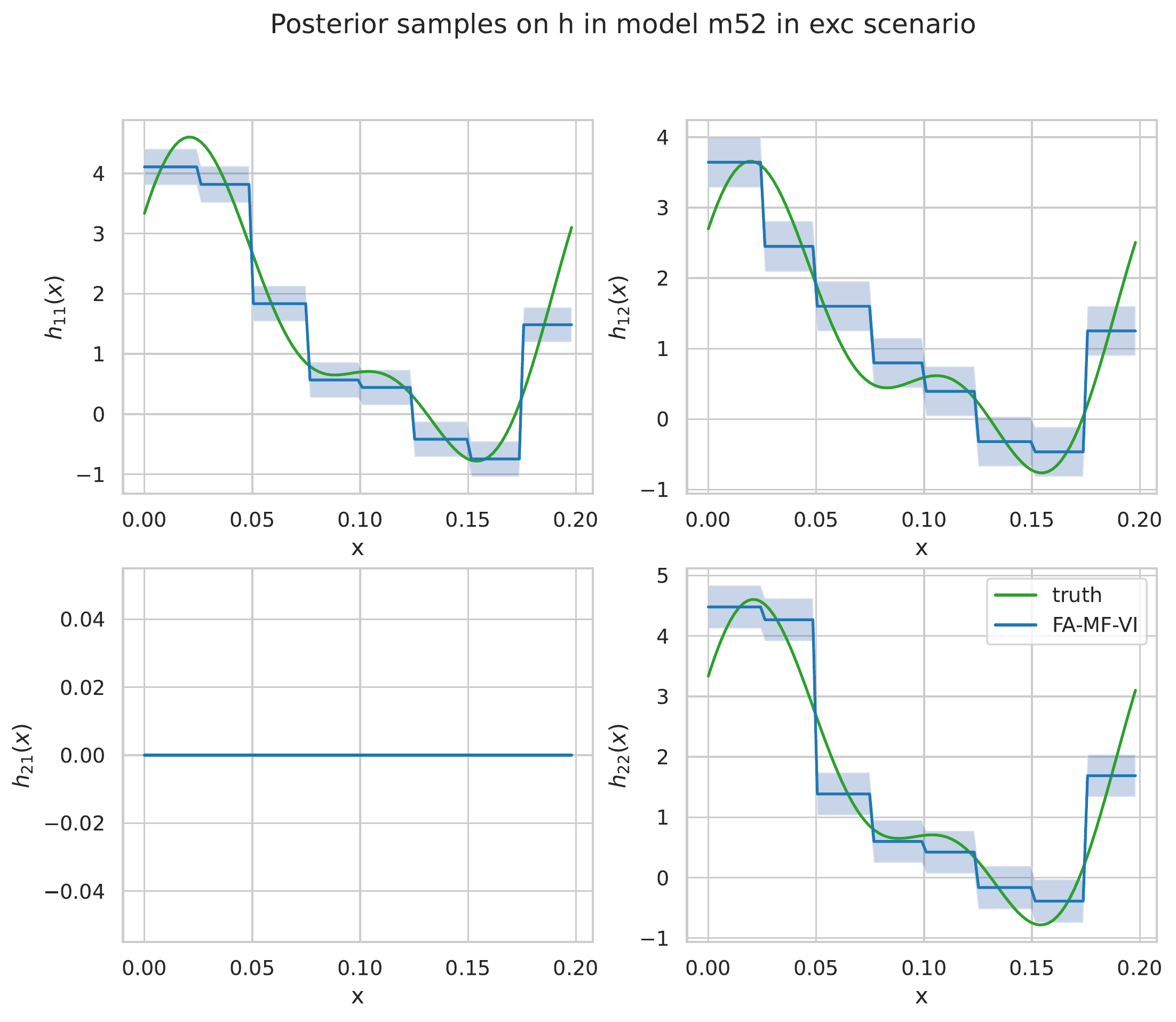}   
\caption{Model-selection variational posterior distributions on $f = (\nu, h)$ in the bivariate sigmoid model, and well-specified and mis-specified settings, and \emph{Excitation} scenario of Simulation 3, computed with the fully-adaptive mean-field variational (FA-MF-VI) algorithm (Algorithm \ref{alg:adapt_cavi}).  The first row correspond two columns correspond to the \emph{Excitation} (left) and \emph{Inhibition} (right) settings. The first row contains the marginal distribution on the background rates $(\nu_1, \nu_2)$, and the second and third rows represent the (variational) posterior mean (solid line) and  95\% credible sets (colored areas) on the four interaction function $h_{11}, h_{12}, h_{21}, h_{22}$.  The true parameter $f_0$ is plotted in dotted green line.}
\label{fig:adaptive_VI_2D_continuous}
\end{figure}

\begin{figure}[hbt!]
    \centering
     \begin{subfigure}[b]{0.49\textwidth}
    \includegraphics[width=\textwidth, trim=0.cm 0.cm 0cm  0.cm,clip]{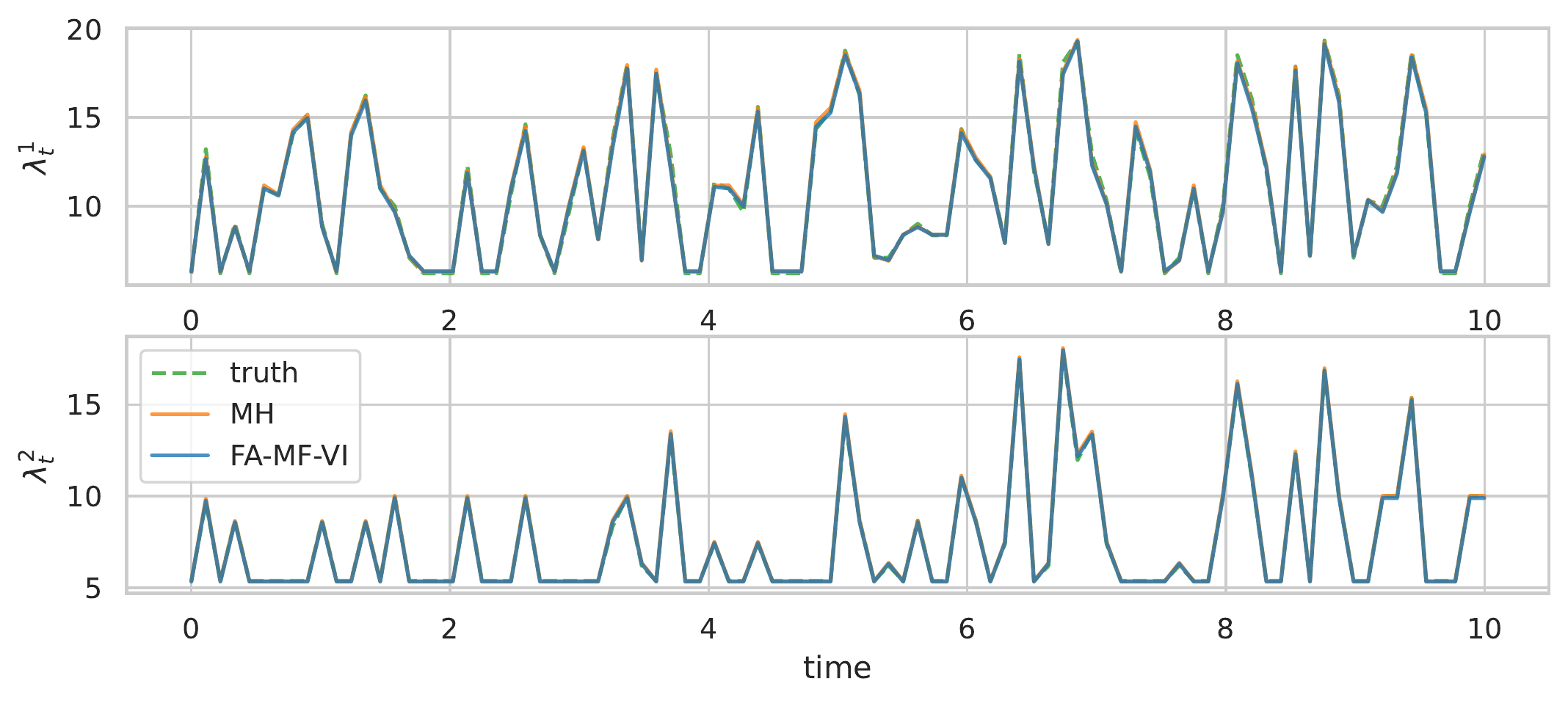}
    \caption{\emph{Excitation} scenario }
    \end{subfigure}
         \begin{subfigure}[b]{0.49\textwidth}
    \includegraphics[width=\textwidth, trim=0.cm 0.cm 0cm  0.cm,clip]{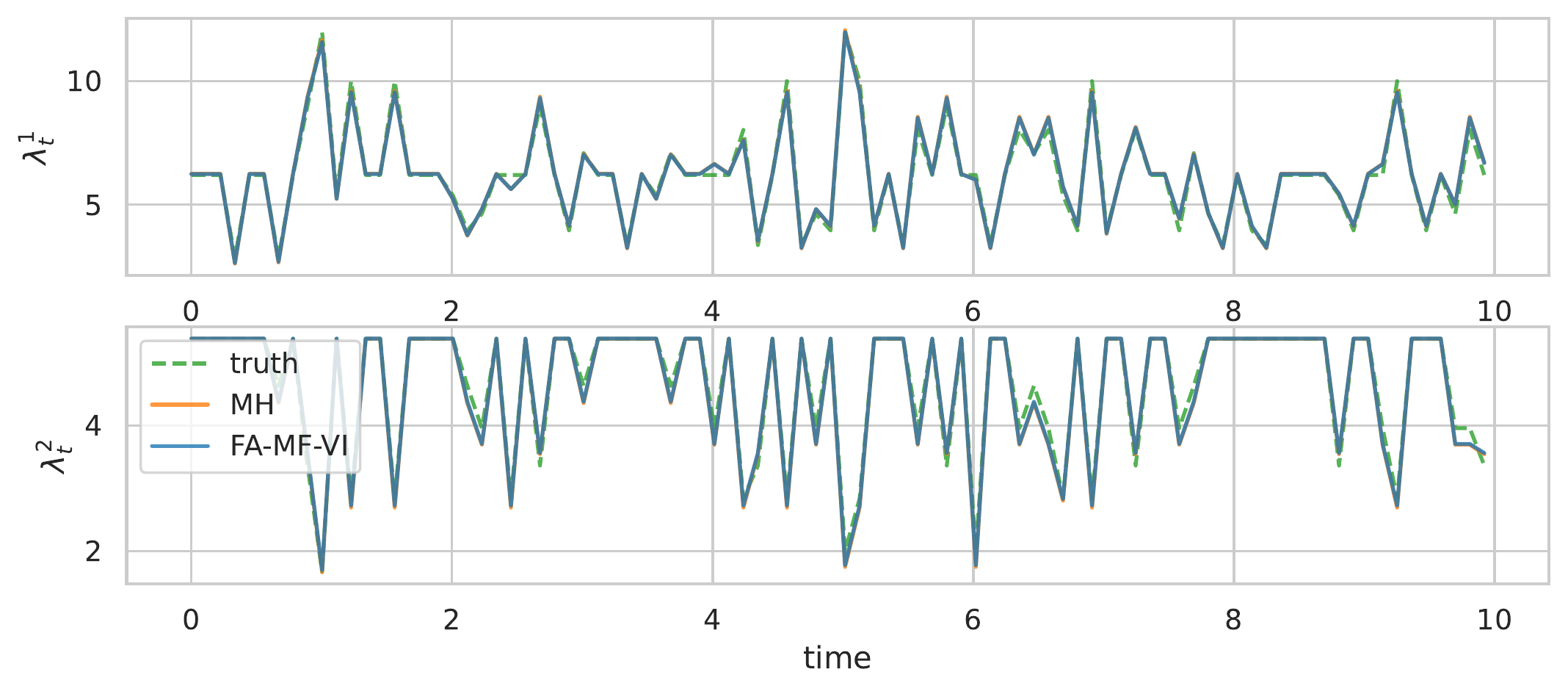}
    \caption{Self-inhibition scenario }
    \end{subfigure}
    \hfill
\caption{Estimated intensity function based on the (variational) posterior mean, in the well-specified and bivariate setting of Simulation 3 on  $[0,10]$,  using the fully-adaptive mean-field variational (FA-MF-VI) algorithm (Algorithm \ref{alg:adapt_cavi}). The true intensity $\lambda_t(f_0)$ is plotted in dotted green line.}
\label{fig:intensity_D2_adaptive}
\end{figure}

The previous observations in the well-specified and mis-specified settings can also be made in the two-dimensional setting. The true connectivity graph and the marginal probabilities in the adaptive variational posterior are plotted in Figure \ref{fig:graph_2D}. We note that in the well-specified case,  $\hat m = m_0$  in both scenarios. % Therefore, both the causality structure and the dimensionality are well recovered in this case.
Moreover, the parameter and the nonlinear intensity are well estimated, as can be seen in Figure \ref{fig:intensity_D2_adaptive} and in Figures \ref{fig:adaptive_VI_2D_selec_exc}, \ref{fig:adaptive_VI_2D_selec_inh} in Appendix \ref{app:details_exp}. Note however that, in the mis-specified setting, the under-coverage phenomenon of the credible regions also occurs (see Figure \ref{fig:adaptive_VI_2D_continuous}).

Finally, we note that our fully-adaptive variational algorithm is more than 10 times faster to compute than the non-adaptive MH sampler, % although here the latter is not adaptive, in particular in the bivariate setting
%\footnote{We further note that the current implementation of our algorithm has not been yet optimised.}, 
as can be seen from the computing times reported in Table \ref{tab:comptime_simu3}. This simulation study therefore shows that our fully-adaptive variational algorithm enjoys several advantages in Bayesian estimation for Hawkes processes: it can infer the dimensionality of the interaction functions $D$, the dependence structure through the graph parameter $\delta$, provides a good approximation of the posterior mean, and is computationally efficient.

%\textcolor{red}{Ju: What are the values of $T$ in the simulations associated to $K=1,2$ below. Is it $2000 $ or $3000$?}

%$f$: $\Pi(f) = \prod_k \mathcal{N}(f_k; \mu, \Sigma)$, with $\mu \in \R^{KB+1}, \Sigma \in \R^{(KB+1) \times (KB+1)} $.

\FloatBarrier

\subsection{Simulation 4: Two-step variational posterior in high-dimensional sigmoid models.}

%\textcolor{red}{ What are the values of $T$ chosen? It is said nowhere . Are they always the same for all $K$? in which case what are they?}

In this section, we test the performance of our two-step variational procedure (Algorithm \ref{alg:2step_adapt_cavi}), first, in sparse settings of the true parameter $h_0$, then, in relatively denser regimes.

\subsubsection{Sparse settings}

\begin{figure}[hbt!]
    \centering
     \begin{subfigure}[b]{0.4\textwidth}
    \includegraphics[width=\textwidth, trim=0.cm 0.cm 0cm  0.cm,clip]{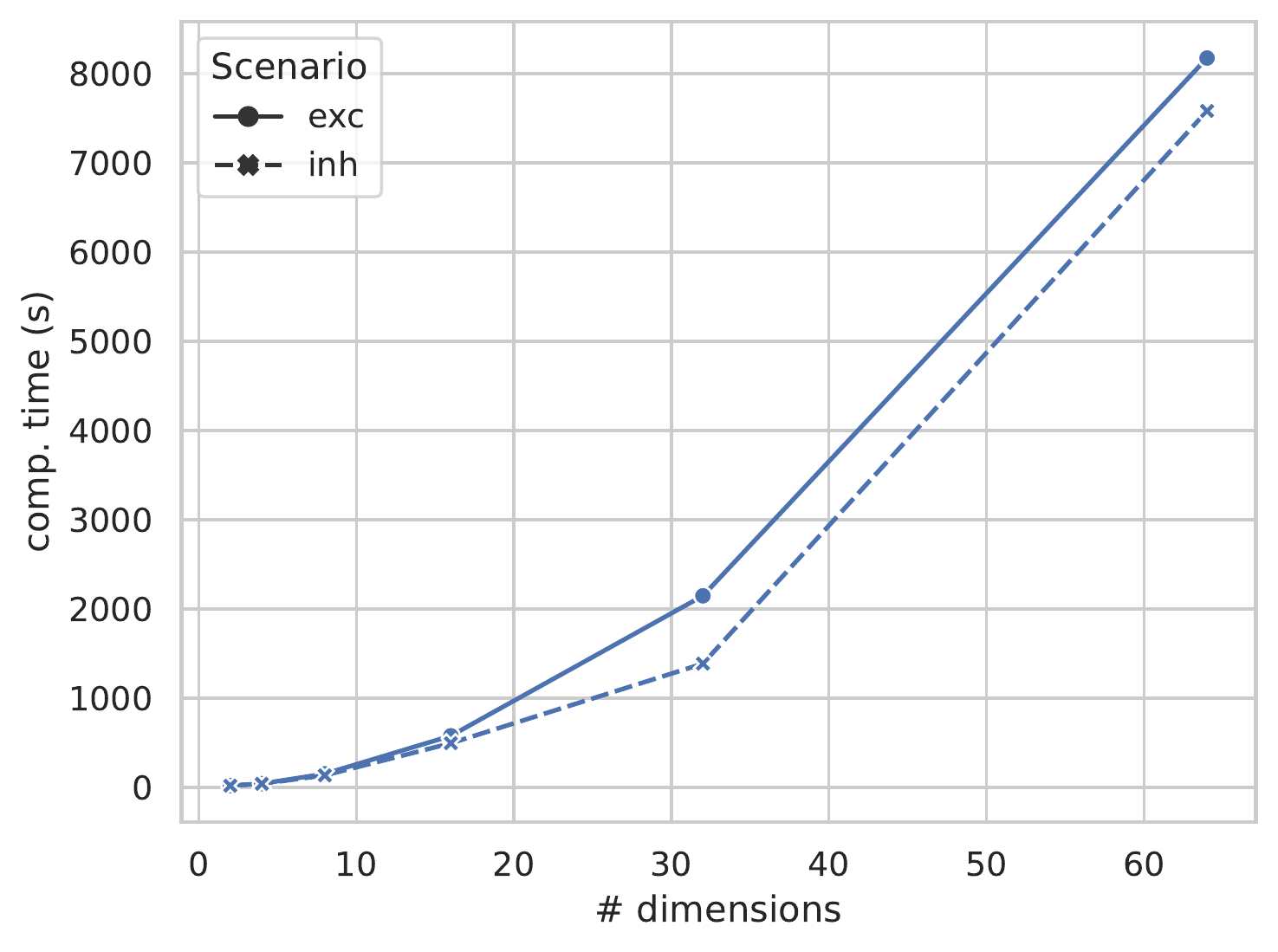}
    \end{subfigure}
    \hfill
\caption{Computational times of our two-step mean-field variational algorithm (Algorithm \ref{alg:2step_adapt_cavi}) in the \emph{Excitation} (exc) and \emph{Inhibition} (inh) scenarios and well-specified setting of Simulation 4, for $K = 2,4,8,16,32, 64$.}
\label{fig:adaptive_VI_comptime}
\end{figure}

\begin{figure}[hbt!]
    \centering
         \begin{subfigure}[b]{0.3\textwidth}
    \includegraphics[width=\textwidth, trim=0.cm .5cm 0cm  0.cm,clip]{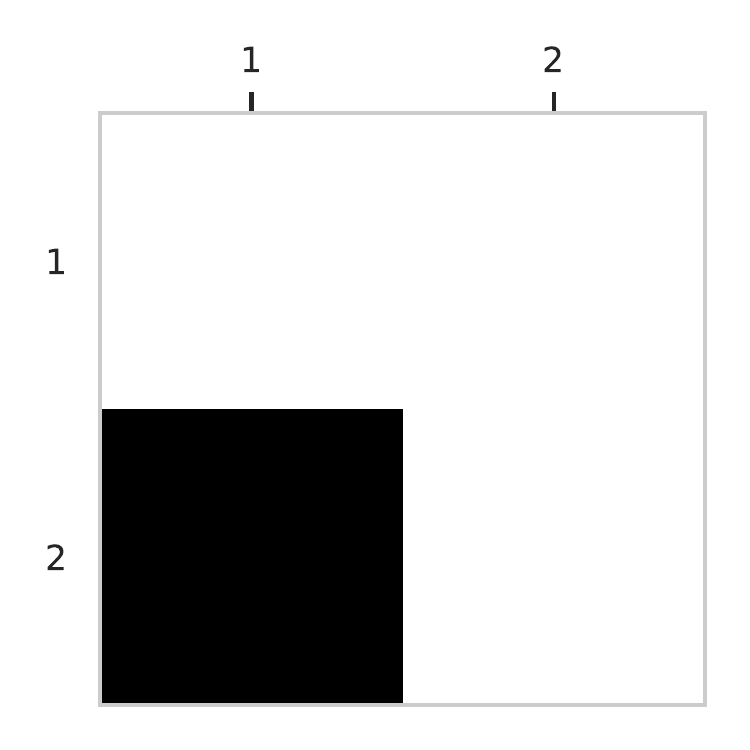}
    \caption{$K=2$}
    \end{subfigure}%
             \begin{subfigure}[b]{0.3\textwidth}
    \includegraphics[width=\textwidth, trim=0.cm .5cm 0cm  0.cm,clip]{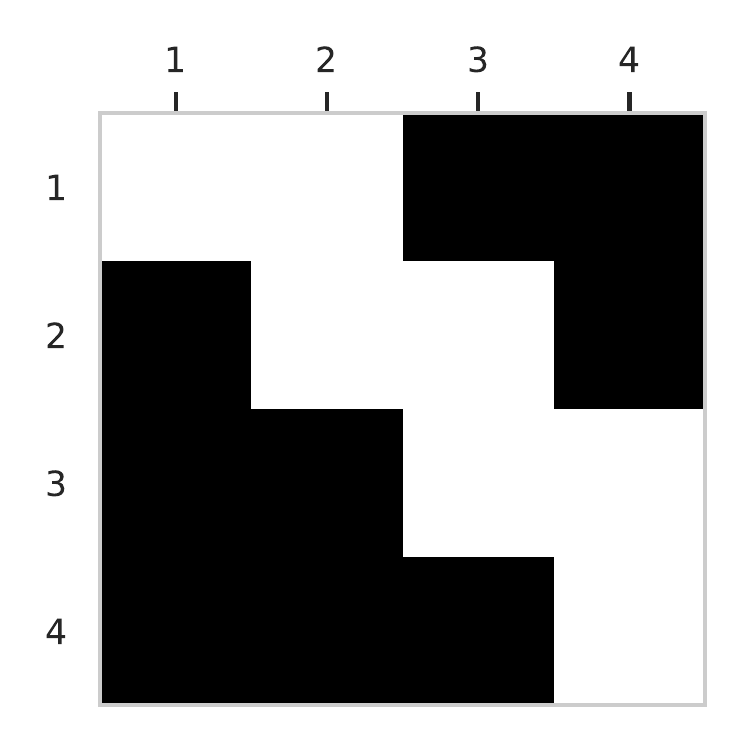}
    \caption{$K=4$}
    \end{subfigure}%
             \begin{subfigure}[b]{0.3\textwidth}
    \includegraphics[width=\textwidth, trim=0.cm .5cm 0cm  0.cm,clip]{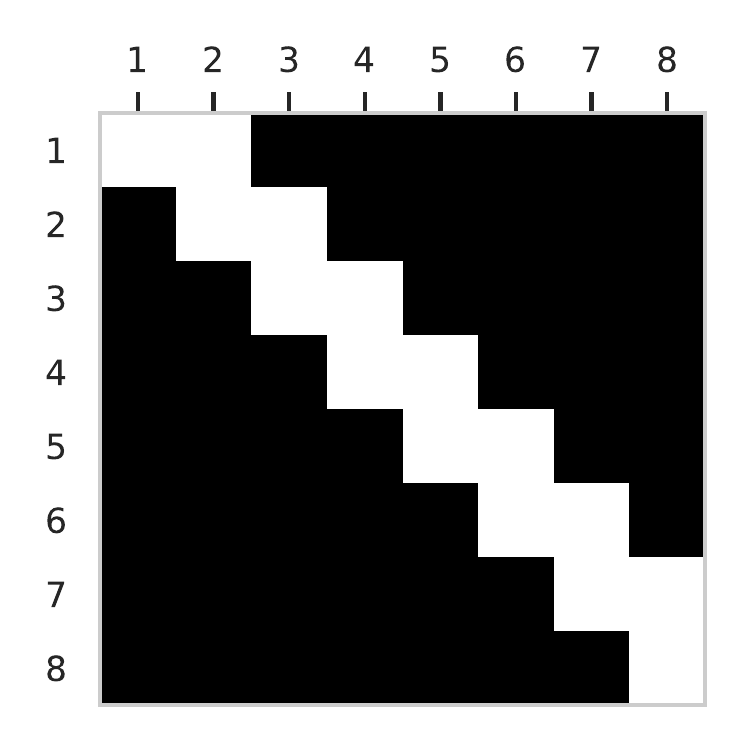}
    \caption{$K=8$}
    \end{subfigure}
        \begin{subfigure}[b]{0.33\textwidth}
    \includegraphics[width=\textwidth, trim=0.cm .5cm 0cm  0.cm,clip]{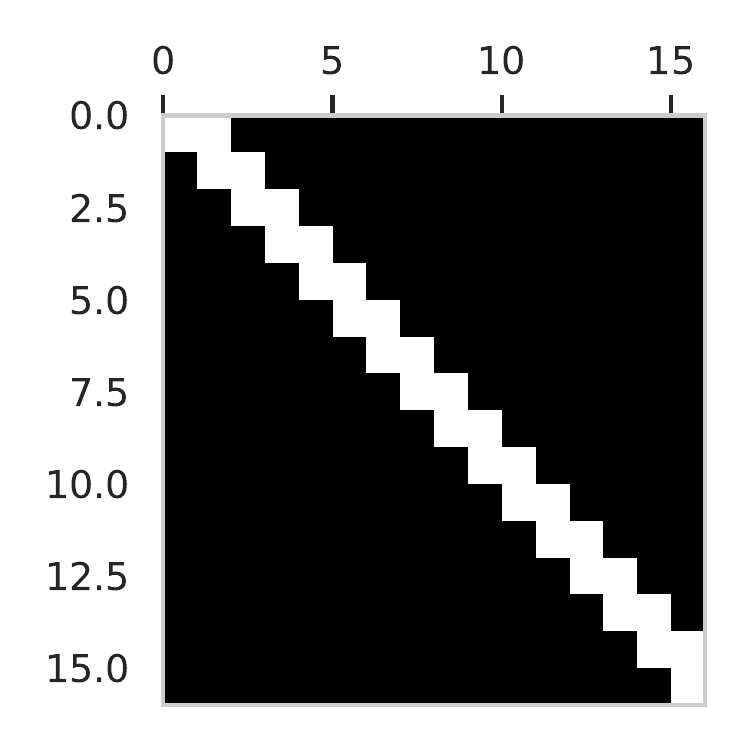}
    \caption{$K=10$}
    \end{subfigure}
    \begin{subfigure}[b]{0.3\textwidth}
    \includegraphics[width=\textwidth, trim=0.cm .5cm 0cm  0.cm,clip]{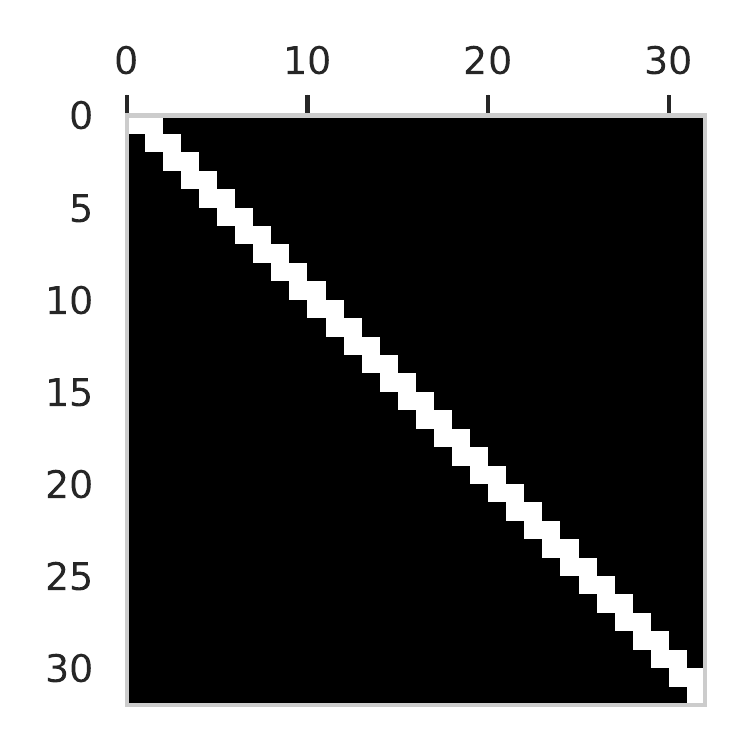}
    \caption{$K=16$}
    \end{subfigure}
        \begin{subfigure}[b]{0.3\textwidth}
    \includegraphics[width=\textwidth, trim=0.cm .5cm 0cm  0.cm,clip]{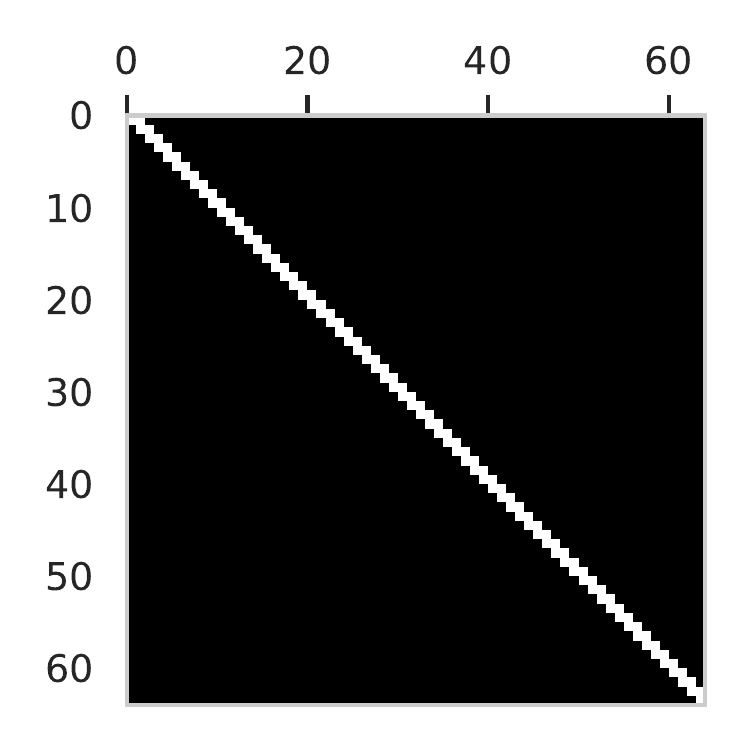}
    \caption{$K=32$}
    \end{subfigure}
\caption{True graph parameter $\delta_0$ (black=0, white=1) in the sparse multivariate settings of Simulations 4 with the number of dimensions $K=2,4,8, 16, 32, 64$.}
\label{fig:graphs}
\end{figure}

In this experiment, we consider sparse multivariate sigmoid models with $K \in \{2,4,8,16, 32, 64\}$ dimensions. We note that to the best of our knowledge, the only Bayesian method that has currently been tested in high-dimensional Hawkes processes is the semi-parametric version of \cite{zhou2021jmlr} where the interaction functions are  also decomposed over a dictionary of functions, but the choice of the number of functions is not driven by a model selection procedure and the graph of interaction is not inferred.
%We compare it to a MCMC method in the fixed model $\{\delta_C, D_{true}\}$,
Here, we construct a well-specified setting with $h_0 \in \mathcal{H}_{hist}^{D_0}$ and $D_0 = 1$, and an  \emph{Excitation} scenario and an \emph{Inhibition} scenario, similar to Simulation 3, and a \emph{sparse} connectivity graph parameter $\delta_0$ with $\sum_{l,k} \delta_{lk}^0 = 2K - 1$, as shown in Figure \ref{fig:graphs}. In Table \ref{tab:simu_4_data}, we report our chosen value of $T$ in each setting and the corresponding  number of events, excursions, and local excursions.
%We only report the results for the former in this section, and the ones for the latter can be found in Appendix \ref{app:additional_exp}.   
 In Table \ref{tab:perf_simu4}, we report  the performance of our method, in terms of the $L_1$-risk of the model-selection variational posterior defined as
\begin{align}\label{eq:risk}
    r_{L_1}(\hat Q) :=  \mathbb{E}_{\hat Q}[\norm{\nu-\nu_0}_{\ell_1}] +  \sum_{l,k} \mathbb{E}_{\hat Q}\left[\norm{h_{lk}-h_{lk}^0}_1\right].
\end{align}
We note that in general, the number of terms in the risk grows with $K$ and the number of non-null interaction functions in $h$ and $h_0$ - which thus can be of order $O(K^2)$ in a \emph{dense} setting.

We first note that for our prior distribution and for the augmented variational posterior distribution $\hat Q$ in a fixed model $m = (\delta, J=(J_k))$, we have that
\begin{align*}
  \mathbb{E}_{\hat Q_1}[\norm{h_{lk}}_1] %&=  \sum_{j=1}^{J_{k}}  \mathbb{E}_{\hat Q_1^{\delta}}[|h_{lk}^j|] \\
    &= \sum_{j=1}^{J_{k}} \sqrt{\frac{2}{\pi} [\Sigma_{lk}^{J_{k}}]_{jj}}  \exp \left \{ - \frac{[\tilde{\mu}_{lk}^{J_{k}}]_{j}^2}{[\Sigma_{lk}^{D_{k}}]_{jj} } \right \} - [\tilde{\mu}_{lk}^{J_{k}}]_{j} \left[ 1 - 2 \Phi\left(- \frac{[\tilde{\mu}_{lk}^{J_{k}}]_{j} }{\sqrt{[\Sigma_{lk}^{J_{k,C}}]_{jj} }}\right) \right].
\end{align*}

We evaluate the accuracy of our algorithm when estimating the graph of interaction and the size $D_k$ at each dimension $k$, defined as
\begin{align*}
    &Acc_{graph}(\hat \delta) = \frac{1}{K^2} \sum_{l,k} \mathds{1}_{\delta_{lk}^0 = \hat \delta_{lk}},
    &Acc_{dim}(\hat D) = \frac{1}{K} \sum_{k} \mathds{1}_{D_{k}^0 = \hat D_{k}},
\end{align*}
where $\hat \delta = (\hat \delta_{lk})_{l,k}$ and $\hat D = (\hat D_{k})_k$ are respectively the estimated graph and the inferred dimensionality of $(h_{.k})_k$ in Algorithm \ref{alg:2step_adapt_cavi}.

\begin{figure}[hbt!]
    \centering
    \begin{subfigure}[b]{0.49\textwidth}
    \includegraphics[width=\textwidth, trim=0.cm 0.cm 0cm  0.cm,clip]{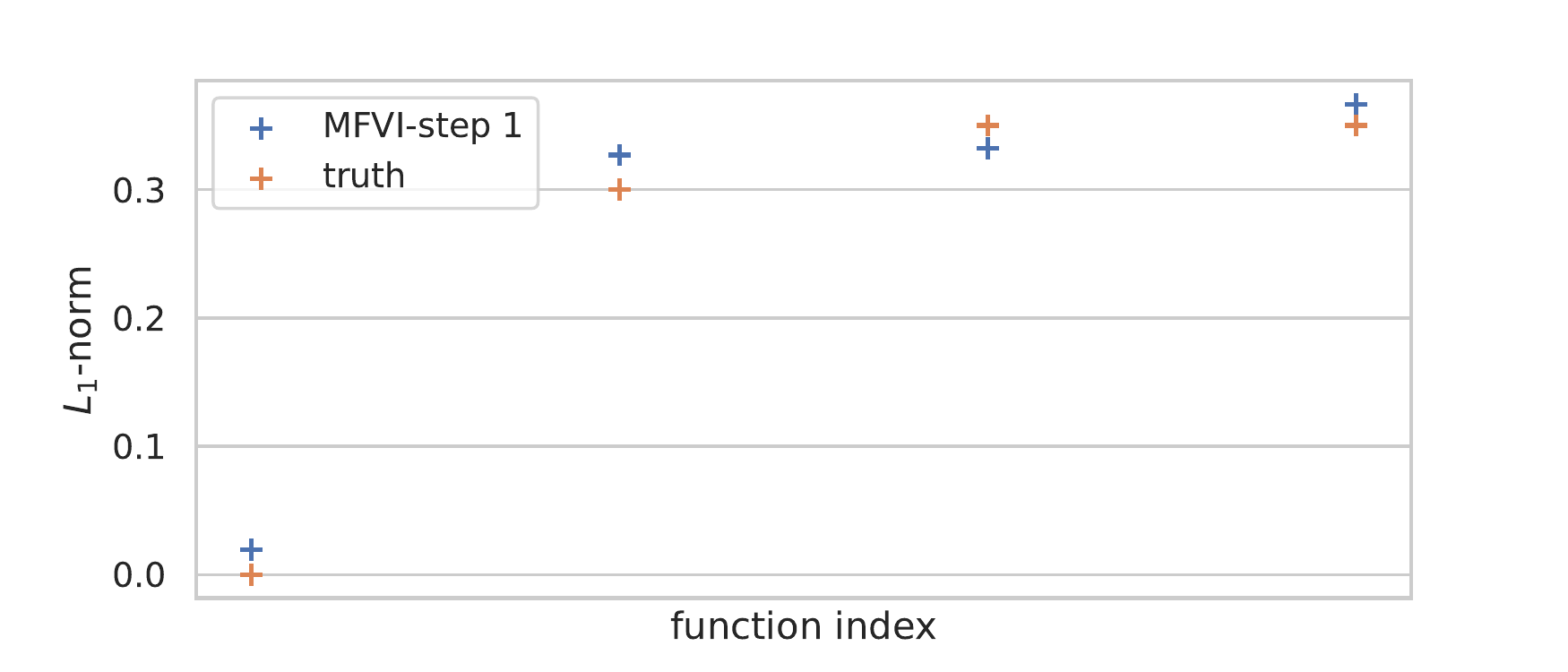}
    \caption{$K=2$}
    \end{subfigure}%
        \begin{subfigure}[b]{0.49\textwidth}
    \includegraphics[width=\textwidth, trim=0.cm 0.cm 0cm  0.cm,clip]{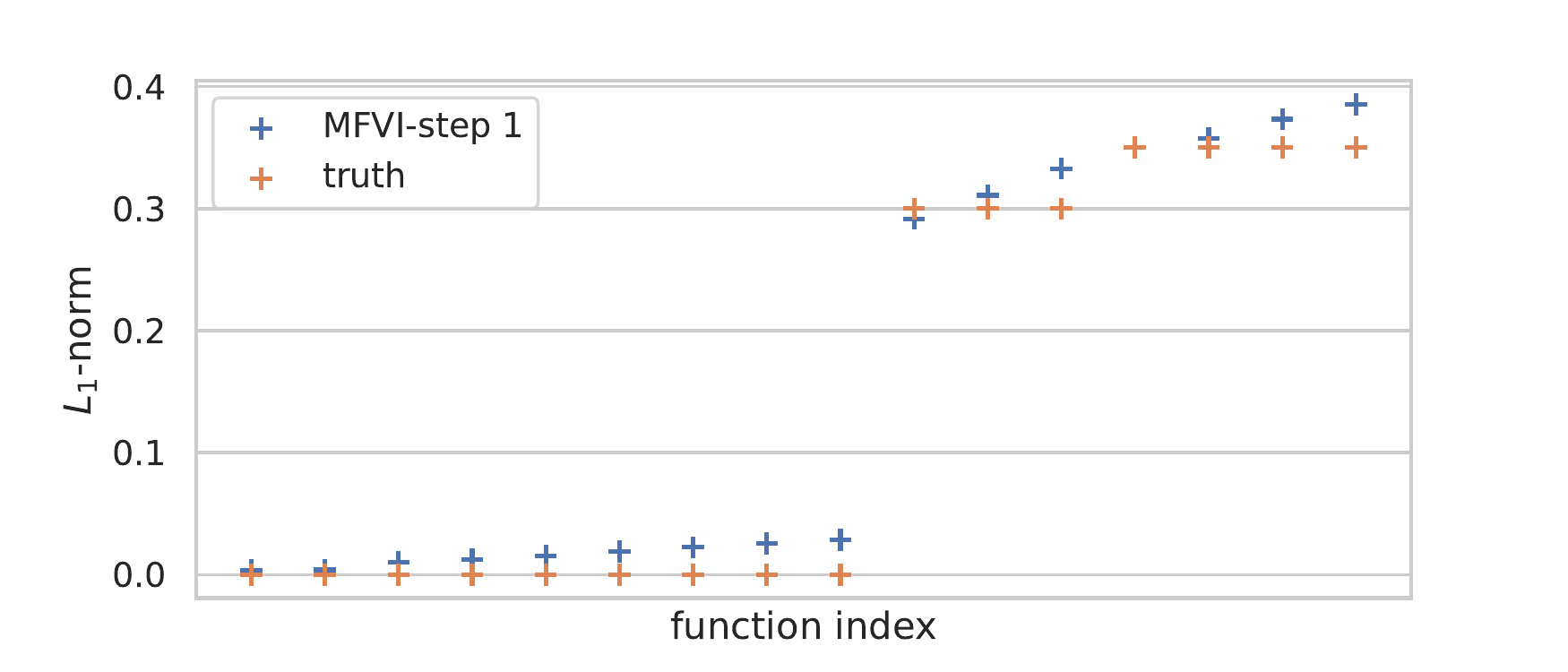}
    \caption{$K=4$}
    \end{subfigure}
     \begin{subfigure}[b]{0.49\textwidth}
    \includegraphics[width=\textwidth, trim=0.cm 0.cm 0cm  0.cm,clip]{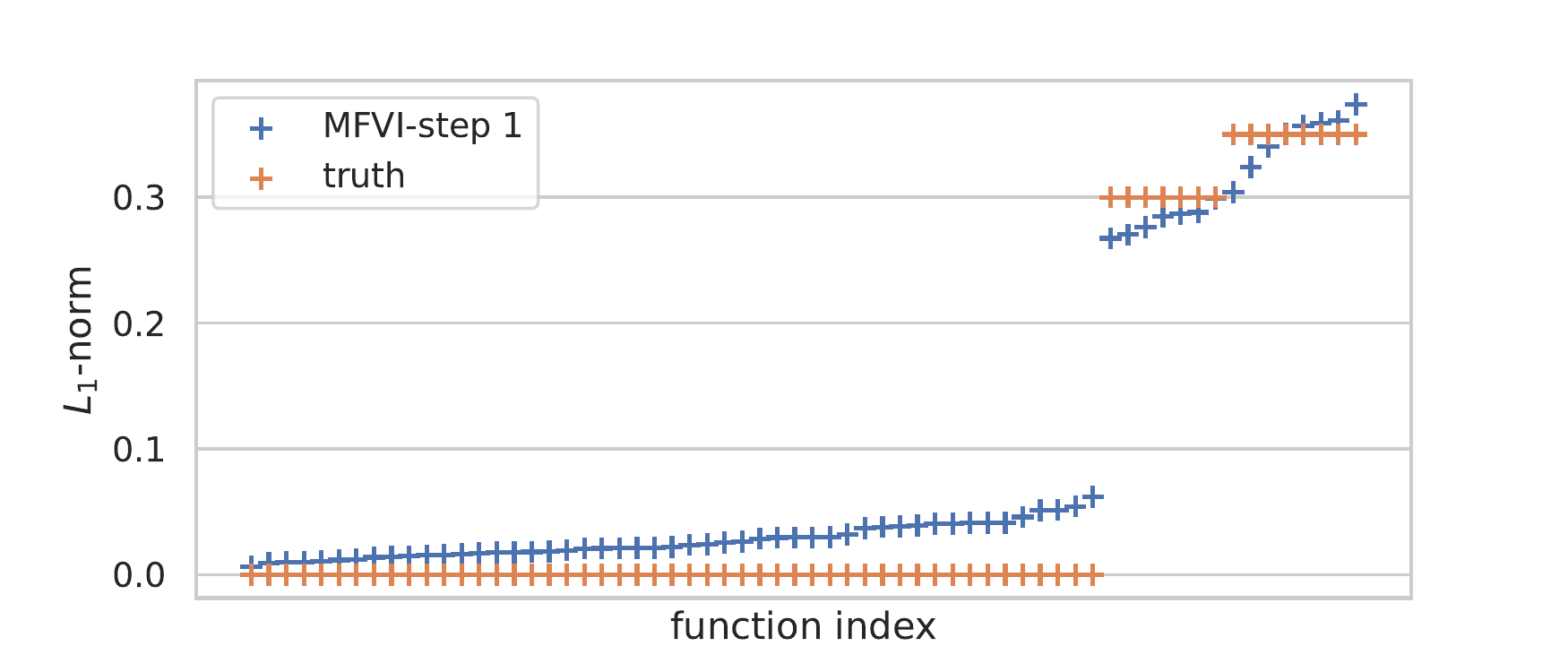}
    \caption{$K=8$}
    \end{subfigure}%
        \begin{subfigure}[b]{0.49\textwidth}
    \includegraphics[width=\textwidth, trim=0.cm 0.cm 0cm  0.cm,clip]{figs/simu_K16_L1_norms_exc.pdf}
    \caption{$K=16$}
    \end{subfigure}
            \begin{subfigure}[b]{0.49\textwidth}
    \includegraphics[width=\textwidth, trim=0.cm 0.cm 0cm  0.cm,clip]{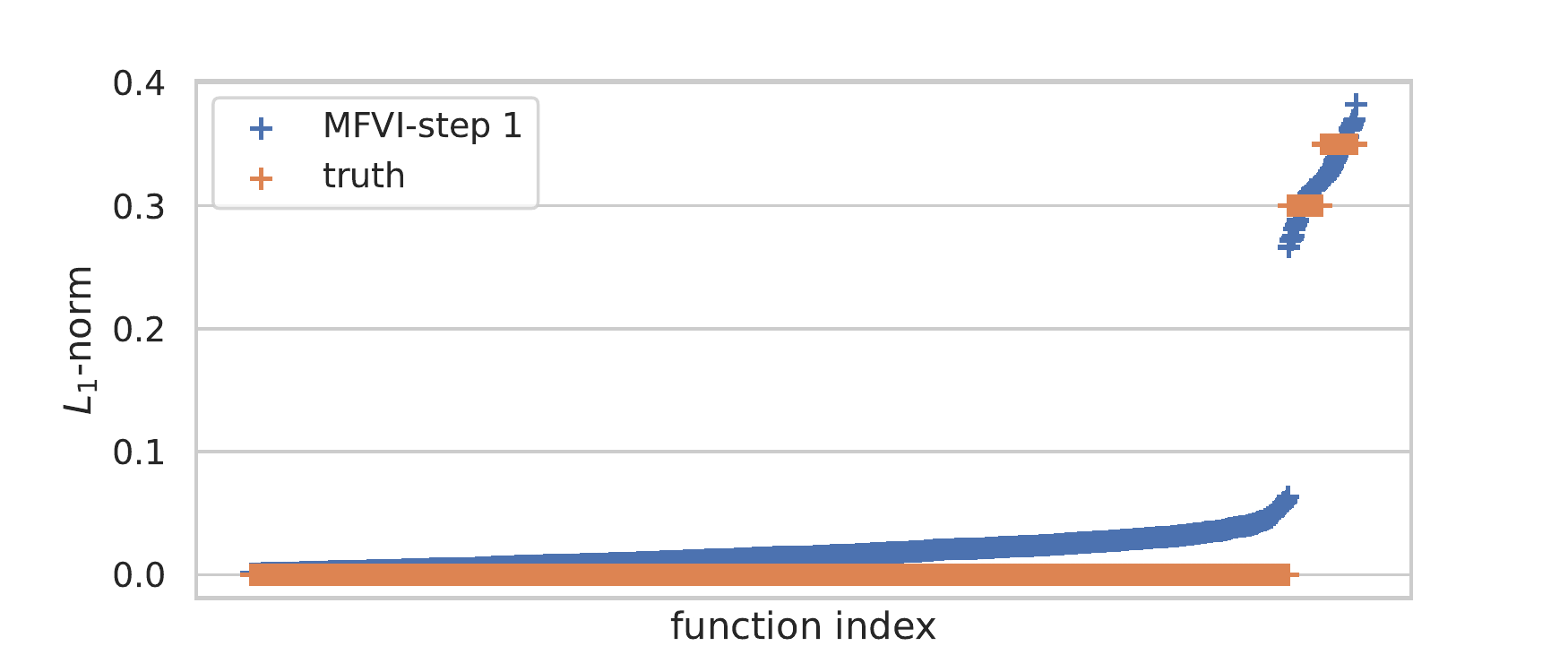}
    \caption{$K=32$}
    \end{subfigure}
     \begin{subfigure}[b]{0.49\textwidth}
    \includegraphics[width=\textwidth, trim=0.cm 0.cm 0cm  0.cm,clip]{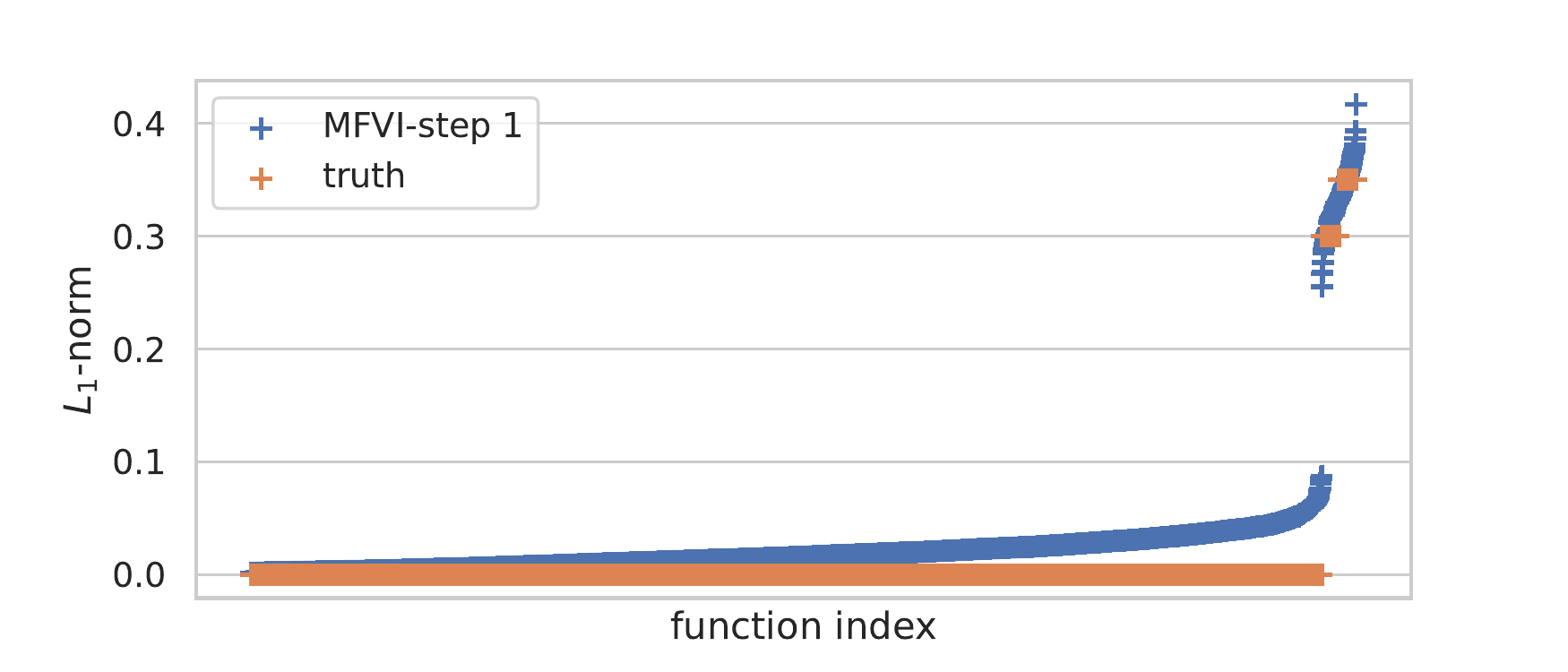}
    \caption{$K=64$}
    \end{subfigure}
\caption{Estimated $L_1$-norms using the model-selection variational posterior obtained after the first step of Algorithm \ref{alg:2step_adapt_cavi}, plotted in increasing order, in the \emph{Excitation} scenario of Simulation 4, for the models with $K=2,4,8,16, 32, 64$. In these settings, our threshold  $\eta_0 = 0.15$ is included in the gap between the estimated norms close to 0 and far from 0, therefore, our gap heuristics allows to recover the true graph parameter (see Section \ref{sec:two_step_mfvi}).}
\label{fig:adaptive_VI_2step_norms_threshold}
\end{figure}

Firstly, we note that, in almost all settings, the accuracy of our algorithm is equal or is very close to 1, therefore, it is able to recover almost perfectly the true graph $\delta_0$ and the dimensionality $D_0$ (the estimated graphs in the \emph{Excitation} and \emph{Inhibition} scenarios are plotted in Figures \ref{fig:graphs_exc} and \ref{fig:graphs_inh} in Appendix). In fact, our gap heuristics for choosing the threshold $\eta_0$ (see Section \ref{sec:two_step_mfvi}) allows to estimate the graph after the first step of Algorithm \ref{alg:2step_adapt_cavi}. In Figure \ref{fig:adaptive_VI_2step_norms_threshold} (and Figure \ref{fig:adaptive_VI_2step_norms_threshold_inhibition} in Appendix in the \emph{Inhibition} scenario), we note that the $L_1$-norms of the interaction functions are well estimated in the first step, leading to a gap between the norms close and far from 0. This gap includes the range $[0.1, 0.2]$ for all $K$'s, therefore, here, we choose $\eta_0 = 0.15$, which allows to discriminate between the true signals and the noise and to recover the true graph parameter.

\begin{figure}[hbt!]
\centering
\setlength{\tempwidth}{.25\linewidth}\centering
\settoheight{\tempheight}{\includegraphics[width=\tempwidth]{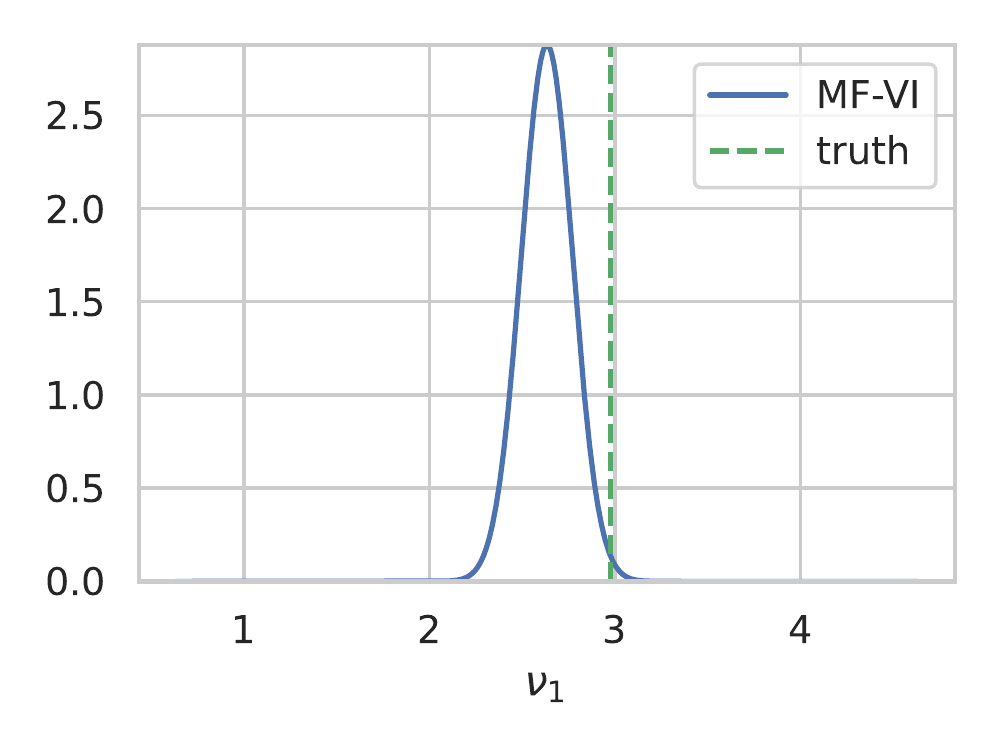}}
\hspace{-30mm} \fbox{\begin{minipage}  {\dimexpr 24mm} \begin{center} \itshape \large  \textbf{Excitation} \end{center} \end{minipage}}
%\hspace{-30mm} 
\columnname{Background $\nu_1$}\hfil
\columnname{Interaction functions $h_{11}$ and $h_{21}$ }\\
\begin{minipage}{\dimexpr 20mm}  \vspace{-35mm} \flushright{\itshape \large  \textbf{$K=2$} } \end{minipage}
    \includegraphics[width=\tempwidth, trim=0.cm 0.cm 0cm 0cm,clip]{figs/simu_K2_estimated_nu_exc.pdf}\hfil
    \includegraphics[width=.6\linewidth, trim=0.cm 0.cm 0cm  1.cm,clip]{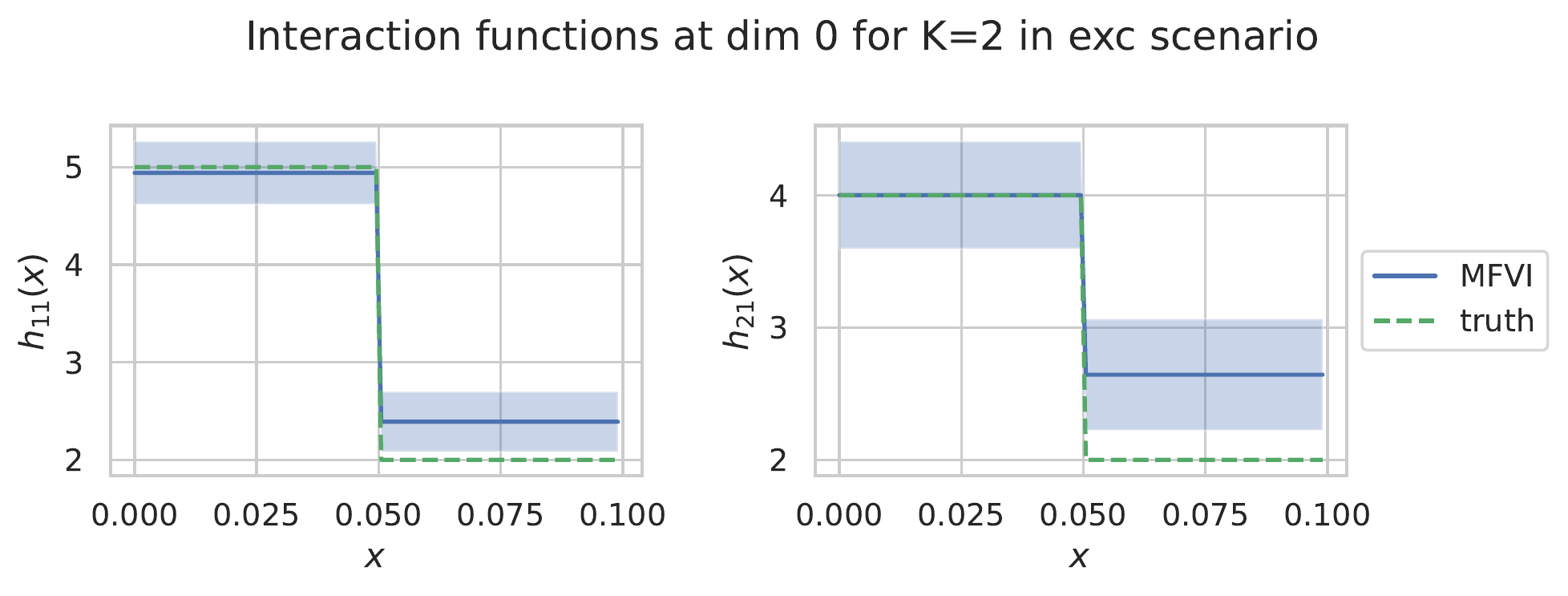}\\
\begin{minipage}{\dimexpr 20mm} \vspace{-35mm} \flushright{\itshape \large  \textbf{$K=4$} } \end{minipage}
    \includegraphics[width=\tempwidth, trim=0.cm 0.cm 0cm 0.cm,clip]{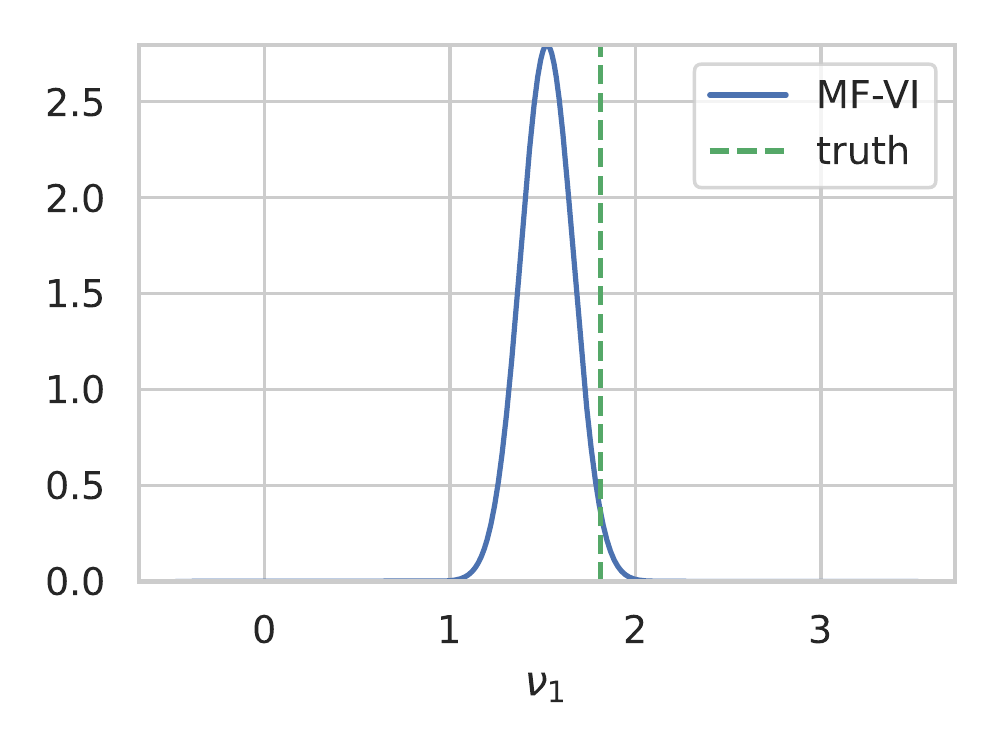}\hfil
    \includegraphics[width=.6\linewidth, trim=0.cm 0.cm 0cm  1.cm,clip]{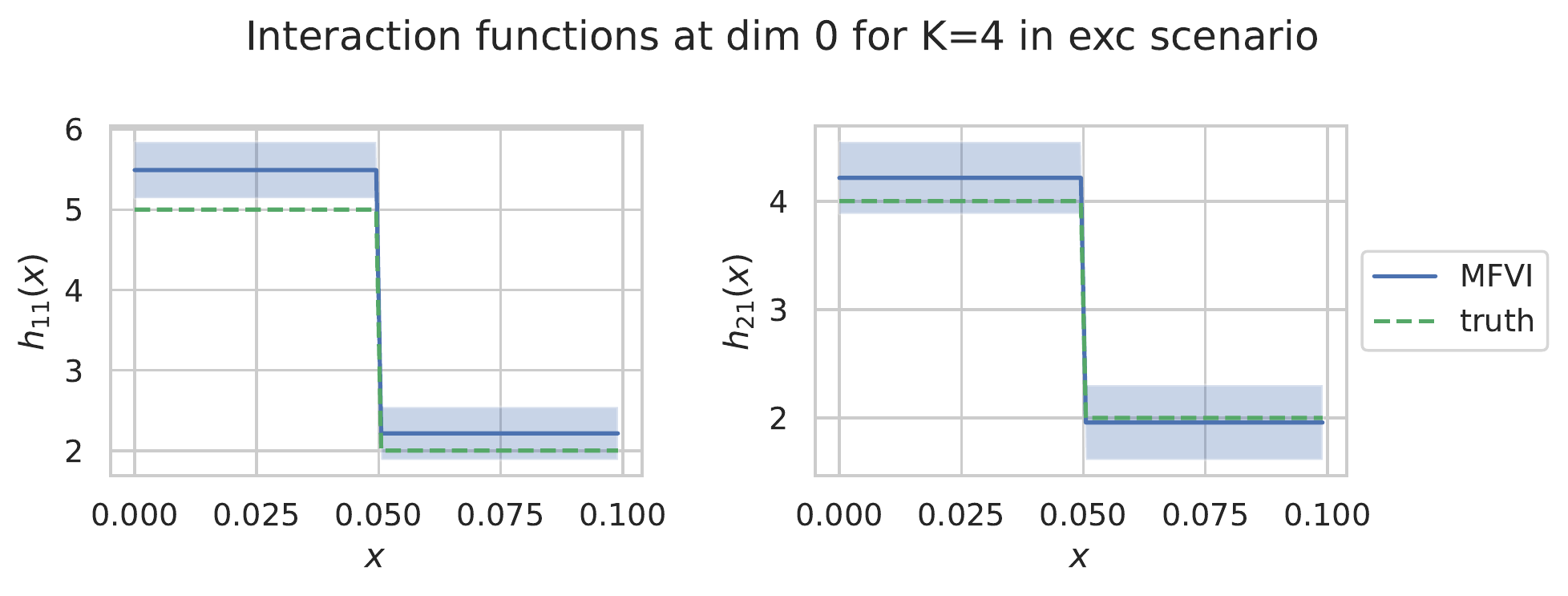}\\
\begin{minipage}{\dimexpr 20mm} \vspace{-35mm} \flushright{\itshape \large  \textbf{$K=8$} } \end{minipage}
    \includegraphics[width=\tempwidth, trim=0.cm 0.cm 0cm 0cm,clip]{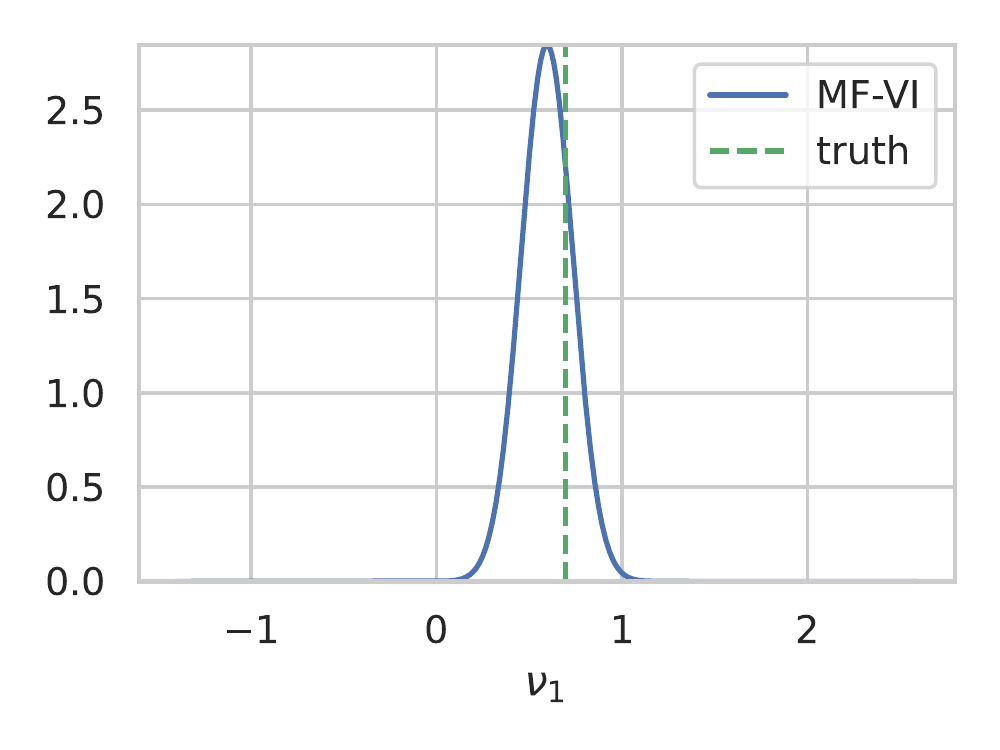}\hfil
    \includegraphics[width=.6\linewidth, trim=0.cm 0.cm 0cm  1.cm,clip]{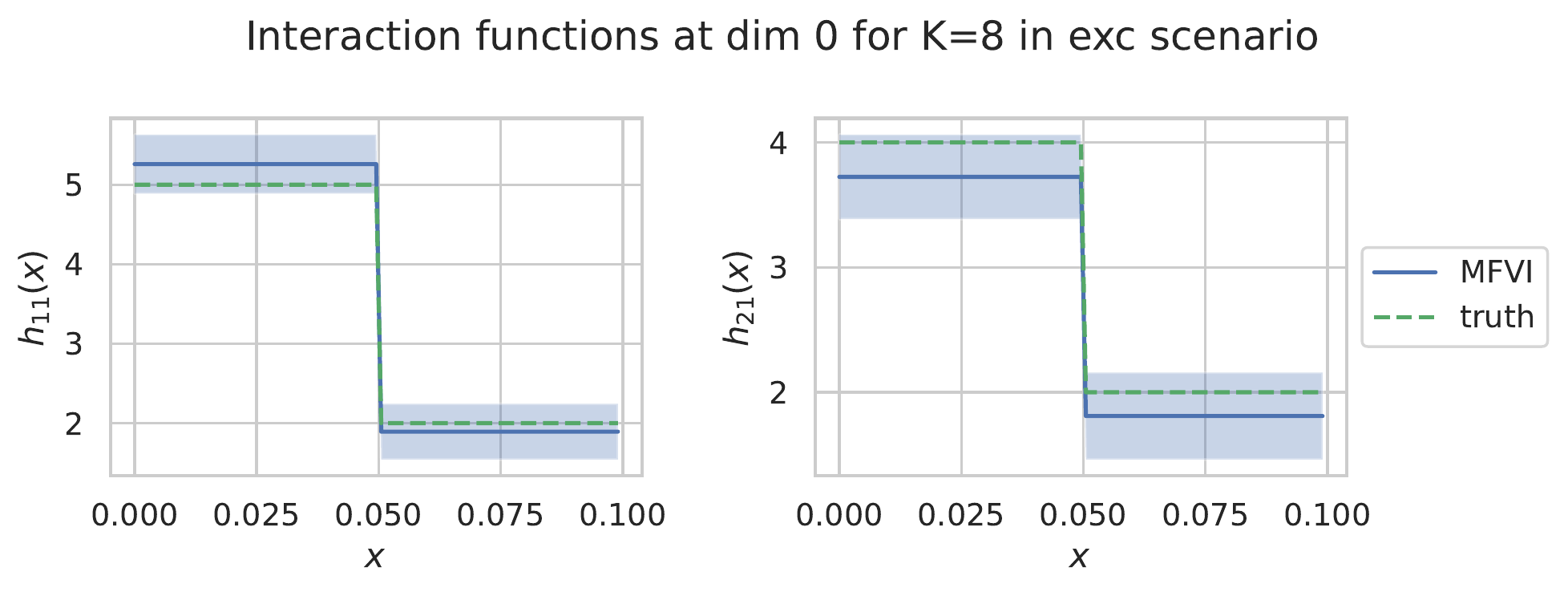}\\
\begin{minipage}{\dimexpr 20mm} \vspace{-35mm} \flushright{\itshape \large  \textbf{$K=16$} } \end{minipage}
    \includegraphics[width=\tempwidth, trim=0.cm 0.cm 0cm 0cm,clip]{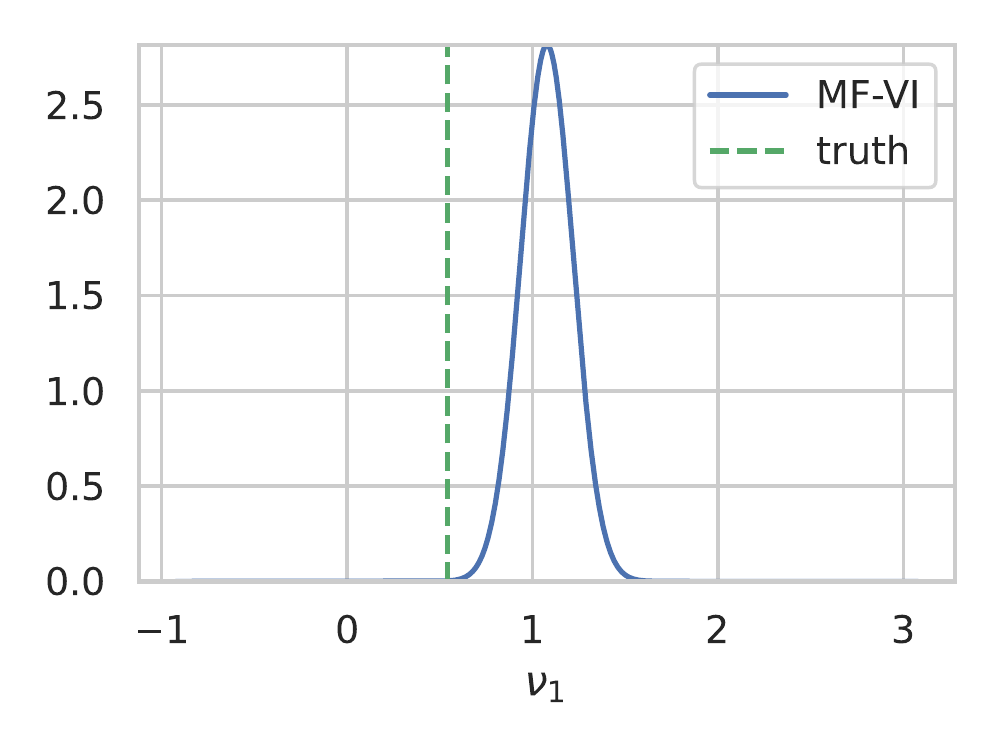}\hfil
    \includegraphics[width=.6\linewidth, trim=0.cm 0.cm 0cm  1.cm,clip]
    {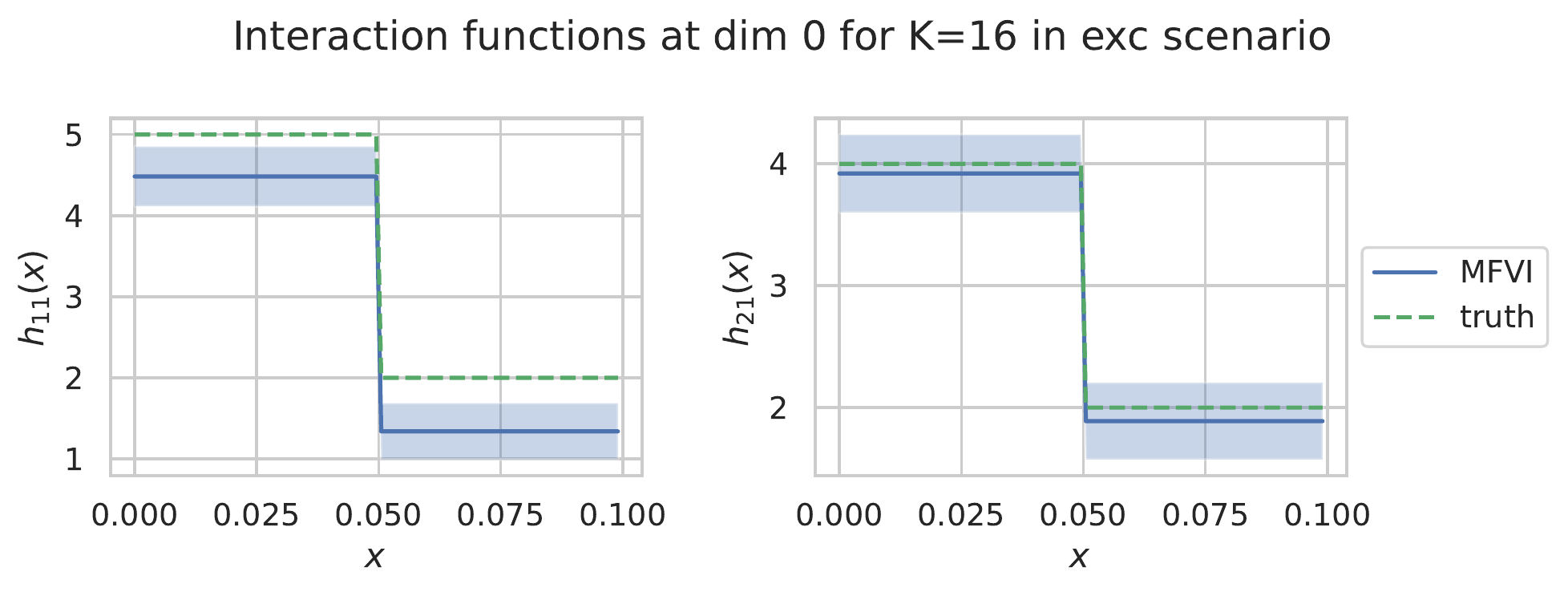}\\
\begin{minipage}{\dimexpr 20mm} \vspace{-35mm} \flushright{\itshape \large  \textbf{$K=32$} } \end{minipage}
    \includegraphics[width=\tempwidth, trim=0.cm 0.cm 0cm 0cm,clip]{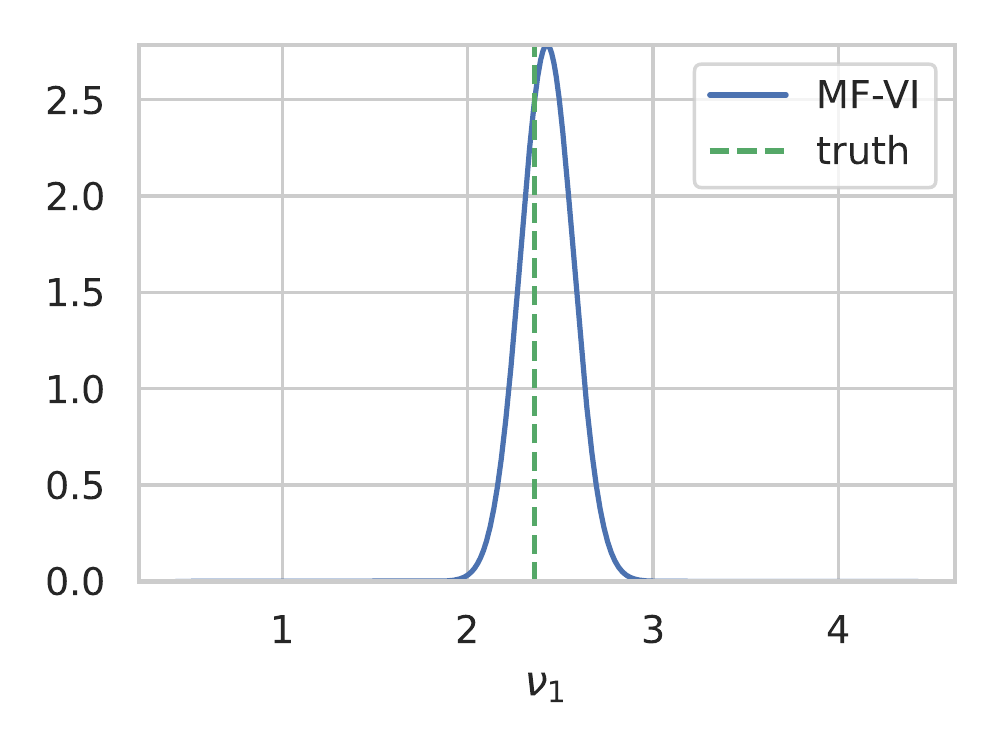}\hfil
    \includegraphics[width=.6\linewidth, trim=0.cm 0.cm 0cm  1.cm,clip]{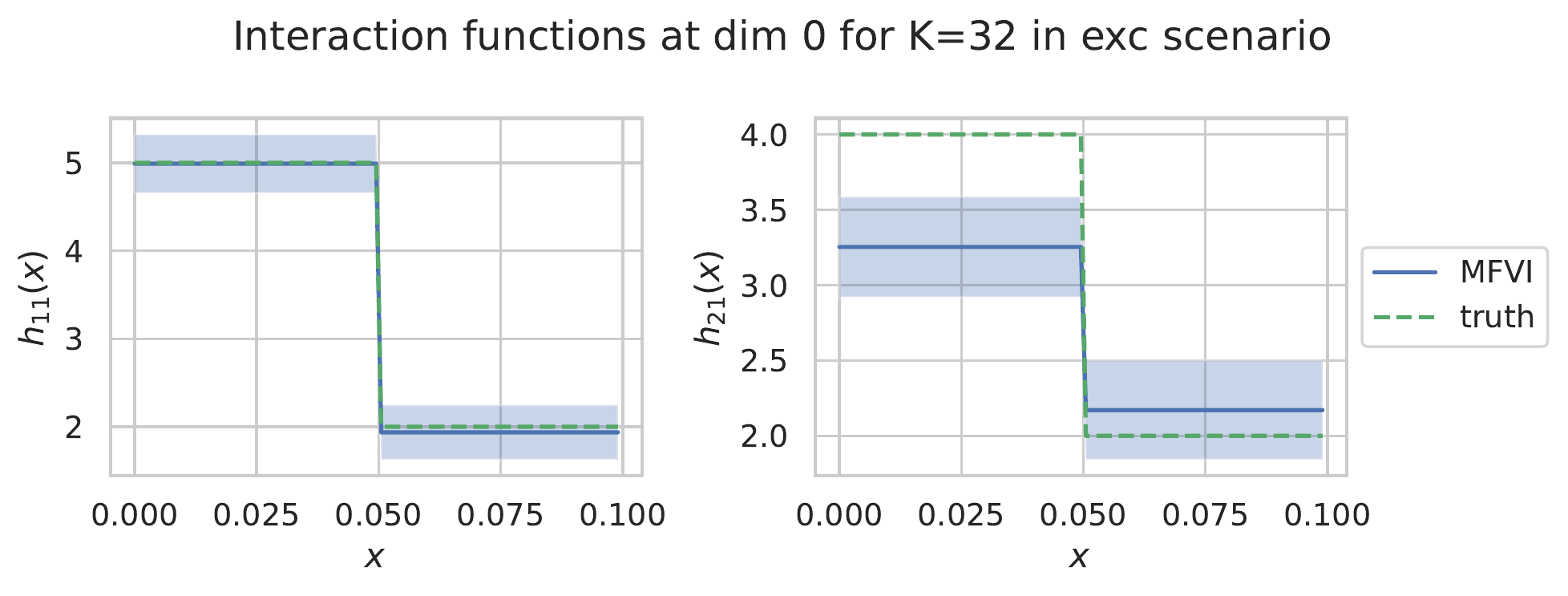}
\begin{minipage}{\dimexpr 20mm} \vspace{-35mm} \flushright{\itshape \large  \textbf{$K=64$} } \end{minipage}
    \includegraphics[width=\tempwidth, trim=0.cm 0.cm 0cm 0cm,clip]{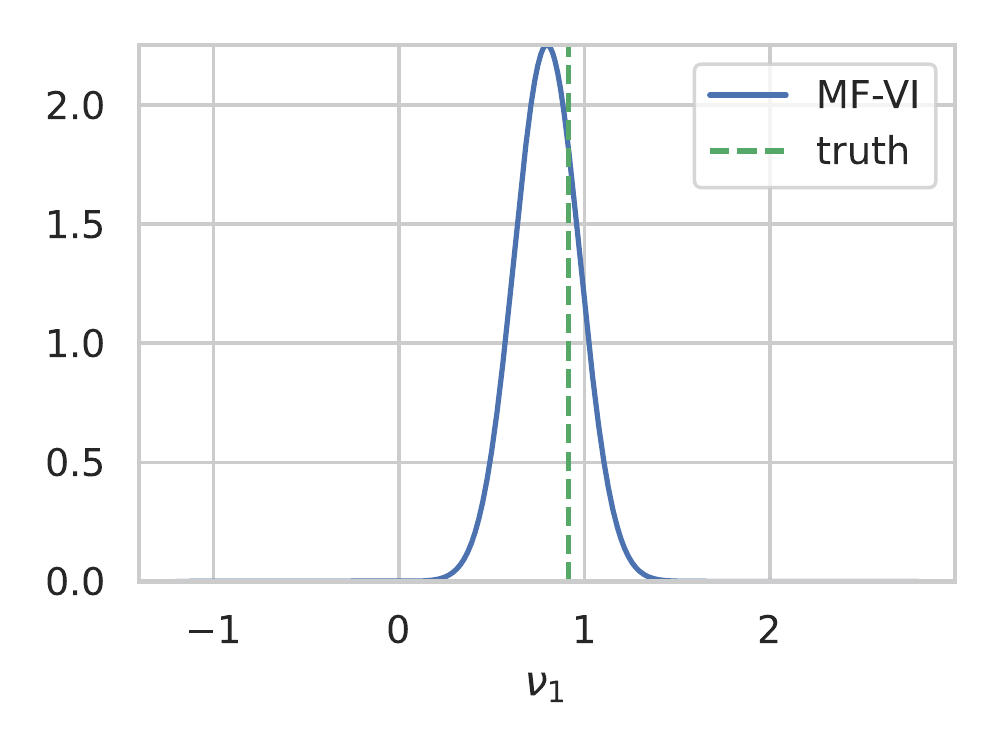}\hfil
    \includegraphics[width=.6\linewidth, trim=0.cm 0.cm 0cm  1.cm,clip]{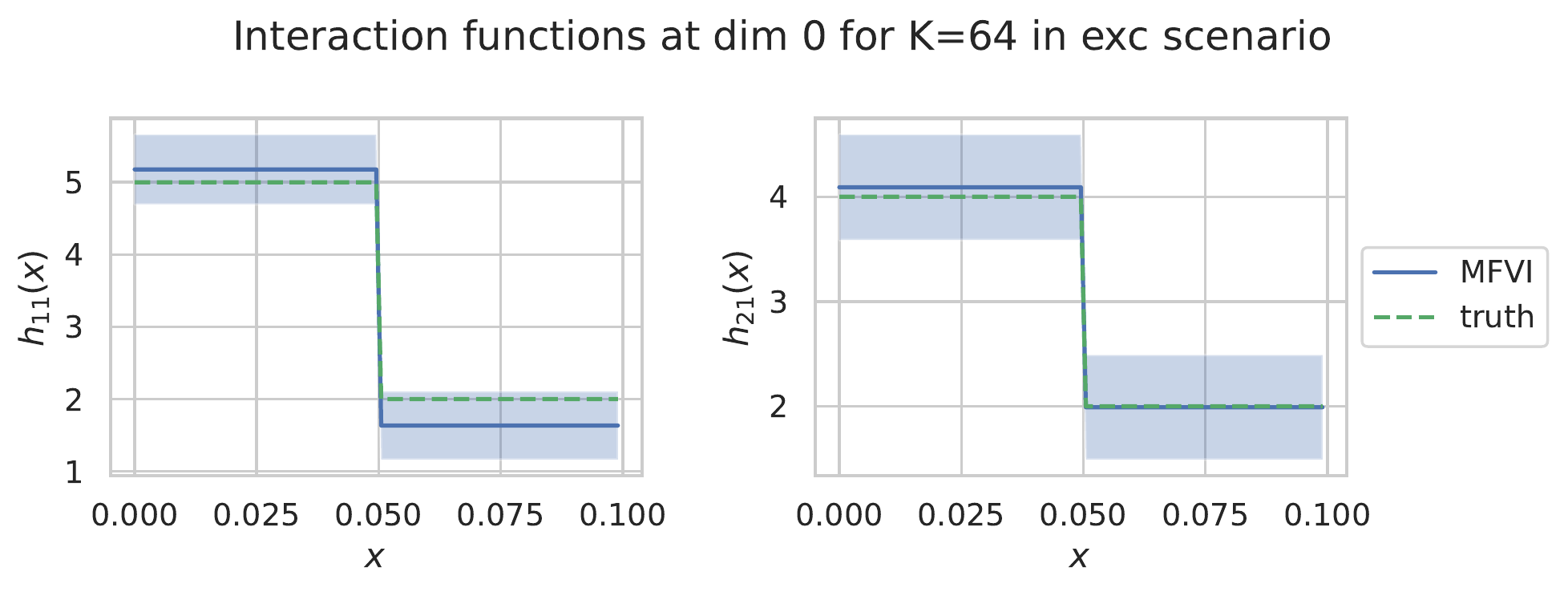}
\caption{Model-selection variational posterior distributions on $\nu_1$ (left column) and interaction functions $h_{11}$ and $ h_{21}$ (second and third columns) in the \emph{Excitation} scenario and multivariate sigmoid models of Simulation 4, computed with our two-step mean-field variational (MF-VI) algorithm (Algorithm \ref{alg:2step_adapt_cavi}). The different rows correspond to different multivariate settings $K=2,4,8,16,32, 64$.}
\label{fig:2step_adaptive_VI_exc_f}
\end{figure}

Secondly, from Table \ref{tab:perf_simu4}, we note that the risk seems to grow linearly with $K$, which indicates that the estimation does not deteriorate with larger $K$. In Figure \ref{fig:adaptive_VI_2step_norms_exc} (and Figure \ref{fig:adaptive_VI_2step_norms_inh} in Appendix, we plot the risk on the $L_1$-norms using the model-selection variational posterior, i.e., $(\mathbb{E}_{\hat Q_{MV}}\left[\norm{h_{lk}-h_{lk}^0}_1\right])_{l,k}$,  in the form of a heatmap compared to the true norms, and note that for all $K$'s, these errors are relatively small. Moreover, our variational algorithm estimates well the parameter, as can be visually checked in Figure \ref{fig:2step_adaptive_VI_exc_f}, where we plot the model-selection variational posterior distribution on a subset of the parameter for each value of $K$, in the \emph{Excitation} scenario (see Figure \ref{fig:2step_adaptive_VI_inh_f} in Appendix for our results in the \emph{Inhibition} scenario). 
%the estimated background rate $\nu_2$ and interaction functions $h_{12}, h_{22}, h_{32}, h_{42}$parameter
%We note that  the 95\% credible bands on the background rate $\nu_2$ become larger for larger $K$, however, this phenomenon does not appear for the  interaction functions. 
%In Figure \ref{fig:adaptive_VI_2step_norms}, we also plot the $L_1$-errors using $\hat Q^{\delta_c}$, i.e., $(\mathbb{E}_{\hat Q^{\delta_c}}\left[\norm{h_{lk}-h_{lk}^0}_1\right])_{l,k}$ in the form of a heatmap compared to the true norms. We recall that $\hat Q^{\delta_c}$ is mode variational distribution obtained after the first step of Algorithm \ref{alg:2step_adapt_cavi}, and is used to estimate the graph for the second step.  We note that in all settings, these errors are relatively small, therefore allowing us to select the true graph parameter for the second step. Indeed, in Figure \ref{fig:adaptive_VI_2step_norms_threshold}, we plot the estimated $L_1$-norms of the interaction functions using $\hat Q^{\delta_c}$, i.e., $(\mathbb{E}_{\hat Q^{\delta_c}}\left[\norm{h_{lk}}_1\right])_{l,k}$ in increasing order of magnitude and observe a gap between the small and larger ``signals". 
%Similar observation can also be made for the Self-inhibition scenario, which results are in Appendix \ref{app:additional_exp},  although the estimation is slightly worse in this scenario. 
Besides, the computing times of our algorithm seem to scale well with $K$ and the number of events in these sparse settings, as can be seen from Table \ref{tab:simu_4_data} and Figure \ref{fig:adaptive_VI_comptime}.  For $K=64$, our algorithm runs in less than 2.5 hours, in spite of  the large number of events (about 133 000). We also note that these experiments have been run using only two processing units. \footnote{The computing time of our algorithm could thus be greatly decreased if it is computed on a machine with more processing units.}

\subsubsection{Testing different graphs and sparsity levels.}

\begin{figure}[hbt!]
    \centering
         \begin{subfigure}[b]{0.3\textwidth}
    \includegraphics[width=\textwidth, trim=0.cm 0cm 0cm  0.cm,clip]{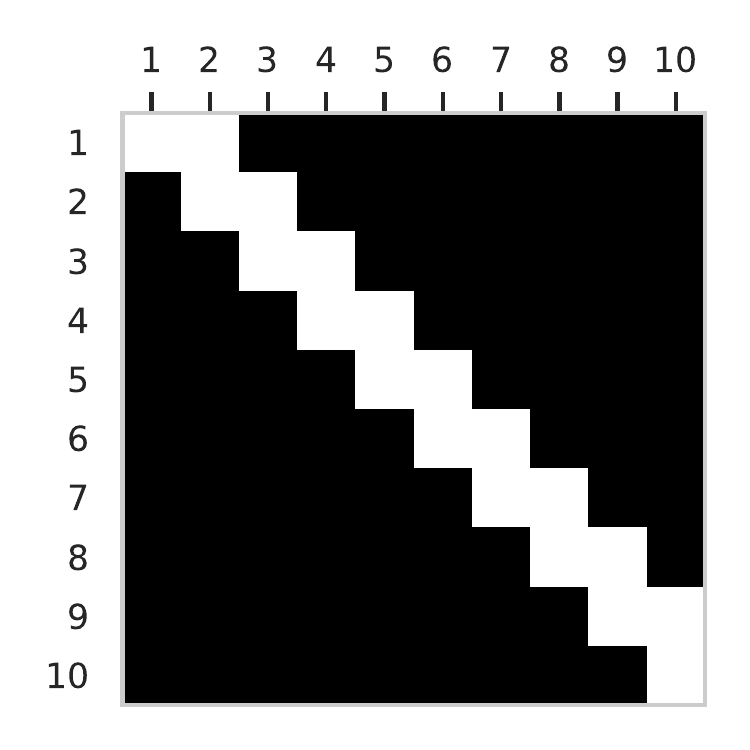}
    \caption{sparse}
    \end{subfigure}%
             \begin{subfigure}[b]{0.31\textwidth}
    \includegraphics[width=\textwidth, trim=1.cm 0cm 0cm  0.cm,clip]{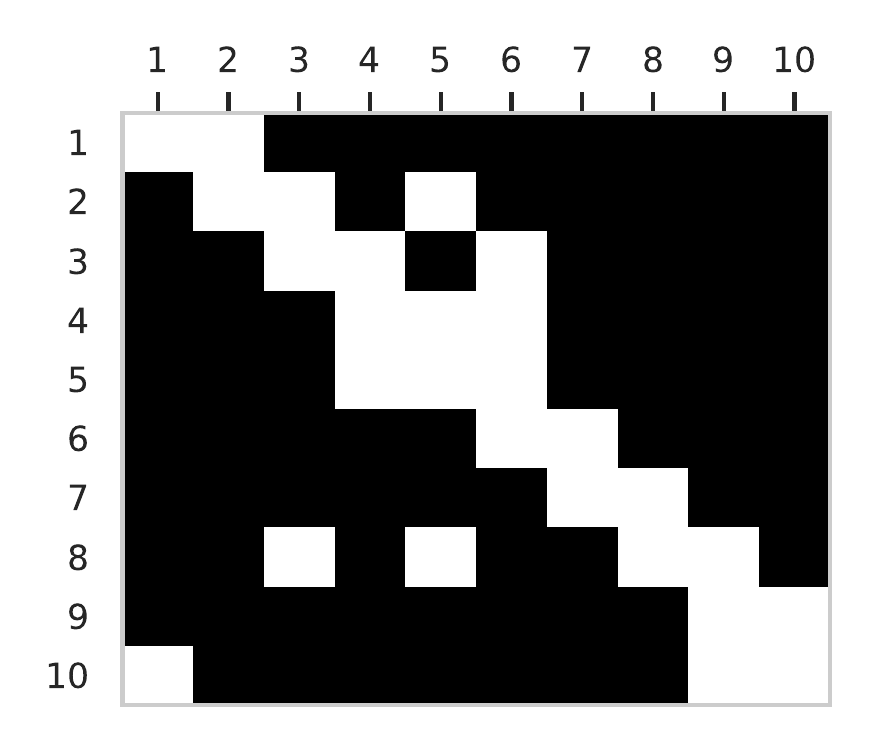}
    \caption{random}
    \end{subfigure}%
             \begin{subfigure}[b]{0.31\textwidth}
    \includegraphics[width=\textwidth, trim=1.cm 0cm 0cm  0.cm,clip]{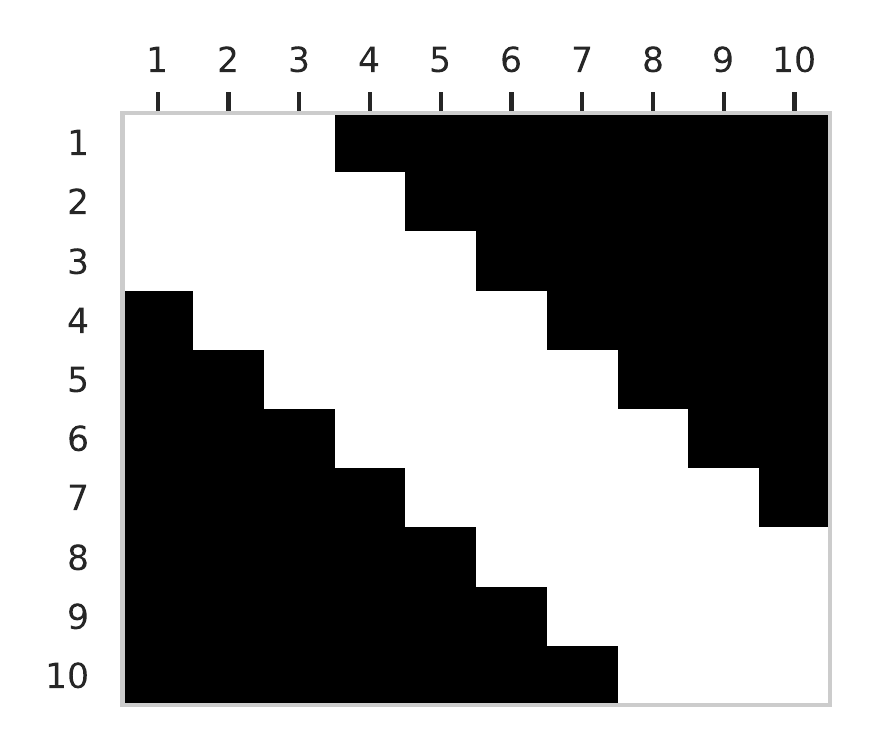}
    \caption{dense}
    \end{subfigure}
\caption{True graph parameter $\delta_0$ (black=0, white=1) in the sparse, random, and dense settings of Simulations 4 with $K=10$ dimensions.}
\label{fig:graphs_K10}
\end{figure}

In this experiment, we evaluate Algorithm \ref{alg:2step_adapt_cavi} on different settings of the graph parameter $\delta_0$, namely a sparse, a random, and a dense settings, illustrate in  Figure \ref{fig:graphs_K10}. The sparse setting is similar to the previous section, while the random setting corresponds to a slightly less sparse regime where additional edges are present in $\delta_0$. Note that these three settings have different numbers of edges in $\delta_0$, therefore, different numbers of non-null interaction functions to estimate. From Table \ref{tab:simu_4_data_graph}, we also note that there are more events and less global excursions in the dense setting that in the two other ones, in particular, in the \emph{Excitation} scenario where this number drops to 2. 

\begin{figure}[hbt!]
    \centering
    \begin{subfigure}[b]{0.49\textwidth}
    \includegraphics[width=\textwidth, trim=0.cm 0.cm 0cm  0.7cm,clip]{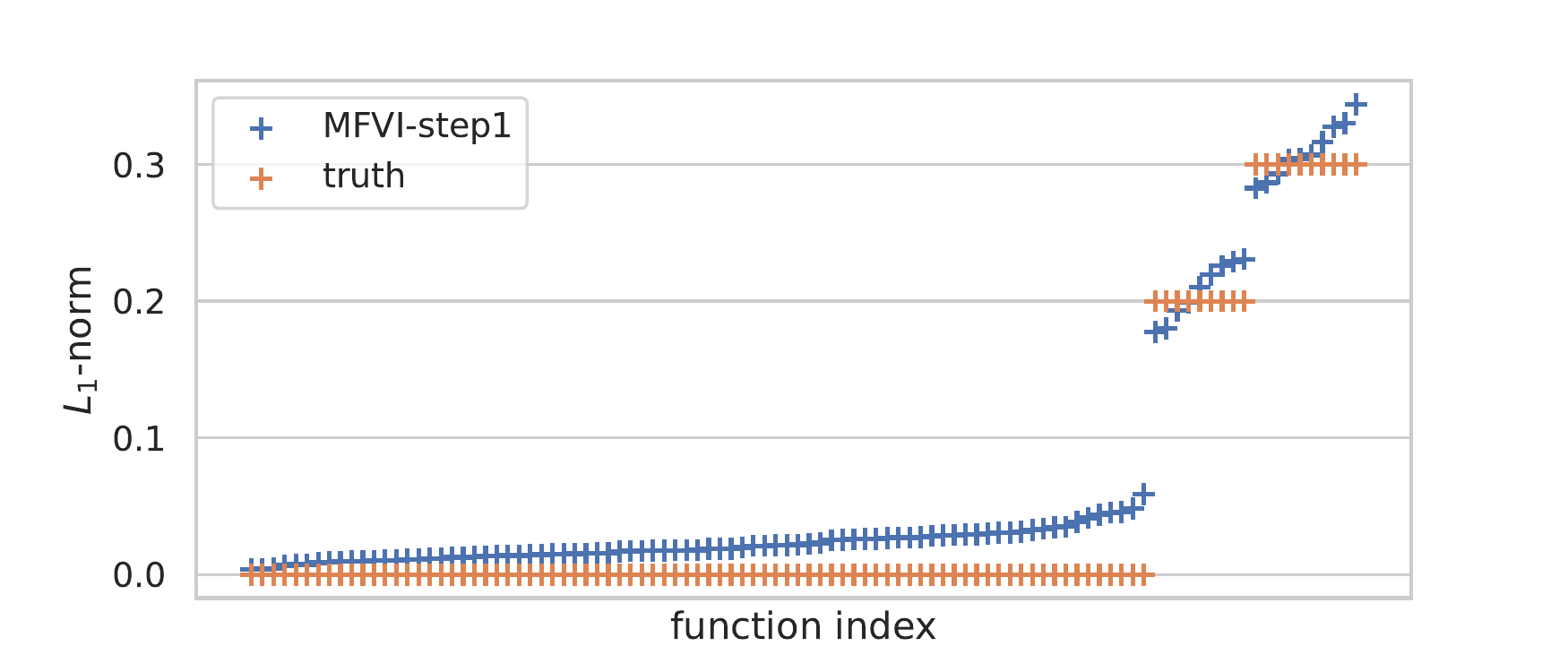}
    \caption{\emph{Excitation} - sparse}
    \end{subfigure}%
        \begin{subfigure}[b]{0.49\textwidth}
    \includegraphics[width=\textwidth, trim=0.cm 0.cm 0cm  0.6cm,clip]{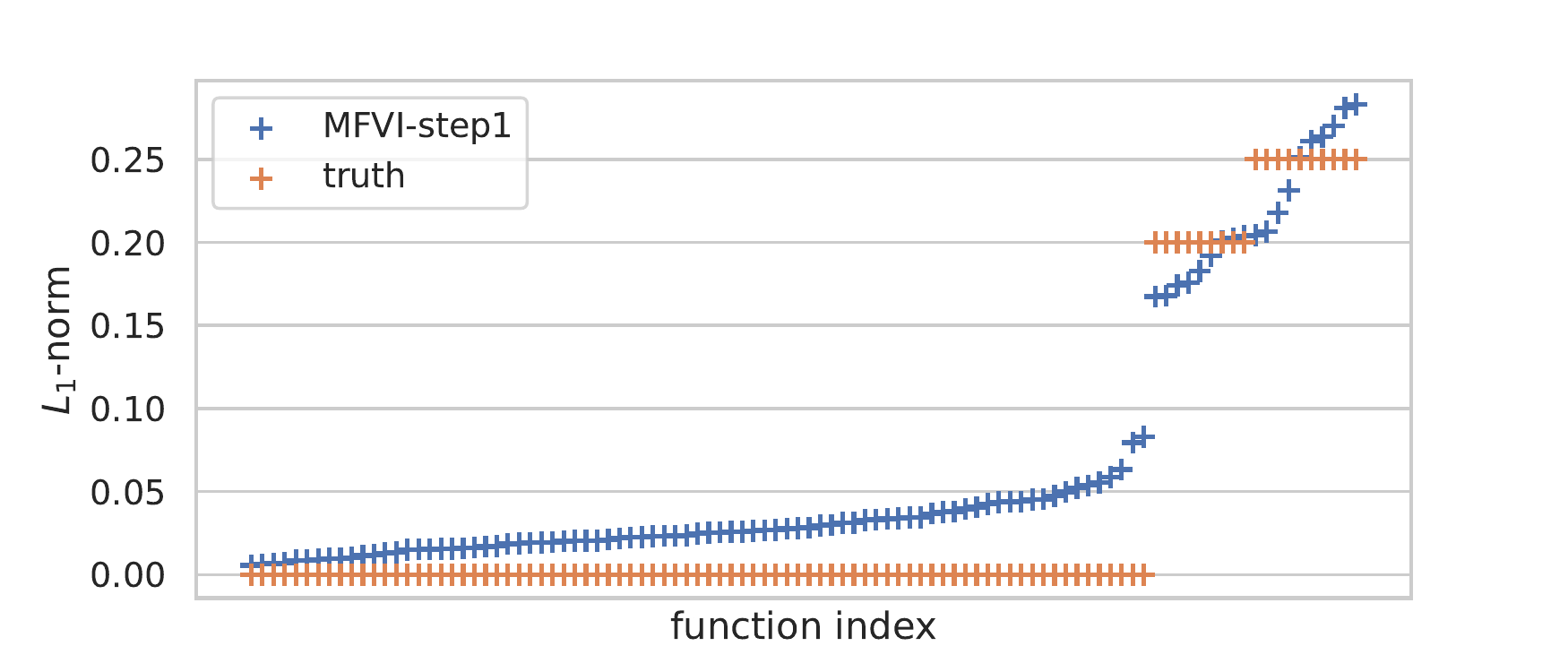}
    \caption{\emph{Inhibition} - sparse}
    \end{subfigure}
        \begin{subfigure}[b]{0.49\textwidth}
    \includegraphics[width=\textwidth, trim=0.cm 0.cm 0cm  .7cm,clip]{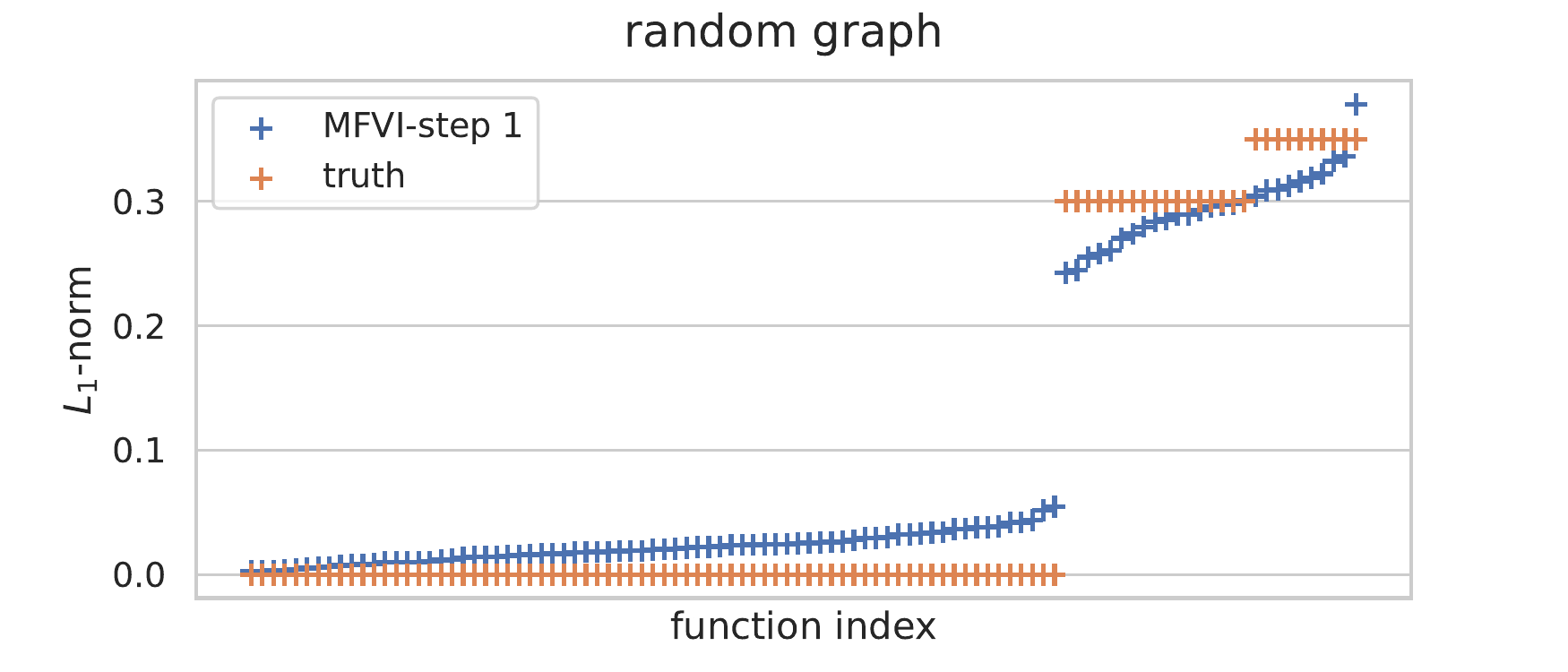}
    \caption{\emph{Excitation} - random}
    \end{subfigure}%
            \begin{subfigure}[b]{0.49\textwidth}
    \includegraphics[width=\textwidth, trim=0.cm 0.cm 0cm  .7cm,clip]{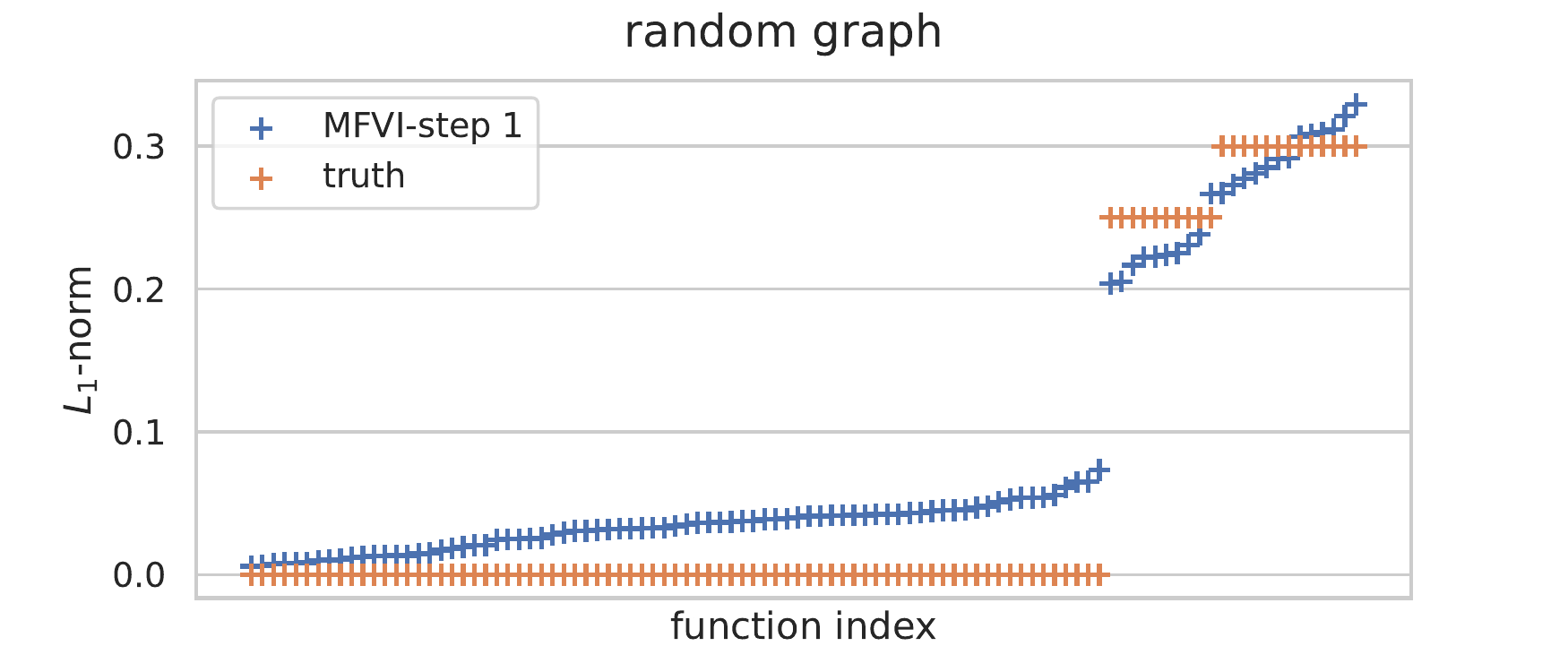}
    \caption{\emph{Inhibition} - random}
    \end{subfigure}
     \begin{subfigure}[b]{0.49\textwidth}
    \includegraphics[width=\textwidth, trim=0.cm 0.cm 0cm  .7cm,clip]{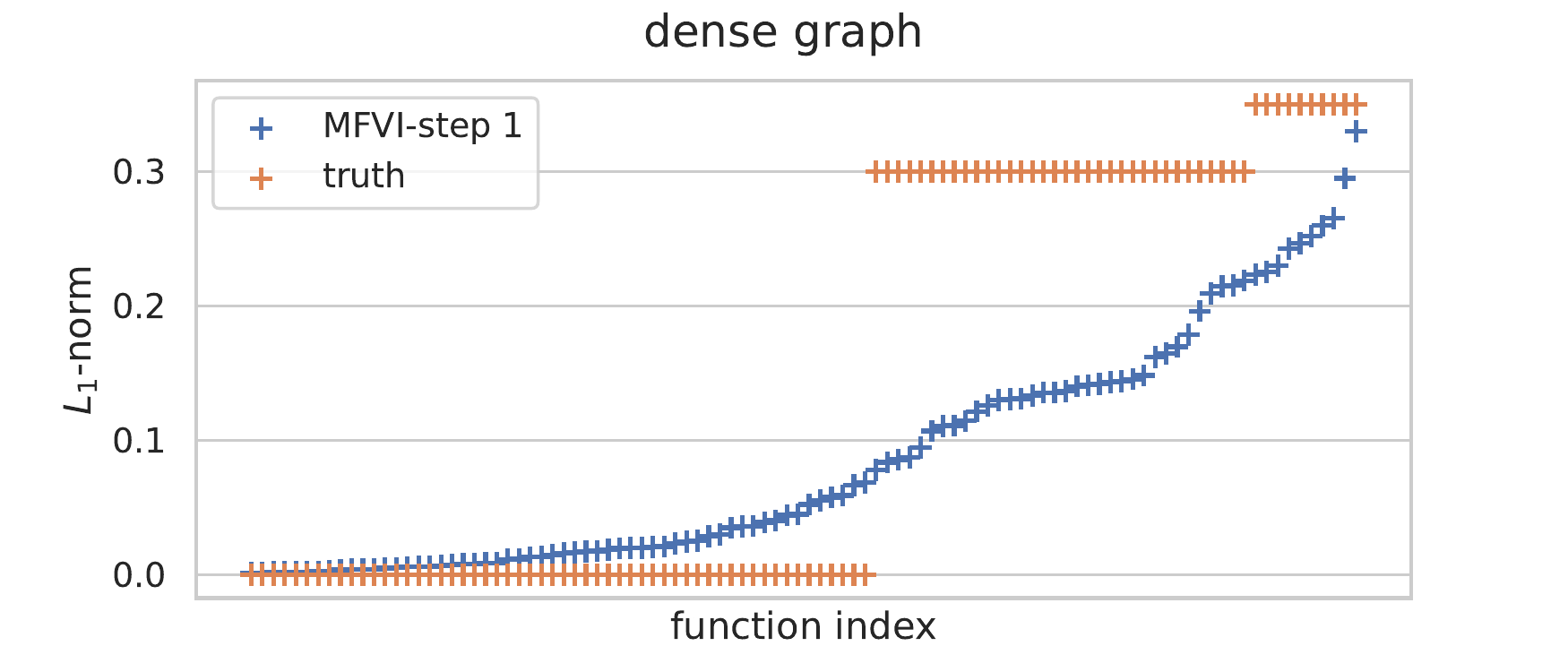}
    \caption{\emph{Excitation} - dense}
    \end{subfigure}%
         \begin{subfigure}[b]{0.49\textwidth}
    \includegraphics[width=\textwidth, trim=0.cm 0.cm 0cm  .7cm,clip]{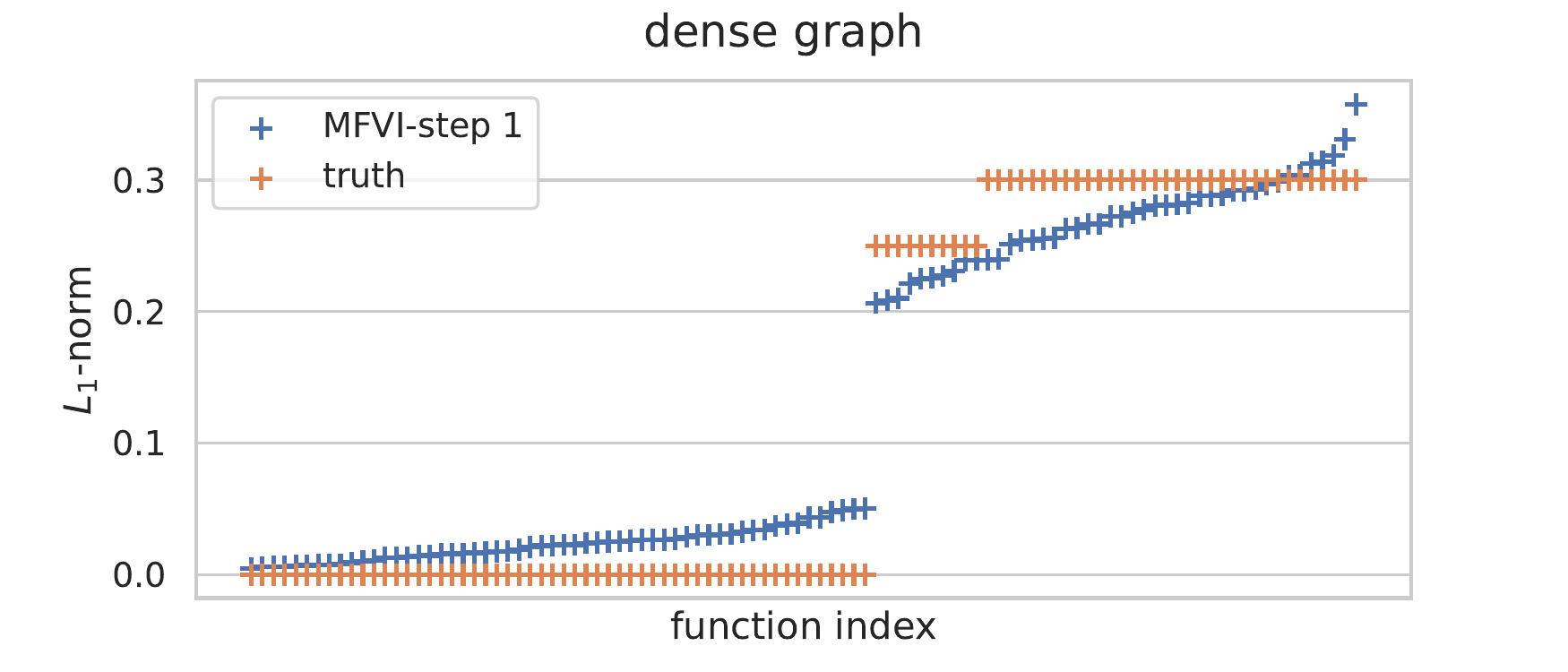}
    \caption{\emph{Inhibition} - dense}
    \end{subfigure}
\caption{Estimated $L_1$-norms using the model-selection variational posterior obtained after the first step of Algorithm \ref{alg:2step_adapt_cavi}, plotted in increasing order, in the different graph settings (sparse, random, and dense $\delta_0$, see Figure \ref{fig:graphs_K10}) and scenarios of Simulation 4 with $K=10$. We note that in the dense graph setting, although the norms are not very well estimated after the first step, the gap heuristics still allows to recover the true graph parameter.}
\label{fig:adaptive_VI_2step_norms_threshold_D10_graphs}
\end{figure}

Our numerical results in Table \ref{tab:perf_simu4_graph} show that in the dense setting, the graph accuracy of our estimator is slightly worse, and the risk of the variational posterior is much higher than in the other settings. We conjecture that this loss of performance is related to the smaller number of global excursions, which leads to a more difficult estimation problem. We can also see from Figure \ref{fig:adaptive_VI_2step_norms_threshold_D10_graphs} that in this particular setting, the estimation of the norms of the interaction functions is deteriorated, and the gap that allows to discriminate between the null and non-null functions is not present anymore. Nonetheless, in the \emph{Inhibition} scenario, for which the number of global excursions is not too small, this phenomenon does not happen and the estimation is almost equivalent in all graph settings. 

\begin{figure}[hbt!]
    \centering    \begin{subfigure}[b]{0.3\textwidth}
    \includegraphics[width=\textwidth, trim=0.cm 0.cm 0cm  0.cm,clip]{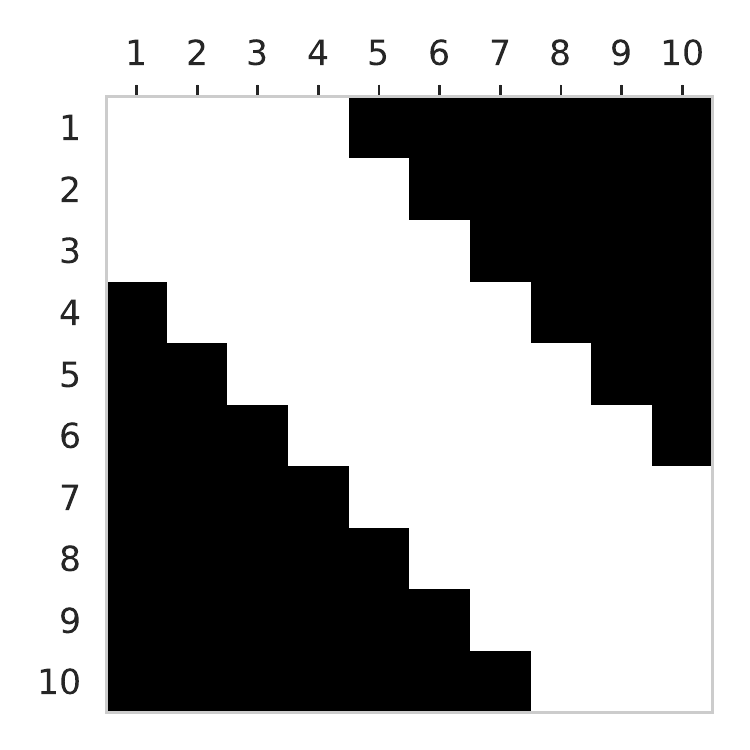}
    \caption{True graph $\delta_0$}
    \end{subfigure}%
    \begin{subfigure}[b]{0.55\textwidth}
    \includegraphics[width=\textwidth, trim=0.cm 0.cm 0cm  0.7cm,clip]{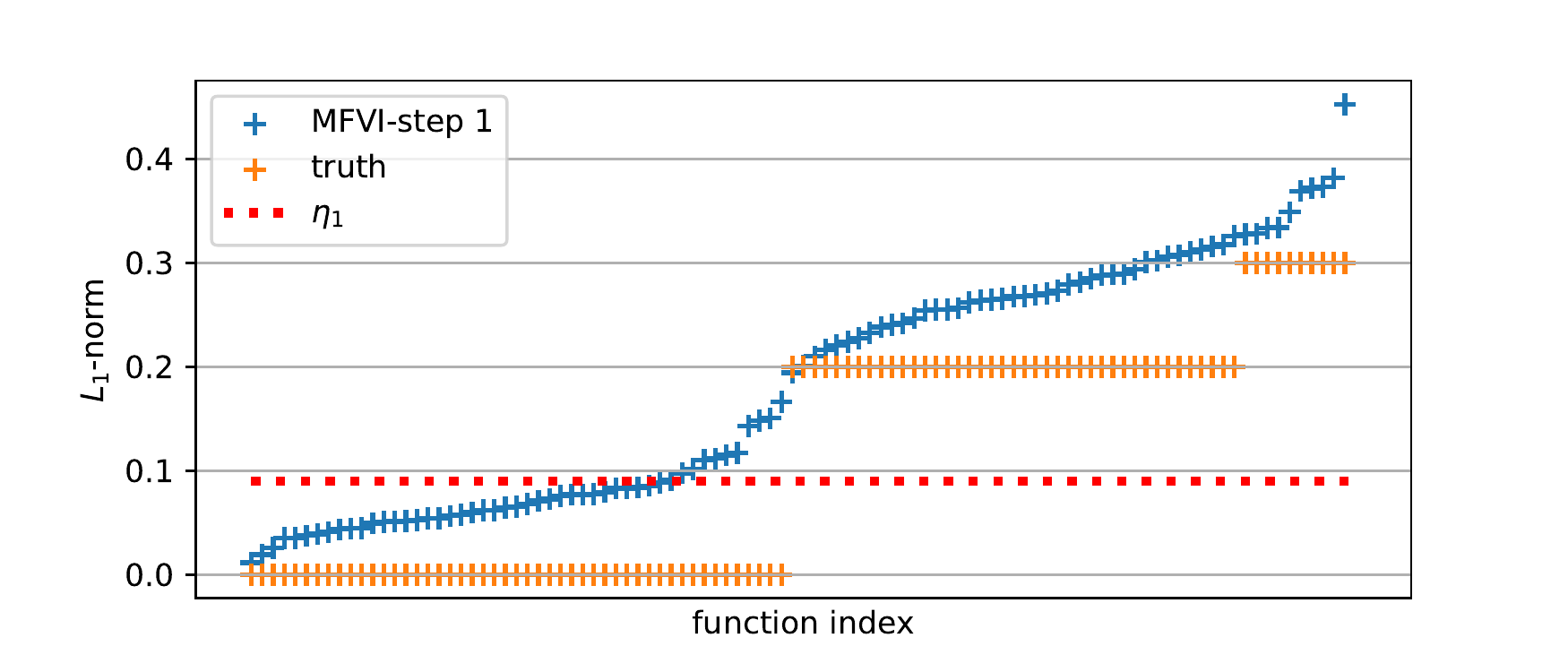}
    \caption{Step 1}
    \end{subfigure}
        \begin{subfigure}[b]{0.55\textwidth}
    \includegraphics[width=\textwidth, trim=0.cm 0.cm 0cm  0.7cm,clip]{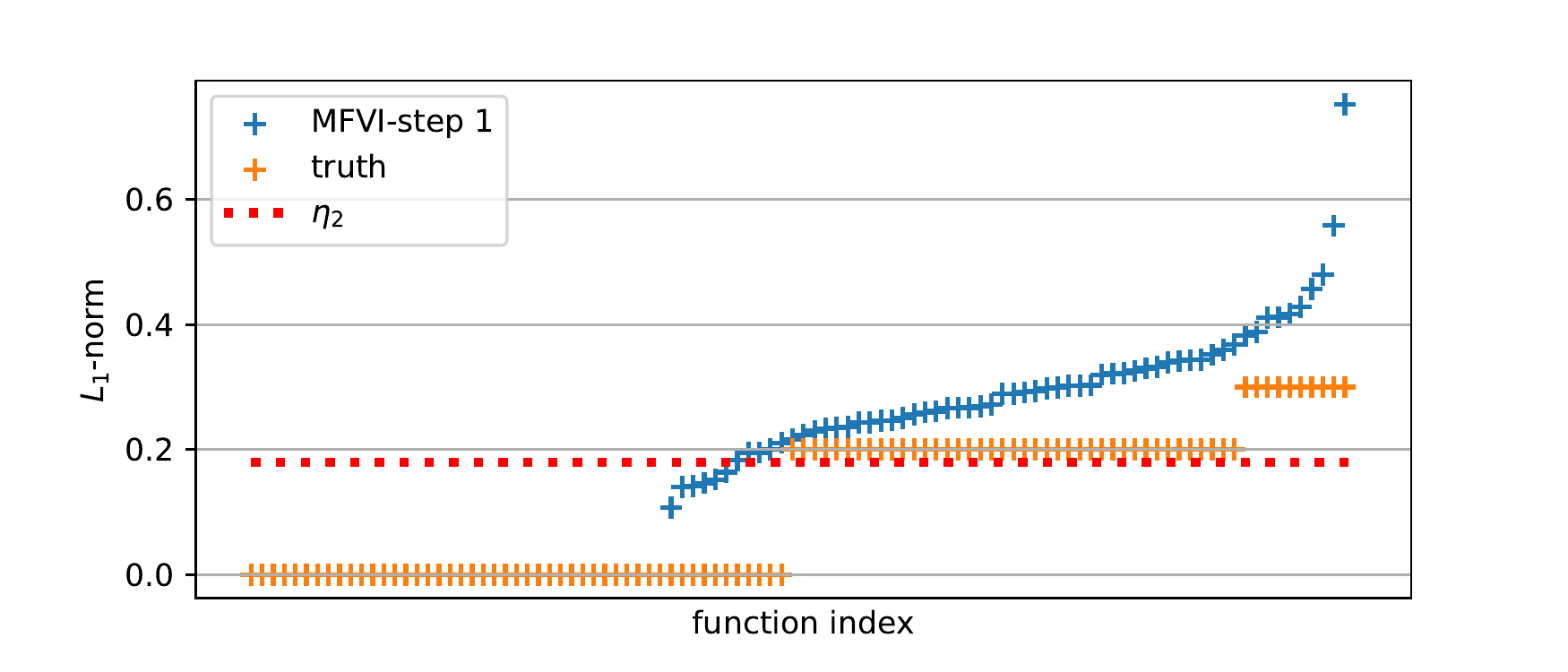}
    \caption{Step 2}
    \end{subfigure}%
        \begin{subfigure}[b]{0.55\textwidth}
    \includegraphics[width=\textwidth, trim=0.cm 0.cm 0cm  0.7cm,clip]{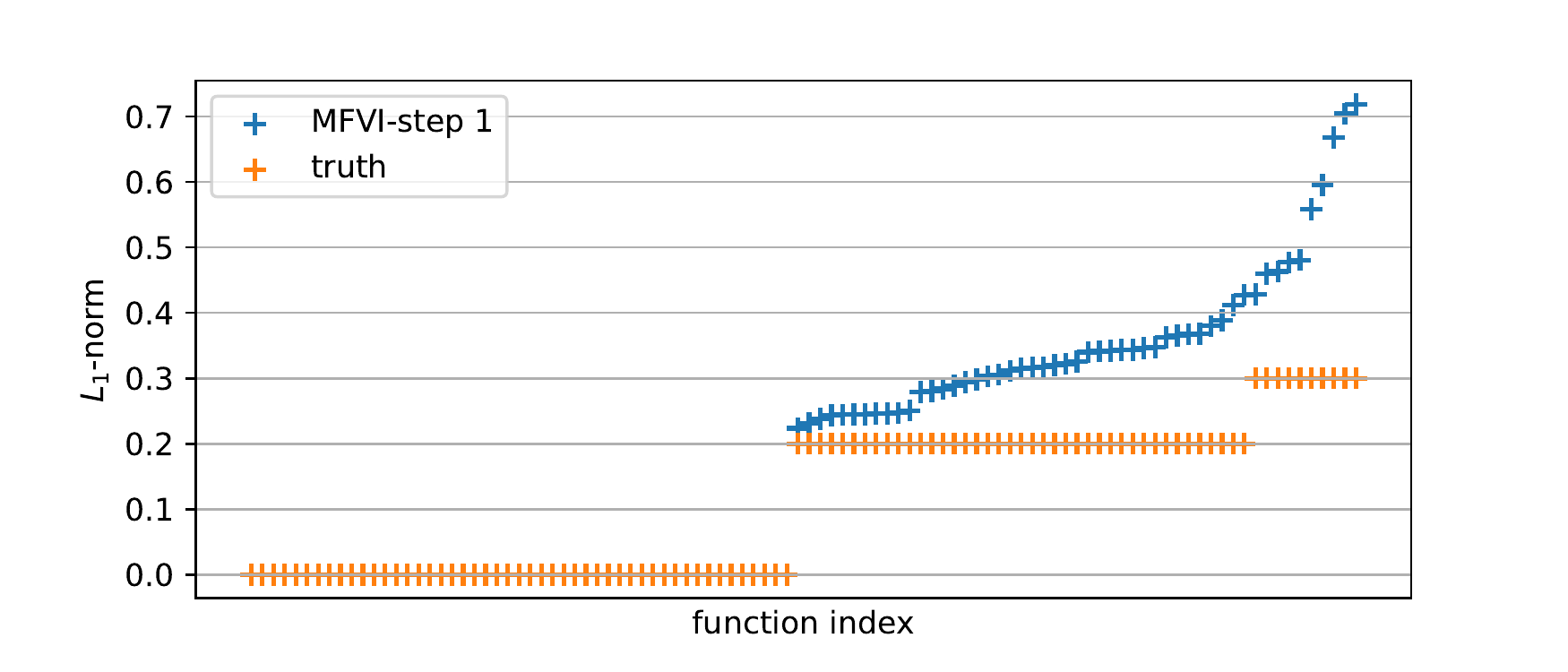}
    \caption{Step 3}
    \end{subfigure}
\caption{True graph $\delta_0$ and estimated $L_1$-norms using the model selection adaptive variational posterior obtained after each step of our three-step procedure, proposed for the dense graph setting of Simulation 4. In Step 1 and Step 2, we plot the data-driven thresholds $\eta_1$ and $\eta_2$, chosen with a ``slope change" heuristics.}
\label{fig:adaptive_VI_3step_norms}
\end{figure}

To further explore the applicability of our thresholding approach in the dense setting, we test the following three-step approach in the \emph{Excitation} scenario, with  $K=10$ and a dense graph $\delta_0$:
\begin{itemize}
    \item The first step is similar to the one of our two-step procedure, i.e., we estimate an adaptive variational posterior distribution within models that contain the complete graph $\delta_C$.
    
    Then, if there is no significant gap in the variational posterior mean estimates of the $L_1$-norms, we look for a (conservative) threshold $\eta_1$ corresponding to the first ``slope change", and estimate a (dense) graph $\hat \delta$.
    \item In a second step, we compute the adaptive variational posterior distribution within models that contain $\hat \delta$ and re-estimate the $L_1$-norms of the functions.
    
    If we now see a significant gap in the norms estimates, we choose a second threshold within that gap; otherwise, we look again for a slope change and pick a conservative threshold $\eta_2$ to compute a second graph estimate  $\hat \delta_2$.
    \item In the third and last step, we repeat the second step with now our second graph estimate, $\hat \delta_2$.
\end{itemize}
In Figure \ref{fig:adaptive_VI_3step_norms}, we plot our estimates of the norms after each step of the previous procedure. In this case, we have chosen visually the threshold $\eta_1=0.09$ and $\eta_2=0.18$ after respectively the first and second step, using the slope change heuristics. We note that the previous method indeed provides a conservative graph estimate in the first step, but in the second step, allows to refine our estimate of the graph and approach the true graph. Besides, we note that the large norms are inflated along the three steps of our procedure. Therefore, our method performs better in sparse settings where a significant gap allows to correctly infer the true graph $\delta_0$.

In conclusion, our simulations in low and high-dimensional settings, with different levels of sparsity in the graph, show that our two-step procedure is able to correctly select the graph parameter and dimensionality of the process in sparse settings, and hence allows to scale up variational Bayes approaches to larger number of dimensions. Nonetheless, from the moderately high-dimensional settings, the estimation of the parameter $f$ becomes sensitive to the difficulty of the problem. In particular, the performance is sensitive to the graph sparsity, tuning the number of non-null functions to estimate, and, as we conjecture, the number of global excursions in the data. Finally, we note that heuristic approaches for the choice of the threshold - needed to  estimate the graph parameter - need to further explored in noisier and denser settings.

\begin{table}[hbt!]
    \centering
\begin{tabular}{c|c|c|c|c|c|c|c}
\toprule
  K & Scenario & T & \# events & \# excursions & \# local excursions & computing time (s) \\
 \midrule
 \multirow{2}{*}{2}  & Excitation & 500 & 5680 &       2416 &           1830 &             19   \\
 & Inhibition & 700 & 4800 &       2416 &           1830 &             18 \\
 \midrule
 \multirow{2}{*}{4}  & Excitation & 500 & 11338 &       2378 &           1878 &             41 \\
 & Inhibition & 700 & 9895 &       2378 &           1878 &             39 \\
 \midrule
 \multirow{2}{*}{8}  & Excitation & 500 & 22514 &       1207 &           1857 &            151 \\
 & Inhibition & 700 & 19746 &       1207 &           1857 &            134 \\
 \midrule
 \multirow{2}{*}{16}  & Excitation & 500 & 51246 &        200 &           1784 &            577 \\
 & Inhibition & 700 & 37166 &        200 &           1784 &            494 \\
 \midrule
 \multirow{2}{*}{32}  & Excitation & 500 & 96803 &          4 &           1824 &           2147 \\
 & Inhibition & 700 & 76106 &          4 &           1824 &           1386 \\
  \midrule
  \multirow{2}{*}{64}  & Excitation & 200 & 117862 &          0 &           1481 &           8176 \\
 & Inhibition & 300 & 133200 &          0 &           1481 &           7583 \\
\bottomrule
\end{tabular}
    \caption{Number of observed events, excursions, and computing times of Algorithm \ref{alg:2step_adapt_cavi} in the multivariate settings of Simulation 4.}
    \label{tab:simu_4_data}
\end{table}

% \begin{table}
% \begin{tabular}{c|c|c|c}
% \toprule
%   K & Scenario & \# points & computing time (s) \\
%  \midrule
%  \multirow{2}{*}{2}  & Excitation & 5680   &           19   \\
%  & Inhibition & 4800 &                   18 \\
%  \midrule
%  \multirow{2}{*}{4}  & Excitation & 11338 &         41 \\
%  & Inhibition & 9895 &                   39 \\
%  \midrule
%  \multirow{2}{*}{8}  & Excitation & 22514 &             151 \\
%  & Inhibition & 19746 &               134 \\
%  \midrule
%  \multirow{2}{*}{16}  & Excitation & 51246 &              577 \\
%  & Inhibition & 37166 &            494 \\
%  \midrule
%  \multirow{2}{*}{32}  & Excitation & 96803 &           2147 \\
%  & Inhibition & 76106 &                1386 \\
%   \midrule
%   \multirow{2}{*}{64}  & Excitation & 117862 &               8176 \\
%  & Inhibition & 133200 &              7583 \\
% \bottomrule
% \end{tabular}
% \end{table}

\begin{table}[hbt!]
    \centering
\begin{tabular}{c|c|c|c|c|c|c}
\toprule
   Scenario & Graph  & \# Edges & \# Events & \# Excursions & \# Local excursions \\
 \midrule
 \multirow{3}{*}{Excitation} &  Sparse & $2K-1$ & 24638 &        431 &           1212 \\
 &  Random & $3K-1$ &  27475    &        398 &           1262 \\
 &  Dense &  $5K-6$ & 90788 &          2 &           1432 \\
 \midrule
 \multirow{2}{*}{Inhibition}  & Sparse  & $2K-1$ &  22683 &        911 &           1778 \\
 & Random & $3K-1$ &  24031 &        884 &           1834 \\
 & Dense &  $5K-6$ &  35291 &        547 &           2170 \\
\bottomrule
\end{tabular}
    \caption{Number of edges, observed events, and excursions in the different graph settings of Simulation 4 ($K=10$).}
    \label{tab:simu_4_data_graph}
\end{table}

\begin{table}[hbt!]
    \centering
\begin{tabular}{c|c|c|c|c}
\toprule
  \# dimensions & Scenario & Graph accuracy & Dimension accuracy & Risk  \\
 \midrule
 \multirow{2}{*}{2}  & Excitation & 1.00 & 1.00 &  0.79 \\
 & Inhibition &1.00 & 1.00 &  0.35 \\
\midrule
 \multirow{2}{*}{4} &  Excitation  & 1.00 & 1.00 &  1.01 \\
 & Inhibition &  1.00 & 1.00 &  0.92 \\
 \midrule
 \multirow{2}{*}{8} &  Excitation  & 1.00 & 1.00 &  2.10 \\
 & Inhibition &   1.00 & 1.00 &  2.12 \\
 \midrule
 \multirow{2}{*}{16} &  Excitation  & 1.00 & 1.00 &  5.77 \\
 & Inhibition &   1.00 & 1.00 &  4.48 \\
 \midrule
 \multirow{2}{*}{32} &  Excitation  & 1.00 & 0.97 & 10.57 \\
 & Inhibition &  1.00 & 1.00 &  8.53 \\
 \midrule
  \multirow{2}{*}{64} &  Excitation  & 1.00 & 1.00 & 23.74 \\
 & Inhibition &  1.00 & 1.00 & 18.43 \\
\bottomrule
\end{tabular}
    \caption{Performance of Algorithm \ref{alg:2step_adapt_cavi} in the multivariate settings of Simulation 4. We report the accuracy of our graph estimate $\hat \delta$ and the selected dimensionality of the interaction functions in the model-selection variational posterior, and the risk on the whole parameter $f$ defined in \eqref{eq:risk}.}
    \label{tab:perf_simu4}
\end{table}

\begin{table}[hbt!]
    \centering
\begin{tabular}{c|c|c|c|c}
\toprule
  Scenario & Graph & Graph accuracy & Dimension accuracy & Risk  \\
 \midrule
 \multirow{3}{*}{Excitation}  & Sparse &  1.00  &  1.00    &     2.91 \\
  & Random &  1.00 &    1.00 &       4.00 \\
  & Dense  & 0.5 &    1.00 &      17.67 \\
\midrule
 \multirow{3}{*}{Inhibition}  & Sparse & 1.00 & 1.00 &  2.62 \\
  & Random &   0.99 &    1.00 &         3.44 \\
  & Dense   &  1.00 &    1.00 &         2.67 \\
\bottomrule
\end{tabular}
    \caption{Performance of Algorithm \ref{alg:2step_adapt_cavi} in the different graph settings of Simulation 4 ($K=10$). We note in that the dense graph setting, there are more parameters to estimate, and therefore non-null terms in the risk metric.}
    \label{tab:perf_simu4_graph}
\end{table}

\begin{figure}[hbt!]
\setlength{\tempwidth}{.20\linewidth}\centering
\settoheight{\tempheight}{\includegraphics[width=\tempwidth, trim=0.cm 0.cm 0.cm  0.cm,clip]{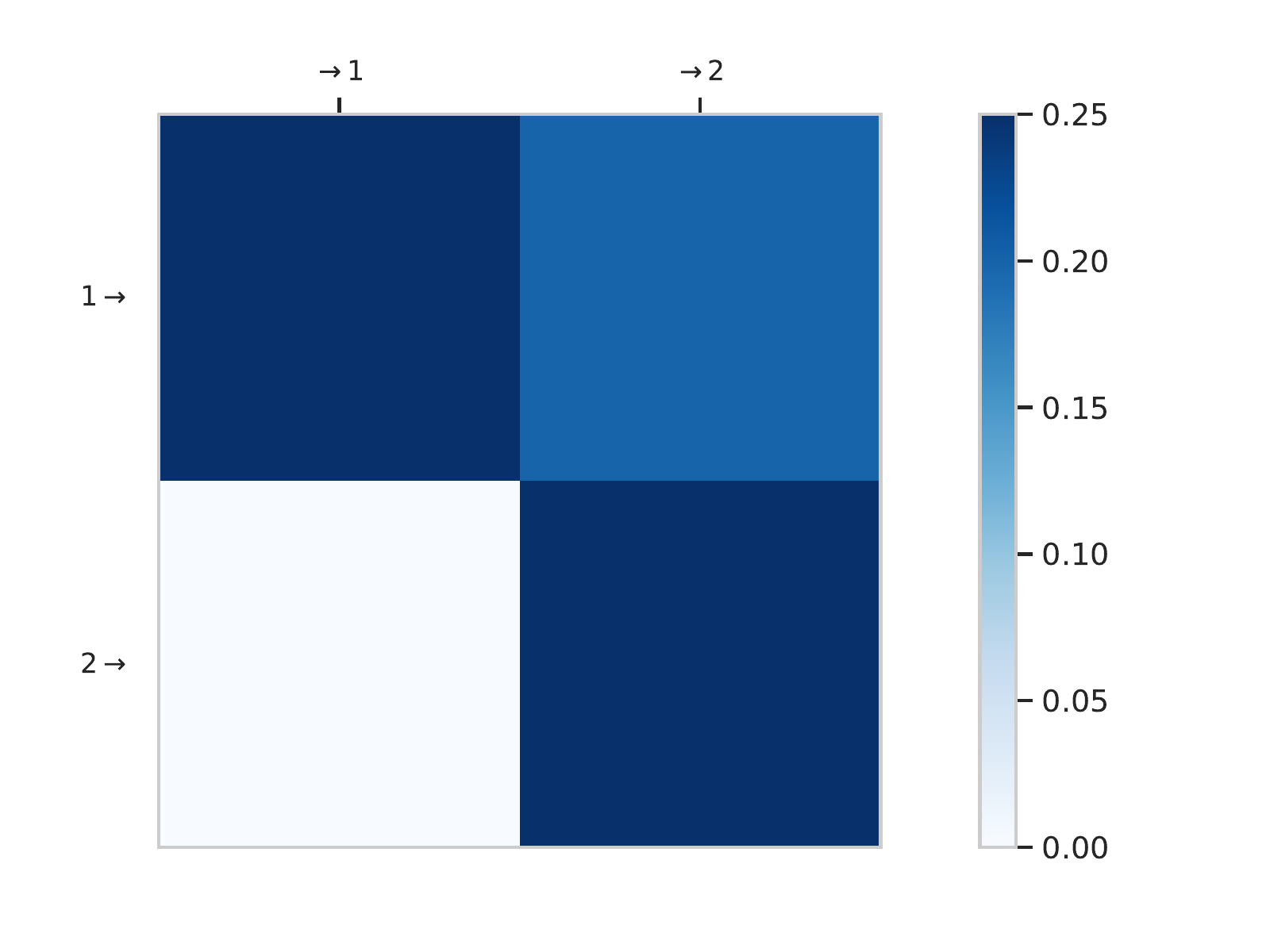}}%
\hspace{-5mm}
\fbox{\begin{minipage}{\dimexpr 24mm} \begin{center} \itshape \large  \textbf{Excitation} \end{center} \end{minipage}}
\hspace{-5mm}
\columnname{$K=2$}
\columnname{$K=4$}
\columnname{$K=16$}
\columnname{$K=64$}\\
\begin{minipage}{\dimexpr 15mm} \vspace{-30mm} \begin{center} \itshape  \textbf{Ground-\\truth} \end{center} \end{minipage}
    \includegraphics[width=\tempwidth, trim=0.cm 0.cm 5.cm  0.cm,clip]{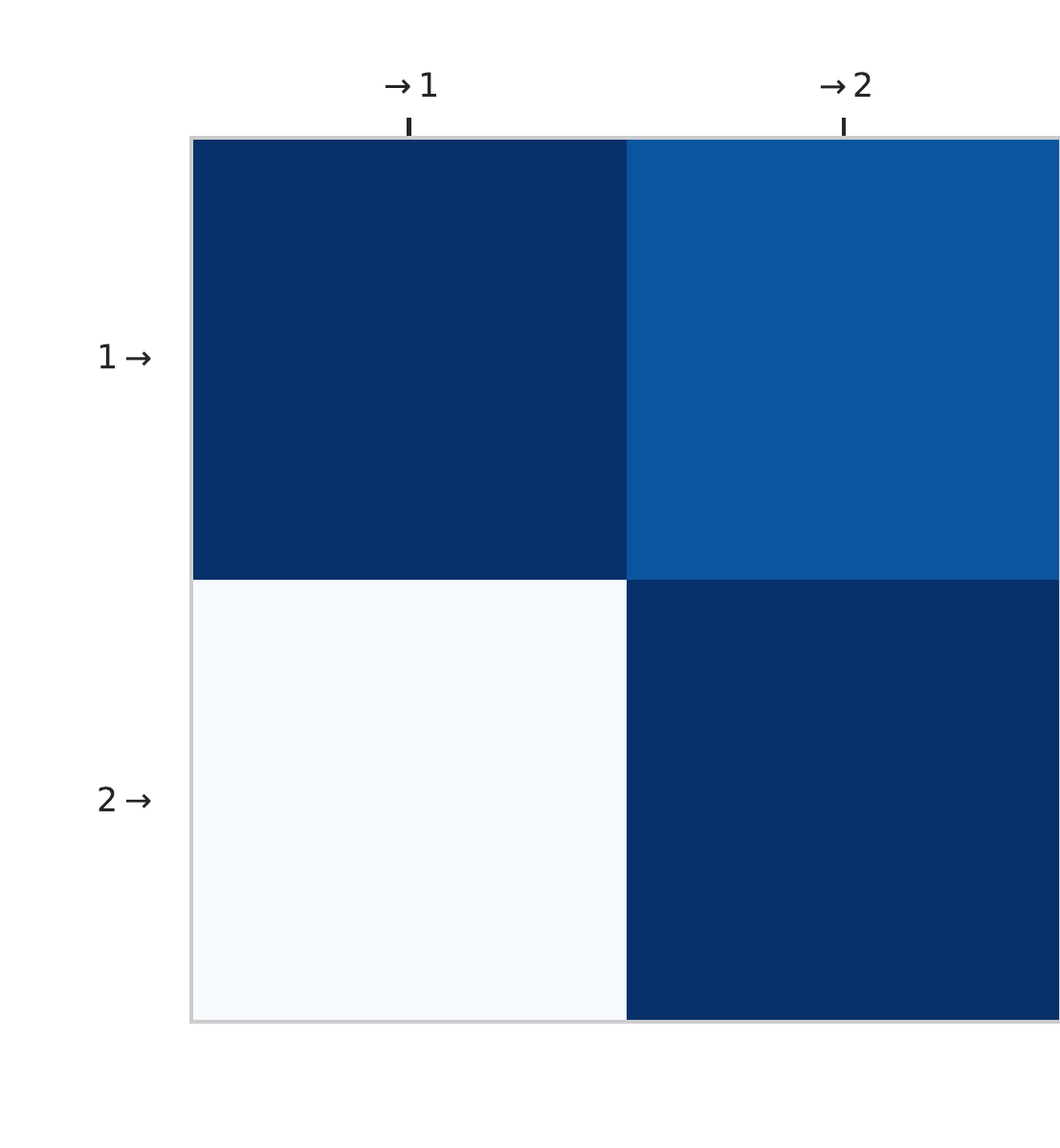}
    \includegraphics[width=\tempwidth, trim=0.cm 0.cm 5cm  0.cm,clip]{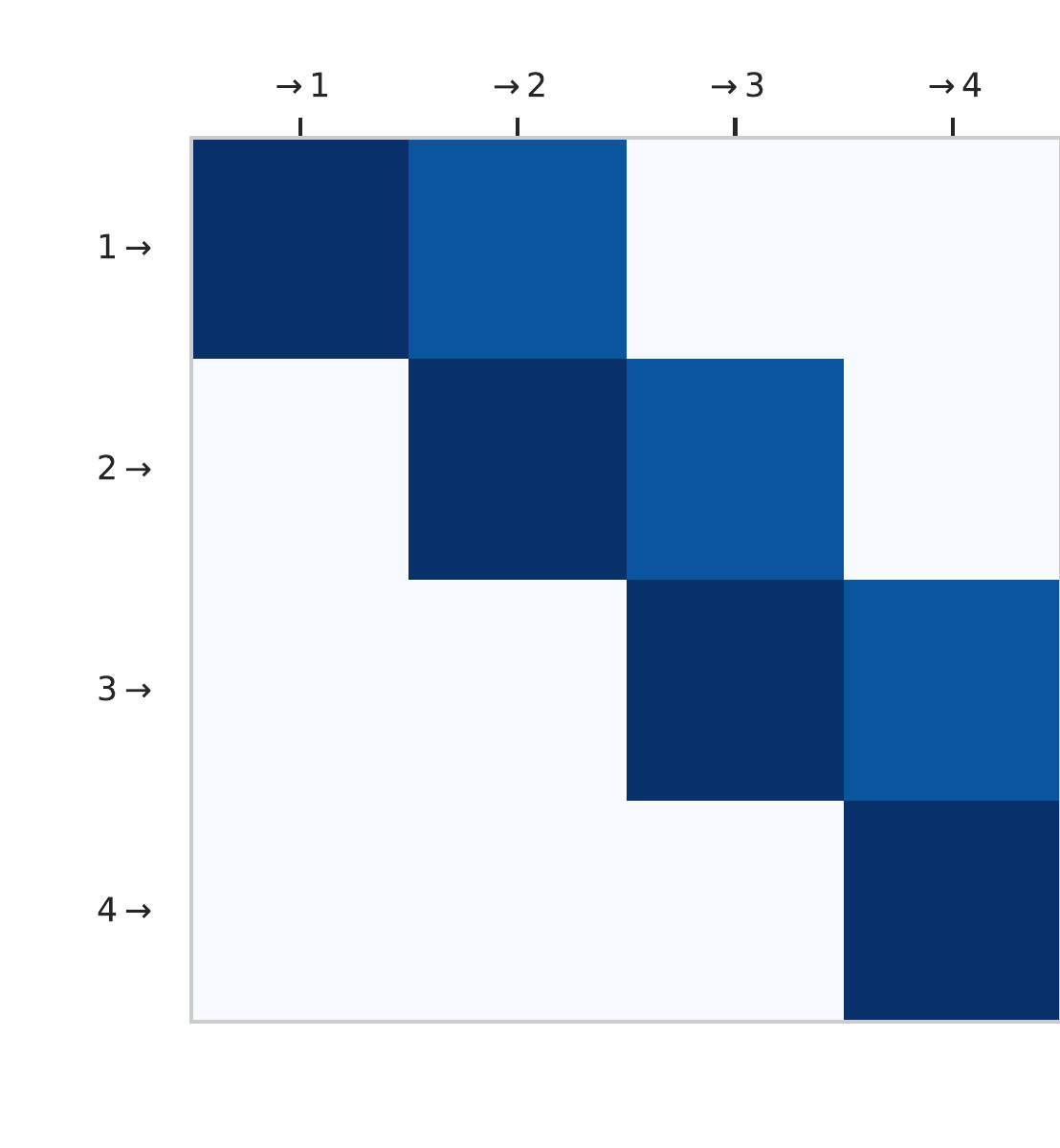}
    \includegraphics[width=\tempwidth, trim=0.cm 0.cm 5cm  0.cm,clip]{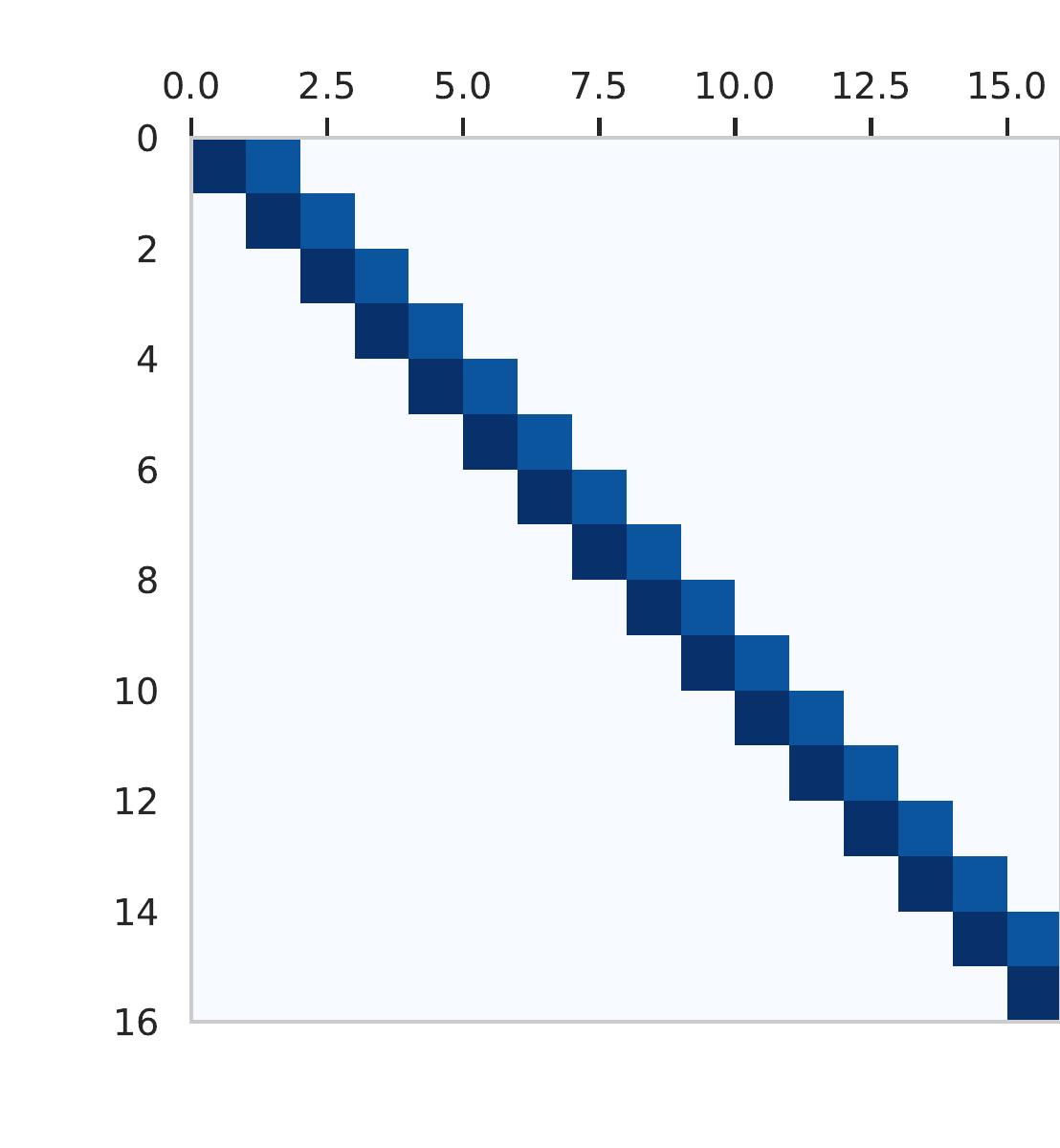}
    \includegraphics[width=1.3\tempwidth, trim=0.cm 0.cm 2cm  0.cm,clip]{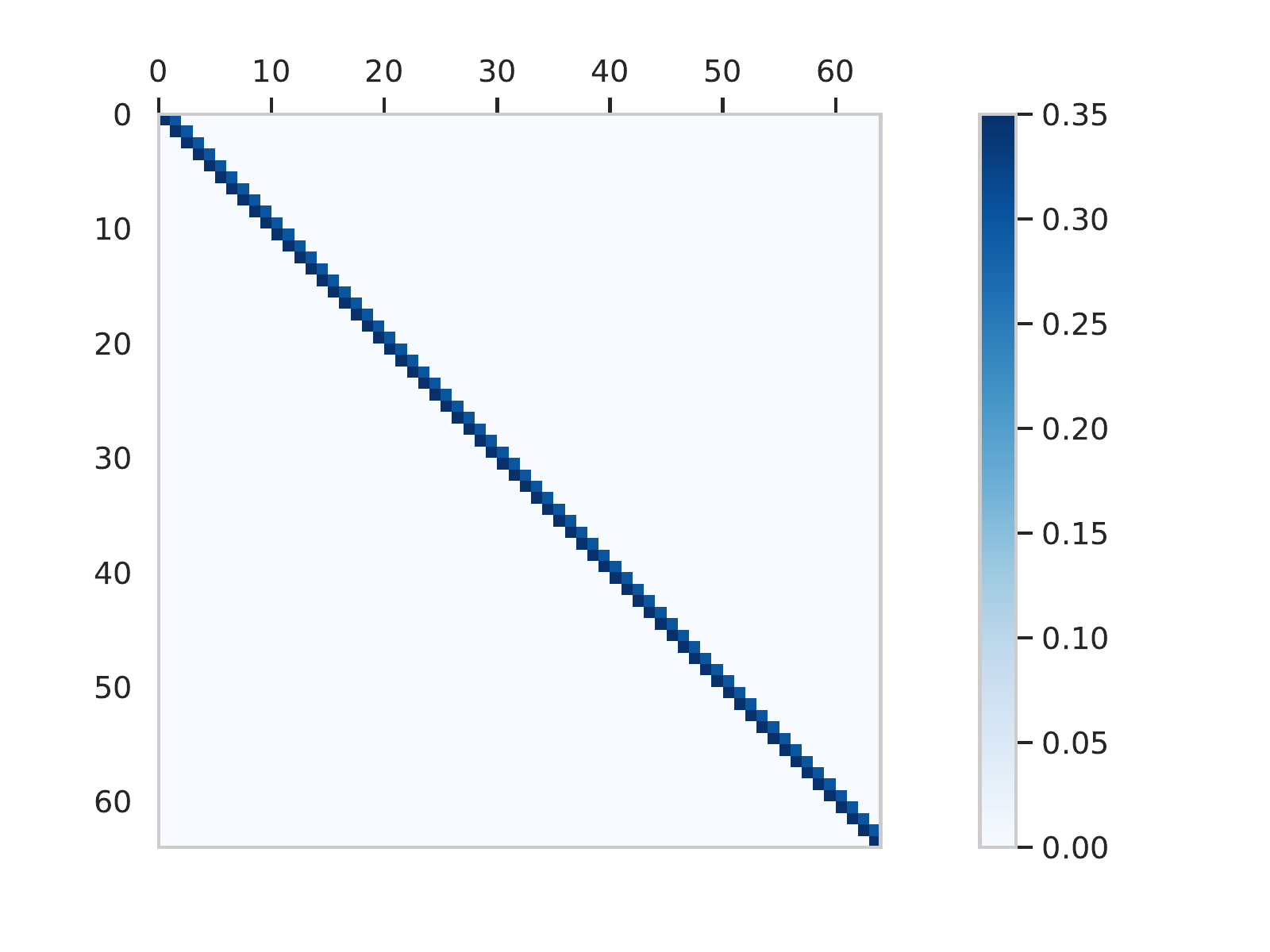}\\
\begin{minipage}{\dimexpr 15mm} \vspace{-30mm} \begin{center} \itshape \textbf{Error} \end{center} \end{minipage}
    \includegraphics[width=\tempwidth, trim=0.cm 0.cm 5cm  0.cm,clip]{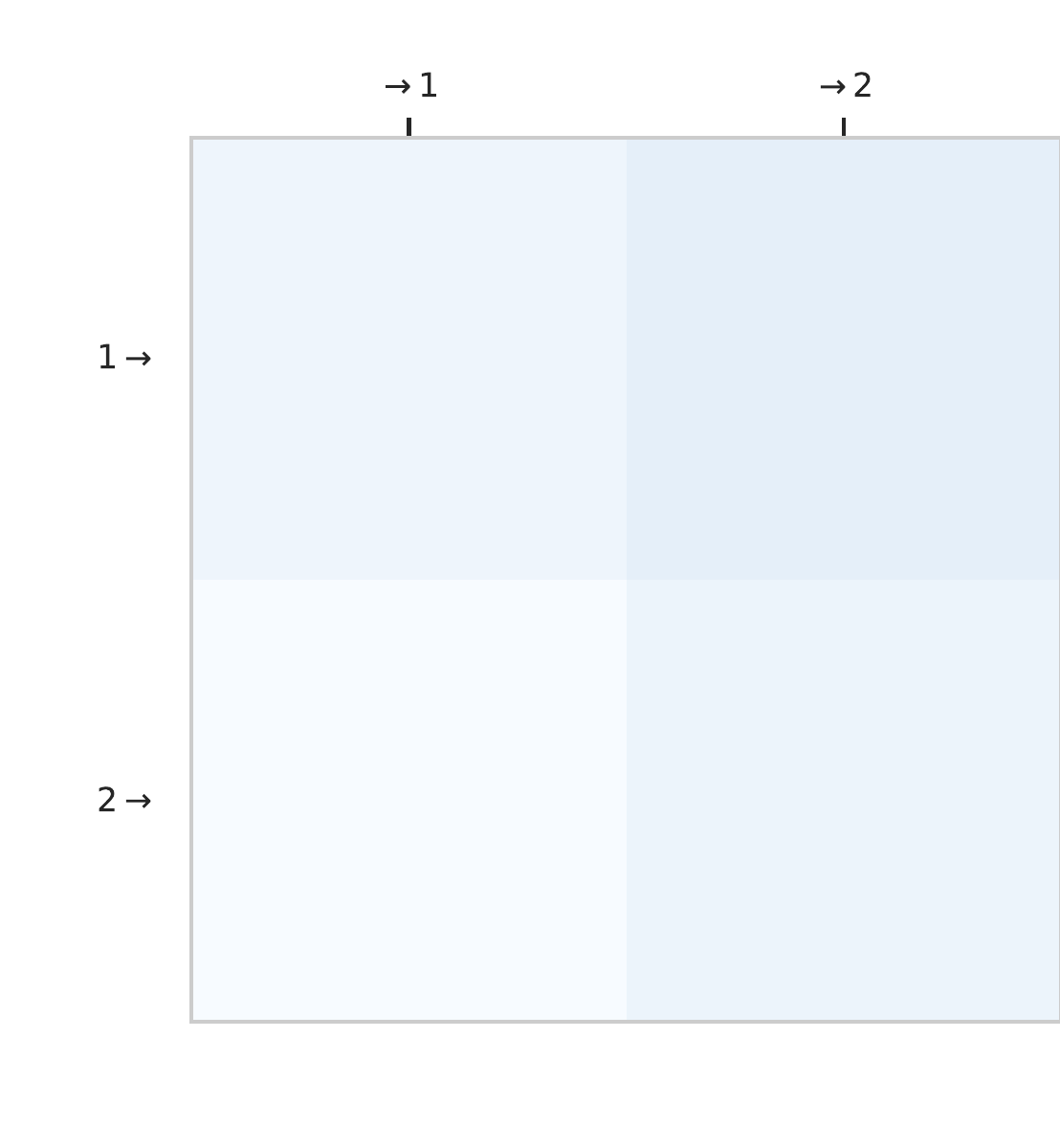}
    \includegraphics[width=\tempwidth, trim=0.cm 0.cm 5cm  0.cm,clip]{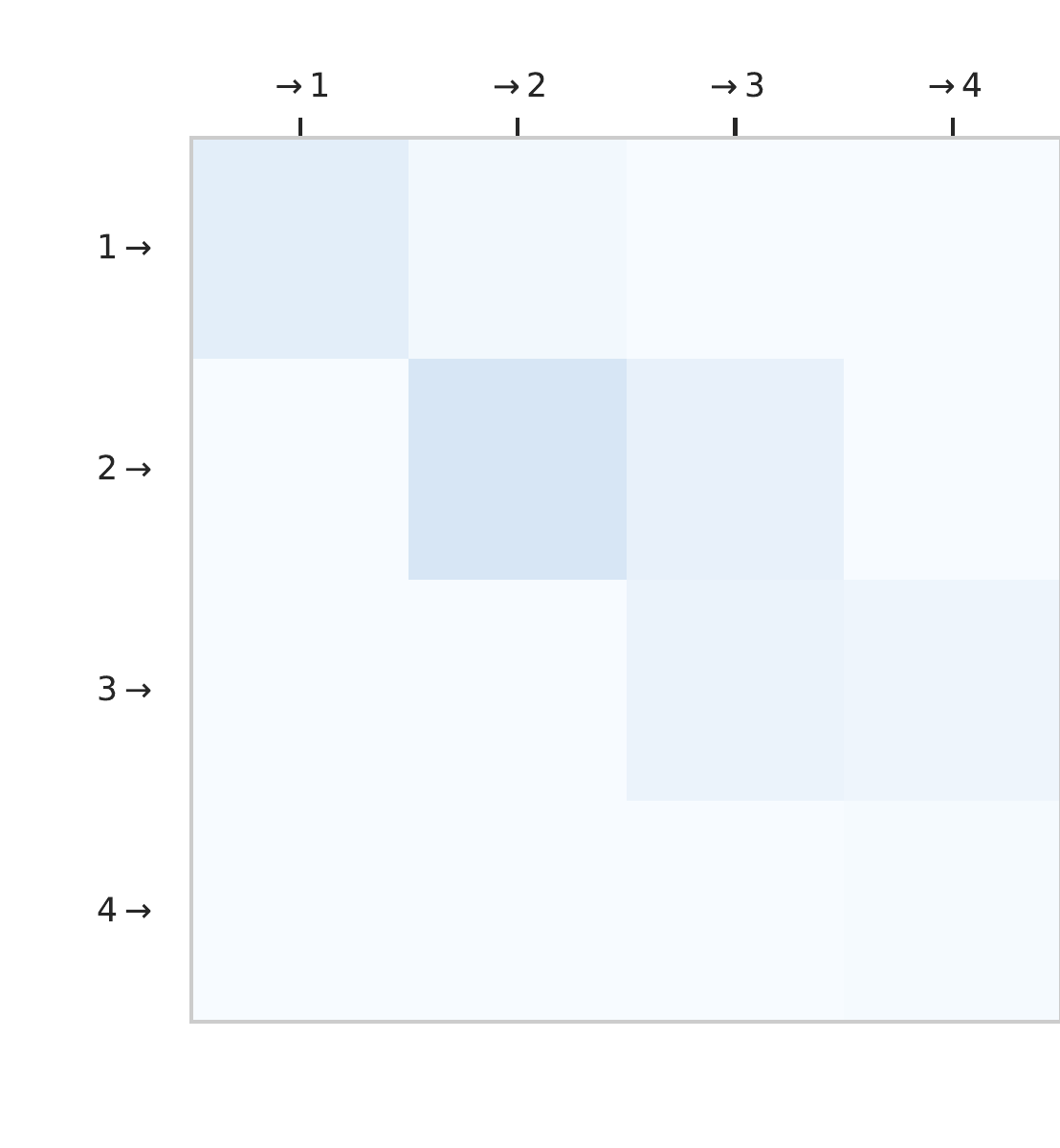}
    \includegraphics[width=\tempwidth, trim=0.cm 0.cm 5cm  0.cm,clip]{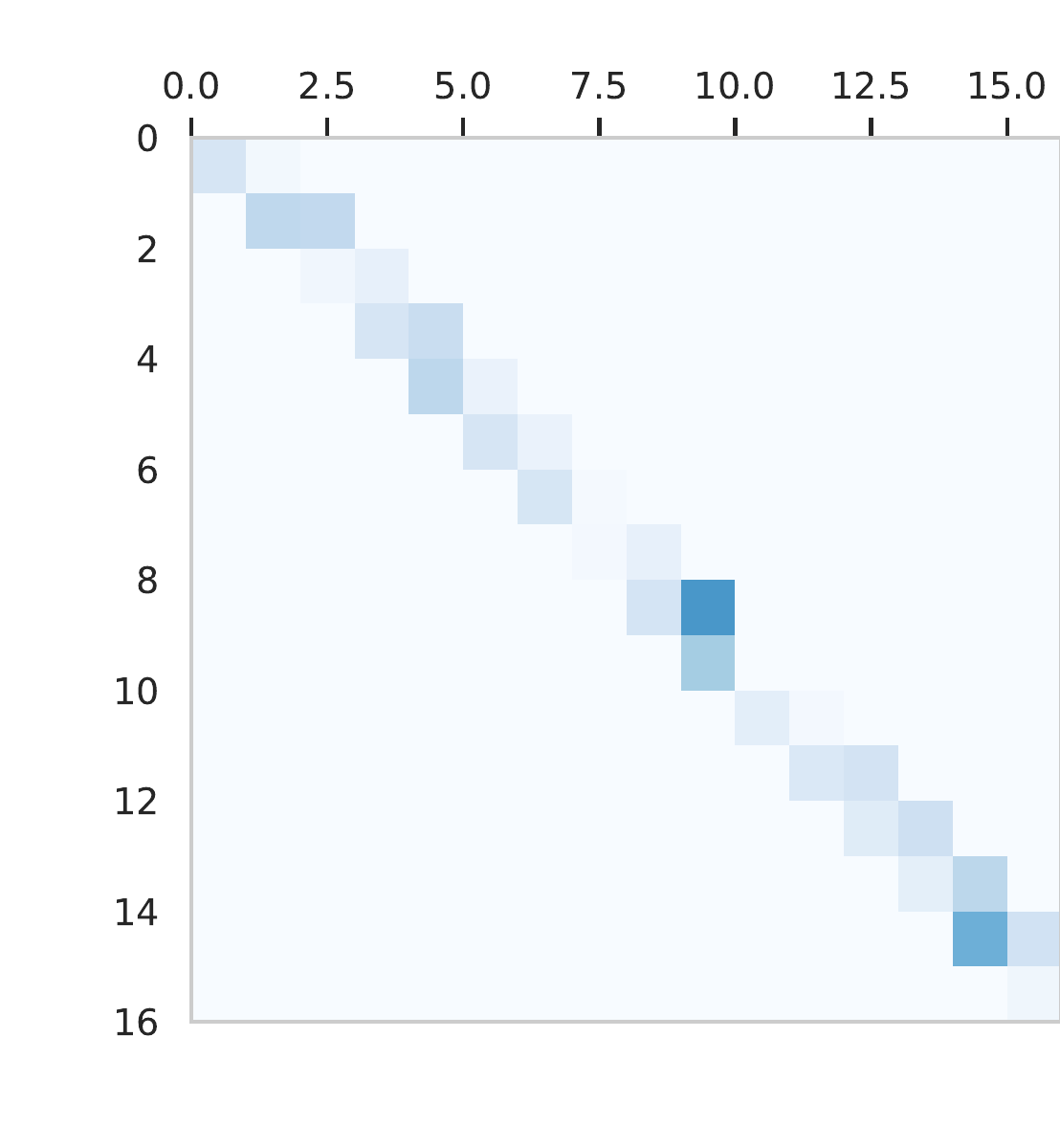} 
    \includegraphics[width=1.3\tempwidth, trim=0.cm 0.cm 2cm  0.cm,clip]{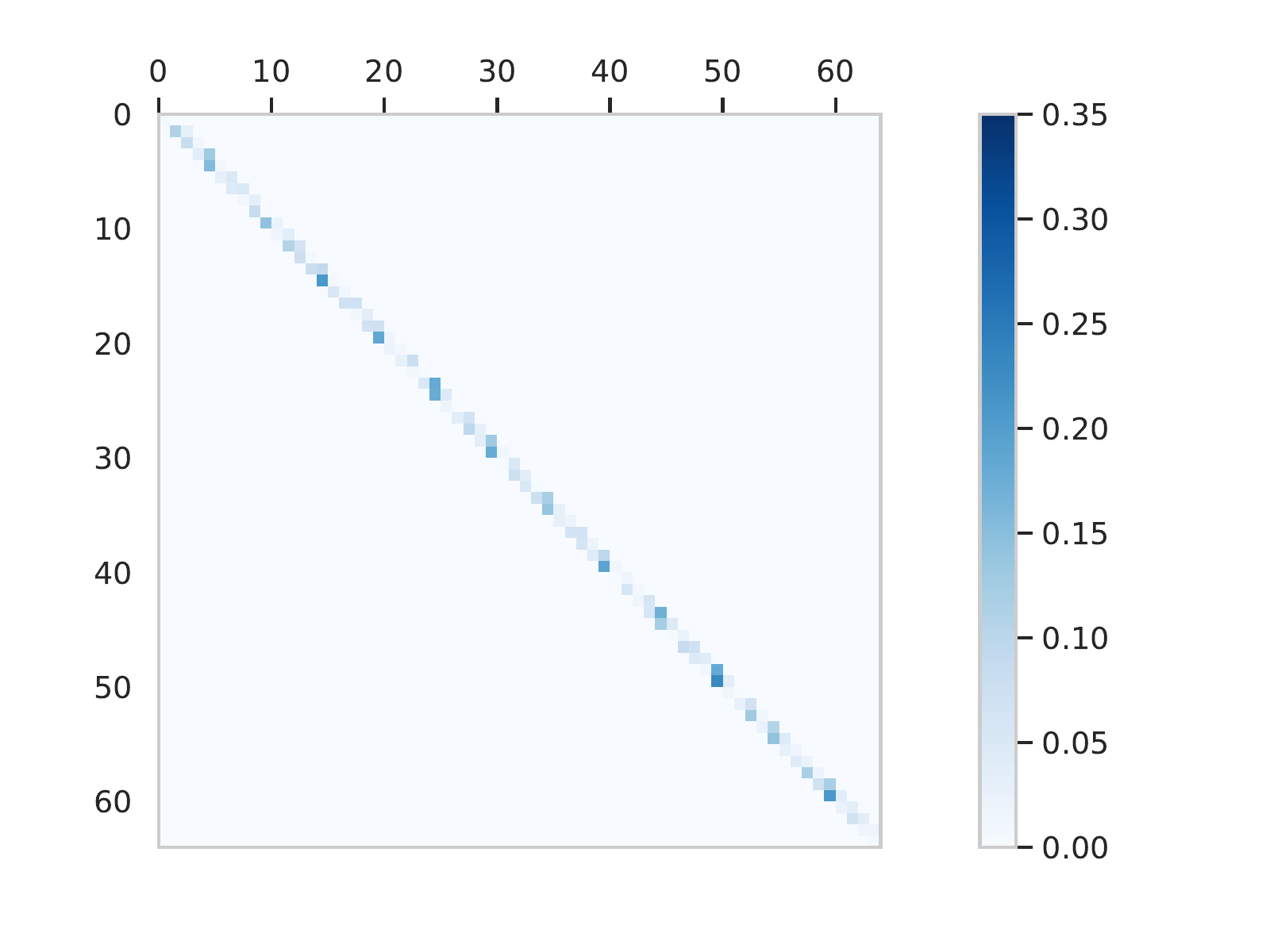}\\
\caption{Heatmaps of the $L_1$-norms of the true parameter $h_0$, i.e., the entries of the matrix $S_0 = (S^0_{lk})_{l,k} = (\norm{h_{lk}^0}_1)_{l,k}$ (left column) and the $L_1$-risk of the model-selection variational posterior obtained with Algorithm \ref{alg:2step_adapt_cavi}, i.e., $(\mathbb{E}^{\hat Q_{MV}}[\norm{h_{lk}^0 - h_{lk}}_1])_{l,k}$ (right column), in the \emph{Excitation} scenario of Simulation 4. The rows correspond to $K=2,4,8,16,32,64$.}
\label{fig:adaptive_VI_2step_norms_exc}
\end{figure}

\FloatBarrier

\subsection{Simulation 5:  Convergence of the two-step variational posterior for varying data set sizes. }

In this experiment, we study the variations of performances of Algorithm \ref{alg:2step_adapt_cavi} with increasing lengths of the observation window, i.e., increasing number of data points. We consider multidimensional data sets with $K=10$, $T \in \{50, 200, 400, 800\}$, the same connectivity graph as in Simulation 4, and an \emph{Excitation} and an \emph{Inhibition} scenarios. The number of events and excursions in each data sets are reported in Table \ref{tab:simu_5_data} in Appendix \ref{app:simu_5}.

\begin{figure}
    \centering
    \begin{subfigure}[b]{0.49\textwidth}
    \includegraphics[width=\textwidth, trim=0.cm 0.cm 0cm  0.7cm,clip]{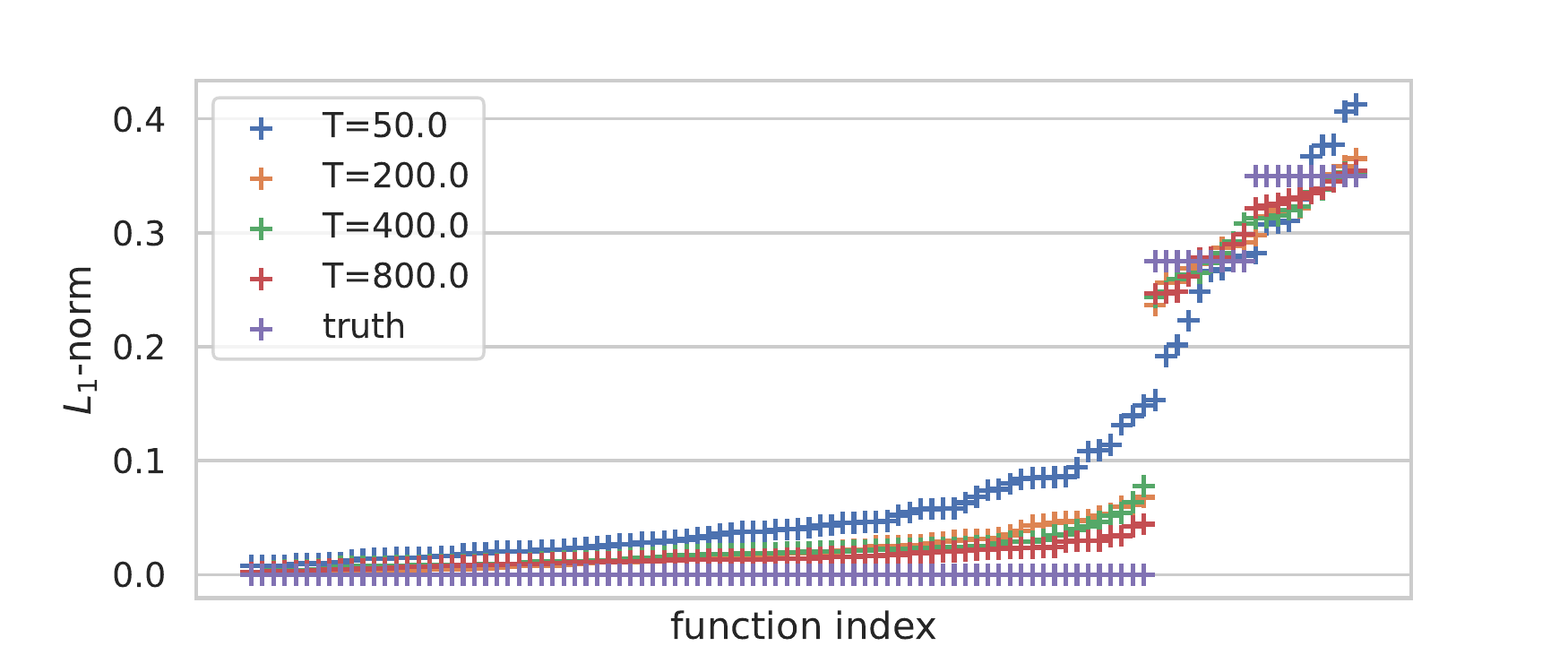}
    \caption{\emph{Excitation} scenario}
    \end{subfigure}%
        \begin{subfigure}[b]{0.49\textwidth}
    \includegraphics[width=\textwidth, trim=0.cm 0.cm 0cm  0.7cm,clip]{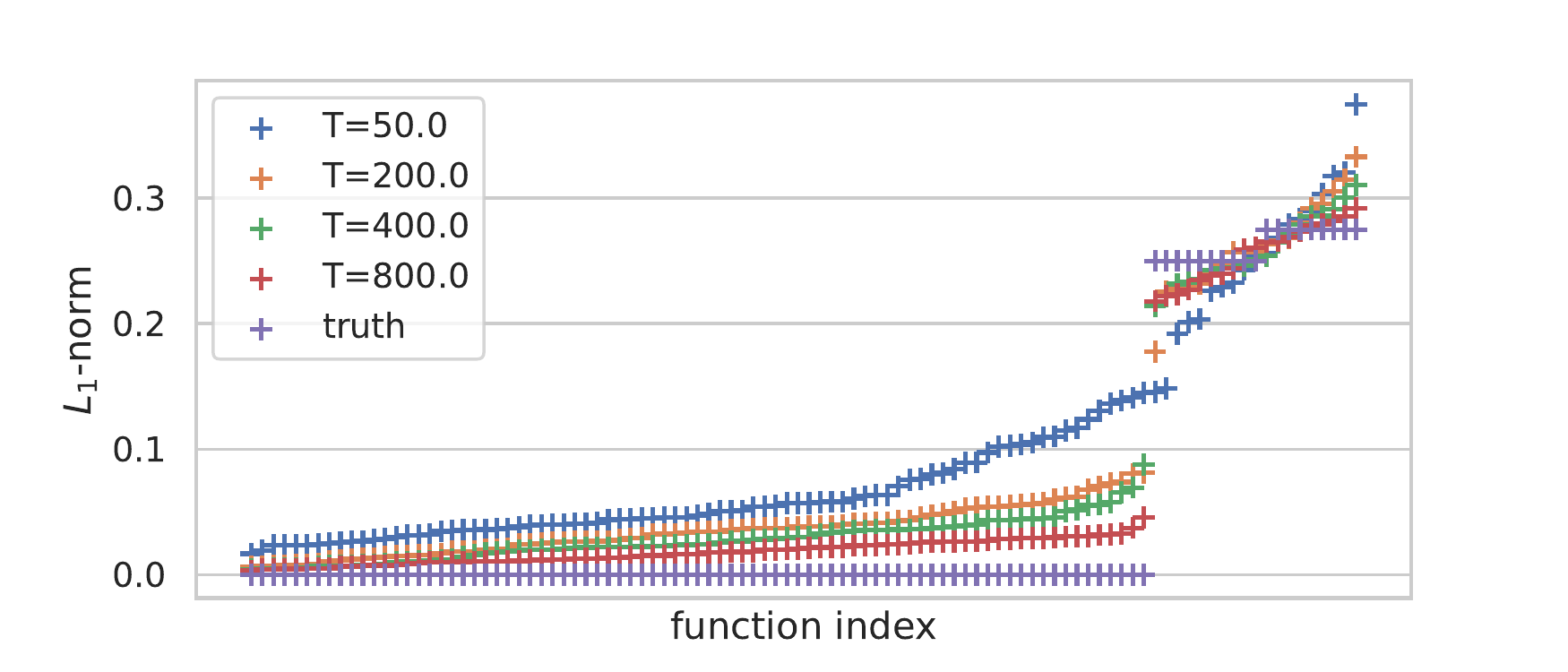}
    \caption{\emph{Inhibition} scenario}
    \end{subfigure}%
    \caption{Estimated $L_1$-norms after the first step of Algorithm \ref{alg:2step_adapt_cavi}, for different observation lengths $T$, in the \emph{Excitation} and \emph{Inhibition} scenarios of Simulation 5 with $K = 10$. % For $T=50$, there are approximately 2600 (resp. 1750) points in the \emph{Excitation} (resp. \emph{Inhibition}) scenario. The number of points and excursions in each setting are reported in Table \ref{tab:simu_5_data}. %$\eta_0$  corresponds to the smallest non-null norm of the true parameter, i.e., $\eta_0 = \min_{l,k} \norm{h_{lk}^0}_1$, and is unknown.
    We note that the norms are better estimated, after the first step of our algorithm, for larger $T$, leading to a larger gap between the small and large estimated norms, in both scenarios.}
    \label{fig:l1norms_T}
\end{figure}

We estimate the parameters using the model-selection variational posterior in Algorithm \ref{alg:2step_adapt_cavi} for each data set. From Table \ref{tab:perf_simu5}, we note that our  graph estimator converges quickly to the true graph and the risk also decreases with the number of observations. We can also see from Figure \ref{fig:l1norms_T} that the estimation of the $L_1$-norms after the first step of the algorithm improves for larger $T$, leading to a bigger gap between the small and large norms. %Besides, the connectivity graph is well recovered for all $T$ in the \emph{Excitation} scenario and almost perfectly estimated in the \emph{Inhibition} scenario for the smaller $T$.
Finally, in Figure \ref{fig:est_nu_T} (and Figure \ref{fig:est_h_T} in Appendix), we plot the model-selection variational posterior and note that its mean gets closer to the ground-truth parameter and its credible set shrinks for larger $T$. % Finally, we note that for all $T$, the number of bins of the interaction functions (here, equal to 2) is well inferred in the second step of Algorithm \ref{alg:2step_adapt_cavi}.

\begin{table}[hbt!]
    \centering
\begin{tabular}{c|c|c|c|c}
\toprule
  Scenario & T & Graph accuracy & Dimension accuracy & Risk  \\
 \midrule
 \multirow{3}{*}{Excitation}  & 50 &      1.00 &    0.40 & 7.06 \\
  &  200 &       1.00 &    1.00 & 5.07 \\
  & 400  & 1.00 &    1.00 & 5.06 \\
  & 800 & 1.00 &    1.00 & 4.01 \\
\midrule
 \multirow{3}{*}{Inhibition}  & 50 &      0.98 &    0.40 & 8.61 \\
  &  200 &      1.00 &    1.00 & 4.30 \\
  & 400  & 1.00 &    1.00 & 3.96 \\
  & 800 & 1.00 &    1.00 & 2.91 \\
\bottomrule
\end{tabular}
    \caption{Performance of Algorithm \ref{alg:2step_adapt_cavi} for the different data set sizes $T \in \{50,200,400,800\}$ in the scenarios of Simulation 5 with $K=10$. We note the graph estimator quickly converges to the true graph $\delta_0$.}
    \label{tab:perf_simu5}
\end{table}

% \begin{figure}
%     \centering
%     \begin{subfigure}[b]{0.95\textwidth}
%     \includegraphics[width=\textwidth, trim=0.cm 0.cm 0cm  0.cm,clip]{figures/estimated_graph_D10_exc.pdf}
%     \caption{\emph{Excitation} scenario}
%     \end{subfigure}
%         \begin{subfigure}[b]{0.95\textwidth}
%     \includegraphics[width=\textwidth, trim=0.cm 0.cm 0cm  0.cm,clip]{figures/estimated_graph_D10_inh.pdf}
%     \caption{\emph{Inhibition} scenario}
%     \end{subfigure}%
%     \caption{Estimated graph using the gap heuristic for different observation lengths $T \in \{50,200,400, 800\}$, in the  \emph{Excitation} and \emph{Inhibition} scenarios of Simulation 5. The true graph has non-null principal and first off-diagonal (see Figure \ref{fig:graphs}\subref{fig:graph_d10}), and is well estimated in all settings, except for a couple of \emph{False Positives} in the Inhibition scenario for the smallest T's.}
%     \label{fig:graphs_T}
% \end{figure}

\begin{figure}
    \centering
    \begin{subfigure}[b]{\textwidth}
    \includegraphics[width=\textwidth, trim=0.cm 0.cm 0cm  0.7cm,clip]{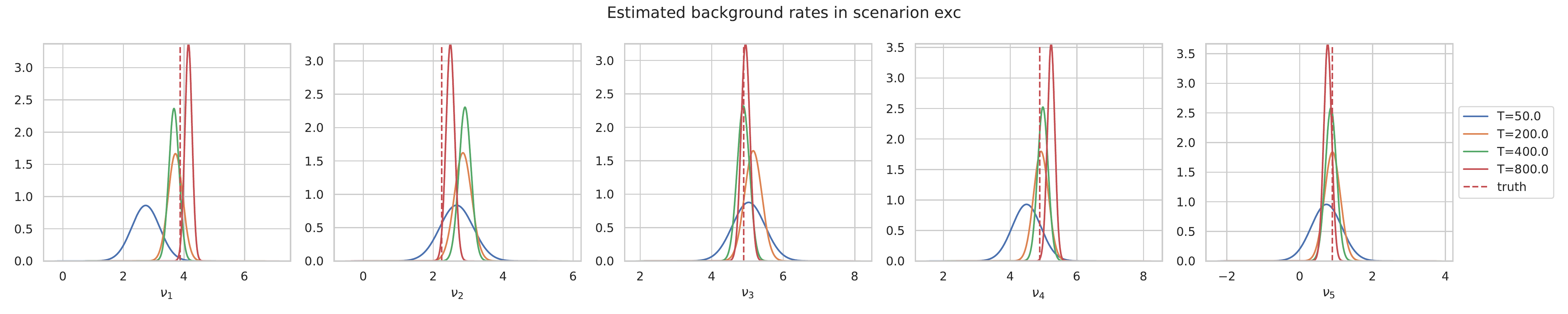}
    \caption{\emph{Excitation} scenario}
    \end{subfigure}
        \begin{subfigure}[b]{\textwidth}
    \includegraphics[width=\textwidth, trim=0.cm 0.cm 0cm  0.7cm,clip]{figs/estimated_nu_D10_exc_T.pdf}
    \caption{\emph{Inhibition} scenario}
    \end{subfigure}%
    \caption{Model-selection adaptive variational posterior on a subset of background rates, $(\nu_1, \dots, \nu_5)$, for different observation lengths $T \in \{50,200,400, 800\}$, in the  \emph{Excitation}  and \emph{Inhibition} scenarios in Simulation 5 with $K=10$. The variational posterior behaves as expected in this simulation: as $T$ increases, its mean gets closer to the ground-truth parameter and its variance decreases.}
    \label{fig:est_nu_T}
\end{figure}

\FloatBarrier

\subsection{Simulation 6: robustness to mis-specification of the link function and the memory parameter} \label{sec:simu6}

In this experiment, we  first test the robustness of our variational method based on the sigmoid model parametrised by \eqref{eq:nonlinearity} with $\xi = (0.0, 20.0, 0.2, 10.0)$ to mis-specification of the nonlinear link functions $(\phi_k)_k$. Specifically, we set $K=10$ and construct synthetic mis-specified data by simulating a Hawkes process where for each $k$, the link $\phi_k$ is chosen as:
\begin{itemize}
    \item ReLU: $\phi_k(x) = (x)_+$;
    \item Softplus: link $\phi_k(x) = \log(1 + e^x)$;
    \item Mis-specified sigmoid, with unknown $\theta_k \overset{i.i.d.}{\sim} U([15,25])$.
\end{itemize}
We also consider \emph{Excitation} and \emph{Inhibition} scenarios. Here,  $T=300$ in all settings. 

\begin{figure}
    \centering
    \begin{subfigure}[b]{0.49\textwidth}
    \includegraphics[width=\textwidth, trim=0.cm 0.cm 0cm  0.7cm,clip]{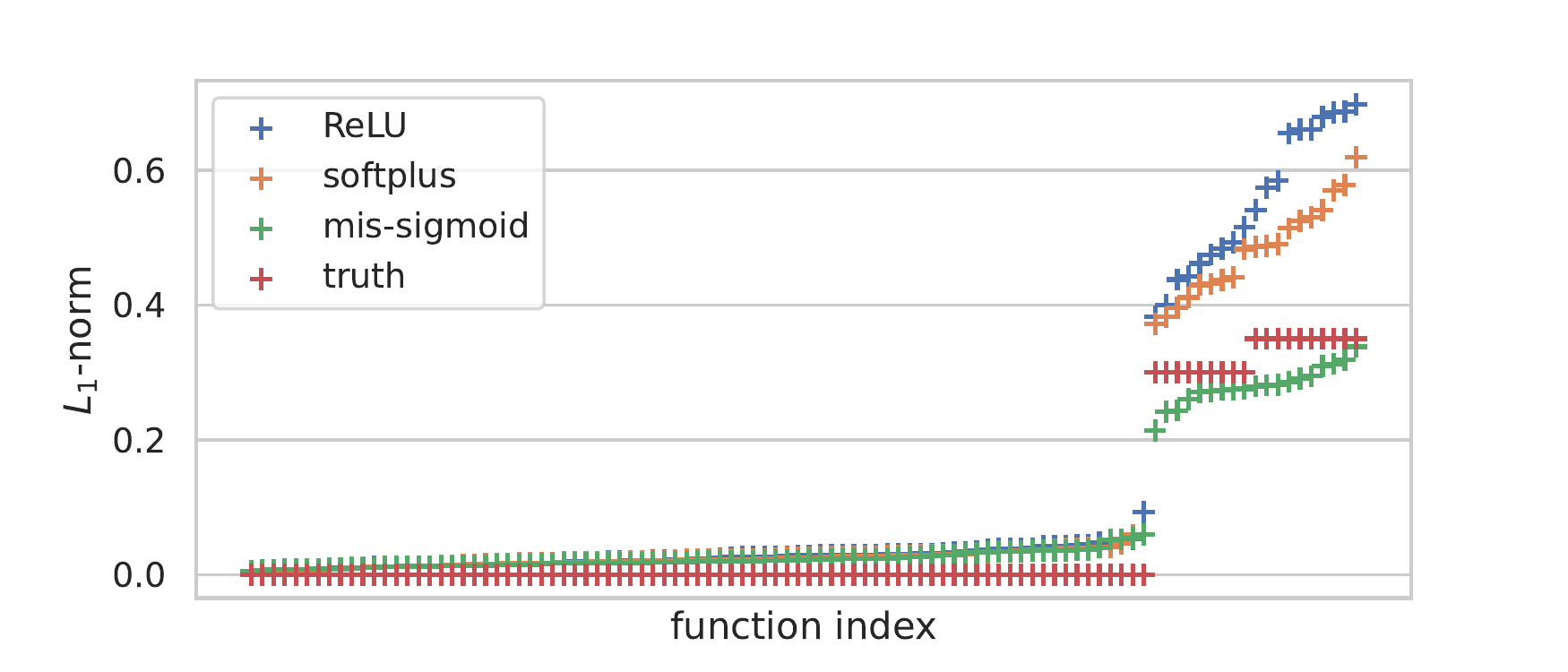}
    \caption{\emph{Excitation} scenario}
    \end{subfigure}%
        \begin{subfigure}[b]{0.49\textwidth}
    \includegraphics[width=\textwidth, trim=0.cm 0.cm 0cm  0.7cm,clip]{figs/L1_norms_D10D_exc_misspecified.pdf}
    \caption{\emph{Inhibition} scenario}
    \end{subfigure}%
    \caption{Estimated $L_1$-norms after the first step of Algorithm \ref{alg:2step_adapt_cavi}, in the mis-specified settings of Simulation 6. In this simulation, the link functions are set to $\phi_k(x) = 20\sigma(0.2(x - 10)), \forall k$, in our algorithm, while the data sets are generated from a Hawkes process with ReLU, softplus, or a mis-specified sigmoid (mis-sigmoid) link functions, in \emph{Excitation} and \emph{Inhibition} scenarios. We note that for the ReLU and softplus link, the norms are not well estimated after the first step, nonetheless, our gap heuristic can still recover the true graph parameter.}
    \label{fig:l1norms_mis}
\end{figure}

In Figure \ref{fig:l1norms_mis}, we plot the estimated $L_1$-norms after the first step of Algorithm \ref{alg:2step_adapt_cavi} and  note that there is still a gap in all settings and scenarios, although the norms are not well estimated in the case of the ReLU and softplus nonlinearities. The gaps allow to estimate well the connectivity graph parameter, but the other parameters cannot be well estimated for these two links, as can be seen from the risks in Table \ref{tab:perf_simu6a}. Nonetheless, %our gap heuristics still works for estimating the connectivity graphs, plotted in Figure \ref{fig:graphs_mis}, and
the sign of the interaction functions is well recovered in all settings.

Then, we test the robustness of our variational method to mis-specification of the memory parameter $A$, assumed to be known in our framework. We recall that $A$ corresponds to the upper bound of the support of the interaction functions. For this experiment, we generate data from the sigmoid Hawkes process with $K=10$ and with ground-truth parameter $A_0=0.1$, in two sets of parameters corresponding to an \emph{Excitation} and an \emph{Inhibition} scenarios. Here, we set $T=500$ and apply our variational method (Algorithm \ref{alg:2step_adapt_cavi}) with $A \in \{0.5, 0.1, 0.2, 0.4\}$.

\begin{figure}
    \centering
    \begin{subfigure}[b]{0.49\textwidth}
    \includegraphics[width=\textwidth, trim=0.cm 0.cm 0cm  0.7cm,clip]{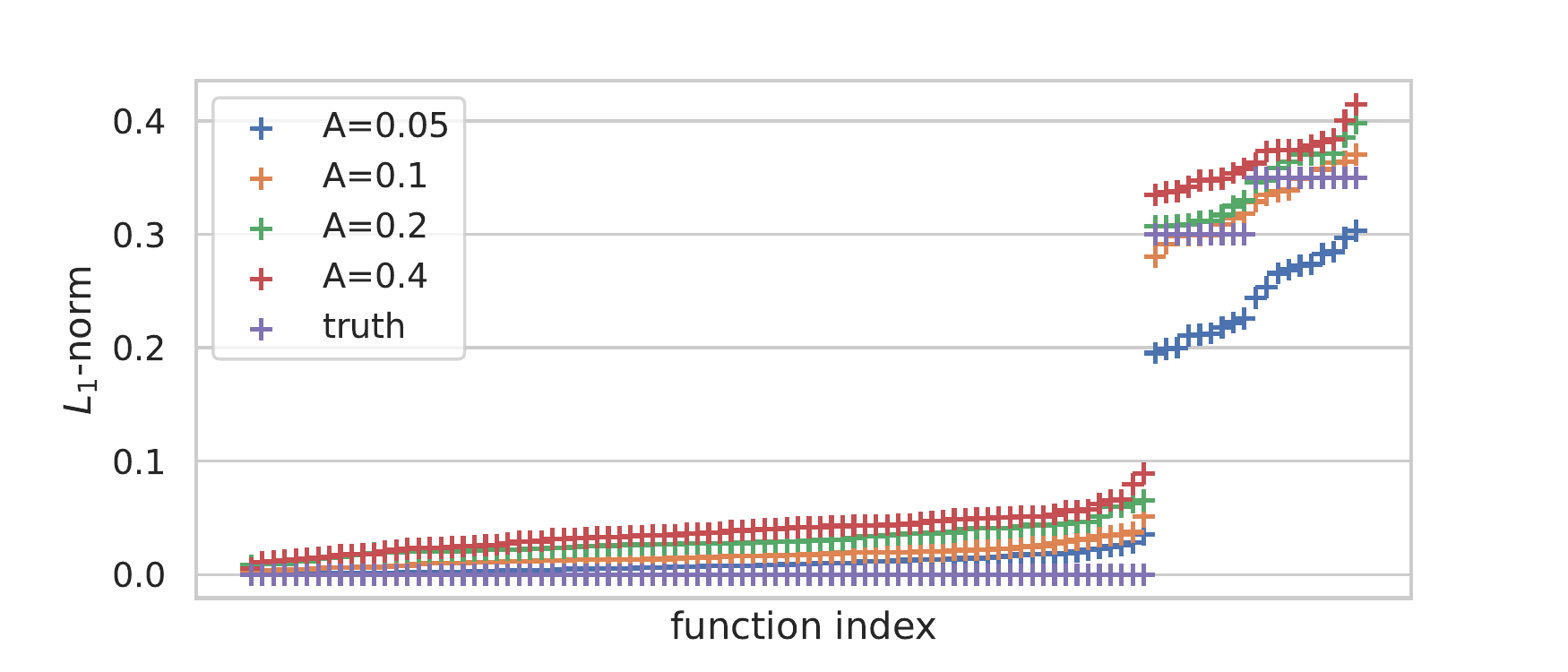}
    \caption{\emph{Excitation} scenario}
    \end{subfigure}%
        \begin{subfigure}[b]{0.49\textwidth}
    \includegraphics[width=\textwidth, trim=0.cm 0.cm 0cm  0.7cm,clip]{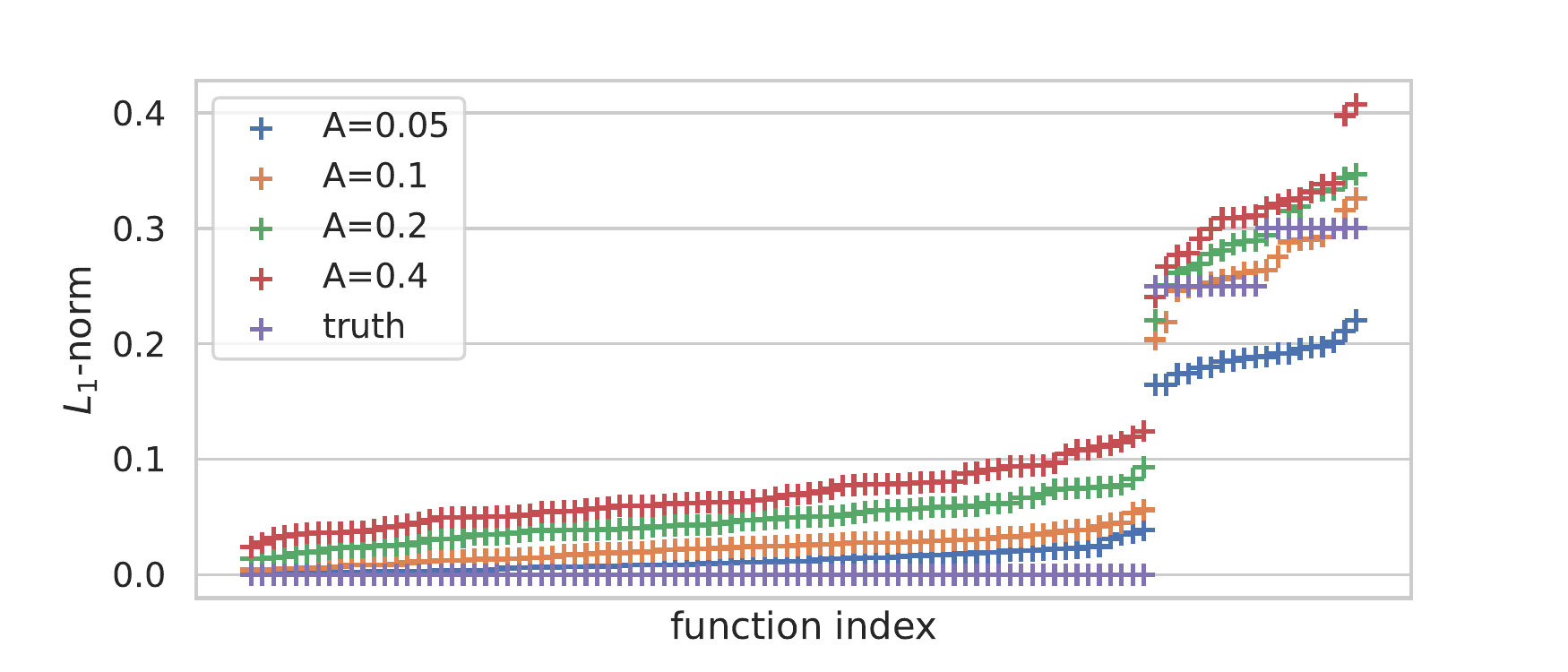}
    \caption{\emph{Inhibition} scenario}
    \end{subfigure}%
    \caption{Estimated $L_1$-norms of the interaction functions after the first step of Algorithm \ref{alg:2step_adapt_cavi} specified with different values of the memory parameter $A=0.05, 0.1, 0.2, 0.4$ containing the true memory parameter $A_0=0.1$, in the scenarios of Simulation 6. In all cases, we still observe a gap, although the norms are under-estimated (resp. over-estimated)} for $A=0.05$ (resp. $A=0.4$)
    \label{fig:l1norms_A}
\end{figure}

\begin{figure}
    \centering
    \begin{subfigure}[b]{\textwidth}
    \centering
    \includegraphics[width=0.8\textwidth, trim=0.cm 0.cm 0cm  0.cm,clip]{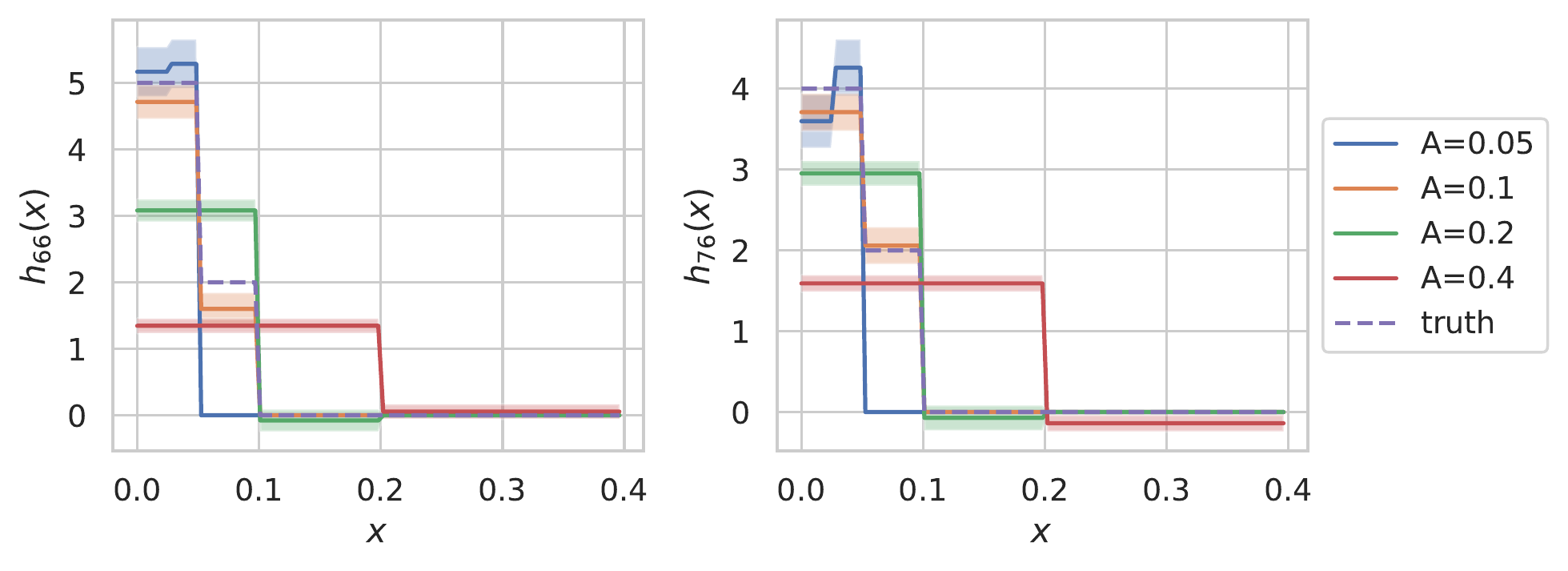}
    \caption{\emph{Excitation} scenario}
    \end{subfigure}
        \begin{subfigure}[b]{\textwidth}
        \centering
    \includegraphics[width=0.8\textwidth, trim=0.cm 0.cm 0cm  0.cm,clip]{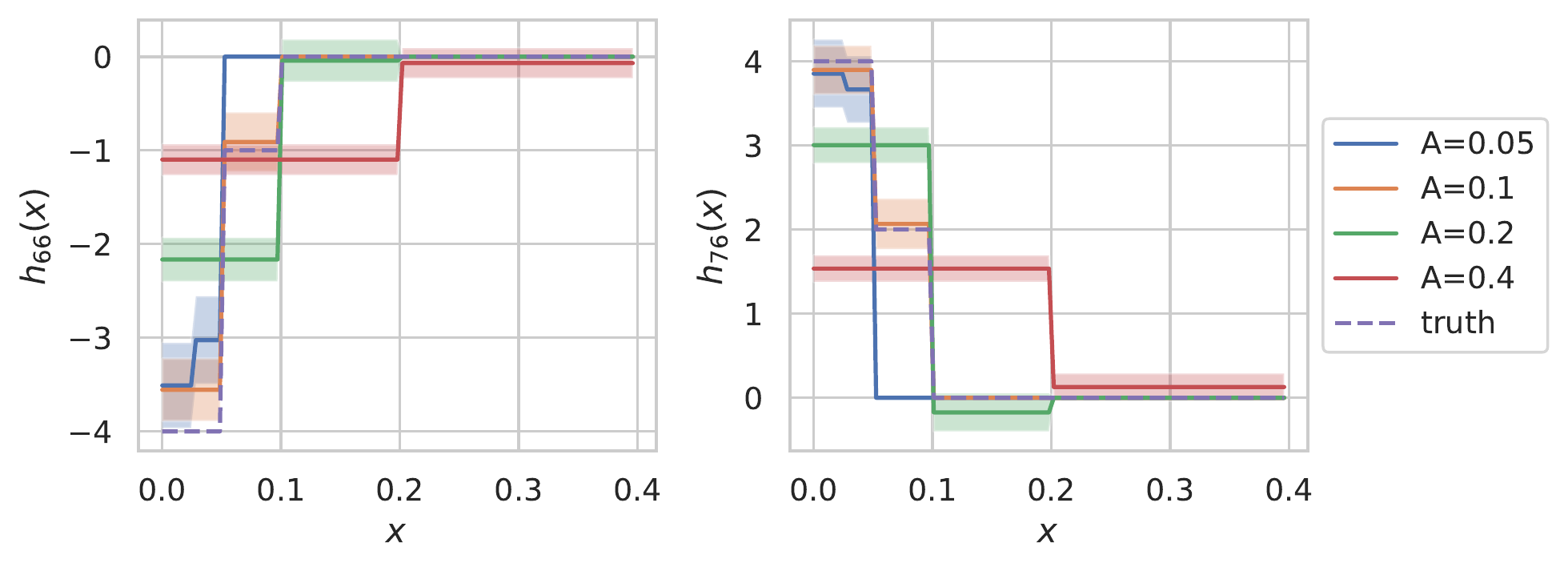}
    \caption{\emph{Inhibition} scenario}
    \end{subfigure}%
    \caption{Model-selection variational posterior on the interaction functions $h_{66}$ and $h_{76}$ obtained with Algorithm \ref{alg:2step_adapt_cavi}, specified with different values of the memory parameter $A=0.05, 0.1, 0.2, 0.4$, in the scenarios of Simulation 6 with $K=10$ and true memory parameter $A_0=0.1$. We note that the estimation of the interaction functions is deteriorated when $A$ is mis-specified, however the signs of the functions are still recovered.}
    \label{fig:est_h_mis}
\end{figure}

In Figure \ref{fig:l1norms_A}, we plot the estimated $L_1$-norms of the interaction functions, after the first step of Algorithm \ref{alg:2step_adapt_cavi}, when using the different values of $A$. We note that when $A$ is smaller than $A_0$, the norms of the non-null functions are underestimated, while if $A$ is larger than $A_0$, the norms are slightly overestimated. We note that, in all settings, the graph can be well estimated with the gap heuristics (see Figure \ref{fig:graphs_A} in Appendix). The model-selection variational posterior on a subset of the interaction functions is plotted in Figure \ref{fig:est_h_mis}. We note that for $A=0.05=A_0/2$, only the first part of the functions can be estimated, while for  $A>A_0$, the mean estimate is close to 0 on the upper part of the support. Nonetheless, in the latter case, the dimensionality of the true functions is not well-recovered.

In conclusion, this experiment shows that our algorithm is robust to the mis-specification of the nonlinear link functions and the memory parameter, for estimating the connectivity graph and the sign of the interaction functions when the latter are either non-negative or non-positive. Nonetheless, the other parameters of the Hawkes model cannot be well recovered.

\begin{table}[hbt!]
    \centering
\begin{tabular}{c|c|c|c|c}
\toprule
  Scenario & Link & Graph accuracy & Dimension accuracy & Risk  \\
 \midrule
 \multirow{3}{*}{Excitation}  & ReLU &    1.00 &   1.00 & 49.58 \\
  &  Softplus &     1.00  &  1.00 & 34.27 \\
  & Mis-specified sigmoid  & 1.00 &    1.00 & 19.69 \\
\midrule
 \multirow{3}{*}{Inhibition} & ReLU &    1.00 &   1.00 & 59.95 \\
  &  Softplus &     1.00  &  1.00 & 33.94 \\
  & Mis-specified sigmoid  & 0.99 &    1.00 & 15.78\\
\bottomrule
\end{tabular}
    \caption{Performance of Algorithm \ref{alg:2step_adapt_cavi} for the different mis-specified settings and scenarios of Simulation 6 ($K=10$). We note that the graph parameter and the dimensionality are still recovered in these cases, although, the other parameters cannot be well estimated, as can be seen from the large risk.}
    \label{tab:perf_simu6a}
\end{table}

\section{Discussion}

In this paper, we proposed a novel adaptive variational Bayes method for sparse and high-dimensional Hawkes processes, and provided a general theoretical analysis of these methods. We notably obtained variational posterior concentration rates, under easily verifiable conditions on the prior and approximating family that we validated  commonly used inference set-ups. Our general theory holds in particular in the sigmoid Hawkes model, for which we developed adaptive variational mean-field algorithms, which improve  existing ones by their ability to infer the graph parameter and the dimensionality of the interaction functions. Moreover, we demonstrated on simulated data that our most computationally efficient algorithm is able to scale up to high-dimensional processes. 
% We tested our approach in an extensive set of simulations, which demonstrates its efficiency and accuracy in the context of Bayesian nonparametric inference. Nonetheless, both 

Nonetheless, our theory does not yet cover the  high-dimensional setting with $K \to \infty$, which is of interest in applications of Hawkes processes to social network analysis and neuroscience. In this limit, previous works have considered sparse models \citep{cai2021latent, bacry2020sparse, chen17b} and mean-field settings \citep{pfaffelhuber2022mean}. We would therefore be interested in extending our results to these models. Moreover, our empirical study shows that the credible sets of variational distributions do not always have good coverage, an observation that sometimes also holds for the posterior distribution. Therefore, it is left for future work to study the property of (variational) posterior credible regions, and potentially design post-processing methods of the latter to improve coverage in practice.
Additionally, %there are still several unknown left around the task of estimating the connectivity graph in Hawkes processes, with provable guarantees.
the thresholding approach for estimating the graph in our two-step adaptive variational procedure could be further explored, in particular, in dense settings.
%However, for this goal, a future line of work could be to design other decision-based or empirical Bayes methods, and to provide some coverage guarantees.

Finally, it would be of practical interest to develop variational algorithms beyond the sigmoid model, e.g., for the ReLU and softplus Hawkes models. While in the sigmoid model, the conjugacy of the mean-field variational posterior using data augmentation leads to particularly efficient algorithms, it is unlikely that such convenient forms could be obtained for more general models. A potential approach for other models could be to parametrise variational families with normalising flows, as it is for instance done for cut posteriors in \cite{carmona22}.

\FloatBarrier

\typeout{}

% \begin{itemize}
% \item Note that figures and tables (such as Figure~\ref{fig:first} and
% Table~\ref{tab:tabone}) should appear in the paper, not at the end or
% in separate files.
% \item In the latex source, near the top of the file the command
% \verb+\newcommand{\blind}{1}+ can be used to hide the authors and
% acknowledgements, producing the required blinded version.
% \item Remember that in the blind version, you should not identify authors
% indirectly in the text.  That is, don't say ``In Smith et. al.  (2009) we
% showed that ...''.  Instead, say ``Smith et. al. (2009) showed that ...''.
% \item These points are only intended to remind you of some requirements.
% Please refer to the instructions for authors
% at \url{http://amstat.tandfonline.com/action/authorSubmission?journalCode=uasa20&page=instructions#.VFkk7fnF_0c}
% \item For more about ASA\ style, please see \url{http://journals.taylorandfrancis.com/amstat/asa-style-guide/}
% \item If you have supplementary material (e.g., software, data, technical
% proofs), identify them in the section below.  In early stages of the
% submission process, you may be unsure what to include as supplementary
% material.  Don't worry---this is something that can be worked out at later stages.
% \end{itemize}

\acks{The project leading to this work has received funding from the European Research Council (ERC) under the European Union’s Horizon 2020 research and innovation programme (grant agreement No 834175). The project is also partially funded by the EPSRC via the CDT OxWaSP.}

% Manual newpage inserted to improve layout of sample file - not
% needed in general before appendices/bibliography.
%\bibliographystyle{agsm}

\bibliography{sample}

\begin{thebibliography}{48}
\providecommand{\natexlab}[1]{#1}
\providecommand{\url}[1]{\texttt{#1}}
\expandafter\ifx\csname urlstyle\endcsname\relax
  \providecommand{\doi}[1]{doi: #1}\else
  \providecommand{\doi}{doi: \begingroup \urlstyle{rm}\Url}\fi

\bibitem[Adams et~al.(2009)Adams, Murray, and MacKay]{adams09}
Ryan~Prescott Adams, Iain Murray, and David J.~C. MacKay.
\newblock Tractable nonparametric bayesian inference in poisson processes with
  gaussian process intensities.
\newblock In \emph{Proceedings of the 26th Annual International Conference on
  Machine Learning}, ICML '09, page 9–16, New York, NY, USA, 2009.
  Association for Computing Machinery.
\newblock ISBN 9781605585161.
\newblock \doi{10.1145/1553374.1553376}.
\newblock URL \url{https://doi.org/10.1145/1553374.1553376}.

\bibitem[Arbel et~al.(2013)Arbel, Gayraud, and Rousseau]{arbeletal:13}
J.~Arbel, G.~Gayraud, and J.~Rousseau.
\newblock Bayesian adaptive optimal estimation using a sieve prior.
\newblock \emph{Scand. J. Statist.}, 40:\penalty0 549--570, 2013.

\bibitem[Bacry and Muzy(2015)]{bacry2015second}
Emmanuel Bacry and Jean-Francois Muzy.
\newblock Second order statistics characterization of hawkes processes and
  non-parametric estimation, 2015.

\bibitem[Bacry et~al.(2020)Bacry, Bompaire, Ga{\"\i}ffas, and
  Muzy]{bacry2020sparse}
Emmanuel Bacry, Martin Bompaire, St{\'e}phane Ga{\"\i}ffas, and Jean-Francois
  Muzy.
\newblock Sparse and low-rank multivariate hawkes processes.
\newblock \emph{Journal of Machine Learning Research}, 21\penalty0
  (50):\penalty0 1--32, 2020.

\bibitem[Bishop(2006)]{bishop2006pattern}
Christopher~M. Bishop.
\newblock \emph{Pattern recognition and machine learning}.
\newblock Information Science and Statistics. Springer, New York, 2006.
\newblock ISBN 978-0387-31073-2; 0-387-31073-8.
\newblock \doi{10.1007/978-0-387-45528-0}.
\newblock URL
  \url{https://doi-org.proxy.bu.dauphine.fr/10.1007/978-0-387-45528-0}.

\bibitem[Bonnet et~al.(2021)Bonnet, Herrera, and Sangnier]{bonnet2021maximum}
Anna Bonnet, Miguel~Martinez Herrera, and Maxime Sangnier.
\newblock Maximum likelihood estimation for hawkes processes with
  self-excitation or inhibition.
\newblock \emph{Statistics \& Probability Letters}, 179:\penalty0 109214, 2021.

\bibitem[Bremaud and Massoulie(1996)]{bremaud96}
Pierre Bremaud and Laurent Massoulie.
\newblock Stability of nonlinear hawkes processes.
\newblock \emph{The Annals of Probability}, 1996.

\bibitem[Cai et~al.(2021)Cai, Zhang, and Guan]{cai2021latent}
Biao Cai, Jingfei Zhang, and Yongtao Guan.
\newblock Latent network structure learning from high dimensional multivariate
  point processes, 2021.

\bibitem[Carmona and Nicholls(2022)]{carmona22}
Chris~U. Carmona and Geoff~K. Nicholls.
\newblock Scalable semi-modular inference with variational meta-posteriors,
  2022.
\newblock URL \url{https://arxiv.org/abs/2204.00296}.

\bibitem[Carstensen et~al.(2010)Carstensen, Sandelin, Winther, and
  Hansen]{carstensen2010multivariate}
Lisbeth Carstensen, Albin Sandelin, Ole Winther, and Niels~R Hansen.
\newblock Multivariate hawkes process models of the occurrence of regulatory
  elements.
\newblock \emph{BMC bioinformatics}, 11\penalty0 (1):\penalty0 1--19, 2010.

\bibitem[Chen et~al.(2017{\natexlab{a}})Chen, Shojaie, Shea-Brown, and
  Witten]{chen17b}
Shizhe Chen, Ali Shojaie, Eric Shea-Brown, and Daniela Witten.
\newblock The multivariate hawkes process in high dimensions: Beyond mutual
  excitation.
\newblock \emph{arXiv:1707.04928v2}, 2017{\natexlab{a}}.

\bibitem[Chen et~al.(2017{\natexlab{b}})Chen, Witten, and Shojaie]{Chen2017}
Shizhe Chen, Daniela Witten, and Ali Shojaie.
\newblock Nearly assumptionless screening for the mutually-exciting
  multivariate {H}awkes process.
\newblock \emph{Electron. J. Stat.}, 11\penalty0 (1):\penalty0 1207--1234,
  2017{\natexlab{b}}.
\newblock ISSN 1935-7524.
\newblock \doi{10.1214/17-EJS1251}.
\newblock URL \url{https://doi.org/10.1214/17-EJS1251}.

\bibitem[Costa et~al.(2020)Costa, Graham, Marsalle, and Tran]{costa18}
Manon Costa, Carl Graham, Laurence Marsalle, and Viet~Chi Tran.
\newblock Renewal in hawkes processes with self-excitation and inhibition.
\newblock \emph{Advances in Applied Probability}, 52\penalty0 (3):\penalty0
  879–915, 2020.
\newblock \doi{10.1017/apr.2020.19}.

\bibitem[Daley and Vere-Jones(2007)]{daley2007introduction}
Daryl~J Daley and David Vere-Jones.
\newblock \emph{An introduction to the theory of point processes: volume II:
  general theory and structure}.
\newblock Springer Science \& Business Media, 2007.

\bibitem[Deutsch and Ross(2022)]{deutsch2022}
Isabella Deutsch and Gordon~J. Ross.
\newblock Bayesian estimation of multivariate hawkes processes with inhibition
  and sparsity, 2022.
\newblock URL \url{https://arxiv.org/abs/2201.05009}.

\bibitem[Donner and Opper(2019)]{donner2019efficient}
Christian Donner and Manfred Opper.
\newblock Efficient bayesian inference of sigmoidal gaussian cox processes,
  2019.

\bibitem[Donnet et~al.(2020)Donnet, Rivoirard, and Rousseau]{donnet18}
Sophie Donnet, Vincent Rivoirard, and Judith Rousseau.
\newblock Nonparametric {B}ayesian estimation for multivariate {H}awkes
  processes.
\newblock \emph{Ann. Statist.}, 48\penalty0 (5):\penalty0 2698--2727, 2020.
\newblock ISSN 0090-5364.
\newblock \doi{10.1214/19-AOS1903}.
\newblock URL \url{https://doi-org.proxy.bu.dauphine.fr/10.1214/19-AOS1903}.

\bibitem[Eichler et~al.(2017)Eichler, Dahlhaus, and
  Dueck]{eichler2016graphical}
Michael Eichler, Rainer Dahlhaus, and Johannes Dueck.
\newblock Graphical modeling for multivariate hawkes processes with
  nonparametric link functions.
\newblock \emph{Journal of Time Series Analysis}, 38\penalty0 (2):\penalty0
  225--242, 2017.

\bibitem[Gerhard et~al.(2017)Gerhard, Deger, and Truccolo]{gerhard17}
Felipe Gerhard, Moritz Deger, and Wilson Truccolo.
\newblock On the stability and dynamics of stochastic spiking neuron models:
  Nonlinear hawkes process and point process glms.
\newblock \emph{PLOS Computational Biology}, 13:\penalty0 1--31, 02 2017.
\newblock \doi{10.1371/journal.pcbi.1005390}.
\newblock URL \url{https://doi.org/10.1371/journal.pcbi.1005390}.

\bibitem[Golub and Welsch(1969)]{golub1969calculation}
Gene~H Golub and John~H Welsch.
\newblock Calculation of gauss quadrature rules.
\newblock \emph{Mathematics of computation}, 23\penalty0 (106):\penalty0
  221--230, 1969.

\bibitem[Hansen et~al.(2015)Hansen, Reynaud-Bouret, and
  Rivoirard]{Hansen:Reynaud:Rivoirard}
Niels~Richard Hansen, Patricia Reynaud-Bouret, and Vincent Rivoirard.
\newblock Lasso and probabilistic inequalities for multivariate point
  processes.
\newblock \emph{Bernoulli}, 21\penalty0 (1):\penalty0 83--143, 2015.
\newblock ISSN 1350-7265.
\newblock \doi{10.3150/13-BEJ562}.
\newblock URL \url{http://dx.doi.org/10.3150/13-BEJ562}.

\bibitem[Hawkes(1971)]{hawkes1971point}
Alan~G Hawkes.
\newblock Point spectra of some mutually exciting point processes.
\newblock \emph{Journal of the Royal Statistical Society: Series B
  (Methodological)}, 33\penalty0 (3):\penalty0 438--443, 1971.

\bibitem[Hawkes(2018)]{hawkes18}
Alan~G. Hawkes.
\newblock Hawkes processes and their applications to finance: a review.
\newblock \emph{Quantitative Finance}, 18\penalty0 (2):\penalty0 193--198,
  2018.
\newblock \doi{10.1080/14697688.2017.1403131}.
\newblock URL \url{https://doi.org/10.1080/14697688.2017.1403131}.

\bibitem[Kingman(1993)]{kingman-poisson-processes}
J.~F.~C. Kingman.
\newblock \emph{Poisson processes}, volume~3 of \emph{Oxford Studies in
  Probability}.
\newblock The Clarendon Press Oxford University Press, New York, 1993.
\newblock ISBN 0-19-853693-3.
\newblock Oxford Science Publications.

\bibitem[Lemonnier and Vayatis(2014)]{lemonnier2014nonparametric}
Remi Lemonnier and Nicolas Vayatis.
\newblock Nonparametric markovian learning of triggering kernels for mutually
  exciting and mutually inhibiting multivariate hawkes processes.
\newblock In \emph{Joint European Conference on Machine Learning and Knowledge
  Discovery in Databases}, pages 161--176. Springer, 2014.

\bibitem[Lu and Abergel(2018)]{lu2018high}
Xiaofei Lu and Fr{\'e}d{\'e}ric Abergel.
\newblock High-dimensional hawkes processes for limit order books: modelling,
  empirical analysis and numerical calibration.
\newblock \emph{Quantitative Finance}, 18\penalty0 (2):\penalty0 249--264,
  2018.

\bibitem[Malem-Shinitski et~al.(2021)Malem-Shinitski, Ojeda, and
  Opper]{malemshinitski2021nonlinear}
Noa Malem-Shinitski, Cesar Ojeda, and Manfred Opper.
\newblock Nonlinear hawkes process with gaussian process self effects, 2021.

\bibitem[Mei and Eisner(2017)]{mei2017neural}
Hongyuan Mei and Jason Eisner.
\newblock The neural hawkes process: A neurally self-modulating multivariate
  point process, 2017.

\bibitem[Mohler et~al.(2011)Mohler, Short, Brantingham, Schoenberg, and
  Tita]{mohler11}
G.~O. Mohler, M.~B. Short, P.~J. Brantingham, F.~P. Schoenberg, and G.~E. Tita.
\newblock Self-exciting point process modeling of crime.
\newblock \emph{Journal of the American Statistical Association}, 106\penalty0
  (493):\penalty0 100--108, 2011.
\newblock \doi{10.1198/jasa.2011.ap09546}.
\newblock URL \url{https://doi.org/10.1198/jasa.2011.ap09546}.

\bibitem[Nieman et~al.(2021)Nieman, Szabo, and van Zanten]{dennis21}
Dennis Nieman, Botond Szabo, and Harry van Zanten.
\newblock Contraction rates for sparse variational approximations in gaussian
  process regression, 2021.
\newblock URL \url{https://arxiv.org/abs/2109.10755}.

\bibitem[Ogata(1999)]{ogata1999seismicity}
Yosihiko Ogata.
\newblock Seismicity analysis through point-process modeling: A review.
\newblock \emph{Seismicity patterns, their statistical significance and
  physical meaning}, pages 471--507, 1999.

\bibitem[Ohn and Lin(2021)]{ohn2021adaptive}
Ilsang Ohn and Lizhen Lin.
\newblock Adaptive variational bayes: Optimality, computation and applications,
  2021.

\bibitem[Olinde and Short(2020)]{olinde20}
Jack Olinde and Martin~B. Short.
\newblock A self-limiting hawkes process: Interpretation, estimation, and use
  in crime modeling.
\newblock In \emph{2020 IEEE International Conference on Big Data (Big Data)},
  pages 3212--3219, 2020.
\newblock \doi{10.1109/BigData50022.2020.9378017}.

\bibitem[Pfaffelhuber et~al.(2022)Pfaffelhuber, Rotter, and
  Stiefel]{pfaffelhuber2022mean}
Peter Pfaffelhuber, Stefan Rotter, and Jakob Stiefel.
\newblock Mean-field limits for non-linear hawkes processes with excitation and
  inhibition.
\newblock \emph{Stochastic Processes and their Applications}, 2022.

\bibitem[Polson et~al.(2012)Polson, Scott, and Windle]{polson2012}
Nicholas~G. Polson, James~G. Scott, and Jesse Windle.
\newblock Bayesian inference for logistic models using polya-gamma latent
  variables, 2012.
\newblock URL \url{https://arxiv.org/abs/1205.0310}.

\bibitem[Ray and Szab{\'{o} }(2021)]{Ray_2021}
Kolyan Ray and Botond Szab{\'{o} }.
\newblock Variational bayes for high-dimensional linear regression with sparse
  priors.
\newblock \emph{Journal of the American Statistical Association}, pages 1--12,
  jan 2021.
\newblock \doi{10.1080/01621459.2020.1847121}.
\newblock URL \url{https://doi.org/10.1080}.

\bibitem[Shen and Ghosal(2015)]{shen:ghosal14}
Weining Shen and Subhashis Ghosal.
\newblock Adaptive bayesian procedures using random series priors.
\newblock \emph{Scandinavian Journal of Statistics}, 42\penalty0 (4):\penalty0
  1194--1213, 2015.
\newblock \doi{https://doi.org/10.1111/sjos.12159}.
\newblock URL \url{https://onlinelibrary.wiley.com/doi/abs/10.1111/sjos.12159}.

\bibitem[Sulem et~al.(2021)Sulem, Rivoirard, and Rousseau]{sulem2021bayesian}
Deborah Sulem, Vincent Rivoirard, and Judith Rousseau.
\newblock Bayesian estimation of nonlinear hawkes process, 2021.

\bibitem[Titsias and L\'{a}zaro-Gredilla(2011)]{titsias_NIPS2011}
Michalis Titsias and Miguel L\'{a}zaro-Gredilla.
\newblock Spike and slab variational inference for multi-task and multiple
  kernel learning.
\newblock In J.~Shawe-Taylor, R.~Zemel, P.~Bartlett, F.~Pereira, and K.Q.
  Weinberger, editors, \emph{Advances in Neural Information Processing
  Systems}, volume~24. Curran Associates, Inc., 2011.
\newblock URL
  \url{https://proceedings.neurips.cc/paper/2011/file/b495ce63ede0f4efc9eec62cb947c162-Paper.pdf}.

\bibitem[van~der Vaart and van Zanten(2009{\natexlab{a}})]{van_der_Vaart_2009}
A.~W. van~der Vaart and J.~H. van Zanten.
\newblock Adaptive bayesian estimation using a gaussian random field with
  inverse gamma bandwidth.
\newblock \emph{The Annals of Statistics}, 37\penalty0 (5B), oct
  2009{\natexlab{a}}.
\newblock \doi{10.1214/08-aos678}.
\newblock URL \url{https://doi.org/10.1214}.

\bibitem[van~der Vaart and van Zanten(2009{\natexlab{b}})]{vzanten:vdv:09}
A.~W. van~der Vaart and J.~H. van Zanten.
\newblock Adaptive {B}ayesian estimation using a {G}aussian random field with
  inverse gamma bandwidth.
\newblock \emph{Ann. Statist.}, 37\penalty0 (5B):\penalty0 2655--2675,
  2009{\natexlab{b}}.
\newblock ISSN 0090-5364.
\newblock \doi{10.1214/08-AOS678}.
\newblock URL \url{https://doi-org.proxy.bu.dauphine.fr/10.1214/08-AOS678}.

\bibitem[Wang et~al.(2016)Wang, Xie, Du, and Song]{wang16}
Yichen Wang, Bo~Xie, Nan Du, and Le~Song.
\newblock Isotonic hawkes processes.
\newblock In \emph{Proceedings of the 33rd International Conference on
  International Conference on Machine Learning - Volume 48}, ICML'16, page
  2226–2234. JMLR.org, 2016.

\bibitem[Zhang and Gao(2020)]{Zhang2017ConvergenceRO}
Fengshuo Zhang and Chao Gao.
\newblock {Convergence rates of variational posterior distributions}.
\newblock \emph{The Annals of Statistics}, 48\penalty0 (4):\penalty0 2180 --
  2207, 2020.

\bibitem[Zhang et~al.(2020)Zhang, Walder, and Rizoiu]{Zhang_2020}
Rui Zhang, Christian Walder, and Marian-Andrei Rizoiu.
\newblock Variational inference for sparse gaussian process modulated hawkes
  process.
\newblock \emph{Proceedings of the AAAI Conference on Artificial Intelligence},
  34\penalty0 (04):\penalty0 6803–6810, Apr 2020.
\newblock ISSN 2159-5399.
\newblock \doi{10.1609/aaai.v34i04.6160}.
\newblock URL \url{http://dx.doi.org/10.1609/aaai.v34i04.6160}.

\bibitem[Zhou et~al.(2020)Zhou, Li, Fan, Wang, Sowmya, and Chen]{zhou2020}
Feng Zhou, Zhidong Li, Xuhui Fan, Yang Wang, Arcot Sowmya, and Fang Chen.
\newblock Efficient inference for nonparametric hawkes processes using
  auxiliary latent variables.
\newblock \emph{Journal of Machine Learning Research}, 21\penalty0
  (241):\penalty0 1--31, 2020.
\newblock URL \url{http://jmlr.org/papers/v21/19-930.html}.

\bibitem[Zhou et~al.(2021{\natexlab{a}})Zhou, Kong, Zhang, Feng, and
  Zhu]{zhou2021nonlinear}
Feng Zhou, Quyu Kong, Yixuan Zhang, Cheng Feng, and Jun Zhu.
\newblock Nonlinear hawkes processes in time-varying system,
  2021{\natexlab{a}}.

\bibitem[Zhou et~al.(2021{\natexlab{b}})Zhou, Zhang, and
  Zhu]{zhou2021efficient}
Feng Zhou, Yixuan Zhang, and Jun Zhu.
\newblock Efficient inference of flexible interaction in spiking-neuron
  networks, 2021{\natexlab{b}}.

\bibitem[Zhou et~al.(2022)Zhou, Kong, Deng, Kan, Zhang, Feng, and
  Zhu]{zhou2021jmlr}
Feng Zhou, Quyu Kong, Zhijie Deng, Jichao Kan, Yixuan Zhang, Cheng Feng, and
  Jun Zhu.
\newblock Efficient inference for dynamic flexible interactions of neural
  populations.
\newblock \emph{Journal of Machine Learning Research}, 23\penalty0
  (211):\penalty0 1--49, 2022.
\newblock URL \url{http://jmlr.org/papers/v23/21-1273.html}.

\end{thebibliography}

\newpage

\appendix

\section{Mean-field and model-selection variational inference}\label{app:mean-field}

In this section,  we first recall some general notions on mean field variational Bayes and model selection variational Bayes, then present additional details on the construction of variational families in the case of multivariate Hawkes processes.

\subsection{Mean-field approximations}

 In a general inference context, when the parameter of interest, say $\vartheta$, is decomposed into $D$ blocks, $\vartheta = (\vartheta_1, \dots, \vartheta_D)$ with $ D > 1$, a common choice of variational class is a mean-field family %\citep{Zhang2017ConvergenceRO, ohn2021adaptive}, 
 that can be defined as
$
    \mathcal{V}_{MF} = \left \{Q ; \: dQ(\vartheta) = \prod_{d=1}^D dQ_d(\vartheta_d) \right \}.
$
In this case, the mean-field variational posterior distribution corresponds to $\hat Q =  \arg \min_{Q \in \mathcal{V}_{MF}} KL \left(Q|| \Pi(.|N)\right) = \prod_{d=1}^D \hat Q_d.$ Note that the mean-field family removes some dependencies between blocks of coordinates of the parameter in the approximated posterior distribution.

Assuming that the mean-field variational posterior distribution has a density with respect to a dominating measure $\mu = \prod_d \mu_d$, with a slight abuse of notation, we denote $\hat Q$ both the distribution and density with respect to $\mu$. An interesting result from \cite{bishop2006pattern} is that the mean-field variational posterior distribution verifies, for each $d \in [D]$, 
\begin{align}\label{eq:var_post_factor}
    \hat Q_d(\vartheta_d) \propto \exp \left \{ \mathbb{E}_{\hat Q_{-d}} [\log p(\vartheta,N ) ] \right \},
\end{align}
where $p(\vartheta, N )$ is the joint density of the observations and the parameter with respect to  $\prod_d \mu_d \times \mu_N$ with $\mu_N$ the data density, and $\hat Q_{-d}  := \prod_{d' \neq d} \hat Q_{d'} $. 
This property \eqref{eq:var_post_factor} can be used to design efficient algorithms for computing the variational posterior, such as the coordinate-ascent variational inference algorithm. 
%In this type of variational class, the approximating distribution factorises into a product of distributions, or \emph{factors}, which apply on a partition of the parameter's coordinates.

In a general setting where the log-likelihood function of the nonlinear Hawkes model can be augmented with some latent variable $z \in \mathcal{Z}$ (see for instance \cite{zhou2021nonlinear, zhou2021jmlr, malemshinitski2021nonlinear}), with $\mathcal{Z}$ the latent parameter space,  the augmented log-likelihood $L_T^A(f,z)$ leads to an \emph{augmented} posterior distribution, defined  as
\begin{align*}%\label{def:aug_posterior_dist}
    \Pi_A(B |N) = \frac{\int_{B} \exp(L_T^A(f,z)) d(\Pi(f) \times \mathbb{P}_{A}(z))}{\int_{\mathcal{F} \times \mathcal{Z}} \exp(L_T^A(f,z)) d(\Pi(f) \times \mathbb{P}_{A})(z)}, \quad B \subset \mathcal{F} \times  \mathcal{Z},
\end{align*}
where $\mathbb{P}_{A}$ is a prior distribution on $z$ which has a density with respect to a dominating measure $\mu_z$. Recalling the mean-field variational from Section \ref{sec:aug-mf-vi} defined as
\begin{align*}
    \mathcal{V}_{AMF} = \left \{Q: \mathcal{F} \times \mathcal{Z} \to [0,1] ; \: Q(f, z) =  Q_1(f)Q_2(z) \right \},
\end{align*}
% Then, an augmented mean-field variational class factorising along the parameter and latent variables can be defined as
% \begin{align*}
%     \mathcal{V}_{MF, A} = \left \{Q ; \: dQ(f, \omega, \bar N) =  dQ_1(f)dQ_2(\omega, \bar N) \right \}.
% \end{align*}
the augmented mean-field variational posterior corresponds to
\begin{align}\label{eq:mean-field-vp}
    \hat Q_{AMF}(f,z) := \arg \min_{Q \in \mathcal{V}_{AMF}} KL \left(Q(f,z) ||  \Pi_A(f,z |N) \right) =: \hat Q_1(f) \hat Q_2(z), 
\end{align}
and, using property \eqref{eq:var_post_factor}, verifies
\begin{align} \label{eq:var_factor_1} 
    \hat Q_{1}(f) \propto \exp \left \{ \mathbb{E}_{\hat Q_{2}} [\log p(f, z, N ) ] \right \},\quad
    \hat Q_{2}(z) \propto \exp \left \{ \mathbb{E}_{\hat Q_{1}} [\log p(f, z,N ) ] \right \}, 
\end{align}
where  $p(f, z,N)$ is the joint density of the parameter, the latent variable, and the observations with respect to the measure $\prod_d \mu_d \times \mu_z \times \mu_N$.

\subsection{Model-selection variational posterior}\label{app:model-selection}

In this section, we present two model-selection variational approaches to approximate the posterior by an adaptive variational posterior distribution. We recall from our construction in Section \ref{sec:adapt-mf-vi} that our parameter $f$ of the Hawkes processes is indexed by a model $m$ of hyperparameters in the form $m=(\delta, J_{lk}, (l,k) \in \mathcal I(\delta))$, where $\mathcal I(\delta) = \{ (l,k);\, \delta_{lk}=1\}$ is the set of non null functions. 

In a model-selection variational approach, one can consider a set of candidate models $\mathcal M$ and for any $m\in\mathcal M$, a class of variational distributions on $f$ with model $m$, denoted $\mathcal{V}^m$. Then, one can define the total variational class as $\mathcal{V} = \cup_{m \in \mathcal M} \{ \{m\}\times \mathcal{V}^m \}$, which contains distributions on $f$ localised on one model. Then, given $\mathcal{V}$ and as shown for instance in \cite{Zhang2017ConvergenceRO}, the variational posterior distribution has the form 
\begin{align*} %\label{eq:ms_var_post}
    \hat Q &:= \hat Q_{\hat m}, \quad  \hat m :=  \arg \max_{m \in \mathcal{M}} ELBO(\hat Q^m),
\end{align*}
 where  $\hat Q^m = \arg \min_{Q \in \mathcal{V}^m} KL(Q|||\Pi(.|N))$ and $ELBO(\cdot)$ is called the \textit{evidence lower bound (ELBO)}, defined as
\begin{align}\label{eq:elbo_1}
    ELBO(Q) &:= \mathbb{E}_{ Q} \left[ \log \frac{p(f,z, N)}{Q(f,z) }\right], \quad Q \in \mathcal{V}.
\end{align}
The ELBO is  a lower bound of the marginal log-likelihood $p(N)$. %$KL(.|| \Pi(.|N))$, 

An alternative model-selection variational approach  consists in constructing a model-averaging  variational posterior, also called \emph{adaptive} in \cite{ohn2021adaptive}, as a mixture of distributions over the different models, i.e.,
\begin{align}\label{eq:adap_var_post}
    \hat Q = \sum_{m \in \mathcal{M}} \hat \gamma_{m} \hat Q_m,
\end{align}
where $ \{\hat \gamma_{m}\}_{m \in \mathcal{M}} $ are marginal probabilities defined as 
\begin{align}
    \hat \gamma_{m} = \frac{ \Pi_m(m) \exp \left \{ ELBO(\hat Q_m)  \right \}}{\sum_{m \in \mathcal{M}}\Pi_m(m) \exp \left \{ ELBO(\hat Q_m)  \right \}}, \quad \forall m \in \mathcal{M}.
\end{align}
In this strategy, the approximating family of distributions corresponds to
\begin{align*}
    \mathcal{V} = \left \{ \sum_{m \in \mathcal{M}} \alpha_{m} Q_m; \sum_{m} \alpha_{m}  = 1, \: \alpha_m \geq 0, \: Q_m \in \mathcal{V}^m, \: \forall m \right\}.
\end{align*}

\section{Data augmentation in the sigmoid Hawkes model}\label{app:data-aug}

In this section, we recall the latent variable augmentation strategy and the definition of the augmented mean-field variational distribution in sigmoid-type Hawkes processes, proposed in previous work \citep{zhou2021jmlr, malemshinitski2021nonlinear}. In our method in Section \ref{sec:adapt-mf-vi}, we use this construction to efficiently compute an approximated posterior distribution on $\mathcal{F}_m \subset \mathcal{F}$, on parameters $f$ within  a model $m=(\delta, J_{lk}; (l,k) \in \mathcal I(\delta))$.

The first data augmentation step consists in re-writing the sigmoid function as a mixture of Polya-Gamma random variables \citep{polson2012}, i.e.,
    \begin{align}\label{eq:sigmoid}
        \sigma(x) =  \mathbb{E}_{\omega \sim p_{PG}(.;1,0)}\left[e^{g(\omega,x)}\right] = \int_0^{+\infty} e^{g(\omega,x)} p_{PG}(\omega;1,0) d\omega, \quad g(\omega,x) = - \frac{\omega x^2}{2} + \frac{x}{2} - \log 2,
    \end{align}
     with $p_{PG}(.;1,0)$ the Polya-Gamma density. We recall that $p_{PG}(.;1,0)$ is the density of the random variable
     \begin{align*}
         \frac{1}{2\pi^2} \sum_{k=1}^\infty \frac{g_k}{(k-1/2)^2}, \quad g_k \overset{\mathrm{i.i.d.}}{\sim} Gamma(1,1), % \quad \forall k \in \mathbb{N}\backslash\{0\}.
     \end{align*}
and that the \emph{tilted} Polya-Gamma distribution is defined as
\begin{align*}
    p_{PG}(\omega;1,c) = \cosh{\left(\frac{c}{2}\right)} \exp \left \{- \frac{c^2 \omega}{2} \right \}  p_{PG}(\omega;1,0), \quad c\geq 0,
\end{align*}
where $\cosh$ denotes the hyperbolic cosine function. %\textcolor{red}{sauf si on utilise la \emph{tilted} Polya-Gamma distribution, cette partie peut etre enlevee}.
%Let $f = (f_k)_{k\in [K]}$ with for any $k \in [K]$, $f_k = (\nu_k, (h_{lk})_{l\in [K]})$.
With a slight abuse of notation, we re-define the linear intensity \eqref{def:lin_intensity} as  
$
    \tilde{\lambda}^k_t(f) =\alpha \left( \nu_k + \sum_{l=1}^{K} \int_{-\infty}^{t^-} h_{lk}(t-s)dN_s^l - \eta \right), 
$
so that we have $\lambda^k_t(f) = \theta_k \sigma(\tilde{\lambda}^k_t(f) ), t\in \R$.
%$$h(t,f_k) = \beta (\tilde{\lambda}^k_t(f) - \eta) = \alpha \left( \nu_k + \sum_{l=1}^{K} \int_{-\infty}^{t^-} h_{lk}(t-s)dN_s^l -\eta\right).$$
For any $k \in [K]$, let $ N_k := N^k[0,T]$ and $T_1^k, \dots, T_{N_k}^k \in [0,T]$ %$N = (t_i^k)_{k \in [K], i \in [N_k]}$ with $ N_k := N^k[0,T]$
be the times of events at component $N^k$.
Now, let $\omega = (\omega_i^k)_{k\in [K], i \in [N_k]}$ be a set of latent variables such that
\begin{align*}
   \omega_i^k \overset{\mathrm{i.i.d.}}{\sim}   p_{PG}(  \cdot;1,0), \quad i \in [N_k], \quad k \in [K].
\end{align*}
Then, using \eqref{eq:sigmoid}, an \emph{augmented} log-likelihood function can be defined as 
\begin{align}\label{eq:aug_loglik}
    L_T(f, \omega ;N) &= \sum_{k \in [K]} \left \{ \sum_{i \in [N_k]} \left( \log \lsup_k  + g(\omega_i^k, \Tilde{\lambda}_{T_i^k}(f)) + \log p_{PG}( \omega_i^k; 1, 0) \right)  - \int_0^T \lsup_k \sigma(\Tilde{\lambda}^k_{t}(f)) dt \right \},
\end{align}
% \begin{align*}
%     \ell_T(N, \omega |f) &\propto \prod_{k} \left \{ \prod_{i} \lsup_k e^{g(\omega_i^k, \Tilde{\lambda}_{t_i^k}(f))}PG( \omega_i^k; 1, 0)  \exp \left \{ \int_0^T \lsup_k \sigma(\Tilde{\lambda}_t^k)) dt\right \} \right \},
% \end{align*}
%Moreover, following \cite{zhou2021nonlinear, malemshinitski2021nonlinear}
and, using that $\sigma(x) = 1 - \sigma(-x)$, the integral term on the RHS in \eqref{eq:aug_loglik} can be re-written as
\begin{align*}
\int_0^T \lsup_k \sigma(\Tilde{\lambda}^k_{t}(f)) dt &= \int_0^T \int_0^\infty \lsup_k \left[ 1 - e^{g(\bar \omega,- \Tilde{\lambda}^k_{t}(f))}\right] p_{PG}(\bar \omega;1,0)  d \bar \omega dt.
\end{align*}
Secondly,  Campbell's theorem \citep{daley2007introduction, kingman-poisson-processes} is applied. We first recall here its general formulation. For a Poisson point process $\bar Z$ on a space $ \mathcal{X}$ with intensity measure $\Lambda: \mathcal{X} \to \R^+$, and for any function  $\zeta: \mathcal{X} \to \R$, it holds true that
\begin{align}\label{eq:campbell}
    \Ex{\prod_{x \in \bar Z} e^{\zeta(x)}} = \exp \left \{ \int (e^{\zeta(x)} - 1) \Lambda(dx)\right \}.
\end{align}
%We note this theorem derives from Fubini's theorem.
% \begin{align*}
%     \Ex{\sum_{x \in \bar Z}\zeta(x)} = \int \zeta(x) \Lambda(dx).
% \end{align*}
Therefore, using that $\sigma(x) = 1 - \sigma(-x)$, and considering for each $k$ a marked Poisson point process $\bar Z^k$ on $\mathcal{X} = ([0,T], \R^+)$ with intensity measure $\Lambda^k(t,  \omega) = \lsup_k p_{PG}( \omega;1,0)$, and distribution $\mathbb{P}_{\bar Z}$, applying Campbell's theorem with  $\zeta(t, \omega ) := g(\omega, -\Tilde{\lambda}^k_{t}(f))$, one obtains that
\begin{align*}
    \Ex{\prod_{(\bar T_j^k, \bar \omega_j^k) \in \bar Z^k}e^{g(\bar \omega_k, -\Tilde{\lambda}^k_{\bar T_j}(f))}} &= \exp \left \{ \int_0^T \int_0^\infty \lsup_k \left(e^{g(\bar \omega, - \Tilde{\lambda}^k_{t}(f))} - 1\right) p_{PG}(\bar \omega;1,0)  d \bar \omega dt \right \}.
    %=  e^{ - \int_0^T \lsup_k \sigma(\Tilde{\lambda}_{t}(f)dt}.
    %&\approx \prod_{(\bar t_j^k, \bar \omega_j^k) \in \bar N^k} e^{g(\omega_j^k, -\Tilde{\lambda}_{\bar t_j^}(f))},
\end{align*}

%and applying the previous result with $\mathcal{X} = (\R^+, \R)$, $x = (\omega, t)$ and $\zeta(x) = g(\omega, -\Tilde{\lambda}_{t}(f)$, one then obtains that
% \begin{align*}
%     \exp \left( - \int_0^T \lsup_k \sigma(\Tilde{\lambda}_{t}(f)dt\right) =  \Ex{\prod_{(\bar t_j^k, \bar \omega_j^k) \in \bar N^k}e^{g(\bar \omega_k, -\Tilde{\lambda}_{\bar t_j^}(f)}} \approx \prod_{(\bar t_j^k, \bar \omega_j^k) \in \bar N^k} e^{g(\omega_j^k, -\Tilde{\lambda}_{\bar t_j^}(f))},
% \end{align*}
% where $\bar N^k$ is a marked PPP with intensity $\Lambda^k(t,\omega) = \lsup_k PG(\omega;1,0)$.
% Since $g(\omega,x)$ is quadratic in $x$, this leads to a quadratic form of the augmented likelihood wrt to the parameters $\nu, h$. In fact, let
% \begin{align*}
%     &N \times \omega = ((t_i^k, \omega_i^k))_{k=1\dots K,i=1\dots N^k} \quad \text{(augmented Hawkes process)} \\
%     &\bar N = ((\bar t_j^k, \bar \omega_j^k)_{k=1\dots K,j=1\dots \bar Z^k} \quad \text{(latent marked Poisson process with intensity)}.
% \end{align*}
Conditionally on $N$, let  $\bar Z := (\bar Z^1, \dots, \bar Z^K)$ be an observation of the previous Poisson point process  on $[0,T]$. For each $k \in [K]$, we denote  $\bar Z_k := \bar Z^k[0,T]$, $(\bar T_1^k, \bar \omega_1^k), \dots,  (\bar T_1^k, \bar \omega_{\bar Z_k}^k) \in [0,T] \times \R_+$ the times and marks of $\bar Z_k$, and $\bar Z = (\bar Z, \omega_i^k, i\leq  N_k, k\leq K)$, the set of augmented variables. Then, replacing the integral term in \eqref{eq:aug_loglik} by a product over the observation $\bar Z$,  the \emph{doubly augmented} log-likelihood function %for an observation of the process $N$ and the latent variables $(\bar Z, \omega)$
corresponds to
%, given a parameter $f = (f_k)_k, f_k \in \R^{KB+1}$,  is given by
\begin{align*}
    L_T(f,\omega, \bar Z; N) &= \sum_{k \in [K]} \left \{ \sum_{i \in [N_k]} \left [ \log \lsup_k  + g(\omega_i^k, \Tilde{\lambda}_{T_i^k}(f)) + \log p_{PG}( \omega_i^k; 1, 0) \right]\right.\\ &\left.\hspace{1cm}+ \sum_{j \in [\bar Z_k]} \left [ \log \lsup_k +  g(\bar \omega_j^k, -\Tilde{\lambda}_{\bar T_j}(f)) + \log p_{PG}(\bar \omega_j^k; 1, 0)  \right] -\lsup_k T\right \}.
\end{align*}
 
%where $\bar Z := (\bar Z^k)_{k \in [K]}$ and $\bar Z_k := \bar Z^k[0,T]$.
The previous augmented log-likelihood function, and the prior distribution $\Pi$ on the parameter and the latent variables distribution $\mathbb{P}_A = p_{PG}(.|1,0) \times \mathbb{P}_{\bar Z}$, allow to construct an \emph{augmented} posterior distribution proportional to
\begin{align}\label{eq:augm_posterior}
    \Pi(f, \omega, \bar Z |N) 
    &\propto \prod_k \left \{  \prod_{i \in [N_k]} \lsup_k  e^{g(\omega_i^k, \Tilde{\lambda}_{T_i^k}(f))}p_{PG}( \omega_i^k; 1, 0)  \times \prod_{j \in [\bar Z^k]} \lsup_k e^{g(\bar \omega_j^k, -
    \Tilde{\lambda}_{\bar T_j}(f))}p_{PG}(\bar \omega_j^k; 1, 0) \right \} \times \Pi(f).
\end{align}

\section{Analytical derivation in the sigmoid Hawkes model}

\subsection{Mean-field updates in a fixed model}\label{app:fixed_updates}

In this section, we derive the analytic forms of the conditional updates in Algorithm \ref{alg:cavi}, the mean-field variational algorithm with fixed dimensionality described in Section \ref{sec:aug-mf-vi}. For ease of exposition, in this section we consider a model $m$ and a dimension $k$ and we drop the indices $k$ and $m$, e.g., we use the notation $Q_1, Q_2$ for the variational factors. In the following computation, we use the notation $c$ to denote a generic constant which value can vary from one line to the other. For simplicity, we also assume that $J := J_1 = \dots = J_K$ and we recall that $\phi_k(x) = \theta_k \sigma(\alpha(x-\eta))$. % with $\theta_k = 20$, $\alpha := 0.1$ and $\eta := 10$.

From the definition of the augmented posterior \eqref{eq:augm_posterior}, we first note that
\begin{small}
\begin{align}\label{eq:joint_dist}
    \log p(f, N, \omega, \bar Z )    &= \log \Pi(f,\omega, \bar Z| N) + \log p(N) = L_T(f,\omega, \bar Z; N) + \log \Pi(f) + \log p(N) + c \nonumber \\
    &= \log p(\omega | f, N) + \log p(\bar Z|f,N) + \log \Pi(f) + \log p(N) + c.
\end{align}
\end{small}
In the previous equality we have used the facts that $p(\omega | f, N, \bar Z) = p(\omega | f, N)$ and $p(\bar Z|f,N, \omega) = p(\bar Z|f,N)$. We recall our notation $H(t) = (H^0(t), H^1(t), \dots, H^K(t)) \in \R^{KJ + 1},  \: t \in \R$, where for $k \in [K]$, $H^k(t) = (H_j^k(t))_{j=1, \dots, J}$ and $H_j^k(t)$ is defined in \eqref{eq:notation_Hjk}. In the following, $H(t)$ denotes $H^k(t)$ for the chosen $k$. We have that
\begin{align*}
    \mathbb{E}_{Q_2} [\log p(\omega|f, N) ] 
&= \mathbb{E}_{Q_2} \left[ \sum_{i \in [N]} g(\omega_i, \Tilde{\lambda}_{T_i}(f)) \right] + c = \mathbb{E}_{Q_2} \left[\sum_{i \in [N]} - \frac{\omega_i \Tilde{\lambda}_{T_i}(f)^2}{2} +  \frac{ \Tilde{\lambda}_{T_i}(f)}{2} \right] + c \\
    &= \mathbb{E}_{Q_2} \left[ \sum_{i \in [N]} - \frac{\omega_i \alpha^2 ( f^T H(T_i) H(T_i)^T f - 2 \eta H(T_i)^T f + \eta^2)}{2} +  \frac{ \alpha H(T_i)^T f }{2} \right] + c \\
        &= \mathbb{E}_{Q_2} \left[ - \frac{1}{2} \sum_{i \in [N]} \left \{ \omega_i \alpha^2  f^T H(T_i) H(T_i)^T f - \alpha  (2 \omega_i \alpha \eta  + 1 ) H(T_i)^T  f + \omega_i \alpha^2 \eta^2 \right \}   \right] + c \\
    &=- \frac{1}{2}  \sum_{i \in [N]} \left \{\mathbb{E}_{Q_2}[\omega_i] \alpha^2  f^T H(T_i)  H(T_i)^T f - \alpha  (2 \mathbb{E}_{Q_2}[\omega_i] \alpha \eta  + 1 ) H(T_i)^T  f + \mathbb{E}_{Q_2}[\omega_i] \alpha^2 \eta^2 \right \}  + c.
        % &=- \frac{1}{2} (f  - \tilde \mu_1)^T \tilde{\Sigma}_1^{-1}(f - \tilde \mu_1)  + c,
\end{align*}
Moreover,  we also have that 
\begin{align*}
    \mathbb{E}_{Q_2} [\log p(\bar Z|f, N) ] %&= \mathbb{E}_{Q_2} \left[ \sum_{j \in [\bar Z]} g(\bar \omega_j, -   h(\bar T_j,f)) \right] + c \\
    %&= \mathbb{E}_{Q_2} \left[\sum_{j} - \frac{\bar \omega_j h(\bar T_j,f)^2}{2} -  \frac{ h(\bar T_j,f)}{2} \right] + c \\
    %&= \mathbb{E}_{Q_2} \left[ \sum_{j} - \frac{\bar \omega_j \alpha^2 ( f^T H(\bar T_j) H(\bar T_j)^T f - 2 \eta H(\bar T_j)^T f + \eta^2)}{2} -  \frac{ \alpha H(\bar T_j)^T f }{2} \right] + c \\
        &= \mathbb{E}_{Q_2} \left[ - \frac{1}{2} \sum_{j \in [\bar Z]}  \left \{ \bar \omega_j \alpha^2  f^T H(\bar T_j) H(\bar T_j)^T f - \alpha  (2 \bar \omega_j \alpha \eta  - 1) H(\bar T_j)^T  f + \bar \omega_j \alpha^2 \eta^2 \right \}   \right] + c \\
    &= \int_{0}^T \int_0^{\infty} \left[  - \frac{1}{2} \left(\bar \omega \alpha^2  f^T H(t) H(t)^T f - \alpha  (2 \bar \omega \alpha \eta  - 1 ) H(t)^T  f + \bar \omega \alpha^2 \eta^2 \right)   \right] \Lambda(t,\bar \omega) d\bar \omega dt + c \\
        &=  - \frac{1}{2} \left[f^T \left(  \alpha^2  \int_0^T \int_0^{\infty} \bar \omega H(t) H(t)^T  \Lambda(t,\bar \omega) d\bar \omega dt \right) f \right.\\ &\hspace{1cm}\left.+ f^T \left( \alpha \int_{0}^T \int_0^{\infty}  (2 \bar \omega \alpha \eta  - 1 ) H(t)^T \Lambda(t,\bar \omega) d\bar \omega dt   \right) \right] + c.
        % &=- \frac{1}{2} (f  - \tilde \mu_2)^T \tilde{\Sigma}_2^{-1} (f - \tilde \mu_2)  + c,
\end{align*}
Besides, we have
$
     \mathbb{E}_{Q_2} [\log \Pi(f) ] = - \frac{1}{2} f^T \Sigma^{-1} f + f^T \Sigma^{-1} \mu + c.
$
Therefore, using \eqref{eq:var_factor_1}, we obtain that
\begin{align*}
    \log Q_1(f)  &=  - \frac{1}{2} \left[f^T \left( \alpha^2   \sum_{i \in [N]} \mathbb{E}_{Q_2}[\omega_i]  H(T_i)  H(T_i)^T  +  \alpha^2  \int_0^T \int_0^{\infty} \bar \omega H(t) H(t)^T  \Lambda(t,\bar \omega) d\bar \omega dt + \Sigma^{-1 }\right) f \right. \\
      &-  \left. f^T \left( \alpha   \sum_{i \in [N]} (2 \mathbb{E}_{Q_2}[\omega_i] \alpha \eta  + 1 ) H(T_i)^T +  \alpha \int_{0}^T \int_0^{\infty}   (2 \bar \omega \alpha \eta  - 1) H(t)^T \Lambda(t,\bar \omega) d\bar \omega dt  + 2\Sigma^{-1 } \mu \right) \right] + c \\
       &=: - \frac{1}{2} (f  - \tilde \mu)^T \tilde{\Sigma}^{-1} (f - \tilde \mu)  + c,
\end{align*}
therefore $Q_1(f)$ is a normal distribution  with mean vector $ \tilde \mu$ and covariance matrix  $\tilde{\Sigma}$ given by
\begin{align}\label{eq:vi_updates}
    &\tilde{\Sigma}^{-1} =  \alpha^2  \sum_{i \in [N]} \mathbb{E}_{Q_2}[\omega_i]  H(T_i)  H(T_i)^T  +  \alpha^2  \int_0^T \int_0^{\infty} \bar \omega H(t) H(t)^T  \Lambda(t,\bar \omega) d\bar \omega dt + \Sigma^{-1 }, \\
    &\tilde{\mu} = \frac{1}{2 } \tilde{\Sigma} \left[ \alpha  \sum_{i \in [N]} (2 \mathbb{E}_{Q_2}[\omega_i] \alpha \eta  + 1 ) H(T_i)^T  +  \alpha \int_{0}^T \int_0^{\infty}  (2 \bar \omega \alpha \eta  - 1 ) H(t)^T \Lambda(t,\bar \omega) d\bar \omega dt + 2\Sigma^{-1 } \mu  \right].
\end{align}

For $Q_2(\omega, \bar Z)$, we first note that using \eqref{eq:var_factor_1} and \eqref{eq:joint_dist}, we have $Q_2(\omega, \bar Z) = Q_{21}(\omega) Q_{22} (\bar Z)$. Using the same computation as \cite{donner2019efficient}) Appendices B and D, one can then show that 
\begin{align*}
    Q_{21}(\omega) &= \prod_{i \in [N]} p_{PG}(\omega_i|1, \underline{\lambda}_{T_i}), \\
      \underline{\lambda}_{t} &= \sqrt{ \mathbb{E}_{Q_1}[\Tilde{\lambda}_t(f)^2]} = \alpha^2 \sqrt{ H(t)^T \Tilde \Sigma H(t) + (H(t)^T \Tilde \mu)^2 - 2 \eta H(t)^T \Tilde \mu + \eta^2 }, \quad \forall t \in [0,T],
\end{align*}
and that $Q_{22}$ is  a marked Poisson point process measure on $[0,T] \times \R^+$ with intensity 
\begin{align*}
     \Lambda(t,\bar \omega) &= \lsup e^{\mathbb{E}_{Q_1}[g(\bar \omega, -\Tilde{\lambda}_t(f)]} p_{PG}(\bar \omega; 1,0) = \lsup \frac{\exp(-\frac{1}{2}\mathbb{E}_{Q_1}[\Tilde{\lambda}_t(f)])}{2\cosh \frac{ \underline{\lambda}_{t}(f)}{2}}  p_{PG}(\bar \omega|1,  \underline{\lambda}_{t}(f)) \\
    &=\lsup \sigma(-  \underline{\lambda}_{t}) \exp \left \{ \frac{1}{2} ( \underline{\lambda}_{t}(f) - \mathbb{E}_{Q_1}[\Tilde{\lambda}_t(f)]) \right \} p_{PG}(\bar \omega|1,  \underline{\lambda}_{t}) \\
    \mathbb{E}_{Q_1}[\Tilde{\lambda}_t(f)] &= \alpha (H(t)^T\tilde{\mu} - \eta).
\end{align*}
Therefore, we have that
\begin{align*}
    \mathbb{E}_{Q_1}[\omega_i] = \frac{1}{2   \underline{\lambda}_{T_i}} \tanh \left( \frac{ \underline{\lambda}_{T_i}}{2} \right), \quad  \forall i \in [N].
\end{align*}

\subsection{Analytic formulas of the ELBO}\label{app:comp_elbo}
In this section, we provide the derivation of the evidence lower bound $(ELBO(\hat Q_k^m))_k$ for a mean-field variational distribution $\hat Q_m(f, \bar Z)= \hat Q_1^m(f)\hat Q_2^m(\bar Z)$ in a fixed model $m = (\delta, D)$. For ease of exposition, we drop the subscript $m$ and $k$. From \eqref{eq:elbo_1}, we have
\begin{align*}
    ELBO(\hat Q) &= \mathbb{E}_{\hat  Q} \left[ \log \frac{p(f,\omega,\bar Z, N)}{\hat Q_{1}(f) \hat Q_{2}(\omega, \bar Z)}\right] \\
    &= \mathbb{E}_{\hat Q_{2}} \left[ - \log \hat Q_2(\omega, \bar Z) \right] +  \mathbb{E}_{\hat Q_{2}} \left[\mathbb{E}_{\hat Q_{1}} \left[ \log p(f,\omega,\bar Z,N) \right]\right]  + \mathbb{E}_{\hat Q_1} [ - \log \hat Q_{1}(f)]. \nonumber
\end{align*}

% then define the variational posterior as the mixture
% \begin{align*}
%     \hat Q = \sum_{m \in \mathcal{M}_T}  \hat \gamma_m \hat Q^m.
% \end{align*}
% Using Theorem 3.6 of \cite{ohn2021adaptive} and Proposition \ref{prop:histo}, we know that the adaptive variational posterior concentrates at the expected rate as soon as $\pi_{\delta}(\delta_0) \geq e^{-c_1 T \epsilon_T^2}$. 
% In the model selection approach of  \cite{Zhang2017ConvergenceRO} for variational Bayes, one would consider a variational $\mathcal{V}^m$ class and a variational posterior $\hat Q^m$ for each model $m \in \mathcal{M}_T$, then select the model maximising the evidence lower bound, i.e.,
% \begin{align*}
%     \hat m &= \arg \max_{m} \int_{\mathcal{F}} L_T(f) d\hat Q^m(f) - KL(\hat Q^m||\Pi_{f|m}) + \log \Pi_{\mathcal{M}}(m) =:  \arg \max_{m} ELBO(\hat Q^m) \\
%     &= \arg \max_m \max_{Q \in \mathcal{V}^m} \int_{\mathcal{F}} L_T(f) dQ(f) - KL(Q||\Pi_{f|m}) + \log \Pi_{\mathcal{M}}(m) = \arg \max_m \max_{Q \in \mathcal{V}^m} ELBO(Q).
% \end{align*}
% From \cite{sulem2021bayesian}, we know that the posterior distribution concentrates as soon as $\Pi_{\mathcal{M}}(m_0) \geq e^{-c_1 T \epsilon_T^2}$, where $m_0$ is the true ``model". Adapting Theorem 4.1 from \cite{Zhang2017ConvergenceRO} to our context, we would also obtain that the ``model selection" variational posterior $\hat Q^{\hat m}$ concentrates at the expected rate. 
%To compute \eqref{eq:adapt_vpost}, we thus need to compute $ELBO(\hat Q)$ for each $\delta$ and $D$. 

Now using the notation of Section \ref{sec:aug-mf-vi}, we first note that defining $K(t) := H(t)H(t)^T$, we have that
\begin{align*}
    &\mathbb{E}_{\hat Q_1}[\tilde{\lambda}_{T_i}(f)^2] = tr(K(t) \tilde{\Sigma}) + \tilde{\mu}^T K(t) \tilde{\mu}  \\
    &\mathbb{E}_{\hat Q_1} \left[ \log \mathcal{N}(f;\mu, \Sigma) \right] = - \frac{1}{2} tr(\Sigma^{-1} \Tilde{\Sigma}) - \frac{1}{2} \Tilde \mu^T \Sigma^{-1} \Tilde \mu + \Tilde \mu^T \Sigma^{-1}  \mu - \frac{1}{2} \mu^T \Sigma^{-1}  \mu - \frac{1}{2} \log |2 \pi \Sigma|.
\end{align*}
% We have that
% \begin{align*}
%     ELBO(\hat Q) &= \mathbb{E}_{\hat Q} \left[ \log \frac{p(f,\omega,\bar Z,N)}{\hat Q_1(f) \hat Q_2(\omega, \bar Z)}\right] \\
%     &= \mathbb{E}_{\hat Q_2} \left[ - \log \hat Q_2(\omega, \bar Z) \right] +  \mathbb{E}_{\hat Q_2} \left[  \mathbb{E}_{\hat Q_1} \left[ \log p(f,\omega,\bar Z|N) \right]\right]  + \mathbb{E}_{\hat Q_1} [ - \log \hat Q_1(f)].
% \end{align*}
Moreover, we have
\begin{align*}
     \mathbb{E}_{\hat Q_1} [ \log \hat Q_1(f)] &=  %\mathbb{E}_{\hat Q_1} [ - \frac{1}{2} (f -\tilde{\mu})^T \tilde{\Sigma}^{-1} (f-\tilde{\mu}) ]
     -\frac{|m|}{2} -\frac{1}{2} \log |2\pi \tilde{\Sigma}|.
    %  &=  - \frac{1}{2} \log |\tilde{\Sigma}| - \frac{|m|}{2} ( 1 +  \log 2 \pi) -\Tilde \mu^T \Tilde \Sigma_D^{-1} \Tilde \mu_D \\
    %  &= - \frac{1}{2} \sum_k \log |\tilde{\Sigma}_k| - \frac{|m|}{2} ( 1 +  \log 2 \pi) -\Tilde \mu_D^T \Tilde \Sigma_D^{-1} \Tilde \mu_D.
\end{align*}
Using that for any $c> 0$,
$
    p_{PG}(\omega;1,c) = e^{-c^2 \omega/2} \cosh{(c/2)}  p_{PG}(\omega;1,0), 
$
we also have
\begin{align*}
    \mathbb{E}_{\hat Q_2} \left[ - \log \hat Q_{2}(\omega, \bar Z) \right] &=  \sum_{i \in [N]}  - \mathbb{E}_{\hat Q_2}[\log p_{PG}(\omega_i,1,0)] + \frac{1}{2}\mathbb{E}_{\hat Q_2} [\omega_i]  \mathbb{E}_{\hat Q_1}[\tilde{\lambda}_{T_i}(f)^2] - \log \cosh{\left(\frac{\underline{\lambda}_{T_i}(f) }{2}\right)} \\ %  tr(K(T_i) \tilde{\Sigma}) - \tilde{\mu}^T K(T_i) \\
    &- \int_{t=0}^T \int_0^{+\infty} [\log \Lambda(t,\bar \omega)]  \Lambda(t,\bar \omega)d \bar \omega dt + \int_{t=0}^T \int_0^{+\infty} \Lambda(t,\bar \omega)d \bar \omega dt\\
    &= \sum_{i \in [N]}  - \mathbb{E}_{\hat Q_2}[\log p_{PG}(\omega_i,1,0)] + \frac{1}{2} \mathbb{E}_{\hat Q_2} [\omega_i]  \mathbb{E}_{\hat Q_{1}}[\tilde{\lambda}_{T_i}(f)^2] - \log \cosh{\left(\frac{\underline{\lambda}_{T_i}(f) }{2}\right)}\\ %  tr(K(T_i) \tilde{\Sigma}) - \tilde{\mu}^T K(T_i) \\
    &- \int_{t=0}^T \int_0^{+\infty} 
    \left[\log \lsup -  \frac{1}{2} \mathbb{E}_{\hat Q_1}[\tilde{\lambda}_{T_i}(f)] - \log 2 - \log \cosh{\left(\frac{\underline{\lambda}_{T_i}(f)}{2} \right)} - \frac{1}{2} \mathbb{E}_{\hat Q_{1}}[\tilde{\lambda}_{T_i}(f)^2] \bar \omega  \right. \\
    &+ \left. \log \cosh{ \left(\frac{1}{2} \underline{\lambda}_{T_i}(f) \right)}  + \log p_{PG}(\bar \omega;1,0) - 1 \right]   \Lambda(t) p_{PG}(\bar \omega;1,\underline{\lambda}_{T_i}(f)) dt d\bar \omega \\
        &= \sum_{i \in [N]}   - \mathbb{E}_{\hat Q_2}[\log p_{PG}(\omega_i,1,0)] +  \frac{1}{2} \mathbb{E}_{\hat Q_2} [\omega_i^k]  \mathbb{E}_{\hat Q_{1}}[\tilde{\lambda}_{T_i}(f)^2] - \log \cosh{\left(\frac{\underline{\lambda}_{T_i}(f) }{2} \right)}\\ %  tr(K(T_i) \tilde{\Sigma}) - \tilde{\mu}^T K(T_i) \\
    &- \int_{t=0}^T 
    \left[\log \lsup -  \frac{1}{2} \mathbb{E}_{\hat Q_1}[\tilde{\lambda}_{T_i}(f)] - \log 2  - \frac{1}{2} \mathbb{E}_{\hat Q_1}[\tilde{\lambda}_{T_i}(f)^2] \mathbb{E}_{\hat Q_2}[\bar \omega ] - 1 \right] \Lambda(t) dt \\
    &- \int_{t=0}^T \int_0^{+\infty} \log p_{PG}(\omega;1,0) \Lambda(t) p_{PG}(\omega;1,\underline{\lambda}_{T_i}(f)) d\omega dt. 
    %      &= \sum_{T_i \in N}  \mathbb{E}_{\hat Q^m_2}[- \log p_{PG}(\omega_i,1,0)]] - \frac{\tanh \sqrt{tr(K(T_i) \tilde{\Sigma}) + \tilde{\mu}^T K(T_i)} /2}{2\sqrt{tr(K(T_i) \tilde{\Sigma}) + \tilde{\mu}^T K(T_i)}} -  tr(K(T_i) \tilde{\Sigma}) - \tilde{\mu}^T K(T_i) \\
    % &- \int_{t=0}^T \int \log \bar \Lambda(t,\bar \omega) \bar \Lambda(t,\bar \omega)d \bar \omega dt
    %  &= \sum_{T_i \in N}  \mathbb{E}_{\hat Q^m_2}[- \log p_{PG}(\omega_i,1,0)]] - \frac{\tanh \sqrt{tr(K(T_i) \tilde{\Sigma}) + \tilde{\mu}^T K(T_i)} /2}{2\sqrt{tr(K(T_i) \tilde{\Sigma}) + \tilde{\mu}^T K(T_i)}} -  tr(K(T_i) \tilde{\Sigma}) - \tilde{\mu}^T K(T_i) \\
    % &- \int_{t=0}^T \int \log \bar \Lambda(t,\bar \omega) \bar \Lambda(t,\bar \omega)d \bar \omega dt
\end{align*}
with $\Lambda(t) = \theta \int_0^\infty \Lambda(t,\bar \omega) d\bar \omega = \frac{e^{-\frac{1}{2} \mathbb{E}_{\hat Q_1}[\tilde{\lambda}_{T_i}(f)] }}{2 \cosh \frac{\underline{\lambda}_{T_i}(f) }{2}}$. Moreover, we have
\begin{align*}
      &\mathbb{E}_{\hat Q_2} \left[   \mathbb{E}_{\hat Q_1} \left[ \log p(f,\omega,\bar Z,N) \right] \right]    =  \sum_{i \in [N]} \left \{ \log \lsup +  \mathbb{E}_{\hat Q_2} \left[  \mathbb{E}_{\hat Q_1} \left[ g(\omega_i, \tilde{\lambda}_{T_i}(f))  \right] + \log p_{PG}(\omega_i;1,0) \right] \right \} \\
      &+   \log \lsup +  \mathbb{E}_{\hat Q_2} \left[\mathbb{E}_{\hat Q_1} \left[ g(\bar \omega_t, - \tilde{\lambda}_{T_i}(f)))  \right] + \log p_{PG}(\bar \omega_t;1,0) \right] +  \mathbb{E}_{\hat Q_1} \left[ \log \mathcal{N}(f;\mu, \Sigma) \right]  \\
          &=  \sum_{i \in [N]} \log \lsup - \log 2 - \frac{1}{2} \mathbb{E}_{\hat Q_1} \left[ \tilde{\lambda}_{T_i}(f)^2  \right] \mathbb{E}_{\hat Q_2} \left[ \omega_i \right]  + \frac{1}{2} \mathbb{E}_{\hat Q_1} \left[\tilde{\lambda}_{T_i}(f) \right] +  \mathbb{E}_{\hat Q_2} \left[   \log p_{PG}(\omega_i;1,0) \right]  \\
          &+ \int_0^T \int_0^{+\infty} \left[ \log \lsup_k - \log 2 - \frac{1}{2} \mathbb{E}_{\hat Q_1} \left[ \tilde{\lambda}_{T_i}(f)^2  \right] \bar \omega - \frac{1}{2} \mathbb{E}_{\hat Q_1} \left[ \tilde{\lambda}_{T_i}(f) \right] +\log p_{PG}(\bar \omega;1,0)  \right] \Lambda^k(t)p_{PG}(\omega;1,\underline{\lambda}_{T_i}(f)) d\omega dt \\
          & %- \int_{t=0}^T \int_0^{+\infty} \Lambda(t,\bar \omega)d \bar \omega dt
          + \mathbb{E}_{\hat Q_1} \left[ \log \mathcal{N}(f;\mu, \Sigma) \right] - \lsup T \\
        &=  \sum_{i \in [N]} \log \lsup_k - \log 2  - \frac{1}{2} \mathbb{E}_{\hat Q_1} \left[ \tilde{\lambda}_{T_i}(f)^2  \right] \mathbb{E}_{\hat Q_2} \left[ \omega_i \right]  + \frac{1}{2} \mathbb{E}_{\hat Q_1} \left[ \tilde{\lambda}_{T_i}(f)  \right] +  \mathbb{E}_{\hat Q_2} \left[   \log p_{PG}(\omega_i;1,0) \right]   \\
          &+ \int_0^T  \left[ \log \lsup - \log 2 - \frac{1}{2} \mathbb{E}_{\hat Q_1} \left[ \tilde{\lambda}_{T_i}(f)^2  \right]\mathbb{E}_{\hat Q_2} \left[ \bar \omega \right] - \frac{1}{2} \mathbb{E}_{\hat Q_1} \left[\tilde{\lambda}_{T_i}(f) \right]\right] \Lambda(t) dt \\
          &+  \int_0^T \int_0^{+\infty} \log p_{PG}(\bar \omega;1,0) \Lambda(t)p_{PG}(\bar \omega;1,\underline{\lambda}_{T_i}(f)) d\bar \omega dt +  \mathbb{E}_{\hat Q_1} \left[ \log \mathcal{N}(f;\mu, \Sigma) \right] - \lsup T.
\end{align*}
Therefore, with $c>0$ a constant that does not depend on the size of the model, with zero mean prior $\mu = 0$,
\begin{align*}
     ELBO(\hat Q) &= \frac{|m|}{2}  +  \frac{1}{2}  \log |2\pi \tilde{\Sigma}| -  \frac{1}{2}  tr(\Sigma^{-1} \Tilde{\Sigma}) - \frac{1}{2} \Tilde \mu^T \Sigma^{-1} \Tilde \mu -  \frac{1}{2}  \log |2\pi\Sigma| \\
     &+  \sum_{i \in [N]} \log \lsup -\log 2 + \frac{\mathbb{E}_{\hat Q_1} \left[ \tilde{\lambda}_{T_i}(f) \right]}{2} - \log \cosh \left(\frac{\tilde{\lambda}_{T_i}(f) }{2} \right) \\
     &+ \int_{t=0}^T \int_0^{+\infty} \Lambda(t,\bar \omega)d \bar \omega dt - \lsup T.
     %+  \log 2 \int_0^T \Lambda^k(t)  dt - \lsup_k T.
    %   &= - \frac{1}{2} tr(\Sigma^{-1} \Tilde{\Sigma}) + \frac{1}{2} \Tilde \mu_D^T \Sigma_D^{-1} \Tilde \mu_D  + \frac{1}{2} \log ( |\tilde{\Sigma}_D| - |\Sigma_D| ) + \frac{|m|}{2}  \\
    %  &+ \sum_k \sum_{i \in [N_k]} \log \lsup_k + \frac{\mathbb{E}_{\hat Q^D_1} \left[ \Tilde{\lambda}_{T_i^k}(f) \right]}{2} + \int_0^T \int_0^\infty \log \left(2 \cosh{ \mathbb{E}_{\hat Q^D_1}[h(t,f_k)]/2} \right) \Lambda^k(t,\bar \omega)  d\bar \omega dt + c.
\end{align*}

\subsection{Gibbs sampler}\label{app:gibbs_sampler}

From the augmented posterior $\Pi_A(f,\omega, \bar |N)$ defined in  \eqref{eq:augm_posterior}  and using the Gaussian prior family described in Section \ref{sec:aug-mf-vi},  similar computation as  Appendix \ref{app:fixed_updates} can provide analytic forms of the conditional posterior distributions $\Pi_A(f|\omega, \bar Z, N), \Pi_A(\omega | N, f)$ and $ \Pi_A(\bar Z|f, N)$ . This allows to design a Gibbs sampler algorithm that sequentially samples the parameter $f$, the latent variables $\omega$ and Poisson process $\bar Z$. With the notation of Appendix \ref{app:fixed_updates}, such procedure can be defined  as

For every $k \in [K]$,
\begin{align*}
    \text{(Sample latent variables)} \: &\omega^k_i|N,f_k \sim p_{PG}(\omega_i^k; 1, \Tilde{\lambda}^k_{T_i^k}(f)), \quad \forall i \in [N_k] \\
    &\text{$\bar Z^k|f_k$, a Poisson process on $[0,T]$ with intensity }\\& \Lambda^k(t, \bar \omega) = \lsup_k \sigma(- \Tilde{\lambda}^k_{t}(f)) p_{PG}(\bar \omega; 1,\Tilde{\lambda}_{t}^k(f))\\
     \text{(Update hyperparameters)} \: &R_k = \bar{N}^k[0,T] \\
     &H_k = [H_{N^k}, H_{\bar Z^k}], \: [H_{N^k}]_{id} = H_j(T_i^k), \\& [H_{\bar Z^k}]_{jd} = H_b(\bar T_j^k), \: d = 0, \dots, KJ, \: i \in [N_k], \: j \in [R_k] \\
     &D_k = Diag([\omega^k_i]_{i \in [N^k]}, [\bar \omega^k_j]_{j \in [R^k]}) \\
     &\tilde \Sigma_{k} = [\beta^2 H_k D_k (H_k)^T + \Sigma^{-1}]^{-1} \\
     &\tilde \mu_{k} = \Tilde \Sigma_{k} \left( H_k \left[\beta v_k + \beta^2 \eta u_k \right] + \Sigma^{-1} \mu \right),\\& \quad v_k = 0.5 [\mathds{1}_{N_k}, - \mathds{1}_{R_k}], \quad u_k =  [[\omega^k_i]_{i \in [N_k]}, [\bar \omega^k_{j}]_{j \in [R_k]}] \\
     \text{(Sample parameter)} %\quad &\pi(\lsup_k|N, \bar Z) = Gamma(\lsup_k; a_0 + N^k + R^k, b_0 + T) \\
     \: &f_{k}|N,\bar Z^k, \omega^k \sim \mathcal{N}(f_k;  \tilde m_{k}, \Tilde \Sigma_{k}).
\end{align*}
These steps are summarised in Algorithm \ref{alg:gibbs} in Appendix. We note that in this algorithm, one does not need to perform a numerical integration, however, sampling the latent Poisson process is computationally intensive. In our numerical experiments, we use the Python package \texttt{polyagamma}\footnote{\url{https://pypi.org/project/polyagamma/}} to sample the Polya-Gamma variables and a thinning algorithm to sample the inhomogeneous Poisson process.

\section{Proofs}\label{sec:proofs}

In this section, we provide the proof of our main theoretical result, namely Theorem \ref{thm:cv_rate_vi}. % and Proposition \ref{prop:conc_rate_sigmoid}. 
We first recall a set of useful lemmas from \cite{sulem2021bayesian}.

\subsection{Technical lemmas}\label{app:main_lemmas}

In the first lemma, we recall the definition of excursions from \cite{sulem2021bayesian}, for stationary nonlinear Hawkes processes verifying conditions (C1) or (C2). Then, Lemma \ref{lem:main_event}, corresponding to Lemma A.1 in \cite{sulem2021bayesian}, provides a control on the main event $\Tilde{\Omega}_T$ considered in the proof of Theorem  \ref{thm:cv_rate_vi}. Finally, Lemma \ref{lem:ef} (Lemma A.4 in \cite{sulem2021bayesian}) is a technical lemma for proving posterior concentration in Hawkes processes.

We also introduce the following notation. For any excursion index $j \in [J_T-1]$, we denote $(U_j^{(1)}, U_j^{(2)})$ the times of the  first two events after the $j$-th renewal time $\tau_j$, and $\xi_j :=  U_j^{(2)}$ if $U_j^{(2)} \in [\tau_j,\tau_{j+1})$ and $\xi_j := \tau_{j+1} $ %= U_j^{(1)} + A$
otherwise. %We note that the interval $[\tau_j, \xi_j]$ corresponds either to the beginning of the $j$-th excursion or to the whole excursion $[\tau_j, \tau_{j+1})$ when the latter contains only one event, implying that $U_j^{(2)} \geq  \tau_{j+1}$.

\begin{lemma}[Lemma 5.1 in \cite{sulem2021bayesian}]\label{lem:excursions}
Let $N$ be a Hawkes process with monotone non-decreasing and Lipschitz link functions $\phi = (\phi_k)_k$ and parameter $f = (\nu, h)$ such that $(\phi, f)$ verify \textbf{(C1)} or \textbf{(C2)}.
% and  such that $\Exz{N[-A,0)} < +\infty$.
Then the point process measure $X_t(.)$ defined as
\begin{equation}\label{eq:pp_measure_x}
    X_t(.) = N|_{(t-A,t]},
\end{equation}
is a strong Markov process with positive recurrent state $\emptyset$. Let $\{\tau_j\}_{j\geq 0}$ be the sequence of random times defined as
\begin{small}
\begin{align*}
    \tau_j = \begin{cases}
    0 & \text{ if } j=0; \\ 
    \inf \left \{t > \tau_{j-1}; \: X_{t^-} \neq \emptyset, \: X_{t} = \emptyset \right \}  = \inf \left \{t > \tau_{j-1}; \: N|_{[t-A,t)} \neq \emptyset, \: N|_{(t-A,t]} = \emptyset \right \} & \text{ if } j\geq 1 .
    \end{cases}
\end{align*}
\end{small}
Then, $\{\tau_j\}_{j\geq 0}$ are stopping times for the process $N$. For $T > 0$, we also define 
\begin{equation}\label{def:J_T}
    J_T=\max\{j\geq 0;\: \tau_j \leq T\}.
\end{equation}
The intervals $\{[\tau_j, \tau_{j+1})\}_{j=0}^{J_{T}-1} \cup [\tau_{J_T}, T]$ form a partition of $[0,T]$. The point process measures $(N|_{[\tau_j, \tau_{j+1})})_{1 \leq j \leq J_T - 1}$ are i.i.d. and independent of $N|_{[0, \tau_1)}$ and $N|_{[\tau_{J_T},T]}$; they are called \emph{excursions} and the stopping times $\{\tau_j\}_{j\geq 1}$ are called \emph{regenerative} or \emph{renewal} times. 
\end{lemma}

\begin{lemma}[Lemma A.1 in \cite{sulem2021bayesian}]\label{lem:main_event}
Let $Q > 0$. We consider $\Tilde{\Omega}_T$ defined in Section~\ref{sec:proof_main_thm}. For any $\beta > 0$, we can choose $C_\beta$ and $c_\beta$ in the definition of $\Tilde{\Omega}_T$ such that
$
	\probz{\Tilde{\Omega}_T^c} \leq T^{-\beta}.
$
Moreover, for any $1 \leq q \leq Q$,
$$
\Exz{\mathds{1}_{\eve^c} \max_l \sup \limits_{t \in [0,T]} \left(N^l[t-A,t)\right)^q} \leq 2 T^{-\beta/2}. 
$$
\end{lemma}

\begin{lemma}[Lemma A.4 in \cite{sulem2021bayesian}]\label{lem:ef}
For any  $f \in \mathcal{F}_T$ and $l \in [K]$, let 
\begin{equation*}
    Z_{1l} = \int_{\tau_1}^{\xi_1} |\lambda^l_t(f) - \lambda^l_t(f_0)|dt.
\end{equation*}
Under the assumptions of Theorem \ref{thm:cv_rate_vi}, for $M_T \to \infty$ such that $M_T > M \sqrt{\kappa_T}$ with $M>0$ and for any $f \in \mathcal{F}_T$ such that $\norm{r-r_0}_1 \leq \max(\norm{r_0}_1, \Tilde{C})$ with $\Tilde{C}>0$,
% $\Tilde d_{1T}(f, f_0) \leq M_T' \epsilon_T$,
there exists $l \in [K]$ such that on $\Tilde{\Omega}_{T}$,
    \begin{equation*}
        \Exf{Z_{1l}} \geq C(f_0) \Big(\norm{r_f - r_0}_1 + \norm{h - h_0}_1\Big),
    \end{equation*}
with $C(f_0) > 0$ a constant that depends only on $f_0$ and $(\phi_k)_k$.
\end{lemma}

\subsection{Proof of Theorem \ref{thm:cv_rate_vi}}\label{sec:proof_main_thm}

We recall that in this result, we consider a general Hawkes model with  known link functions $(\phi_k)_k$. Let  $r_0 = (r_1^0, \dots, r_K^0)$ with $r_k^0 = \phi_k(\nu_k^0)$. With $C_\beta, c_\beta > 0$, we first define $\eve \in \mathcal{G}_T$ as
\begin{align*}
    \eve &= \Omega_N \cap \Omega_J \cap \Omega_U, \\
   \Omega_N &= \left \{ \max \limits_{k \in [K]} \sup \limits_{t \in [0,T]} N^k[t-A,t) \leq C_\beta \log T  \right \} \cap \left \{ \sum_{k=1}^K \left|\frac{N^k[-A,T]}{T} - \mu_k^0\right| \leq \delta_T  \right \}, \\
        \Omega_{J} &= \left\{ J_T \in \mathcal{J}_T \right \}, \quad \Omega_{U} =  \left\{ \sum_{j=1}^{J_T-1} (U_j^{(1)} - \tau_j) \geq 
     \frac{T}{\mathbb{E}_0[\Delta \tau_1] \|r_0\|_1} \left(1 - 2c_\beta\sqrt{\frac{\log T }{T}}\right) \right \}, \\
     \mathcal{J}_T &= \left \{ J \in \N; \: \left|\frac{J-1}{T} - \frac{1}{\mathbb{E}_0[\Delta \tau_1]} \right| \leq c_\beta \sqrt{\frac{\log T}{T}} \right \},
\end{align*}
with $J_T$ the number of excursions as defined in \eqref{def:J_T}, $\mu_k^0 := \Exz{\lambda_t^k(f_0)}, \forall k$, $\delta_T = \delta_0 \sqrt{\frac{\log T}{T}}, \: \delta_0 > 0$ and $\{U_j^{(1)}\}_{j=1, \dots, J_T-1}$ denoting the first events of each excursion (see Lemma \ref{lem:excursions} for a precise definition). Secondly, we define $A_T' \in \mathcal{G}_T$ as
\begin{align*}
   A_T' = \left \{\int e^{L_T(f) - L_T(f_0)} d\widetilde{\Pi}(f) >  e^{- C_1 T \e_T^2} \right \}, \quad \widetilde{\Pi}(B) = \frac{\Pi(B \cap K_T)}{\Pi(K_T)}, \quad K_T \subset \mathcal{F},
\end{align*}
with $C_1 > 0$ and $\e_T, M_T$ positive sequences such that $T\e_T^2 \to \infty$ and  $M_T \to \infty$. From Lemma \ref{lem:main_event}, we have that $\Probz{\eve^c} = o(1)$.  Thus, with $D_T$ defined in \eqref{def:pposterior_dist}, $A_T = \eve \cap A_T'$, $K_T = B_\infty(\epsilon_T)$, and $\e_T = \sqrt{\kappa_T } \epsilon_T$,  we  obtain that
\begin{align*}
    \Probz{A_T^c} &\leq \Probz{\eve^c} + \Probz{A_T'^c \cap \eve}\\  
    &= o(1) + \Probz{ \left \{\int_{K_T} e^{L_T(f) - L_T(f_0)} d\Pi(f) \leq \Pi(K_T) e^{- C_1 T \e_T^2} \right\} \cap \eve} \\
    &\leq o(1) + \Probz{ \left \{ D_T \leq \Pi(K_T) e^{- C_1 T \e_T^2} \right\} \cap \eve} = o(1),
\end{align*}
with $C_1 > 1$, using (A0), i.e., $\Pi(K_T) \geq e^{-c_1 T \e_T^2}$, and the following intermediate result from the proof of Theorem 3.2 in \cite{sulem2021bayesian}
\begin{align*}
    \Probz{\left \{ D_T \leq \Pi(B_\infty(\epsilon_T)) e^{- \kappa_T T \e_T^2} \right \} \cap \eve} = o(1).
\end{align*}
Therefore, we can conclude that
$$\Probz{A_T} \xrightarrow[T \to \infty]{} 1.$$
We now define the stochastic distance $\Tilde{d}_{1T}$ and stochastic neighborhoods around $f_0$ as
\begin{align}\label{def:stoch_dist}
        &\Tilde{d}_{1T}(f,f') = \frac{1}{T} \sum_{k=1}^K \int_0^T \mathds{1}_{A_{2}(T)}(t) |\lambda_{t}^k(f) - \lambda_{t}^k(f')| dt, \quad A_2(T) = \bigcup_{j=1}^{J_T - 1} [\tau_j, \xi_j] \\
           &A_{d_1}(\e) = \left \{f \in \mathcal{F}; \: \Tilde{d}_{1T}(f,f_0) \leq \e \right \}, \quad \e > 0, \nonumber
\end{align}
where for each $j \in [J_T]$, $ U_j^{(2)}$ is the first event after $ U_j^{(1)}$, and $\xi_j :=  U_j^{(2)}$ if $U_j^{(2)} \in [\tau_j,\tau_{j+1})$ and $\xi_j := \tau_{j+1} $ otherwise. Let $\eta_T$ be a positive sequence and $\hat Q$ be the variational posterior as defined in \eqref{eq:var_posterior}. We have
% \begin{align*}
%      \Exz{\hat Q( A_{d_1}(\eta_T)^c)} &\leq \Probz{A_T^c} + \frac{\Exz{\mathds{1}_{A_T }\hat{Q}(\Tilde{d}_{1T}(f,f_0))}}{\eta_T}
% \end{align*}
\begin{align}\label{eq:q_decomp}
     \Exz{\hat Q( A_{d_1}(\eta_T)^c)} &\leq \Probz{A_T^c} + \Exz{\hat Q( A_{d_1}(\eta_T)^c) \mathds{1}_{A_T}}.
\end{align}
We first bound the second term on the RHS of \eqref{eq:q_decomp} using the following technical lemma, which is an adaptation of Theorem 5 of \cite{Ray_2021} and Lemma 13 in \cite{dennis21}.

\begin{lemma}\label{lem:var_inequality}
 Let $B_T \subset \mathcal{F}$, $A_T \in \mathcal{G}_T$, and $Q$ be a distribution on $\mathcal{F}$. If there exist $C, u_T > 0$ such that 
 \begin{align}\label{eq:hyp_post}
     \Exz{\Pi(B_T|N) \mathds{1}_{A_T}} \leq C e^{-u_T},
 \end{align}
 then, we have that
 \begin{align*}
     \Exz{Q(B_T) \mathds{1}_{A_T}} \leq \frac{2}{u_T} \left( \Exz{KL(Q||\Pi(.|N)) \mathds{1}_{A_T}} + C e^{-u_T/2} \right).
 \end{align*}
\end{lemma}
\begin{proof}
We follow the proof of \cite{Ray_2021} and use the fact that, for any $g: \mathcal{F} \to \R$ such that $\int_{\mathcal{F}} e^{g(f)} d\Pi(f|N) < +\infty$, it holds true that
\begin{align*}
     \int_{\mathcal{F}}g(f) dQ(f) \leq KL(Q||\Pi(.|N)) + \log \int_{\mathcal{F}} e^{g(f)}\Pi(f|N).
\end{align*}
Applying the latter inequality with $g = \frac{1}{2} u_T \mathds{1}_{B_T}$, we obtain
\begin{align*}
  \frac{1}{2} u_T Q(B_T) &\leq  KL(Q||\Pi(.|N)) + \log (1 + e^{\frac{1}{2} u_T} \Pi(B_T|N)) \\
  &\leq KL(Q||\Pi(.|N)) + e^{\frac{1}{2} u_T} \Pi(B_T|N).
\end{align*}
Then, multiplying both sides of the previous inequality by $\mathds{1}_{A_T}$ and taking expectation w.r.t. to $\mathbb{P}_0$, using \eqref{eq:hyp_post}, we finally obtain
\begin{align*}
   \frac{1}{2} u_T   \Exz{Q(B_T) \mathds{1}_{A_T}} \leq  \Exz{KL(Q||\Pi(.|N)) \mathds{1}_{A_T}} + C e^{-\frac{1}{2}u_T}.
\end{align*}
\end{proof}
We thus apply Lemma \ref{lem:var_inequality} with $B_T = A_{d_1}(\eta_T)^c$, $\eta_T = M_T' \e_T$, $Q = \hat Q$, and $u_T = M_T T \e_T^2$ with $M_T' \to \infty$. We first check that \eqref{eq:hyp_post} holds, i.e., we show that there exist $C, M_T,M_T' > 0$ such that
\begin{align}\label{eq:cond_post}
     \Exz{\mathds{1}_{A_T }\Pi[ \Tilde{d}_{1T}(f,f_0) > M_T' \e_T |N]} \leq C \exp(-M_T T\e_T^2).
\end{align}
For any test $\phi$, we have the following decomposition
\begin{align*}
     \Exz{\mathds{1}_{A_T}\Pi[ \Tilde{d}_{1T}(f,f_0)   > M_T' \e_T |N]}  \leq  \underbrace{\Exz{\phi \mathds{1}_{A_T} ] }}_{(I)} +  \underbrace{\Exz{(1 - \phi)\mathds{1}_{A_T }\Pi[A_{d_1}(M_T' \e_T)^c|N]}}_{(II)}. %+  \underbrace{\Exz{\mathds{1}_{A_T}\Pi[A_{L_1}(C\e)^c \cap A_{d_1}(M \e) |N]}}_{(III)}
\end{align*}
Note that we have
\begin{align}\label{eq:upper_bound_error2}
    (II) = \Exz{(1 - \phi)\mathds{1}_{A_T }\Pi[A_{d_1}(M_T' \e_T)^c|N]} &= \Exz{\int_{ A_{d_1}(M_T' \e_T)^c } \mathds{1}_{A_T} (1-\phi) \frac{e^{L_T(f) - L_T(f_0)}}{D_T} d\Pi(f)} \nonumber \\
    &\leq \frac{e^{C_1 T \e_T^2}}{\Pi(K_T)} \Exz{ \sup_{f \in \mathcal{F}_T} \Exf{\mathds{1}_{A_{d_1}(M_T'\e_T)^c} \mathds{1}_{A_T} (1-\phi)|\mathcal{G}_0}},
\end{align}
since on $A_T, D_T \geq \Pi(K_T)e^{-C_1 T \e_T^2} $. Using the proof of Theorem 5.5 in \cite{sulem2021bayesian}, we can directly obtain that for $T$ large enough, there exist $x_1, M, M' > 0$ such that
\begin{align*}
    &(I) \leq  2(2K+1) e^{-x_1 {M'_T}^2 T \e_T^2} \\
    &(II) \leq 2 (2K+1) e^{-x_1 {M'_T}^2 T \e_T^2 /2},
\end{align*}
which implies that
\begin{align*}
    \Exz{\mathds{1}_{A_T }\Pi[\Tilde{d}_{1T}(f,f_0)   > M_T' \e_T |N]}  &\leq  4 (2K+1)   e^{-x_1 M_T'^2 T \e_T^2 /2},
\end{align*}
and \eqref{eq:cond_post} with $M_T = x_1 M_T'^2/2$ and $C = 4(2K+1)$. Applying Lemma \ref{lem:var_inequality} thus leads to
\begin{align*}
     \Exz{\hat Q( A_{d_1}(\eta_T)^c) \mathds{1}_{A_T}} \leq 2 \frac{KL(\hat Q || \Pi(.|N)) + Ce^{-M_T T\e_T^2/2}}{M_T T\e_T^2} \leq 2C e^{-M_T T \e_T^2/2} + 2 \frac{KL(\hat Q || \Pi(.|N))}{M_T T\e_T^2}.
\end{align*}
Moreover, from (A2) and the remark following Theorem~\ref{thm:cv_rate_vi}, it holds that $KL(\hat Q || \Pi(.|N)) = O(T\e_T^2) $, therefore we obtain the following intermediate result
\begin{align*}
      \Exz{\hat Q( A_{d_1}(\eta_T)^c) } = o(1).
\end{align*}
Now,  with $M_T >  M_T'$, we note that
\begin{align*}
    \Exz{\hat Q (\norm{f-f_0}_1 > M_T \e_T)} &= \Exz{\hat Q (\Tilde{d}_{1T}(f,f_0) > M_T' \e_T)}\\ &\hspace{0.5cm}+ \Exz{\hat Q (\norm{f-f_0}_1 > M_T \e_T ,\Tilde{d}_{1T}(f,f_0) < M_T' \e_T) \mathds{1}_{A_T}} + \probz{A_T^c}.
\end{align*}
Therefore, it remains to show that
\begin{small}
\begin{align*}
    \Exz{\hat Q (\norm{f-f_0}_1 > M_T \e_T ,\Tilde{d}_{1T}(f,f_0) < M_T' \e_T) \mathds{1}_{A_T}} =  \Exz{\hat Q( A_{L_1}( M_T \epsilon_T)^c \cap A_{d_1}( M_T' \e_T)) \mathds{1}_{A_T}} = o(1).
\end{align*}
\end{small}
For this, we apply again Lemma \ref{lem:var_inequality} with $B_T = A_{L_1}( M_T \e_T)^c \cap A_{d_1}( M_T' \e_T)$ and $u_T = T  M_T^2 \e_T^2$. We have
\begin{align*}
     \Exz{\mathds{1}_{A_T} \Pi(A_{L_1}( M_T \e_T)^c \cap A_{d_1}( M_T' \e_T)|N)}  &\leq \frac{e^{C_1 T \e_T^2}}{\Pi(K_T)} \Exz{\Exf{ \mathds{1}_{A_T}\mathds{1}_{A_{L_1}(M_T \e_T)^c \cap  A_{d_1}(M_T' \epsilon_T) }| \mathcal{G}_0}}.
\end{align*}
Let $f \in A_{L_1}(M_T \e_T)^c \cap A_{d_1}(M_T' \e_T)$. For any $j \in [J_T-1]$ and $l \in [K]$, let
\begin{align}\label{def:zj}
    Z_{jl} = \int_{\tau_j}^{\xi_j} |\lambda^l_t(f) - \lambda^l_t(f_0)|dt, \quad j \in [J_T-1], \quad  l \in [K].
\end{align}

% We first state a technical  lemma corresponding to Lemma A.4 in \cite{sulem2021bayesian}.
% \begin{lemma}\label{lem:A4}
% For any $f \in A_{d_1}(M_T' \epsilon_T)$, there exists $l \in [K]$ such that
% \begin{align*}
%     \Exf{Z_{1l}} \geq C(f_0) \norm{f-f_0}_1.
% \end{align*}
% \end{lemma}

Using Lemma \ref{lem:ef} and the integer $l$ introduced in this lemma, for any $f \in A_{L_1}(M_T \epsilon_T)^c$,  we have
\begin{align*}
    \Exf{ \mathds{1}_{A_T}\mathds{1}_{  A_{d_1}(M_T' \e_T) } | \mathcal{G}_0} &\leq  \Probf{\sum_{j=1}^{J_T-1} Z_{jl} \leq T M_T' \e_T | \mathcal{G}_0} \\
    &\leq \sum_{J \in \mathcal{J}_T} \Probf{\sum_{j=1}^{J-1} Z_{jl} - \Exf{Z_{jl}} \leq T M_T' \epsilon_T - \frac{T}{2\Exz{\Delta \tau_1}} C(f_0) M_T \epsilon_T | \mathcal{G}_0} \\
    &\leq \sum_{J \in \mathcal{J}_T} \Probf{\sum_{j=1}^{J-1} Z_{jl} - \Exf{Z_{jl}} \leq - \frac{T}{4\Exz{\Delta \tau_1}}  C(f_0) M_T \e_T | \mathcal{G}_0},
\end{align*}
for any $M_T \geq 4\Exz{\Delta \tau_1} M_T'$. Similarly to the proof of Theorem 3.2 in \cite{sulem2021bayesian}), we apply Bernstein's inequality for each $J \in \mathcal{J}_T$ and obtain that
\begin{align*}
    \Exf{ \mathds{1}_{A_T}\mathds{1}_{  A_{d_1}(M_T' \e_T) } | \mathcal{G}_0} \leq \exp\{-c(f_0)' T\}, \quad \forall f \in A_{L_1}(M_T \e_T)^c, 
\end{align*}
for $c(f_0)'$ a positive constant. Therefore, we can conclude that
\begin{align*}
     \Exz{\hat Q \left( A_{L_1}( M_T \e_T)^c \cap A_{d_1}( M_T' \e_T)\right) \mathds{1}_{A_T}} \leq \frac{2}{M_T T \e_T^2} \Exz{KL(\hat Q||\Pi(.|N)) } + o(1) = o(1),
\end{align*}
since $ \Exz{KL(\hat Q||\Pi(.|N)) } = O(T \e_T^2)$ by assumption (A2). This leads to our final conclusion
\begin{align*}
     \Exz{\hat Q \left( \norm{f-f_0}_1 > M_T \e_T \right) } = o(1).
\end{align*}

\section{Gibbs sampler in the sigmoid Hawkes model}\label{app:gibbs}

In this section, we describe a non-adaptive Gibbs sampler that computes the posterior distribution in the sigmoid Hawkes model, using the data augmentation scheme of Section \ref{sec:DA_sigmoid} (see also Remark \ref{rem:gibbs}).

\begin{algorithm}
\caption{Gibbs sampler in the sigmoid Hawkes model with data augmentation}\label{alg:gibbs} 
%\begin{algorithmic}
\KwIn{ $N = (N^1, \dots, N^K)$, $n_{iter}$,  $\mu, \Sigma$.}
\KwOut{ Samples $S = (f_i)_{i\in [n_{iter}]}$ from the posterior distribution $\Pi_A(f|N)$.}
Precompute $(H_k(T_i^k))_i, k \in [K]$. \\
Initialise $f \sim \mathcal{N}(f,\mu, \Sigma)$ and $S = []$. \\
\For{$t \gets 1$ to $n_{iter}$}{ 
    \For{$k \gets 1$ to $K$}{
        \For{$i \gets 1$ to $N_k$}{
        Sample $\omega_i^k \sim p_{PG}(\omega_i^k; 1, \Tilde{\lambda}^k_{T_i^k}(f))$
        }
        Sample $(\bar T_j^k)_{j=1,R_k}$ a Poisson temporal point process on $[0,T]$ with intensity $\lsup_k \sigma(- \Tilde{\lambda}^k_{t}(f))$ \\
        \For{$j \gets 1$ to $R_k$}{
        Sample $\bar \omega_j^k \sim p_{PG}(\omega; 1,\Tilde{\lambda}^k_{\bar T_j^k}(f))$
        }
        Update $\tilde \Sigma_{k} = [\beta^2 H_k D_k (H_k)^T + \Sigma^{-1}]^{-1}$     \\
         Update $\tilde \mu_{k} = \Tilde \Sigma_{k} \left( H_k \left[\beta v_k + \beta^2 \eta u_k \right] + \Sigma^{-1} \mu \right)$ \\
         Sample $f_k \sim \mathcal{N}(f_k;  \tilde \mu_{k}, \Tilde \Sigma_{k})$
    }
    Add $f = (f_k)_k$ to $S$.
}
%\end{algorithmic}
\end{algorithm}

\section{Additional results from our numerical experiments}\label{app:details_exp}

%\subsection{Hyperparameters}

%We approximate integrals using the Gaussian quadrature method with  $n_{GQ} = 2^{D+1}T/A$ points in the univariate settings (Simulation 1,2 and 3). In Simulation 4, we set $n_{GQ}$ to reduce the computational time.

In this section, we report results from our simulation study in Section \ref{sec:numerical} that were not added to the main text for conciseness purposes. Each of the following sub-sections corresponds to one of the simulation set-up.

\subsection{Simulation 1}\label{app:simu1}

\begin{figure}
\setlength{\tempwidth}{.3\linewidth}\centering
\settoheight{\tempheight}{\includegraphics[width=\tempwidth, trim=0.cm 0.cm 0cm  0.65cm,clip]{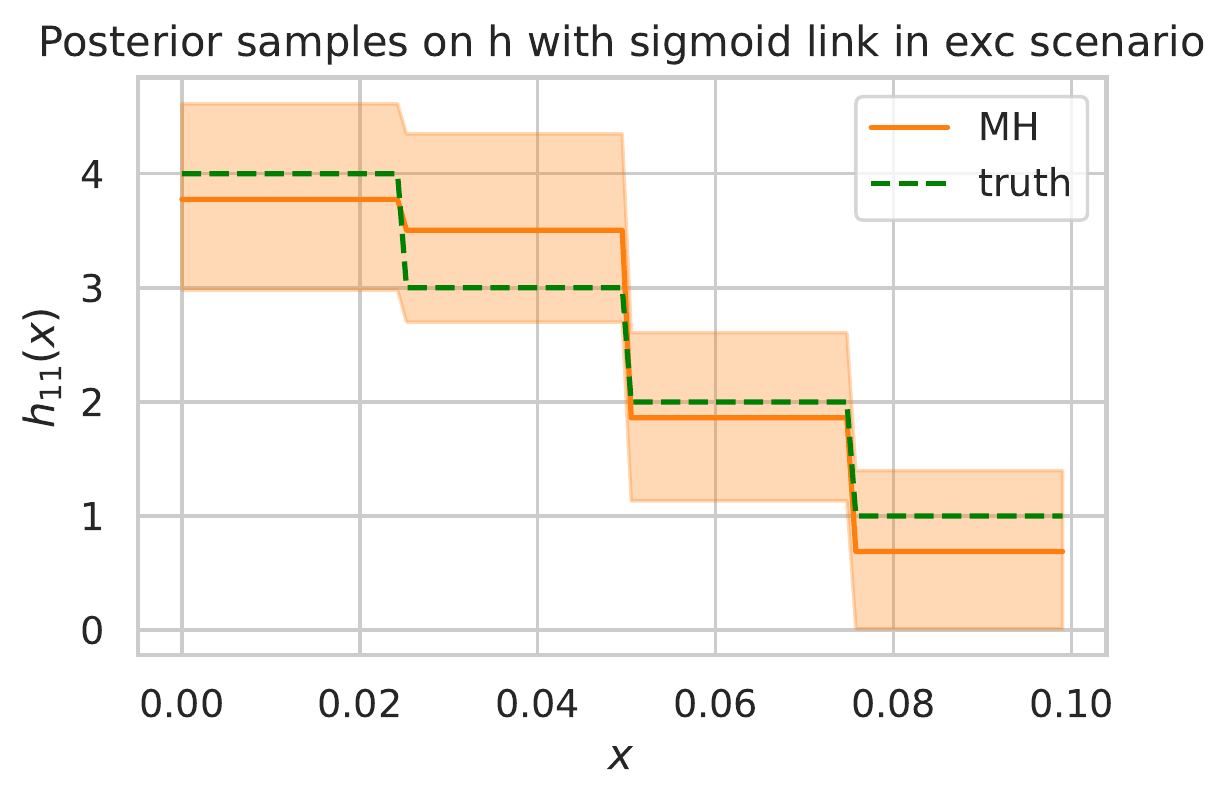}}%
\hspace{-5mm}
\fbox{\begin{minipage}{\dimexpr 15mm} \begin{center} \itshape \large  \textbf{Sigmoid} \end{center} \end{minipage}}
\hspace{-5mm}
\columnname{Excitation}\hfil
\columnname{Mixed}\hfil
\columnname{Inhibition}\\
\rowname{Background}
\includegraphics[width=\tempwidth, trim=0.cm 0.cm 0cm  0.65cm,clip]{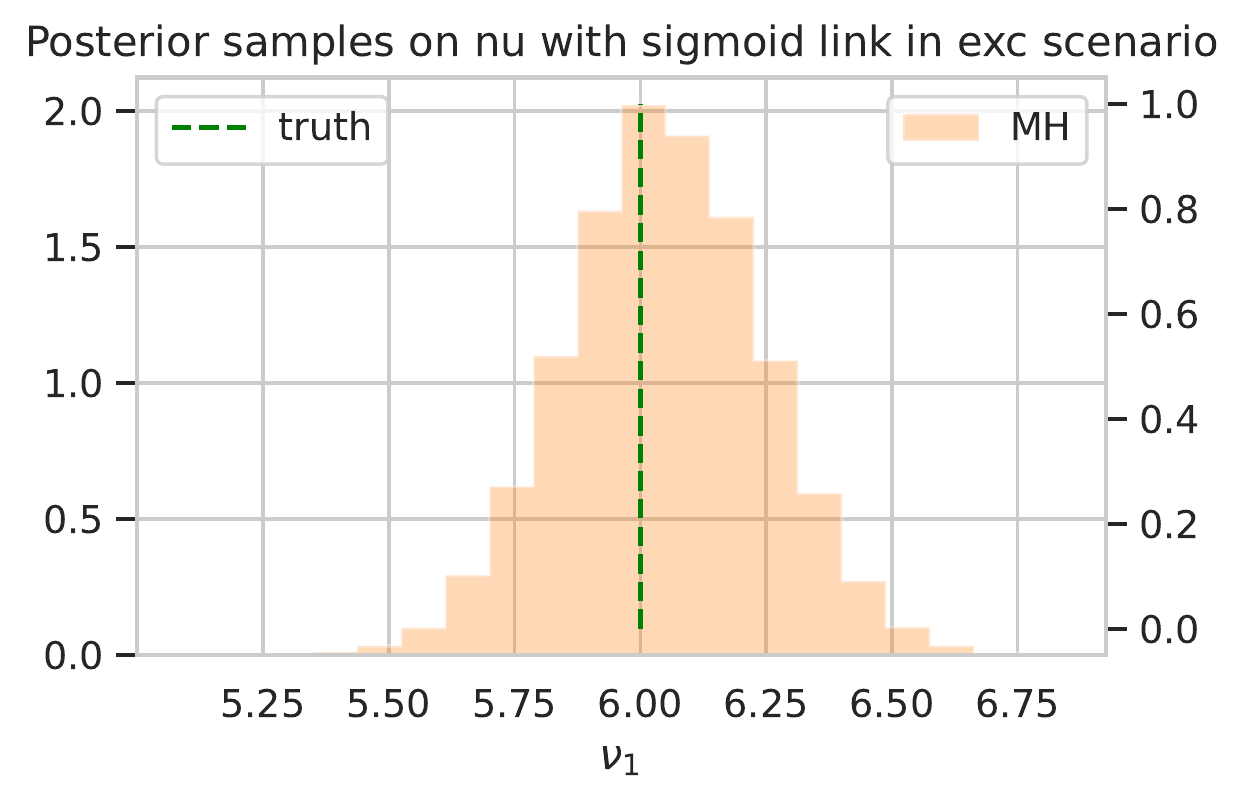}\hfil
\includegraphics[width=\tempwidth, trim=0.cm 0.cm 0cm  0.65cm,clip]{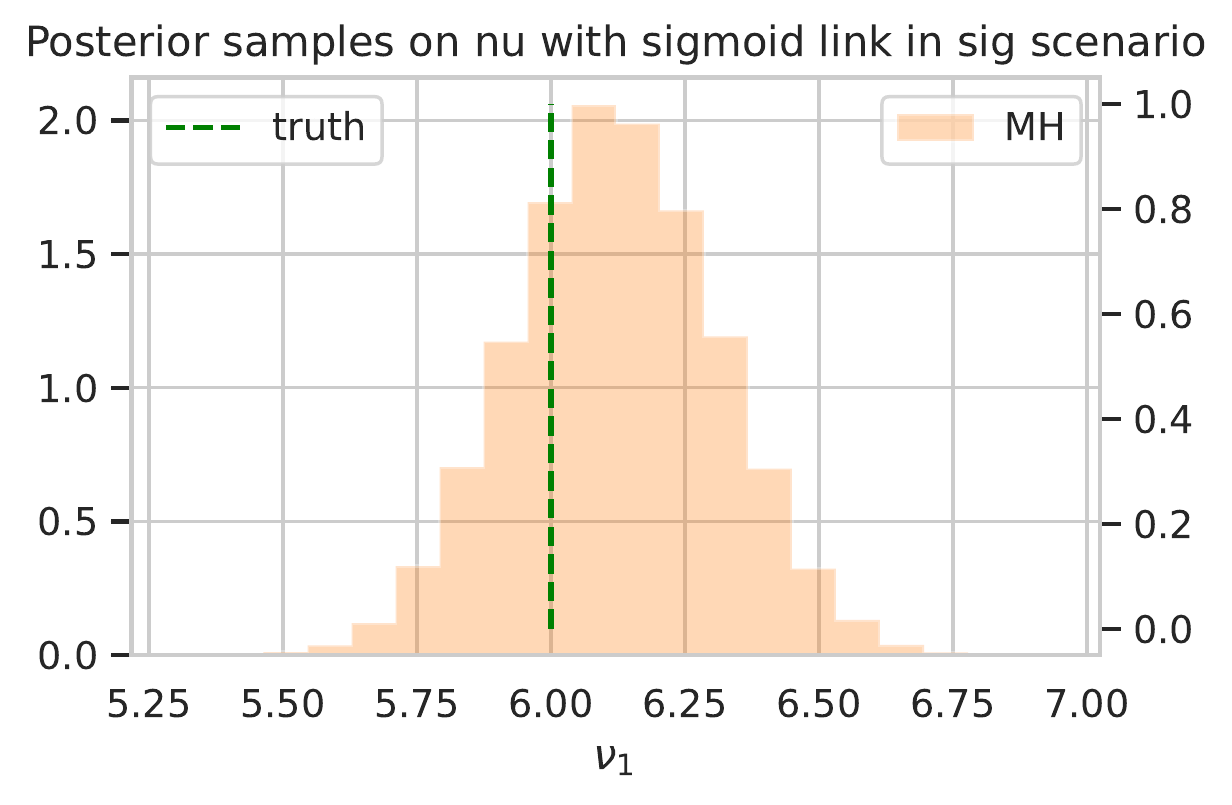}\hfil
\includegraphics[width=\tempwidth, trim=0.cm 0.cm 0cm 0.65cm,clip]{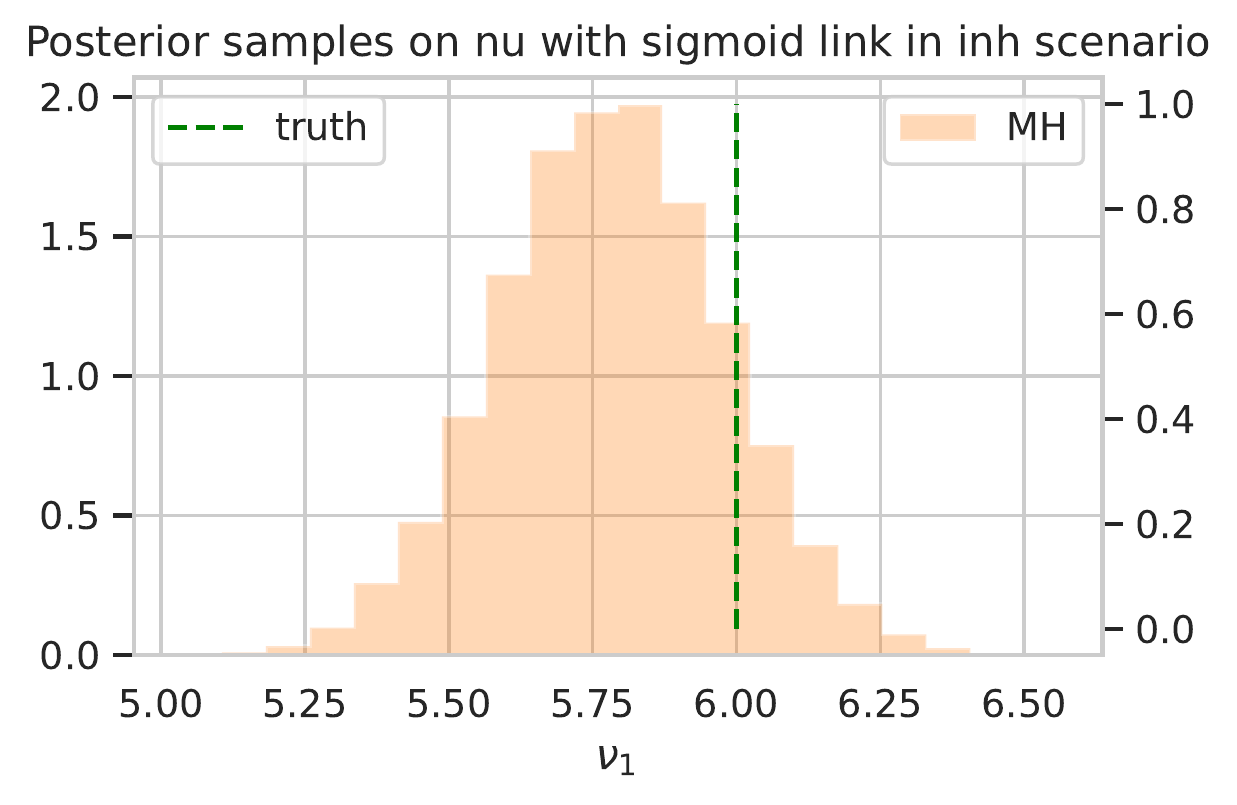}\\
\rowname{Interaction}
\includegraphics[width=\tempwidth, trim=0.cm 0.cm 0cm  0.65cm,clip]{figs/sigmoid_d4_exc_h_mh_True_hmc_False.pdf}\hfil
\includegraphics[width=\tempwidth, trim=0.cm 0.cm 0cm 0.65cm,clip]{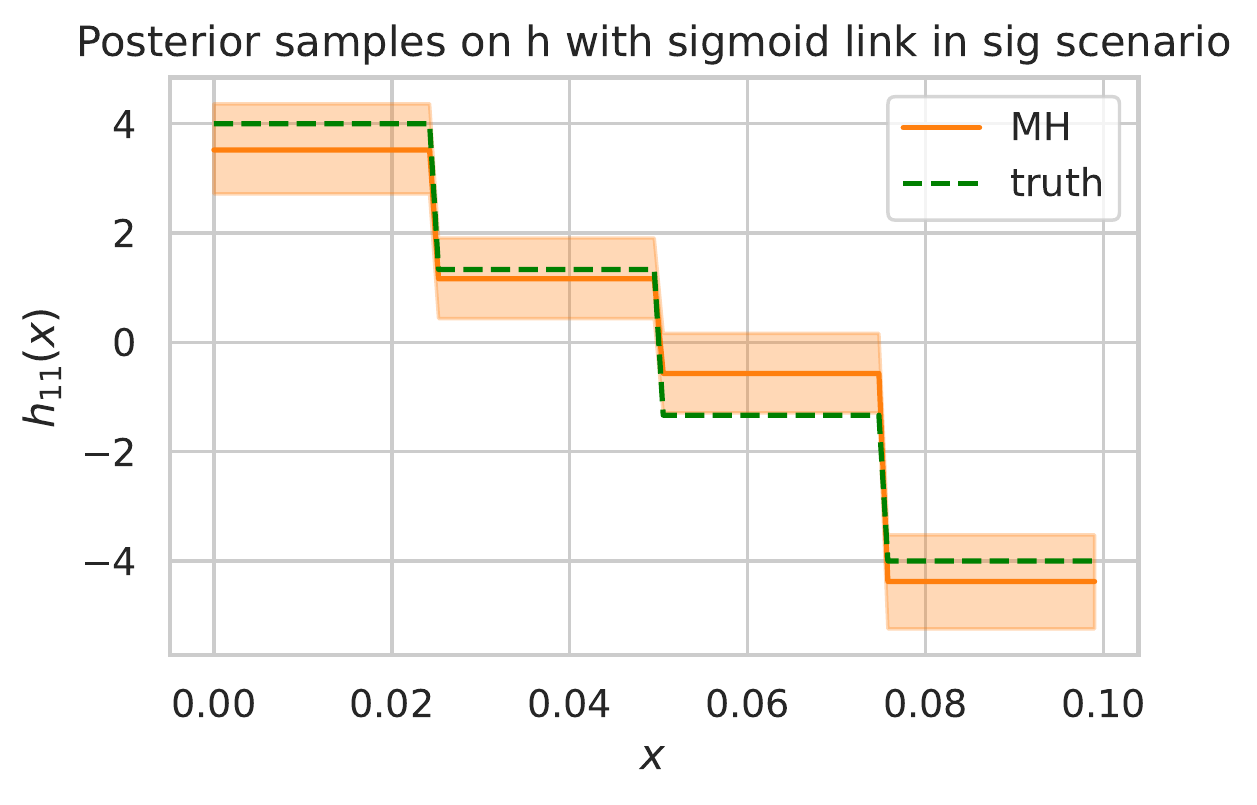}\hfil
\includegraphics[width=\tempwidth, trim=0.cm 0.cm 0cm  0.65cm,clip]{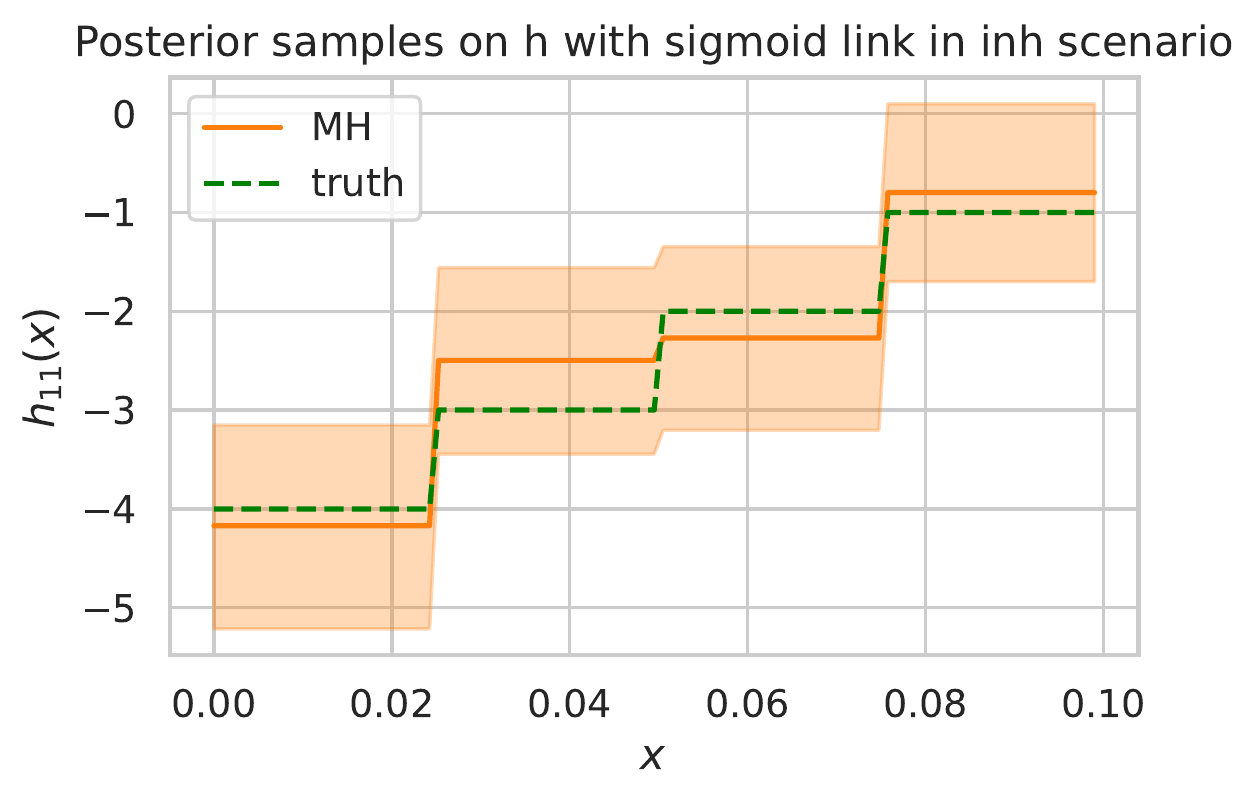}\\
\caption{Posterior distribution on $f = (\nu_1, h_{11})$ obtained with the MH sampler in the sigmoid model, in the three scenarios of Simulation 1 ($K=1$). The three columns correspond to the \emph{Excitation only} (left), \emph{Mixed effect} (center), and \emph{Inhibition only} (right) scenarios. The first row contains the marginal distribution on the background rate $\nu_1$, and the second row represents the posterior mean (solid orange line) and  95\% credible sets (orange areas) on the (self) interaction function $h_{11}$. The true parameter $f_0$ is plotted in dotted green line.}
\label{fig:sigmoid_mcmc_D4}
\end{figure}

\begin{figure}[hbt!]
\setlength{\tempwidth}{.3\linewidth}\centering
\settoheight{\tempheight}{\includegraphics[width=\tempwidth, trim=0.cm 0.cm 0cm  0.65cm,clip]{figs/sigmoid_d4_exc_h_mh_True_hmc_False.pdf}}%
\hspace{-5mm}
\fbox{\begin{minipage}{\dimexpr 15mm} \begin{center} \itshape \large  \textbf{Softplus} \end{center} \end{minipage}}
\hspace{-5mm}
\columnname{Excitation}\hfil
\columnname{Mixed}\hfil
\columnname{Inhibition}\\
\rowname{Background}
    \includegraphics[width=\tempwidth, trim=0.cm 0.cm 0cm  0.65cm,clip]{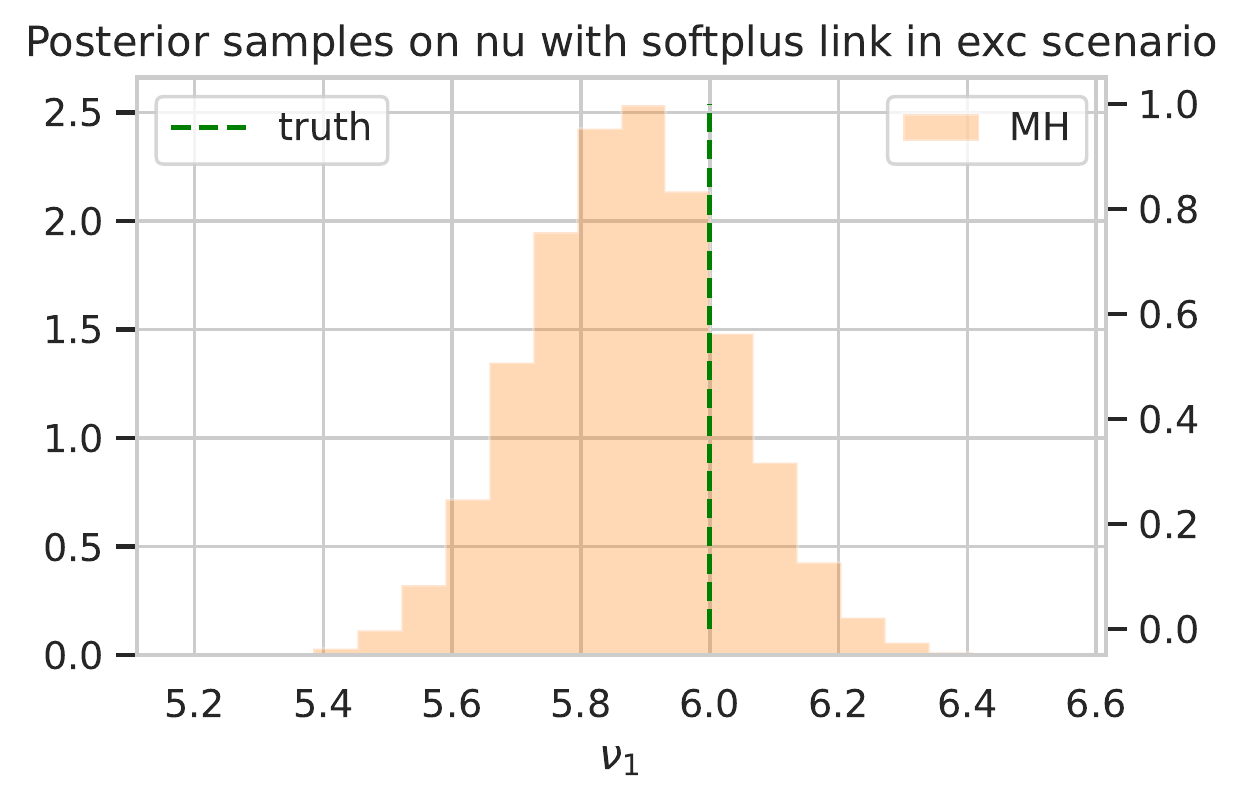}\hfil
    \includegraphics[width=\tempwidth, trim=0.cm 0.cm 0cm  0.65cm,clip]{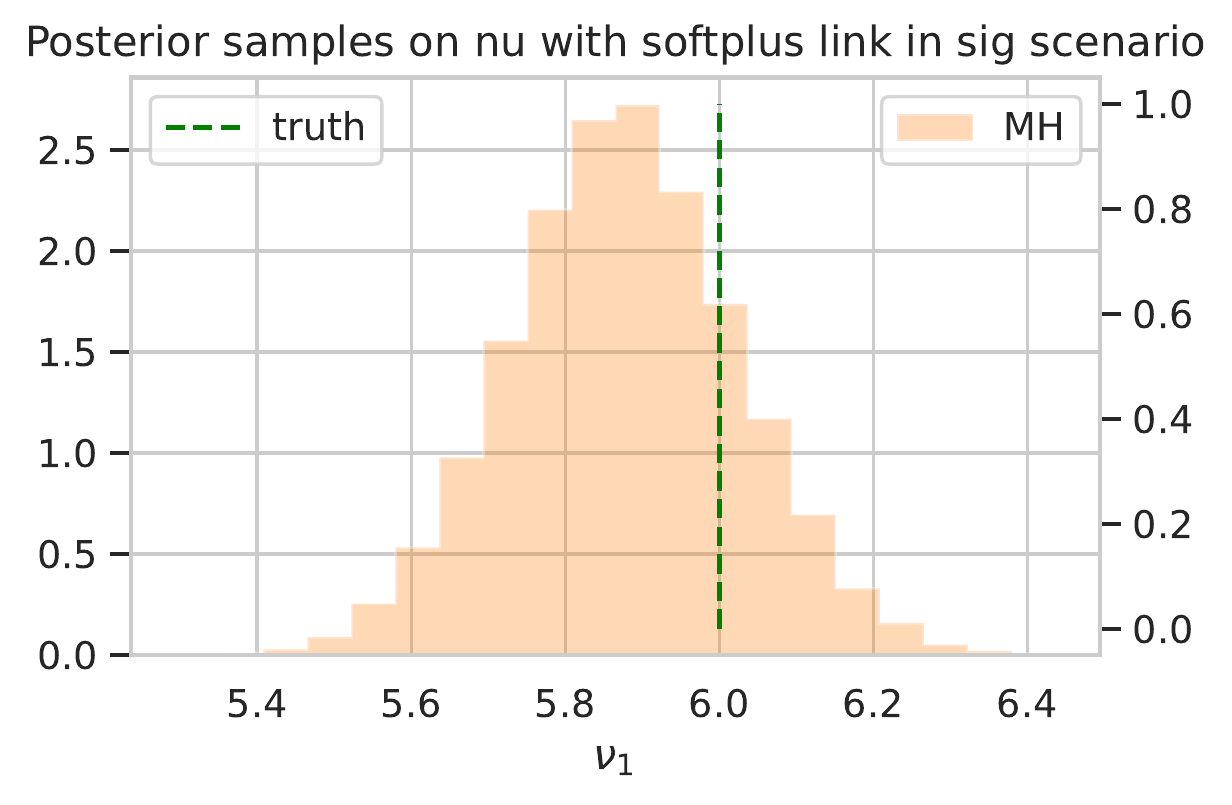}\hfil
    \includegraphics[width=\tempwidth, trim=0.cm 0.cm 0cm  0.65cm,clip]{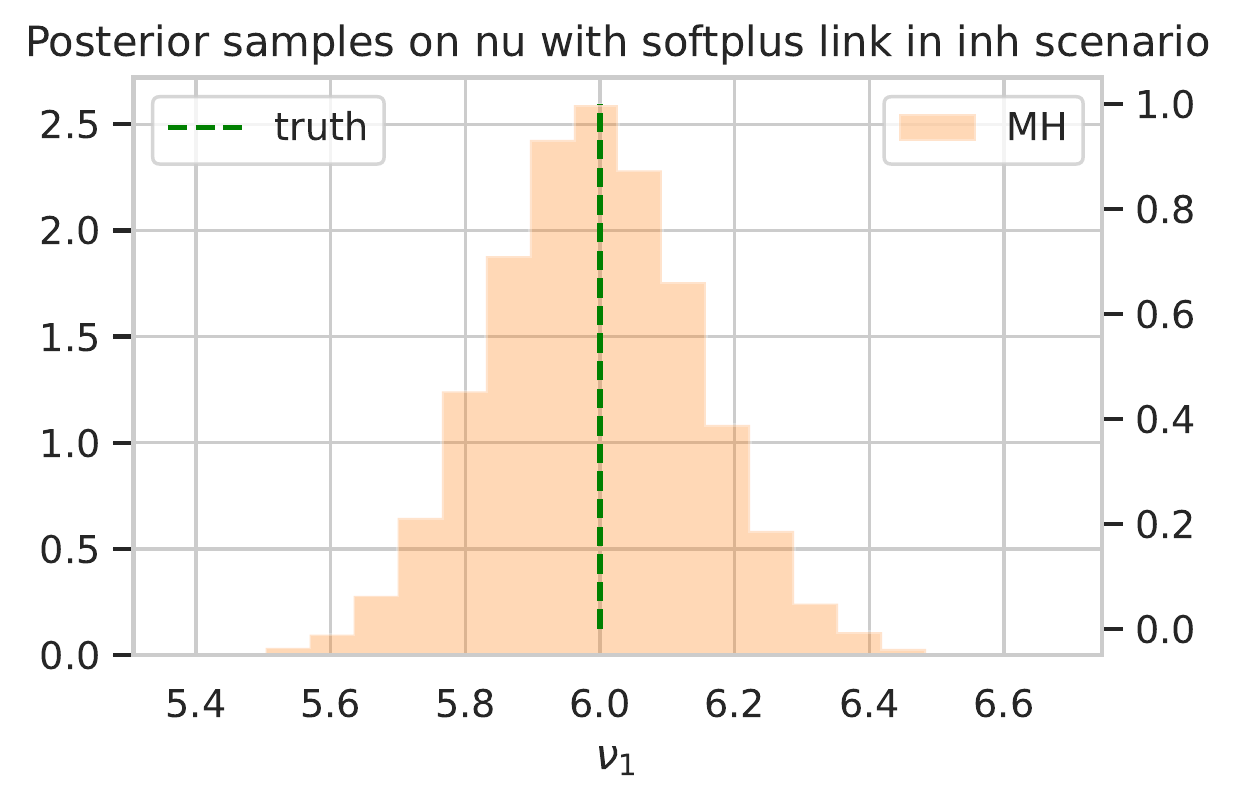}\\
    \includegraphics[width=\tempwidth, trim=0.cm 0.cm 0cm  0.65cm,clip]{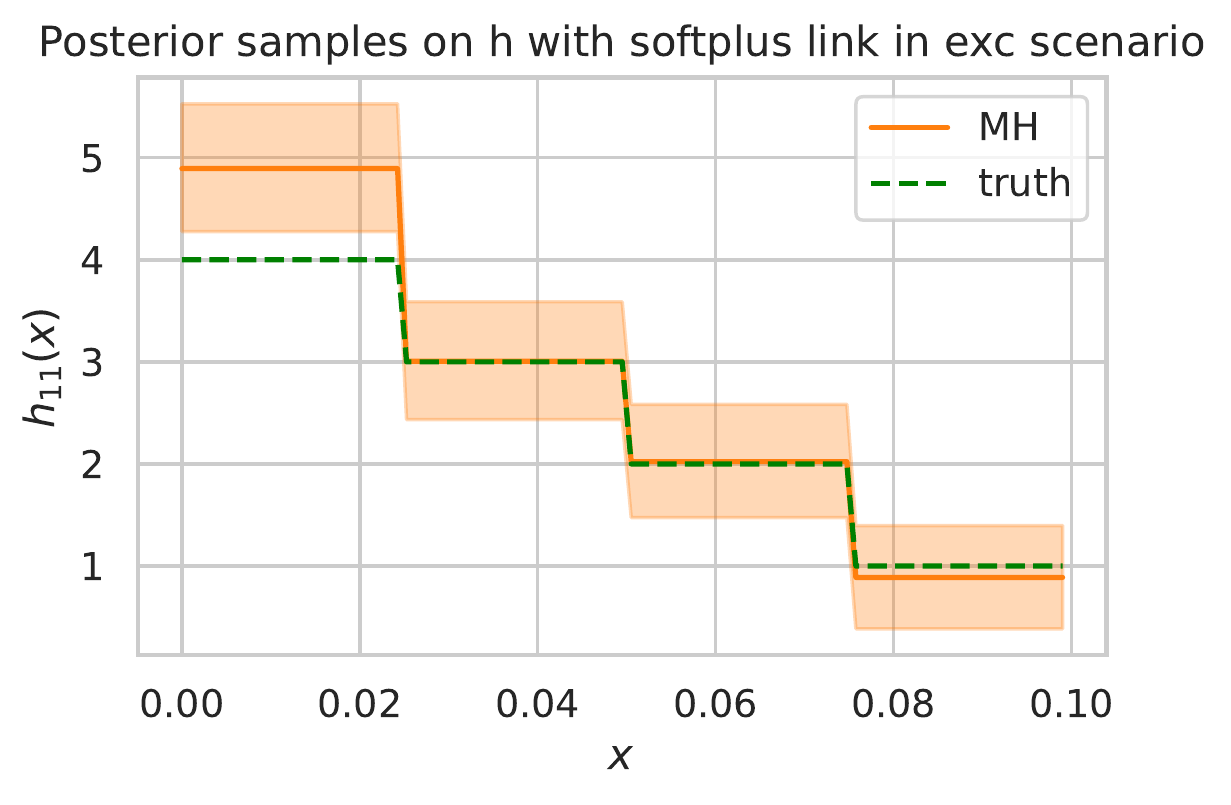}\hfil
    \includegraphics[width=\tempwidth, trim=0.cm 0.cm 0cm  0.65cm,clip]{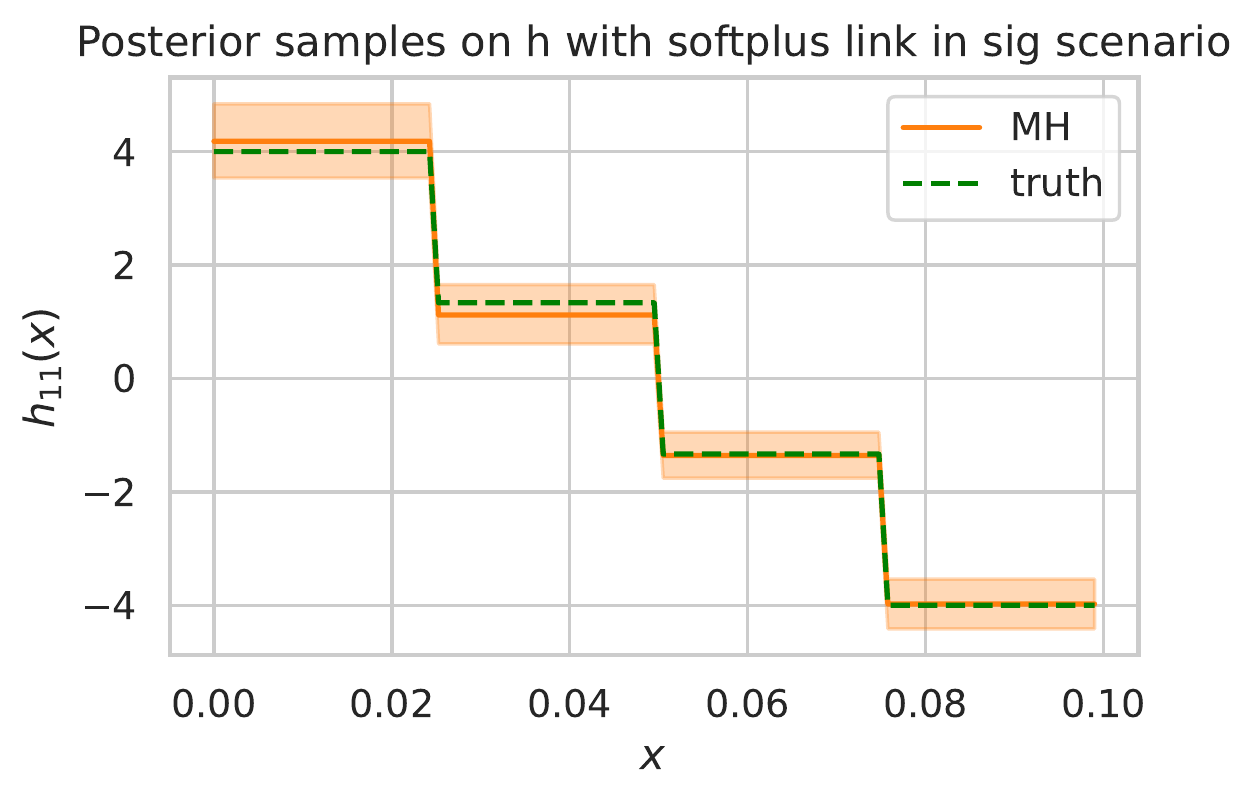}\hfil
    \includegraphics[width=\tempwidth, trim=0.cm 0.cm 0cm 0.65cm,clip]{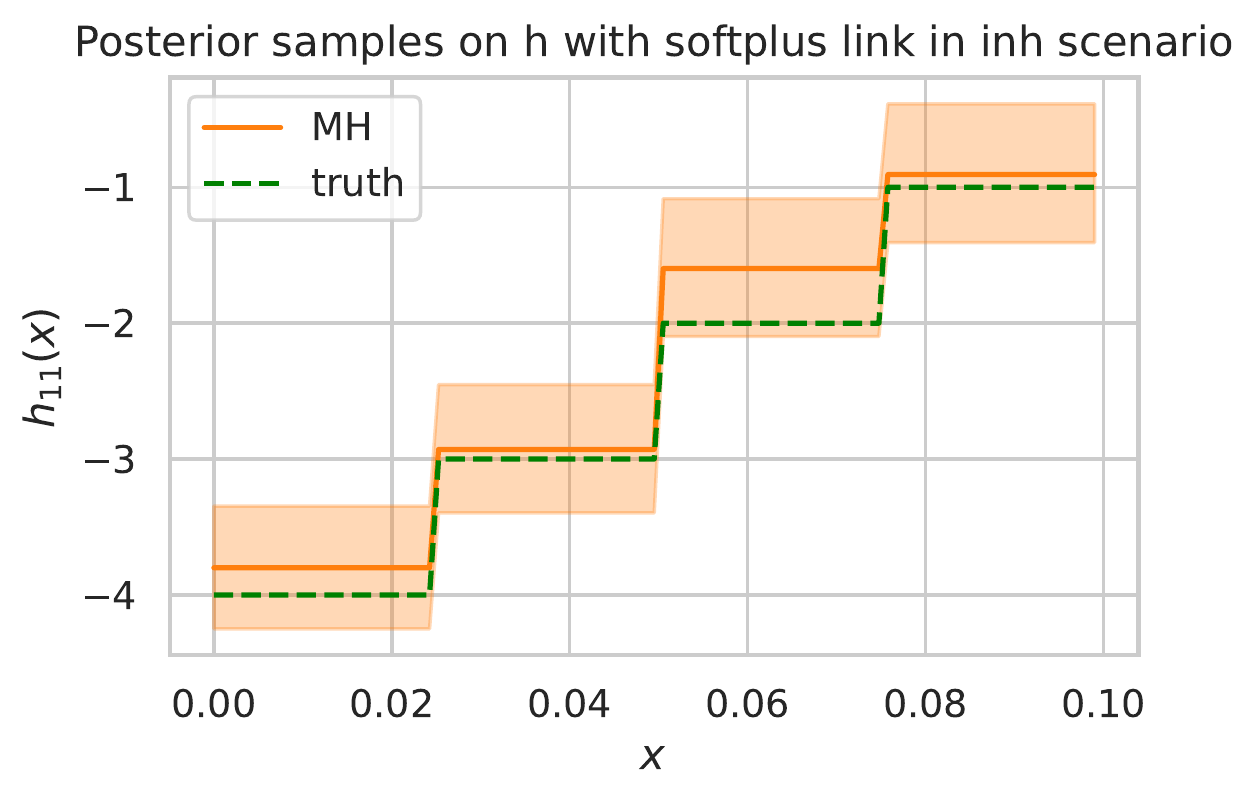}\\
\caption{Posterior distribution on $f = (\nu_1, h_{11})$ obtained with the MH sampler in the softplus model, in the three scenarios of Simulation 1 ($K=1$). The three columns correspond to the \emph{Excitation only} (left), \emph{Mixed effect} (center), and \emph{Inhibition only} (right) scenarios. The first row contains the marginal distribution on the background rate $\nu_1$, and the second row represents the posterior mean (solid orange line) and  95\% credible sets (orange  areas) on the (self) interaction function $h_{11}$. The true parameter $f_0$ is plotted in dotted green line.}
\label{fig:logit_mcmc_D4}
\end{figure}

This section contains our results for the MH sampler, in the univariate settings of Simulation 1 with sigmoid and softplus link functions (see Figures \ref{fig:sigmoid_mcmc_D4} and \ref{fig:logit_mcmc_D4}).

\subsection{Simulation 3}\label{app:simu_3}

\begin{figure}[hbt!]
    \centering
     \begin{subfigure}[b]{0.6\textwidth}
    \includegraphics[width=\textwidth, trim=0.cm 0.cm 0cm  0.cm,clip]{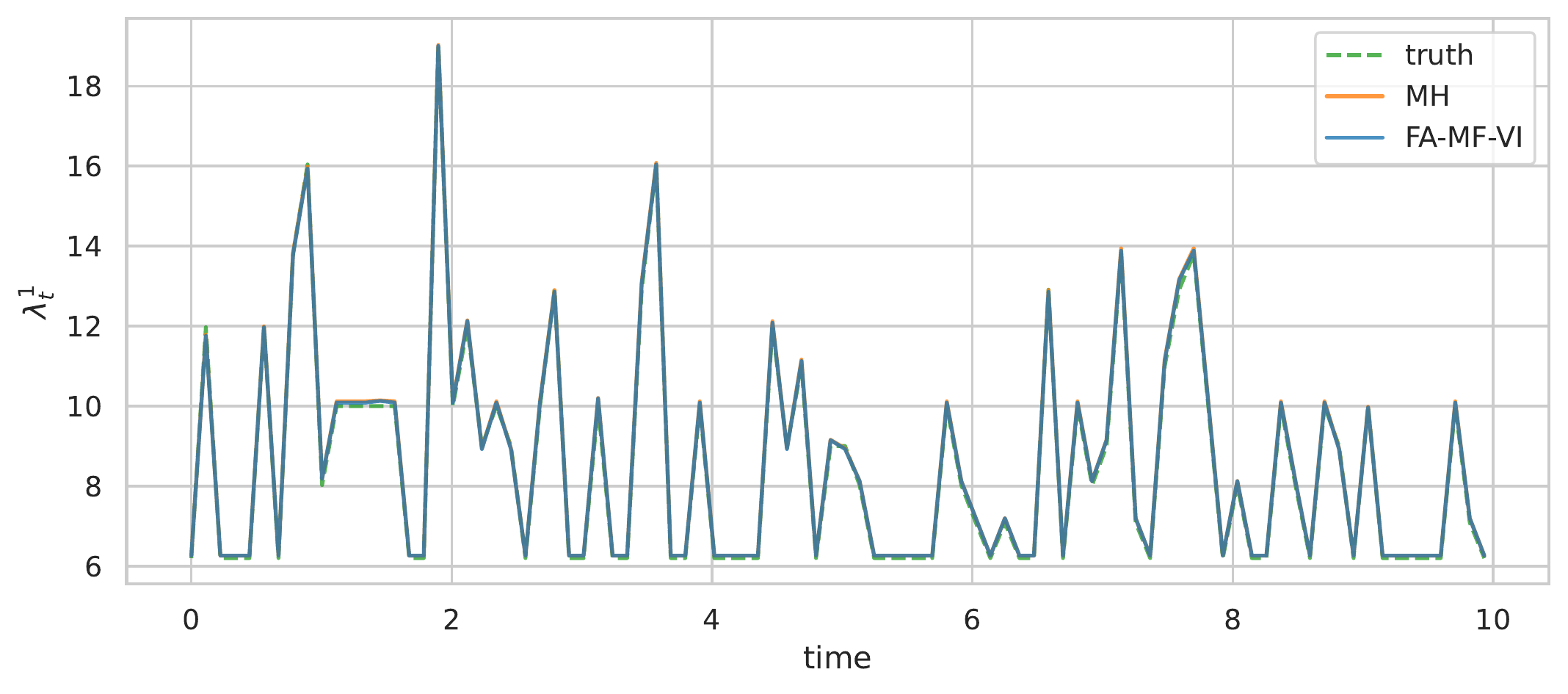}
    \caption{Excitation scenario }
    \end{subfigure}

         \begin{subfigure}[b]{0.6\textwidth}
    \includegraphics[width=\textwidth, trim=0.cm 0.cm 0cm  0.cm,clip]{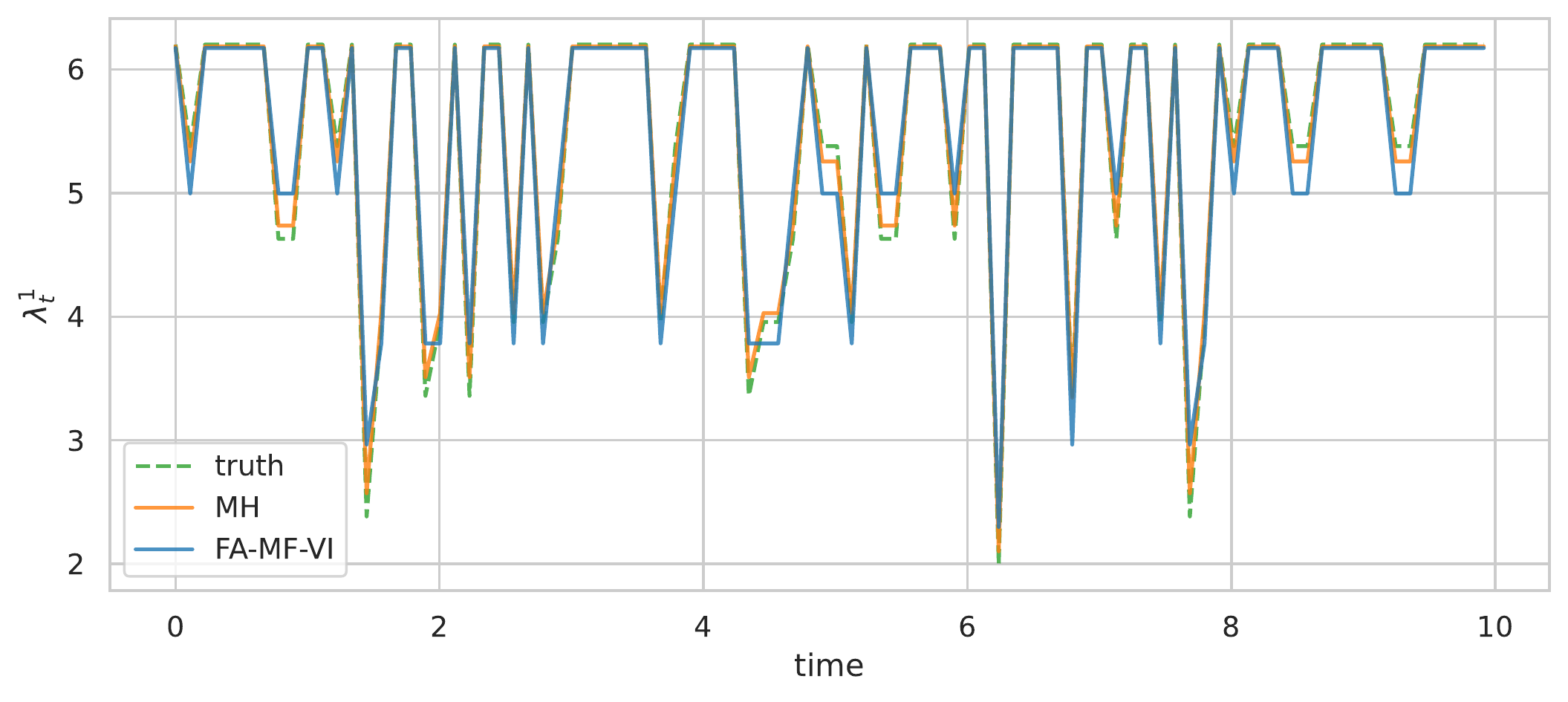}
    \caption{Inhibition scenario }
    \end{subfigure}
    \hfill
\caption{Intensity function on a subwindow of the observation window estimated via the variational posterior mean and via the posterior mean computed with the MH sampler, in the well-specified setting of Simulation 3 on  $[0,10]$,  using the fully-adaptive mean-field variational (FA-MF-VI) algorithm (Algorithm \ref{alg:adapt_cavi}). The true intensity $\lambda_t^1(f_0)$ is plotted in dotted green line. }
\label{fig:intensity_D1_adaptive}
\end{figure}

\begin{figure}[hbt!]
\setlength{\tempwidth}{.35\linewidth}\centering
\settoheight{\tempheight}{\includegraphics[width=\tempwidth, trim=0.cm 0.cm 0cm  1.cm,clip]{figs/adaptive_vi_1D_histogram_exc_nu_smod.pdf}}%
\hspace{-5mm}
\fbox{\begin{minipage}{\dimexpr 25mm} \begin{center} \itshape \large  \textbf{$K = 1$\\ Mis-specified} \end{center} \end{minipage}}
\columnname{Excitation}\hfil
\columnname{Inhibition}\\
\hspace{10mm}
\rowname{Background}
    \includegraphics[width=\tempwidth, trim=0.cm 0.cm 1cm  1.cm,clip]{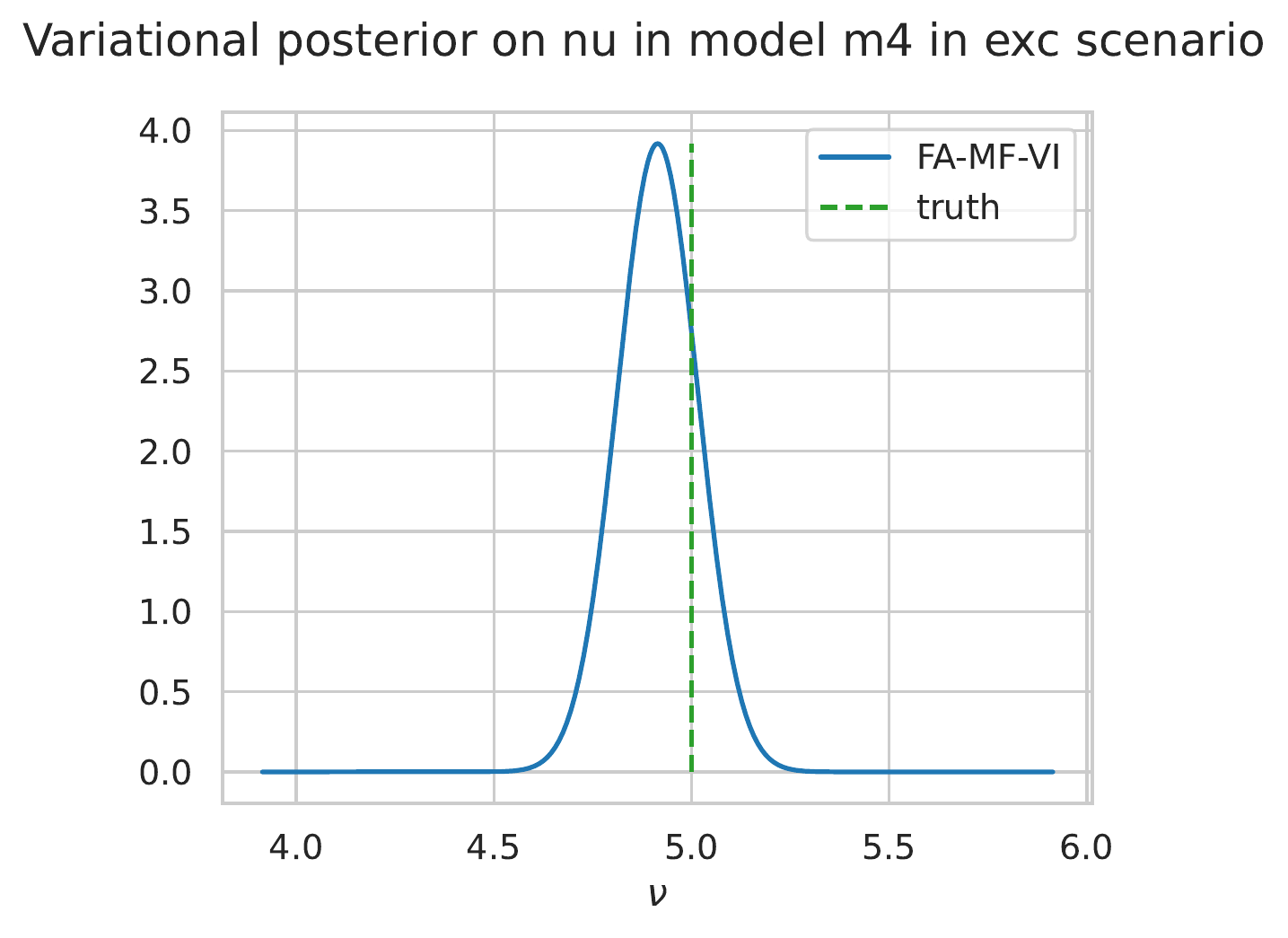}\hfil
    \includegraphics[width=\tempwidth, trim=0.cm 0.cm 1cm  1.cm,clip]{figs/adaptive_vi_1D_inh_continuous_nu.pdf}\\
    \hspace{10mm}
\rowname{Interaction}
    \includegraphics[width=0.95\tempwidth, trim=0.cm 0.cm 1.5cm  1.cm,clip]{figs/adaptive_vi_1D_exc_continuous_h.pdf}\hfil
    \includegraphics[width=0.95\tempwidth, trim=0.cm 0.cm 1.5cm  1.cm,clip]{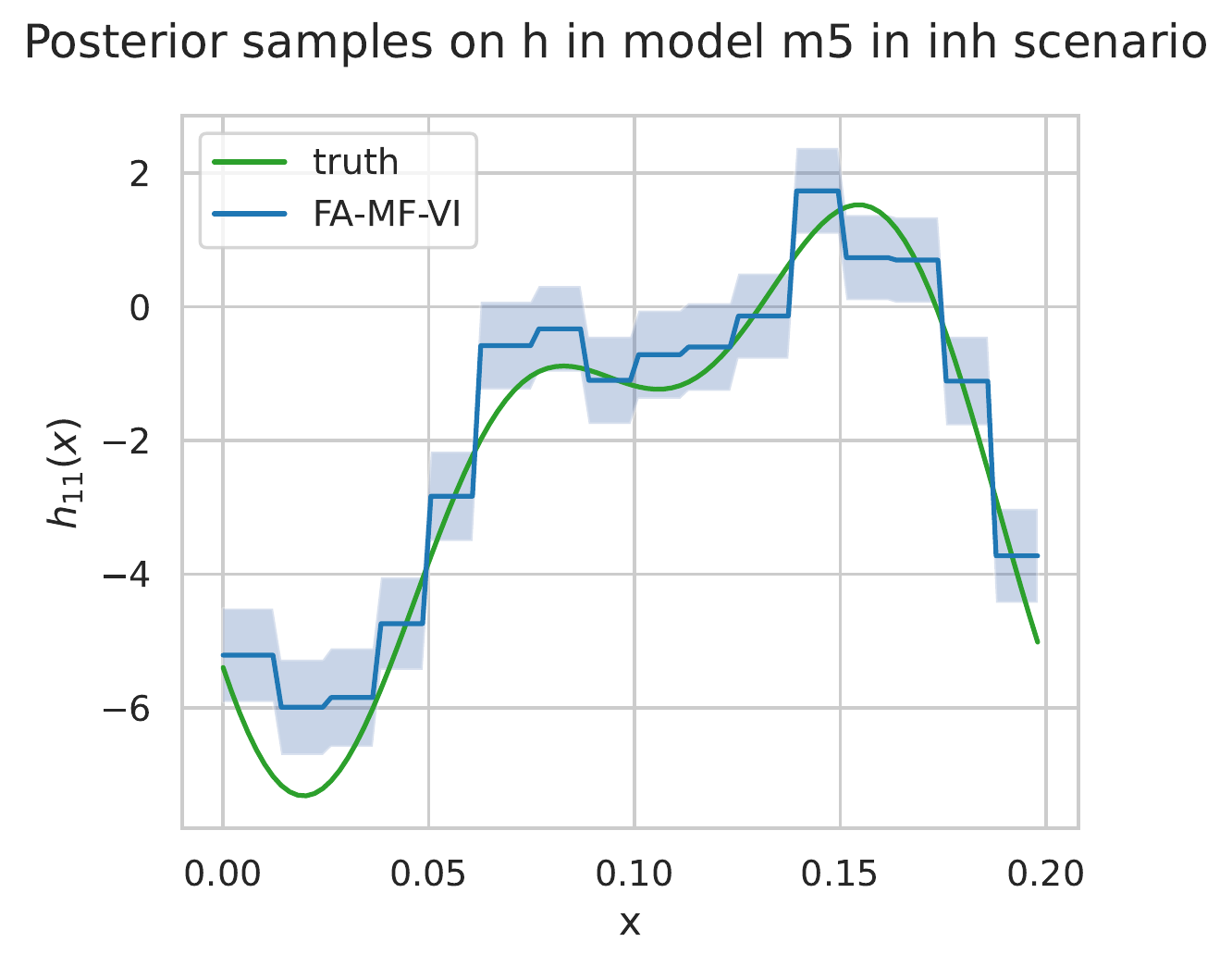}\\
\caption{Model-selection variational posterior distributions on $f = (\nu_1, h_{11})$ in the univariate sigmoid model and mis-specified setting of Simulation 3, evaluated by the fully-adaptive mean-field variational (FA-MF-VI) algorithm (Algorithm \ref{alg:adapt_cavi}). The two columns correspond to a (mostly) \emph{Excitation} (left) and a (mostly) \emph{Inhibition} (right) settings. The first row contains the marginal distribution on the background rate $\nu_1$, and the second row represents the variational posterior mean (solid line) and  95\% credible sets (colored areas) on the (self) interaction function $h_{11}$.  The true parameter $f_0$ is plotted in dotted green line.}
\label{fig:adaptive_VI_1D_fourier}
\end{figure}

\begin{figure}[hbt!]
\setlength{\tempwidth}{.7\linewidth}\centering
\settoheight{\tempheight}{\includegraphics[width=0.5\tempwidth, trim=0.cm 0.cm 0cm  1.cm,clip]{figs/adaptive_vi_1D_histogram_exc_nu_smod.pdf}}%
\hspace{-5mm}
\fbox{\begin{minipage}{\dimexpr 25mm} \begin{center} \itshape \large  \textbf{$K = 2$\\ Inhibition} \end{center} \end{minipage}}
\columnname{}\\
\rowname{\centering Background}
    \includegraphics[width=\tempwidth, trim=0.cm 0.cm 0cm  0.8cm,clip]{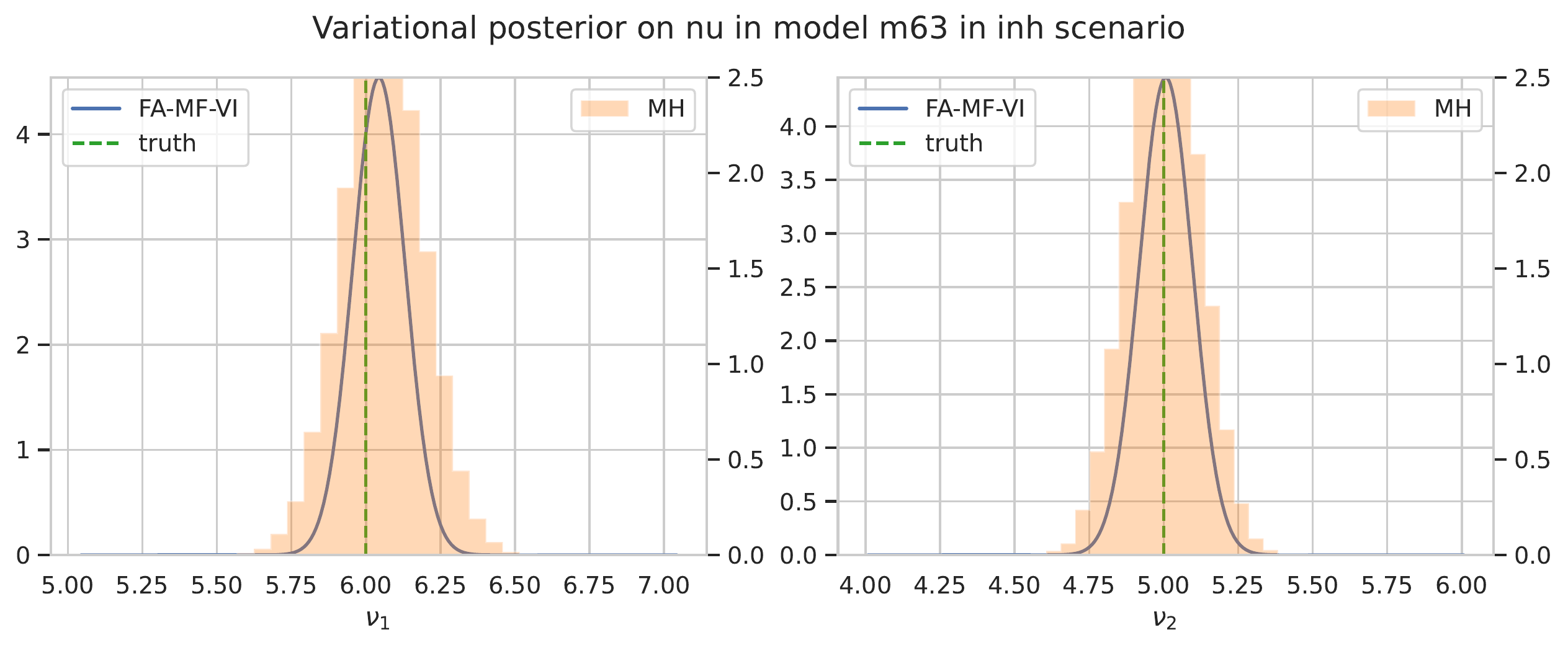}\\
\hspace{-10mm}
\rowname{\hspace{60mm} Interaction functions}
    \includegraphics[width=0.93\tempwidth, trim=0.cm 0.cm 0cm  1.5cm,clip]{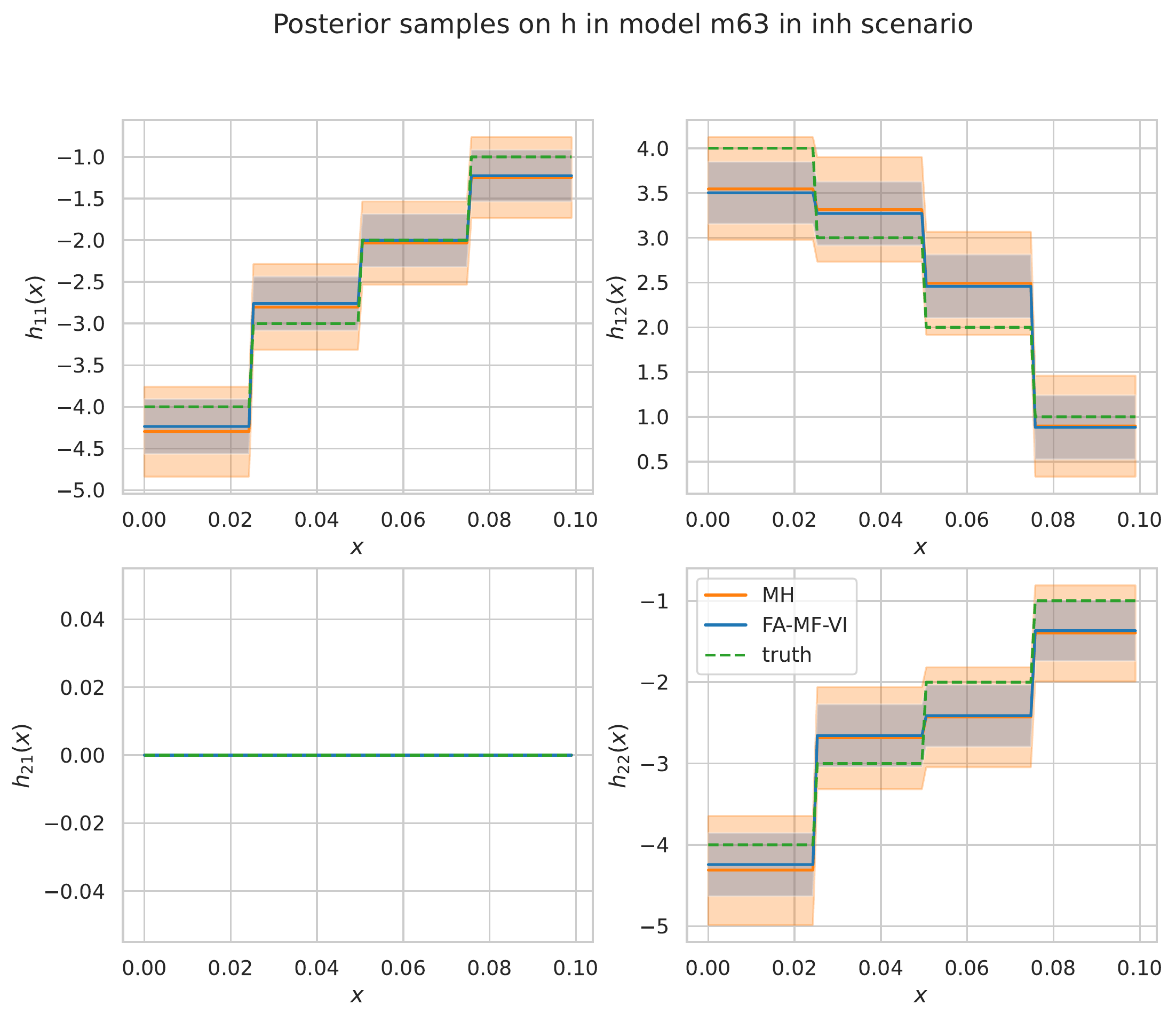}
\caption{Posterior and model-selection variational posterior distributions on $f = (\nu, h)$ in the bivariate sigmoid model, well-specified setting, and Inhibition setting of Simulation 3, evaluated by the non-adaptive MH sampler and the fully-adaptive mean-field variational (FA-MF-VI) algorithm (Algorithm \ref{alg:adapt_cavi}).  The first row contains the marginal distribution on the background rates $(\nu_1, \nu_2)$, and the second and third rows represent the (variational) posterior mean (solid line) and  95\% credible sets (colored areas) on the four interaction function $h_{11}, h_{12}, h_{21}, h_{22}$.  The true parameter $f_0$ is plotted in dotted green line.}
\label{fig:adaptive_VI_2D_selec_inh}
\end{figure}

This section contains our results regarding the estimated intensity function in the univariate and well-specified settings in Simulation 3 (see Figure \ref{fig:intensity_D1_adaptive}), the estimated parameter in the mis-specified settings (see Figure \ref{fig:adaptive_VI_1D_fourier}), and the estimated interaction functions in the bivariate settings (see Figures \ref{fig:adaptive_VI_2D_selec_exc} and \ref{fig:adaptive_VI_2D_selec_inh}).

\begin{figure}[hbt!]
\setlength{\tempwidth}{.7\linewidth}\centering
\settoheight{\tempheight}{\includegraphics[width=0.5\tempwidth, trim=0.cm 0.cm 0cm  1.cm,clip]{figs/adaptive_vi_1D_histogram_exc_nu_smod.pdf}}%
\hspace{-5mm}
\fbox{\begin{minipage}{\dimexpr 25mm} \begin{center} \itshape \large  \textbf{$K = 2$\\ Excitation} \end{center} \end{minipage}}
\columnname{}\\
\rowname{\centering Background}
    \includegraphics[width=\tempwidth, trim=0.cm 0.cm 0cm  0.8cm,clip]{figs/adaptive_vi_2D_histogram_exc_nu_smod.pdf}\\
\hspace{-10mm}
\rowname{\hspace{60mm} Interaction functions}
    \includegraphics[width=0.93\tempwidth, trim=0.cm 0.cm 0cm  1.5cm,clip]{figs/adaptive_vi_2D_4_histogram_exc_h_smod.pdf}
\caption{Posterior and model-selection variational posterior distributions on $f = (\nu, h)$ in the bivariate sigmoid model, well-specified setting, and Excitation setting of Simulation 3, evaluated by the non-adaptive MH sampler and the fully-adaptive mean-field variational (FA-MF-VI) algorithm (Algorithm \ref{alg:adapt_cavi}).  The first row contains the marginal distribution on the background rates $(\nu_1, \nu_2)$, and the second and third rows represent the (variational) posterior mean (solid line) and  95\% credible sets (colored areas) on the four interaction function $h_{11}, h_{12}, h_{21}, h_{22}$.  The true parameter $f_0$ is plotted in dotted green line.}
\label{fig:adaptive_VI_2D_selec_exc}
\end{figure}

\FloatBarrier

\subsection{Simulation 4}\label{app:simu_4}

\begin{figure}[hbt!]
    \centering
         \begin{subfigure}[b]{0.3\textwidth}
    \includegraphics[width=\textwidth, trim=0.cm 0cm 0cm  1.cm,clip]{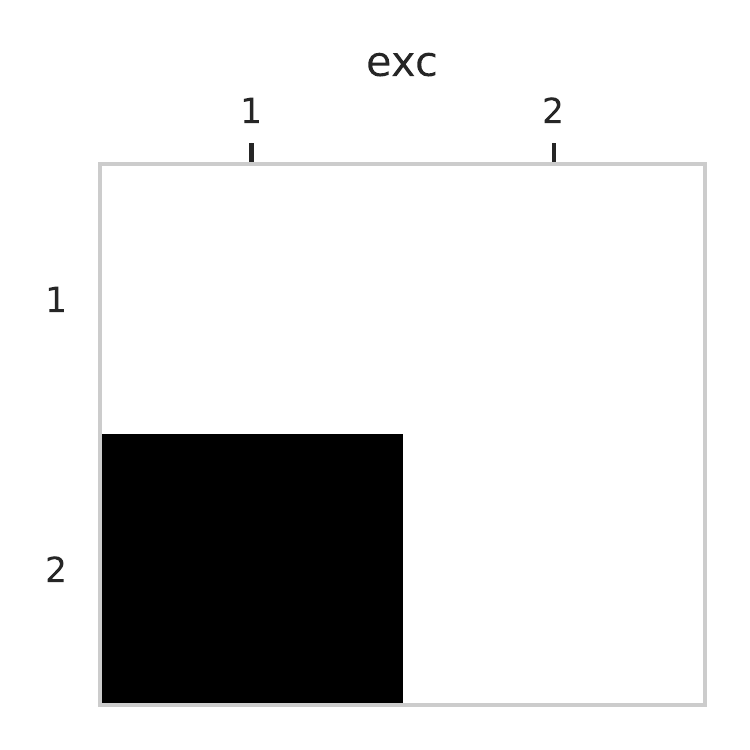}
    \caption{$K=2$}
    \end{subfigure}%
             \begin{subfigure}[b]{0.3\textwidth}
    \includegraphics[width=\textwidth, trim=0.cm 0cm 0cm  1.cm,clip]{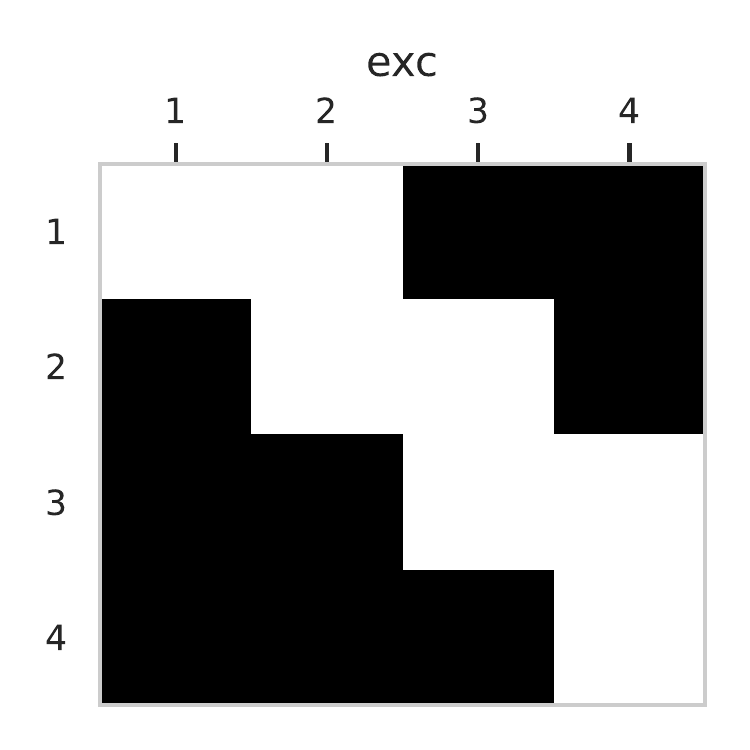}
    \caption{$K=4$}
    \end{subfigure}%
             \begin{subfigure}[b]{0.3\textwidth}
    \includegraphics[width=\textwidth, trim=0.cm 0cm 0cm  1.cm,clip]{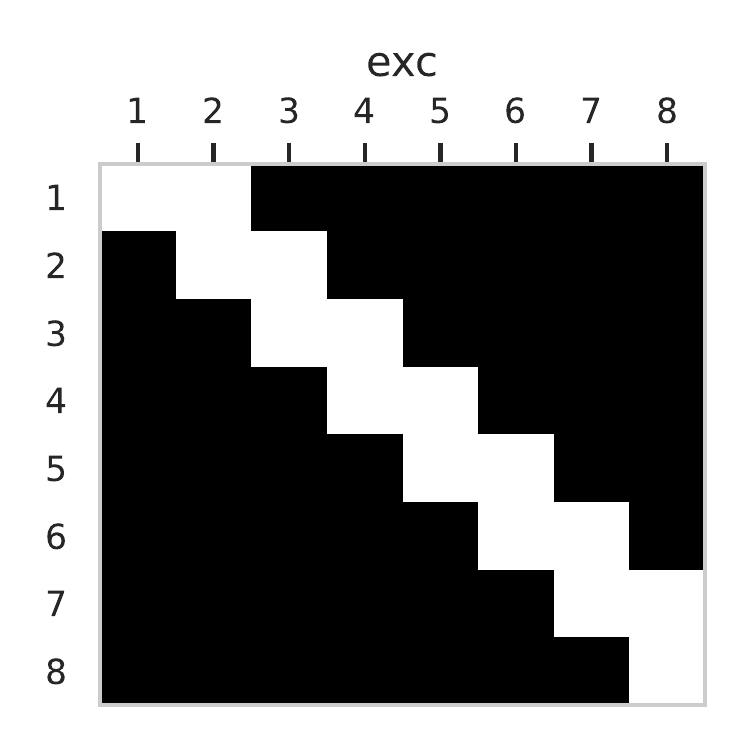}
    \caption{$K=8$}
    \end{subfigure}
        \begin{subfigure}[b]{0.33\textwidth}
    \includegraphics[width=\textwidth, trim=0.cm 0cm 0cm  1.cm,clip]{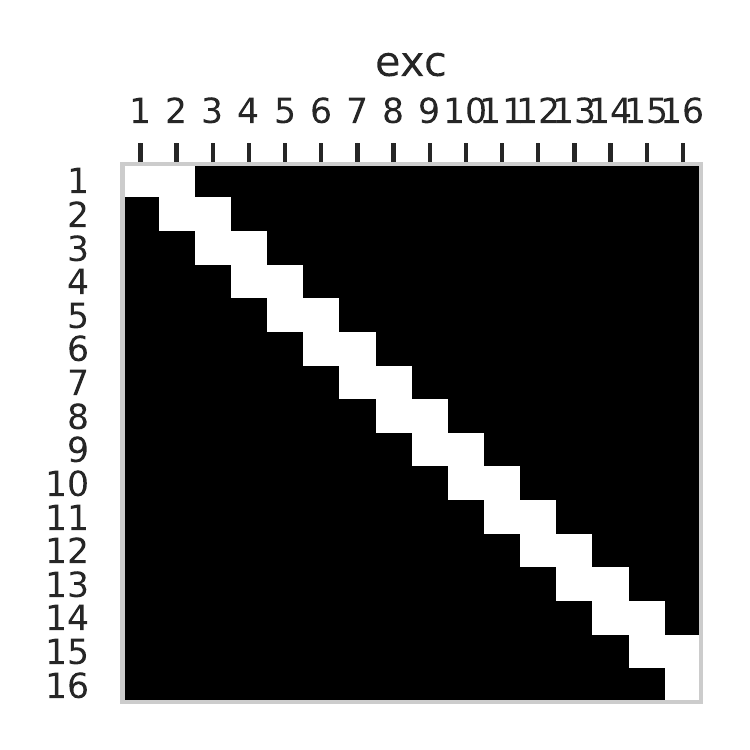}
    \caption{$K=10$}
    \end{subfigure}
    \begin{subfigure}[b]{0.3\textwidth}
    \includegraphics[width=\textwidth, trim=0.cm 0cm 0cm  1.cm,clip]{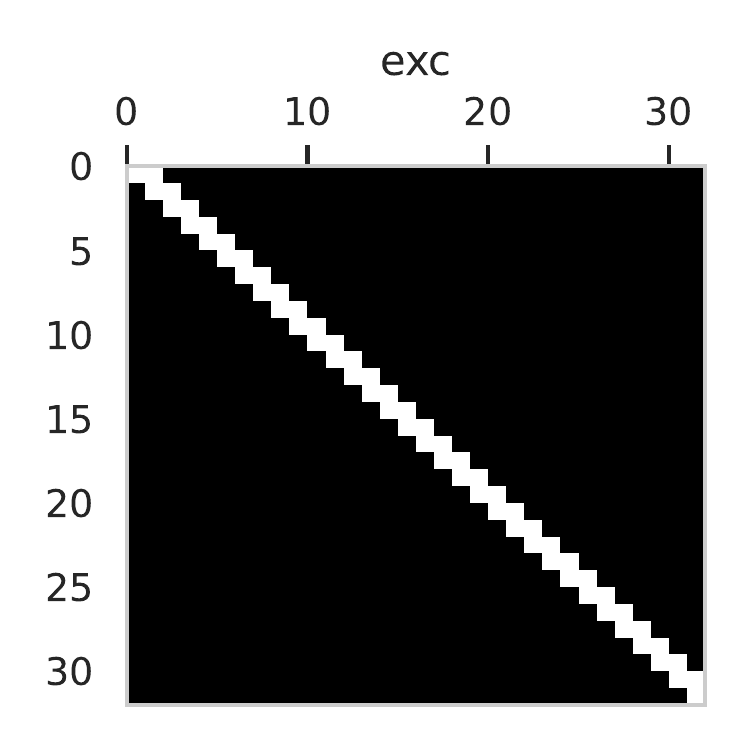}
    \caption{$K=16$}
    \end{subfigure}
        \begin{subfigure}[b]{0.3\textwidth}
    \includegraphics[width=\textwidth, trim=0.cm 0cm 0cm  1.cm,clip]{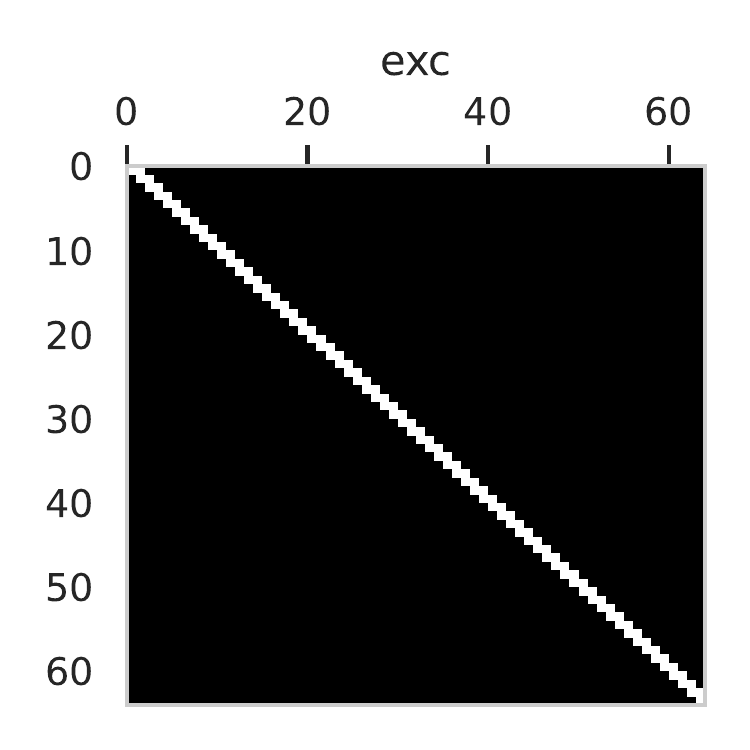}
    \caption{$K=32$}
    \end{subfigure}
\caption{Estimated graph parameter $\hat \delta$ (black=0, white=1) for $K=2,4,8,16,32,64$ in the Excitation scenario of Simulation 4.}
\label{fig:graphs_exc}
\end{figure}

\begin{figure}[hbt!]
    \centering
         \begin{subfigure}[b]{0.3\textwidth}
    \includegraphics[width=\textwidth, trim=0.cm 0cm 0cm  1.cm,clip]{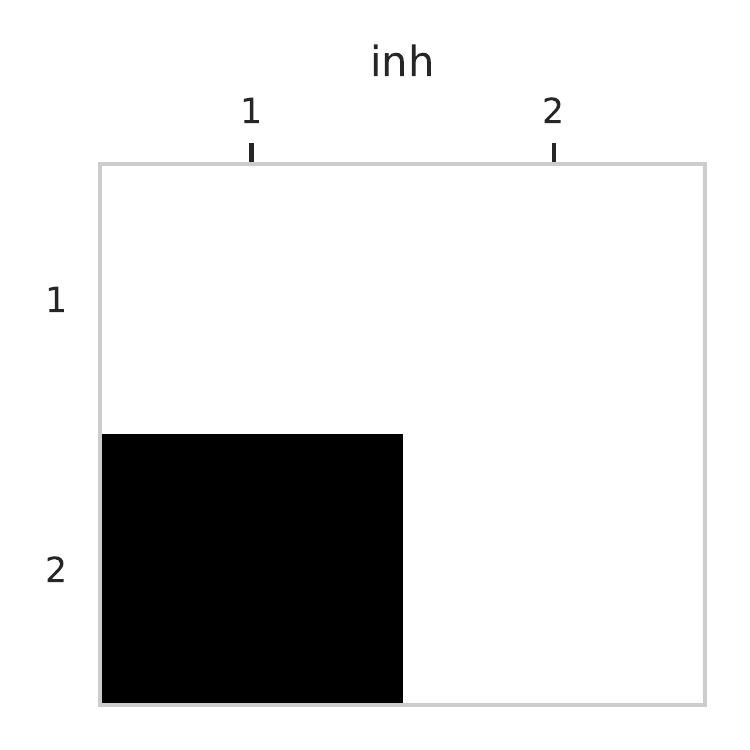}
    \caption{$K=2$}
    \end{subfigure}%
             \begin{subfigure}[b]{0.3\textwidth}
    \includegraphics[width=\textwidth, trim=0.cm 0cm 0cm  1.cm,clip]{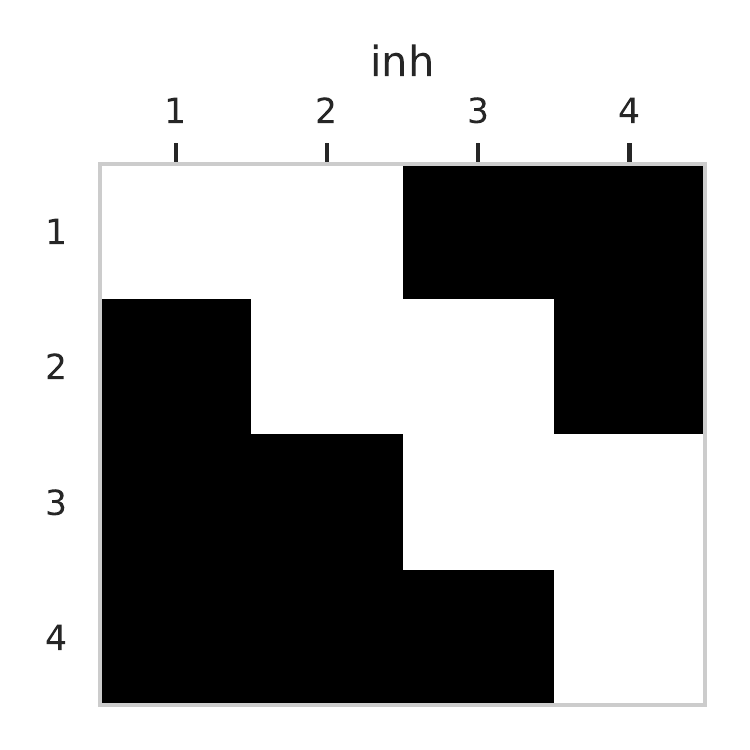}
    \caption{$K=4$}
    \end{subfigure}%
             \begin{subfigure}[b]{0.3\textwidth}
    \includegraphics[width=\textwidth, trim=0.cm 0cm 0cm  1.cm,clip]{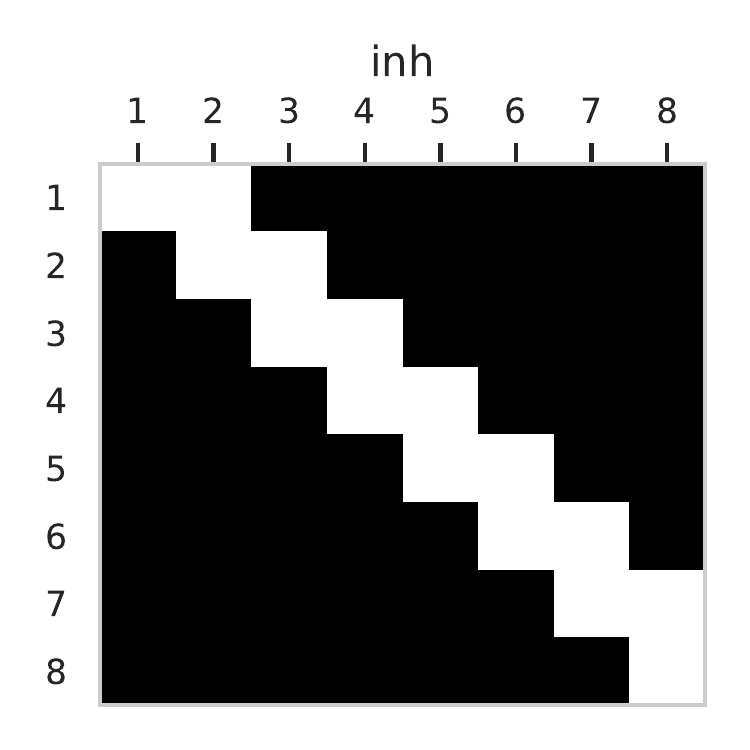}
    \caption{$K=8$}
    \end{subfigure}
        \begin{subfigure}[b]{0.33\textwidth}
    \includegraphics[width=\textwidth, trim=0.cm 0cm 0cm  1.cm,clip]{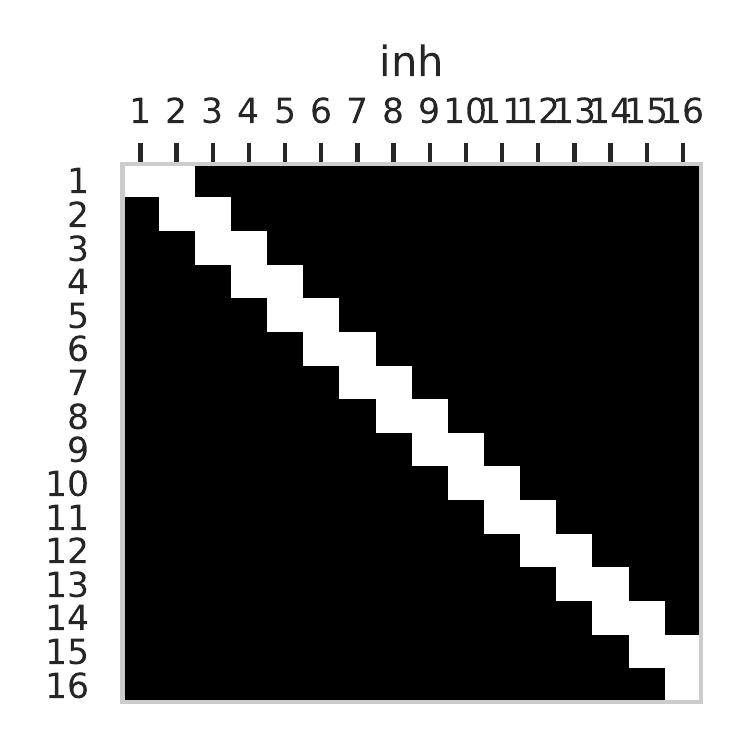}
    \caption{$K=16$}
    \end{subfigure}
    \begin{subfigure}[b]{0.3\textwidth}
    \includegraphics[width=\textwidth, trim=0.cm 0cm 0cm  1.cm,clip]{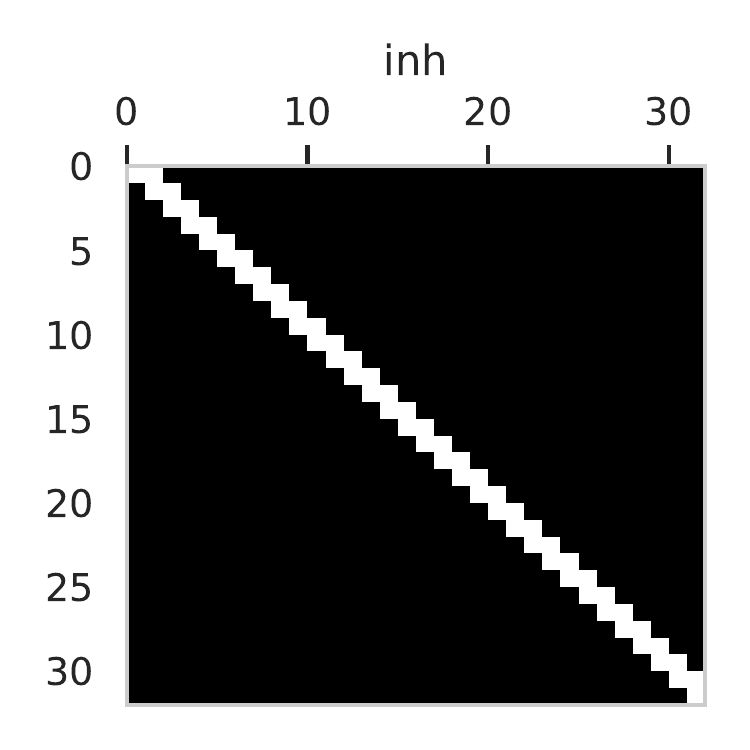}
    \caption{$K=32$}
    \end{subfigure}
        \begin{subfigure}[b]{0.3\textwidth}
    \includegraphics[width=\textwidth, trim=0.cm 0cm 0cm  1.cm,clip]{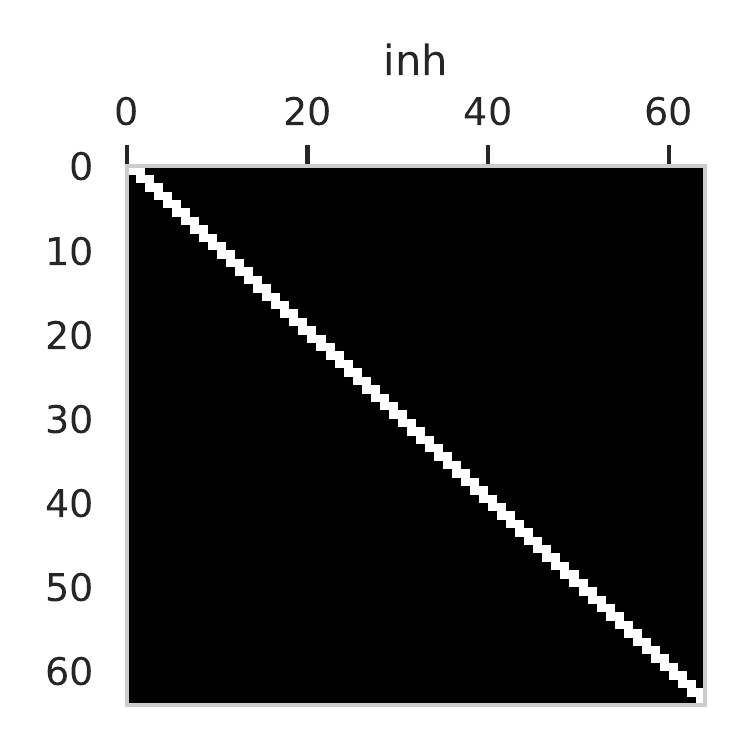}
    \caption{$K=64$}
    \end{subfigure}
\caption{Estimated graph parameter $\hat \delta$ (black=0, white=1) for $K=2,4,8,16,32,64$ in the Inhibition scenario of Simulation 4.}
\label{fig:graphs_inh}
\end{figure}

\begin{figure}[hbt!]
\setlength{\tempwidth}{.3\linewidth}\centering
\settoheight{\tempheight}{\includegraphics[width=\tempwidth, trim=0.cm 0.cm 0.cm  0.cm,clip]{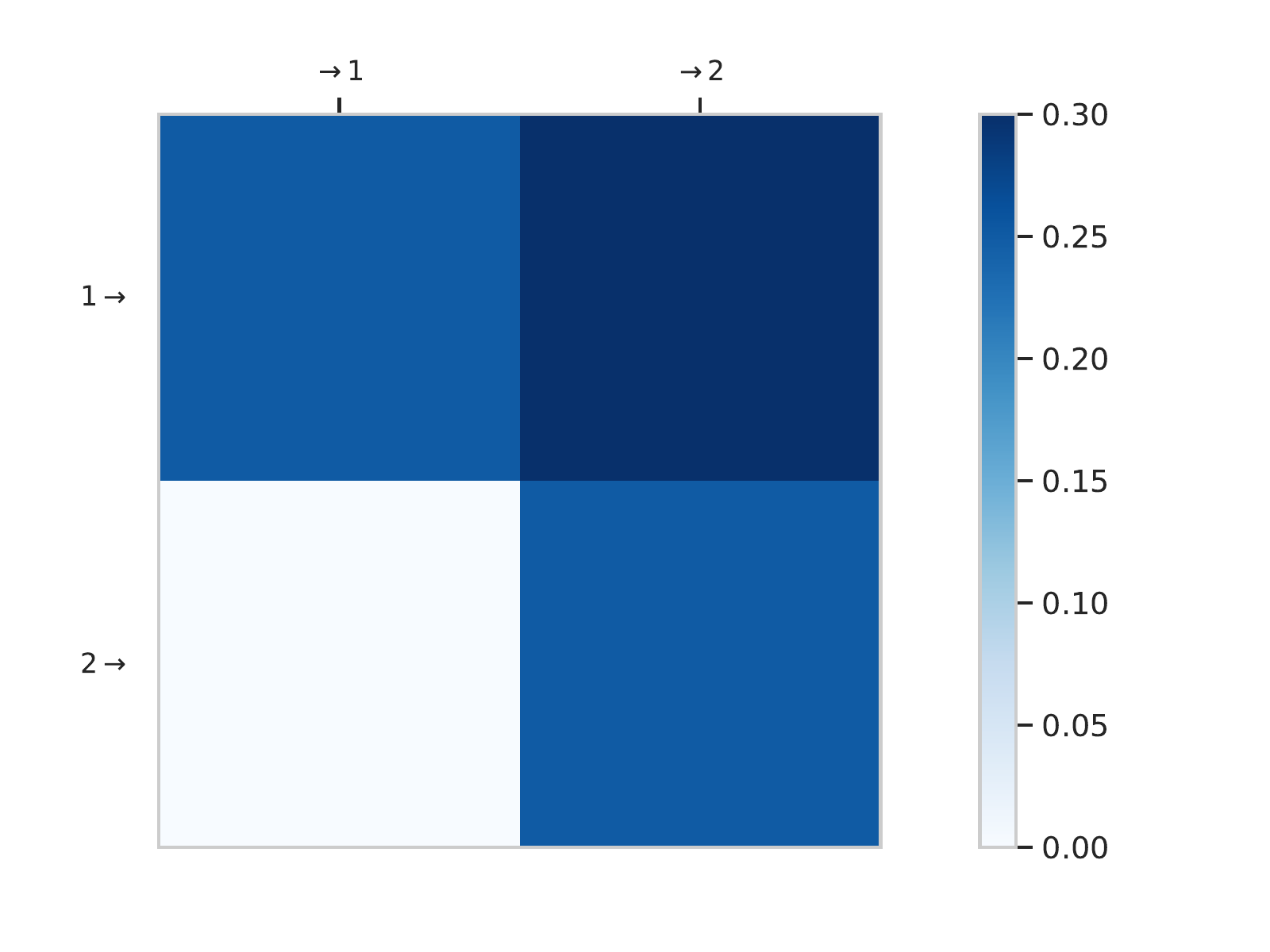}}%
\hspace{-5mm}
\fbox{\begin{minipage}{\dimexpr 40mm} \begin{center} \itshape \large  \textbf{Function norms\\Inhibition} \end{center} \end{minipage}}
\hspace{-5mm}
\columnname{Ground-truth}\hfil
\columnname{Error}\\
\begin{minipage}{\dimexpr 40mm} \vspace{-35mm} \begin{center} \itshape \large  \textbf{$K=2$} \end{center} \end{minipage}
    \includegraphics[width=\tempwidth, trim=0.cm 0.cm 0.cm  0.cm,clip]{figs/simu_K2_true_norms_inh.pdf}\hfil
    \includegraphics[width=\tempwidth, trim=0.cm 0.cm 0cm  0.cm,clip]{figs/simu_K2_err_norms_exc.pdf}\\
\begin{minipage}{\dimexpr 40mm} \vspace{-35mm} \begin{center} \itshape \large  \textbf{$K=4$} \end{center} \end{minipage}
    \includegraphics[width=\tempwidth, trim=0.cm 0.cm 0cm  0.cm,clip]{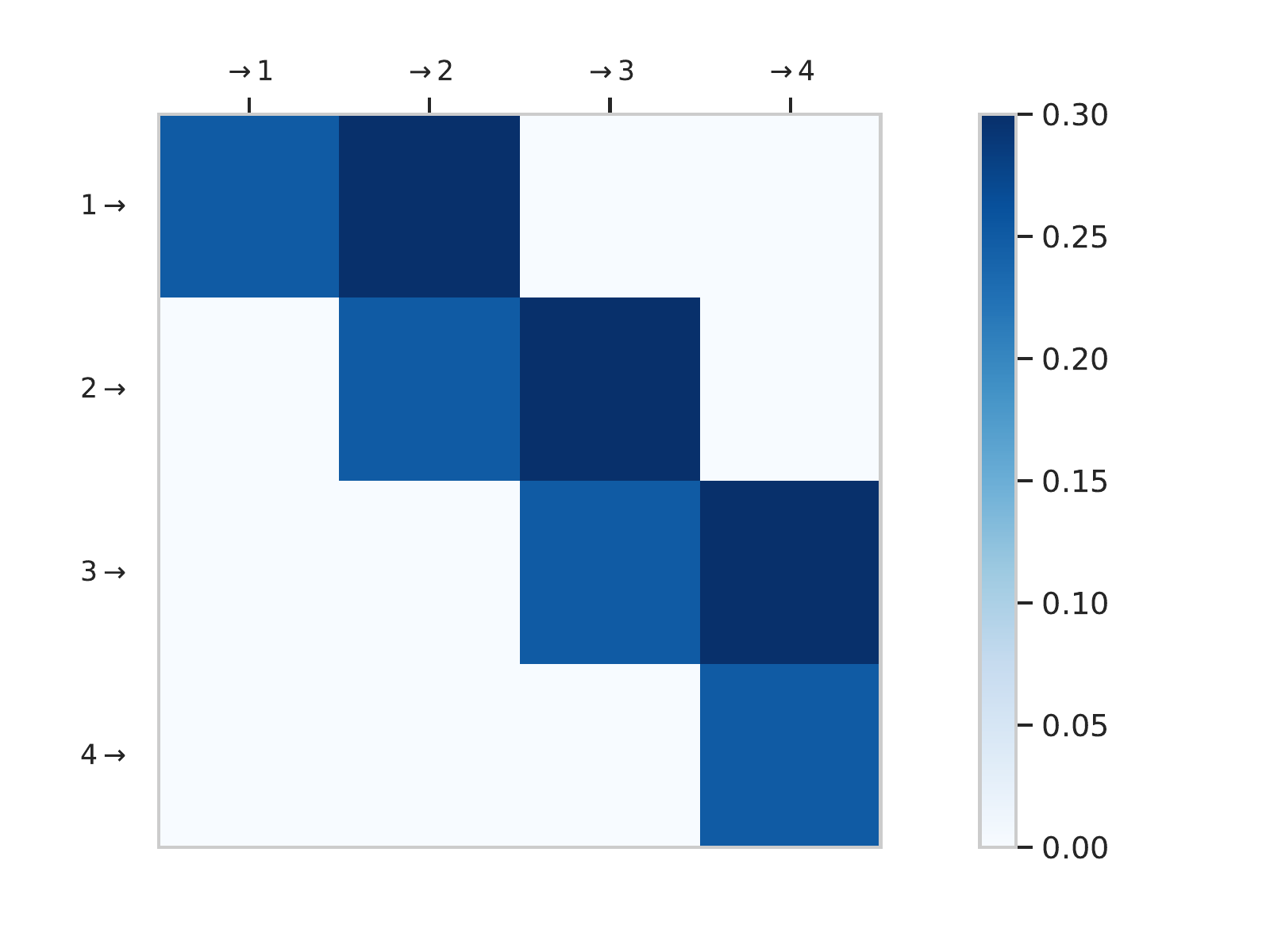}\hfil
    \includegraphics[width=\tempwidth, trim=0.cm 0.cm 0cm  0.cm,clip]{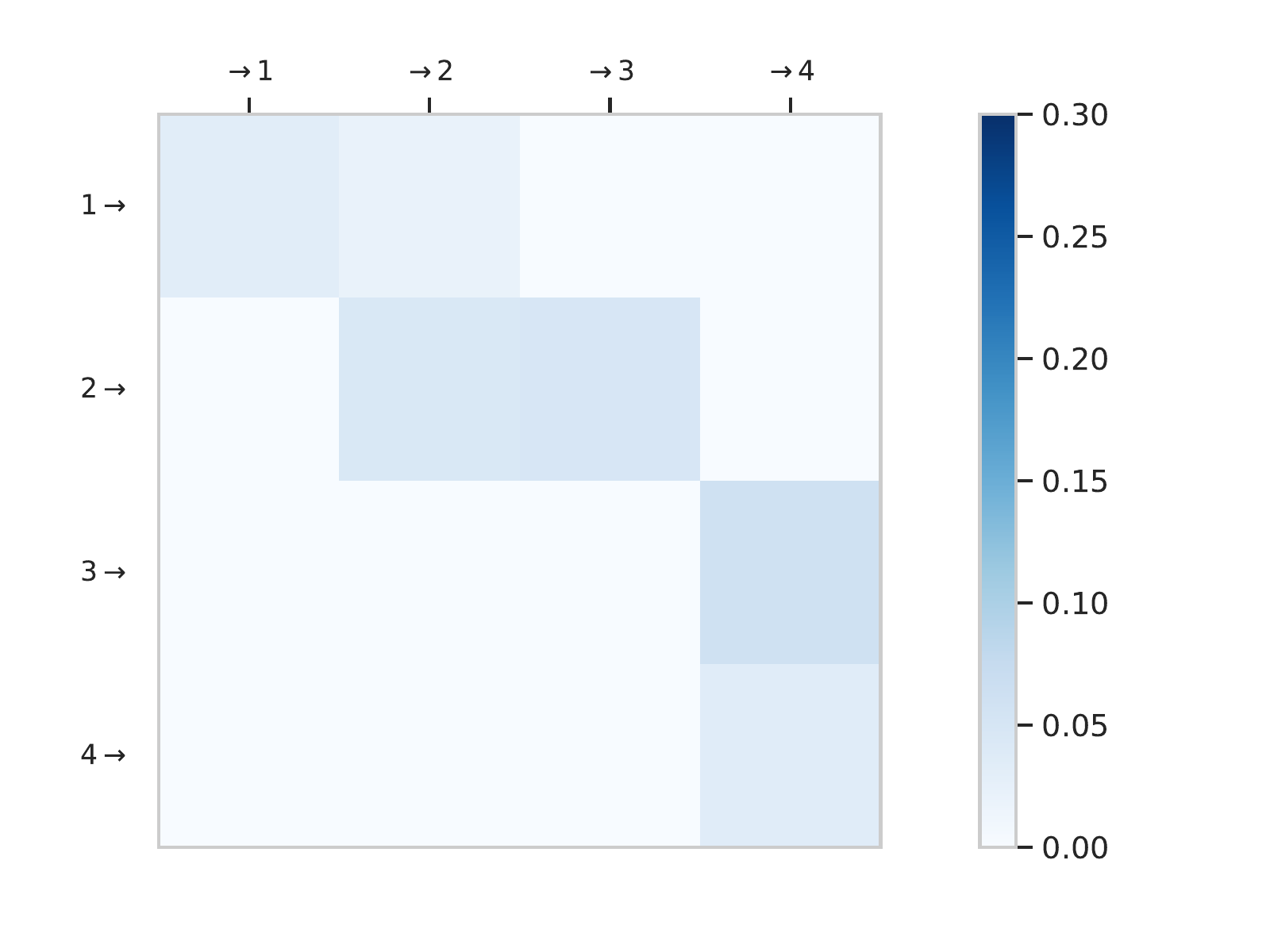}\\
\begin{minipage}{\dimexpr 40mm} \vspace{-30mm} \begin{center} \itshape \large  \textbf{$K=8$} \end{center} \end{minipage}
    \includegraphics[width=\tempwidth, trim=0.cm 0.cm 0cm  0.cm,clip]{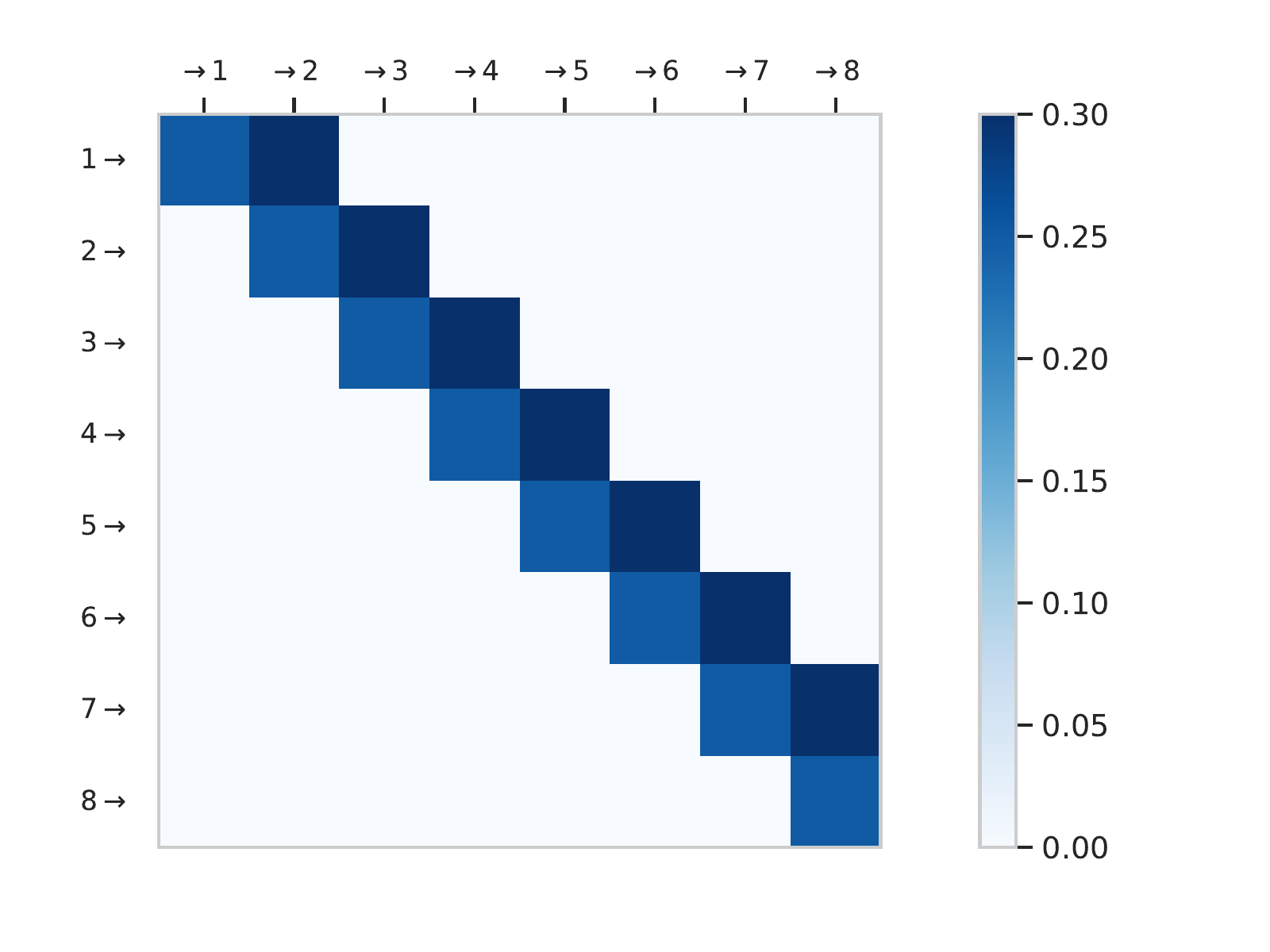}\hfil
    \includegraphics[width=\tempwidth, trim=0.cm 0.cm 0cm  0.cm,clip]{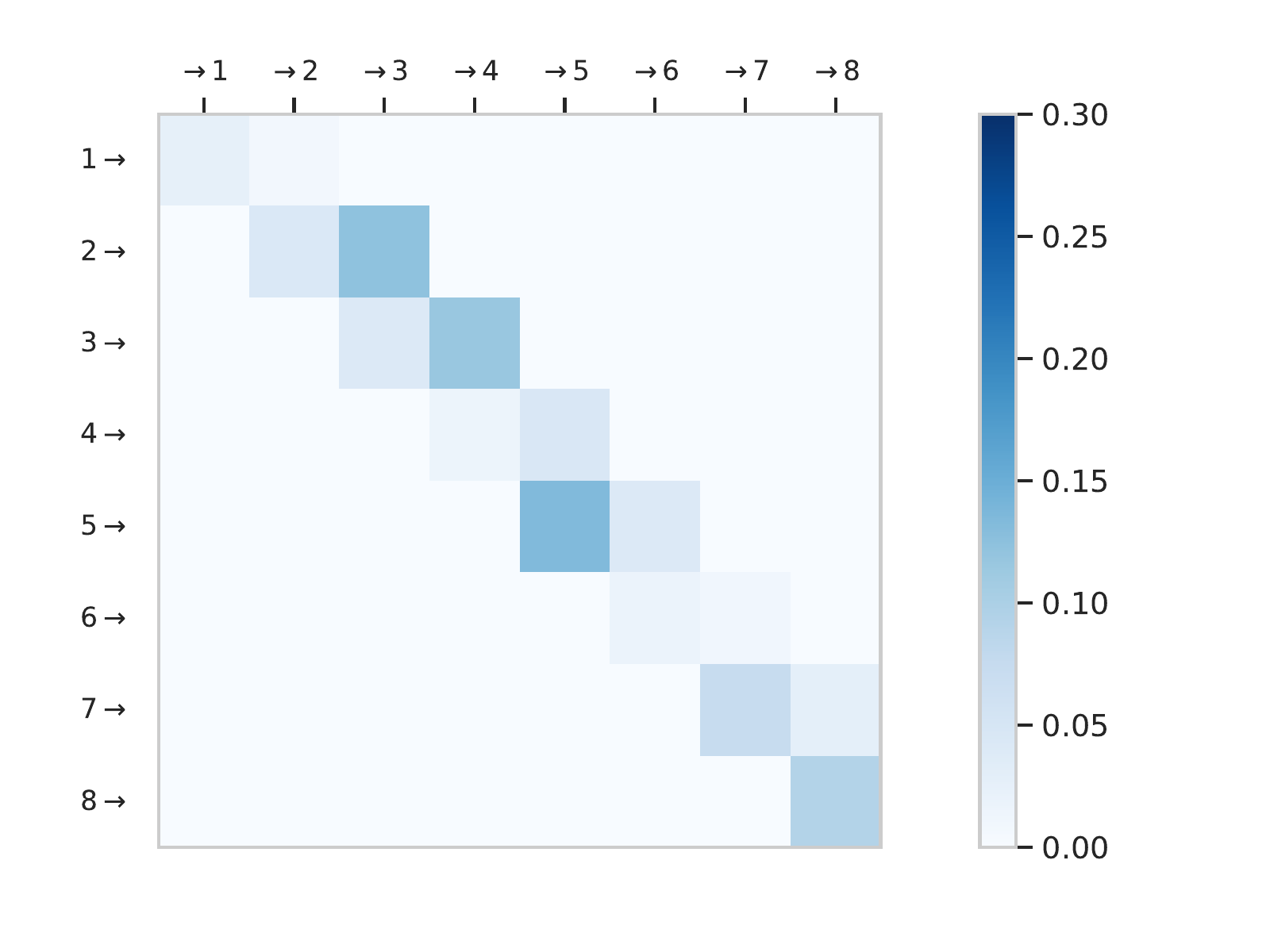}\\
\begin{minipage}{\dimexpr 40mm} \vspace{-30mm} \begin{center} \itshape \large  \textbf{$K=16$} \end{center} \end{minipage}
    \includegraphics[width=\tempwidth, trim=0.cm 0.cm 0cm  0.cm,clip]{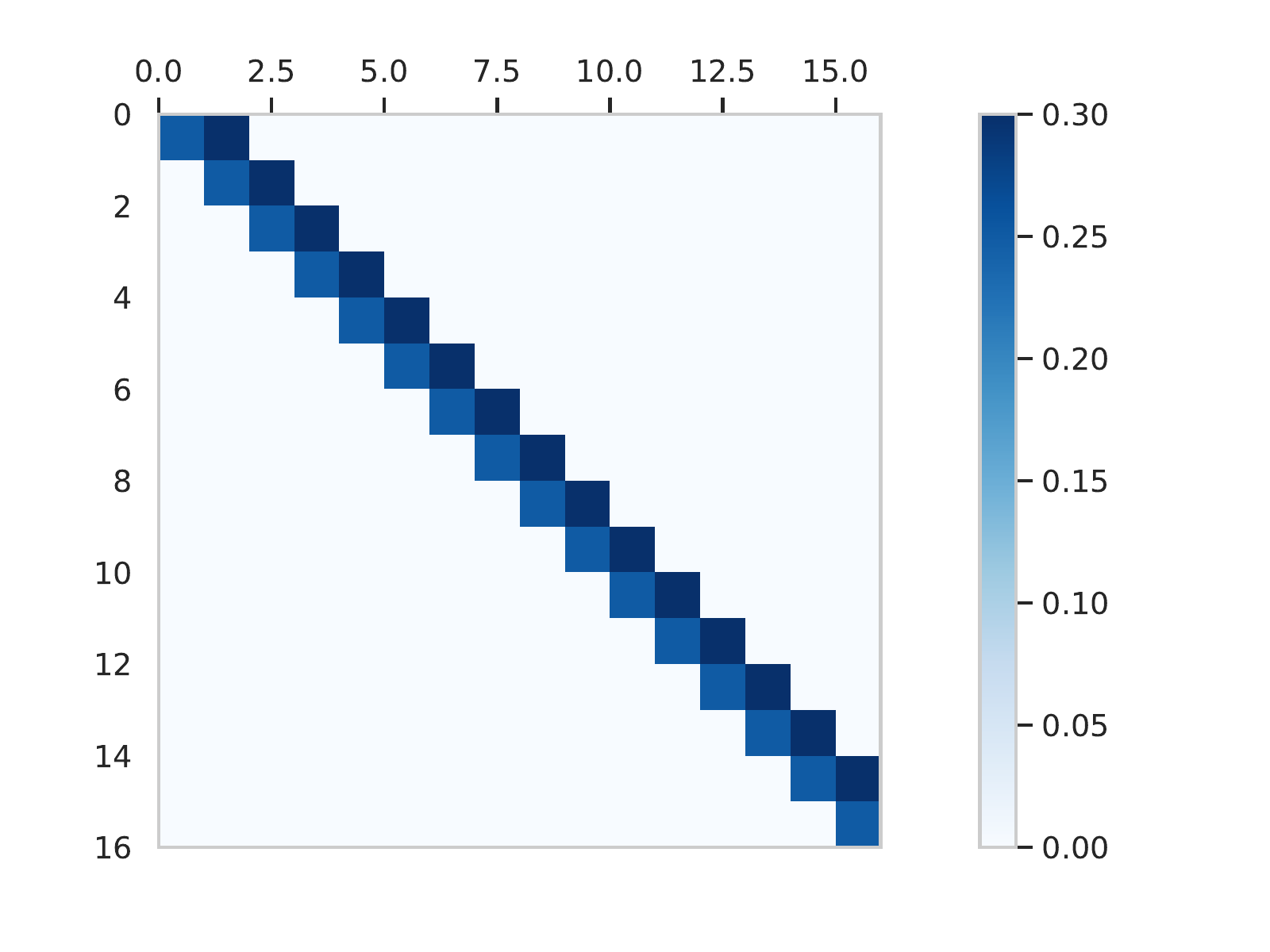}\hfil
    \includegraphics[width=\tempwidth, trim=0.cm 0.cm 0cm  0.cm,clip]{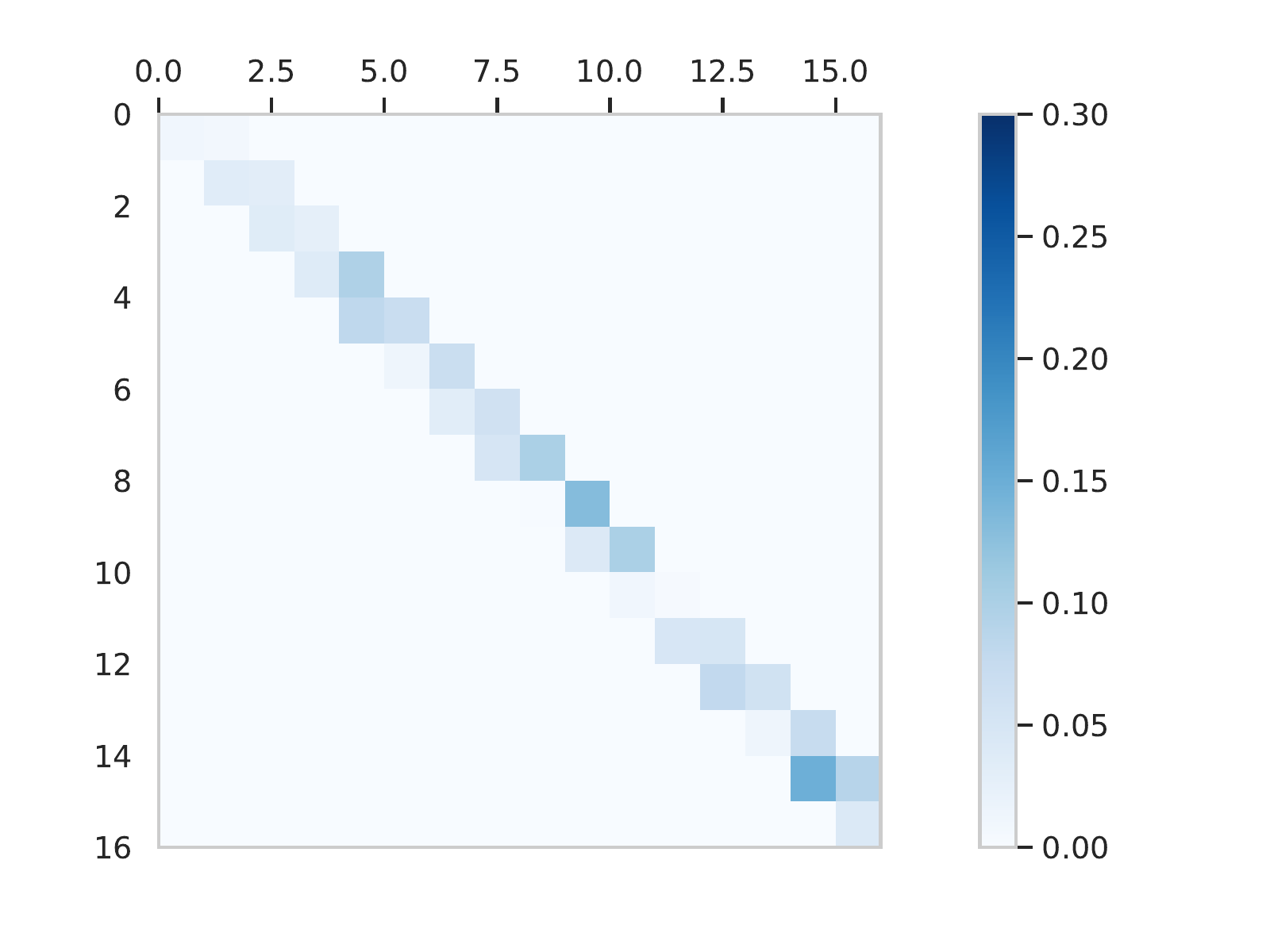}\\
    \begin{minipage}{\dimexpr 40mm} \vspace{-30mm} \begin{center} \itshape \large  \textbf{$K=32$} \end{center} \end{minipage}
    \includegraphics[width=\tempwidth, trim=0.cm 0.cm 0cm  0.cm,clip]{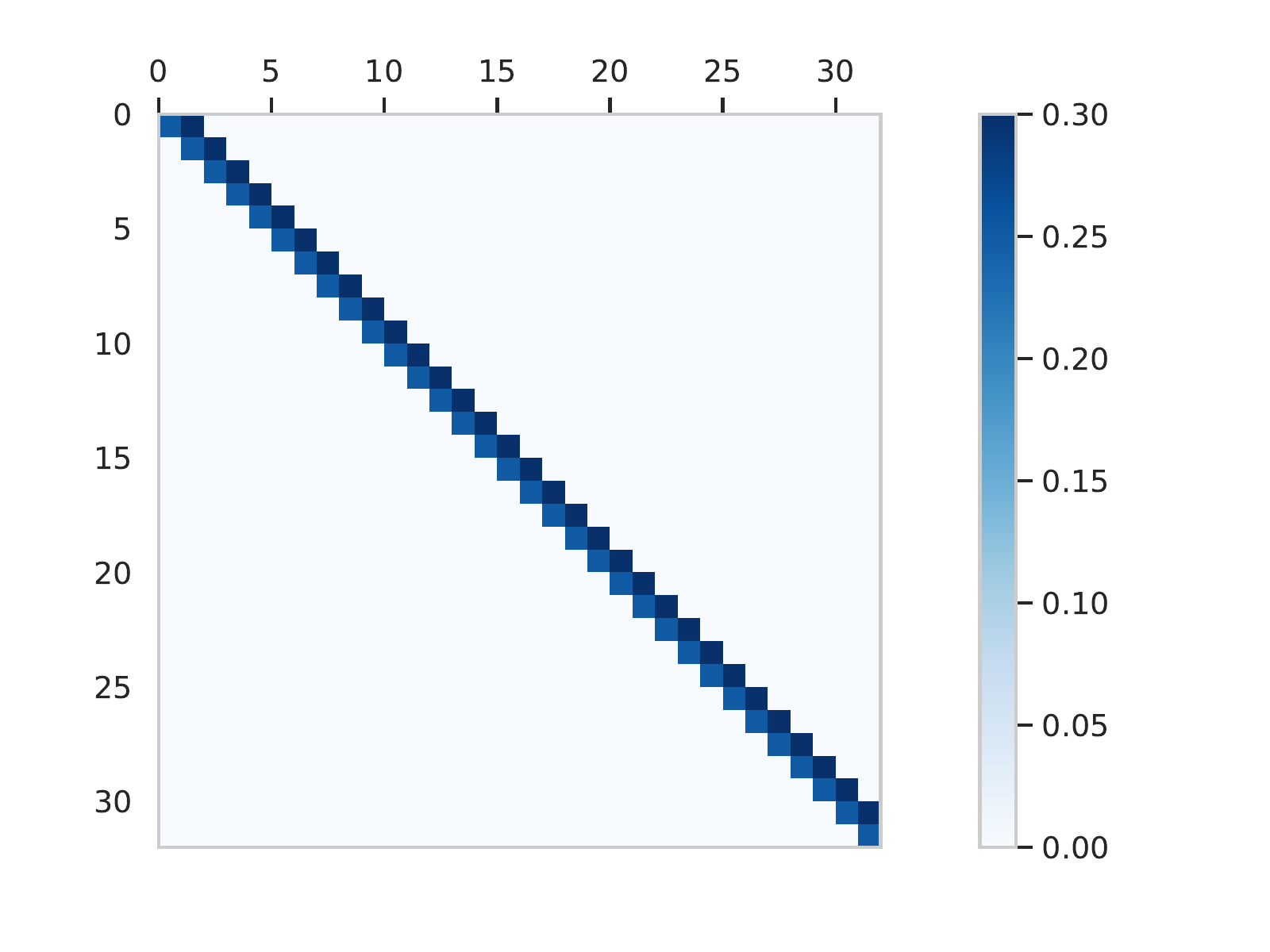}\hfil
    \includegraphics[width=\tempwidth, trim=0.cm 0.cm 0cm  0.cm,clip]{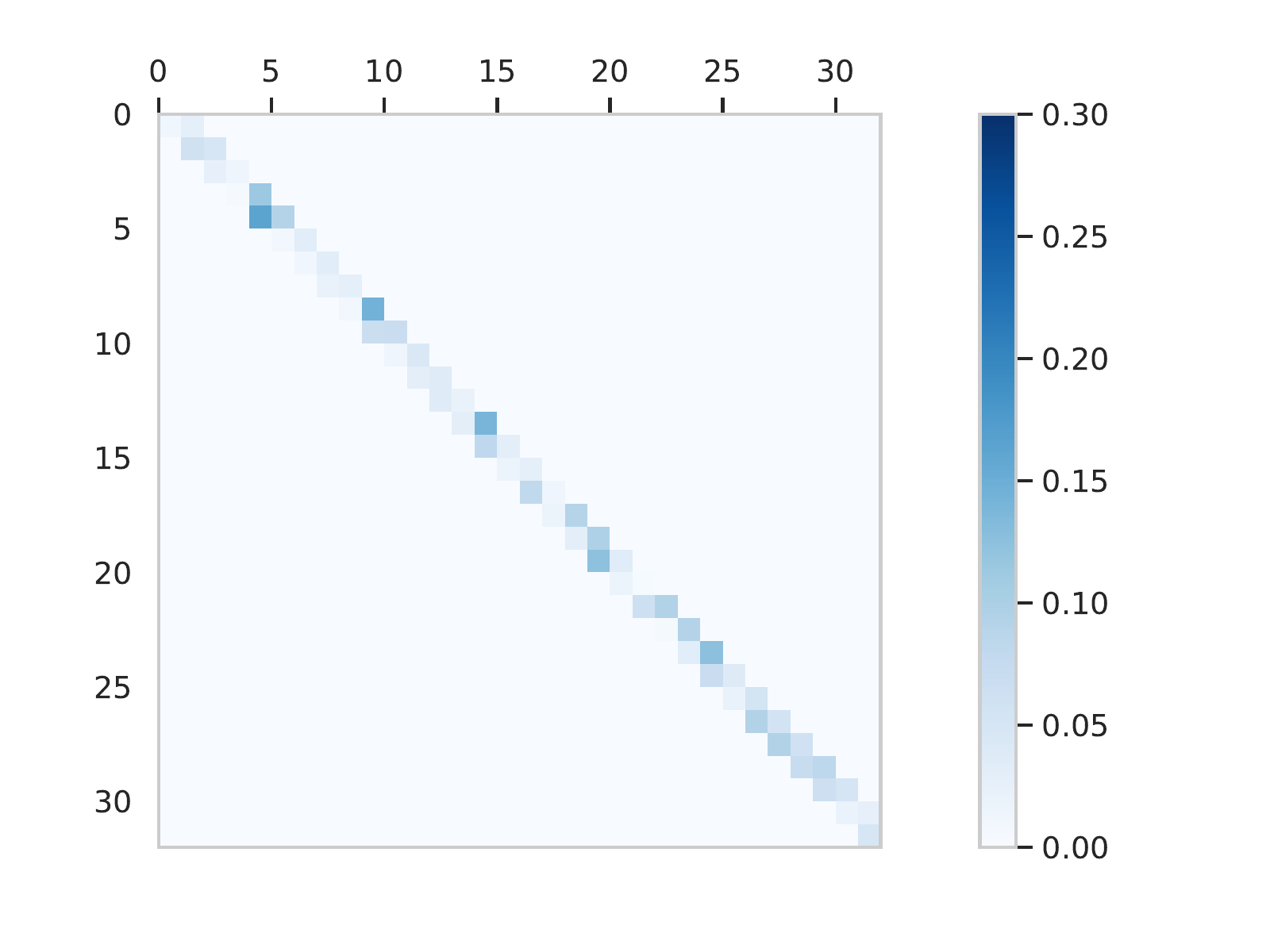}\\
        \begin{minipage}{\dimexpr 40mm} \vspace{-30mm} \begin{center} \itshape \large  \textbf{$K=64$} \end{center} \end{minipage}
    \includegraphics[width=\tempwidth, trim=0.cm 0.cm 0cm  0.cm,clip]{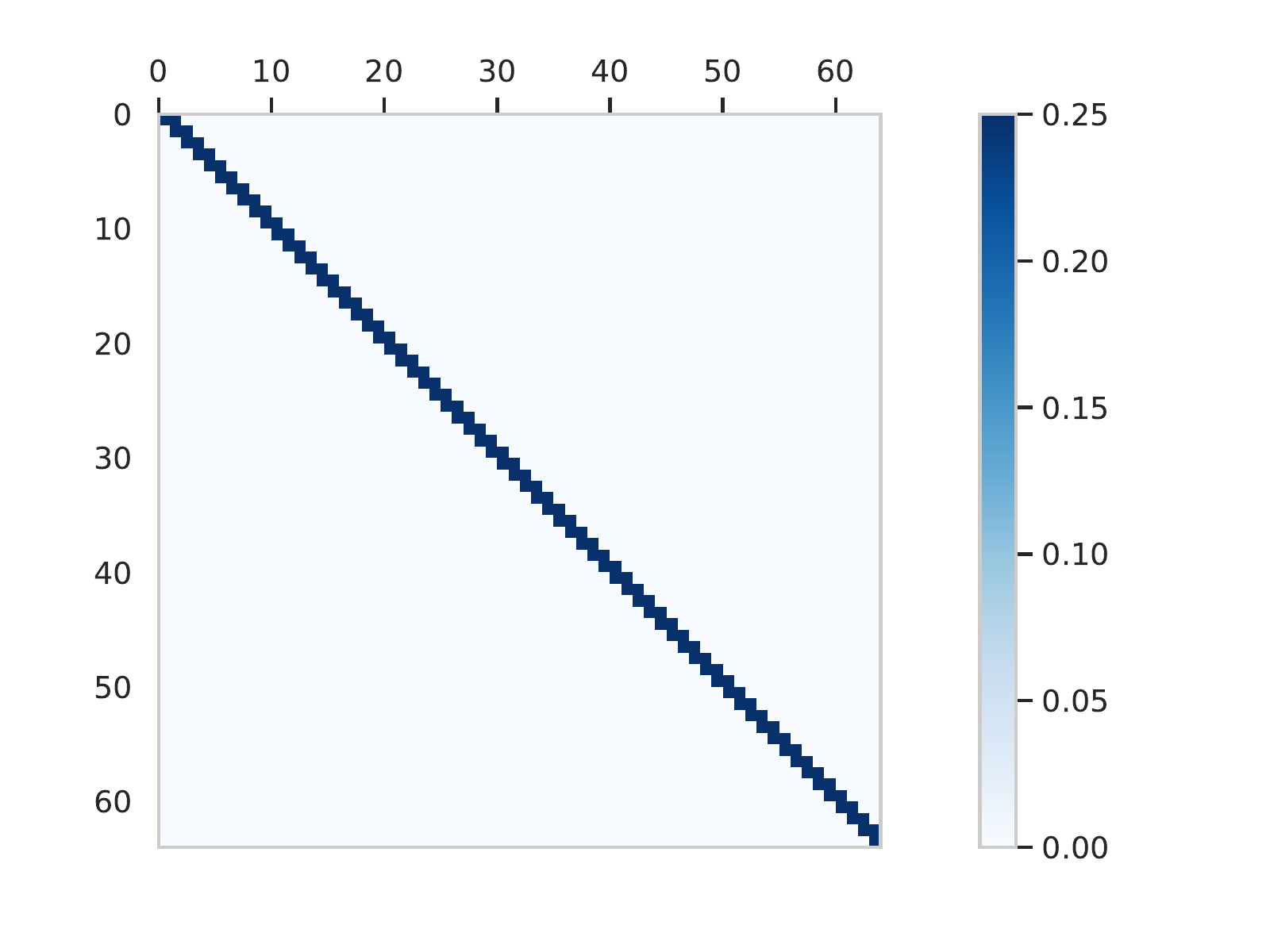}\hfil
    \includegraphics[width=\tempwidth, trim=0.cm 0.cm 0cm  0.cm,clip]{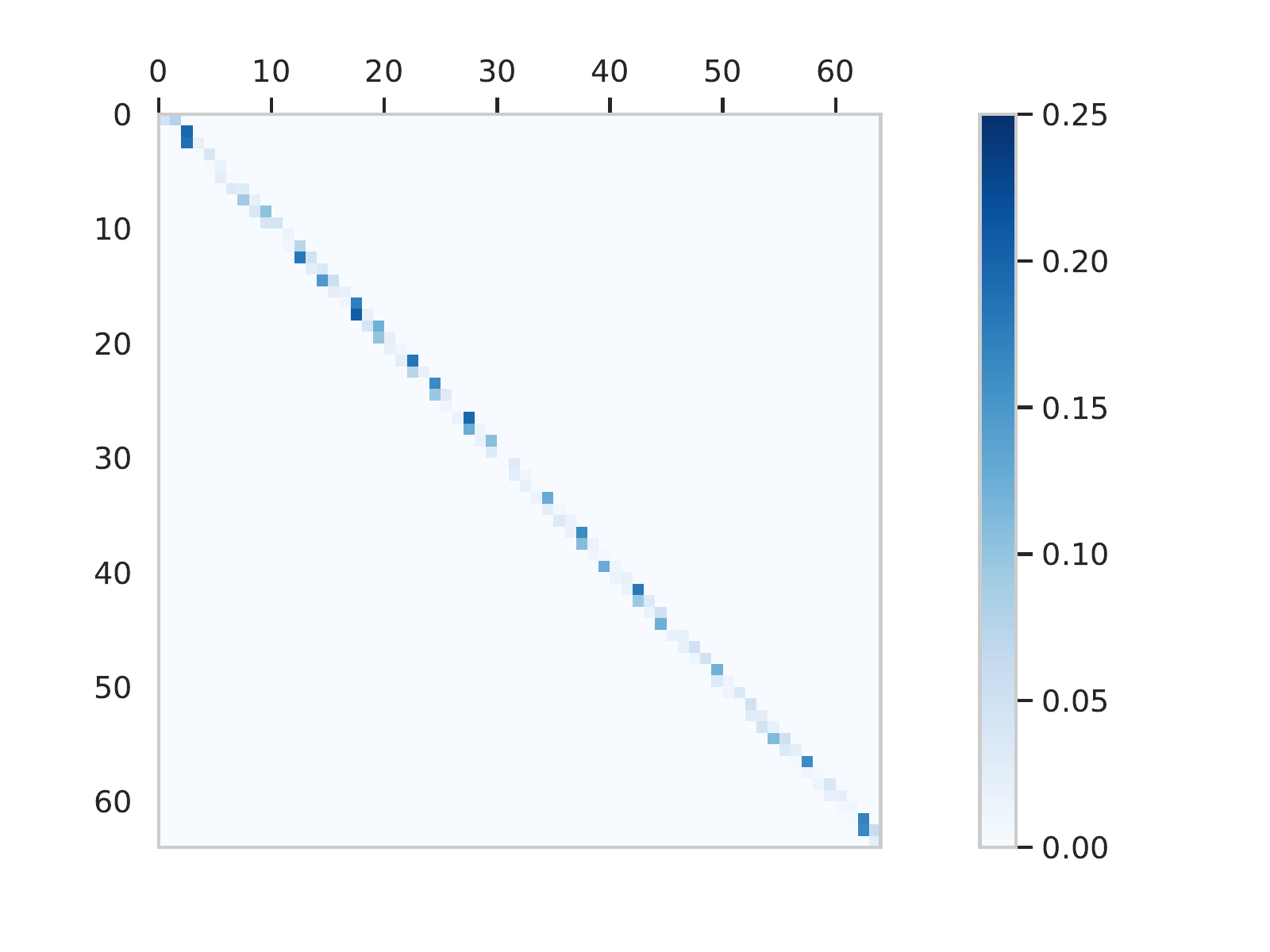}\\
\caption{Heatmaps of the $L_1$-norms of the true parameter $h_0$, i.e., the entries of the matrix $S_0 = (S^0_{lk})_{l,k} = (\norm{h_{lk}^0}_1)_{l,k}$ (left column) and $L_1$-risk, i.e., $(\mathbb{E}^{Q}[\norm{h_{lk}^0 - h_{lk}}_1])_{l,k}$ (right column) after the first step of Algorithm \ref{alg:2step_adapt_cavi}, in the Inhibition scenario of Simulation 4. The rows correspond to $K=2,4,8,16,32,64$.}
\label{fig:adaptive_VI_2step_norms_inh}
\end{figure}

\begin{figure}[hbt!]
    \centering
    \begin{subfigure}[b]{0.49\textwidth}
    \includegraphics[width=\textwidth, trim=0.cm 0.cm 0cm  0.cm,clip]{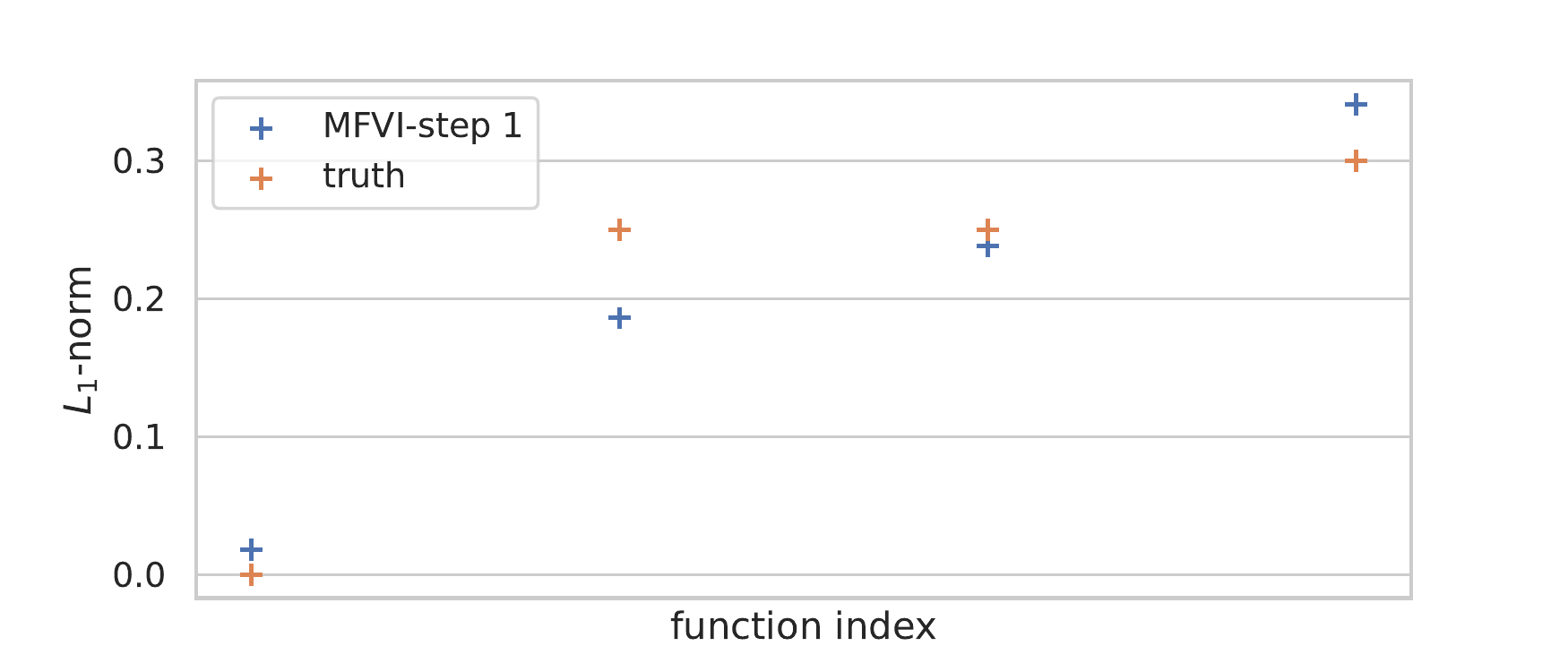}
    \caption{$K=2$}
    \end{subfigure}%
        \begin{subfigure}[b]{0.49\textwidth}
    \includegraphics[width=\textwidth, trim=0.cm 0.cm 0cm  0.cm,clip]{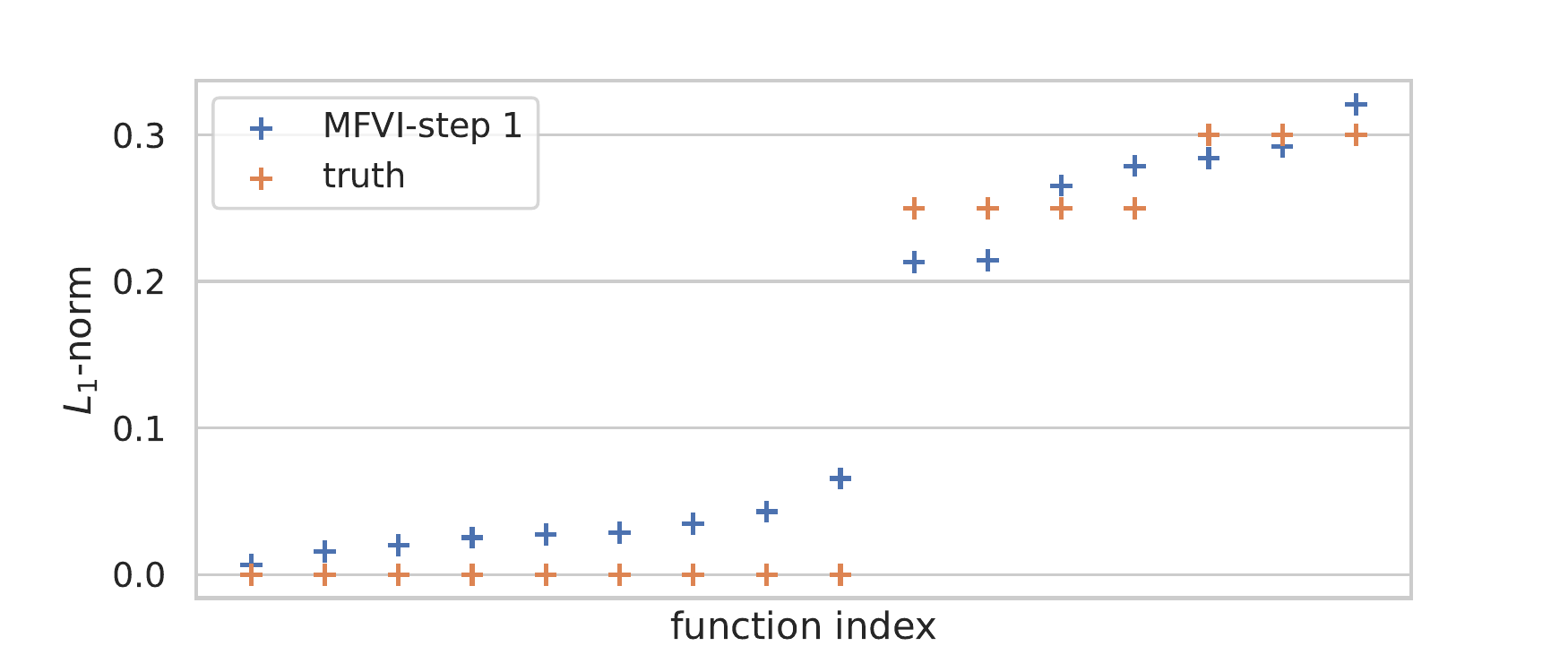}
    \caption{$K=4$}
    \end{subfigure}
     \begin{subfigure}[b]{0.49\textwidth}
    \includegraphics[width=\textwidth, trim=0.cm 0.cm 0cm  0.cm,clip]{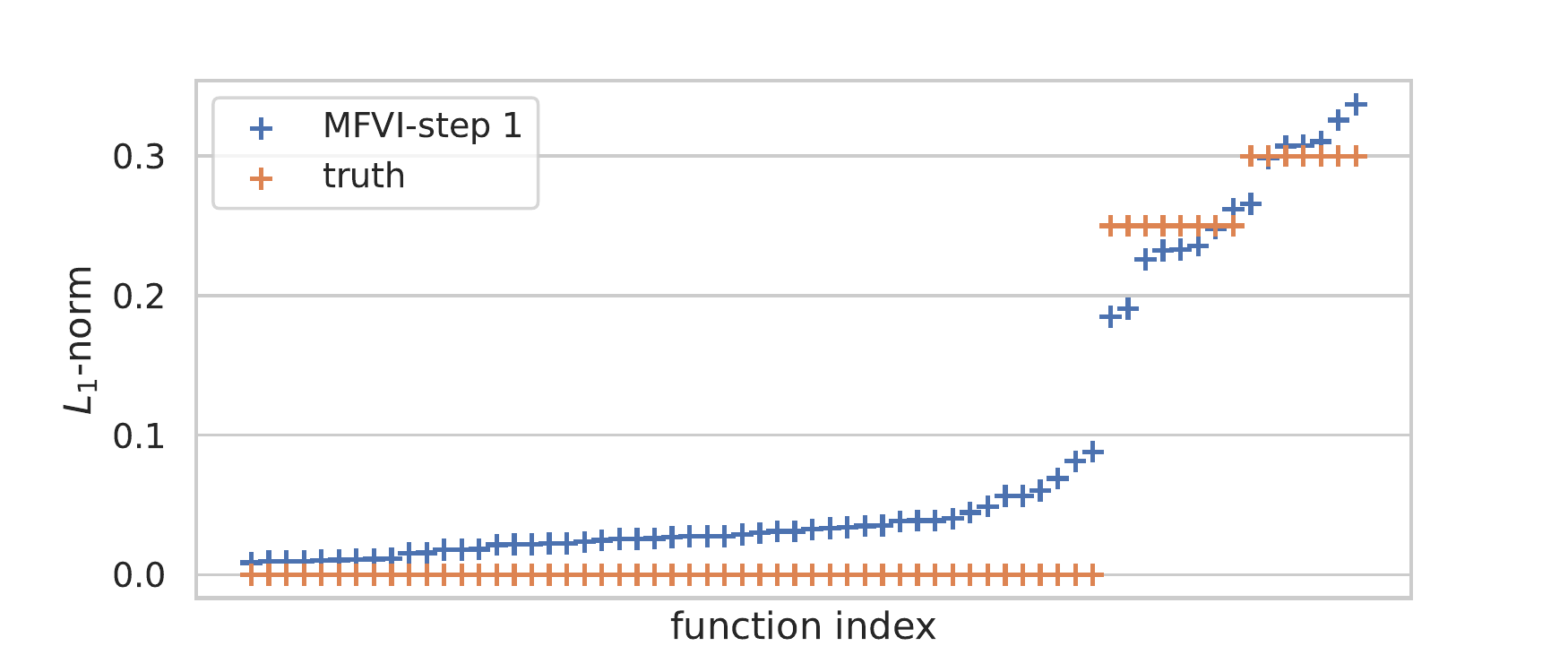}
    \caption{$K=8$}
    \end{subfigure}%
        \begin{subfigure}[b]{0.49\textwidth}
    \includegraphics[width=\textwidth, trim=0.cm 0.cm 0cm  0.cm,clip]{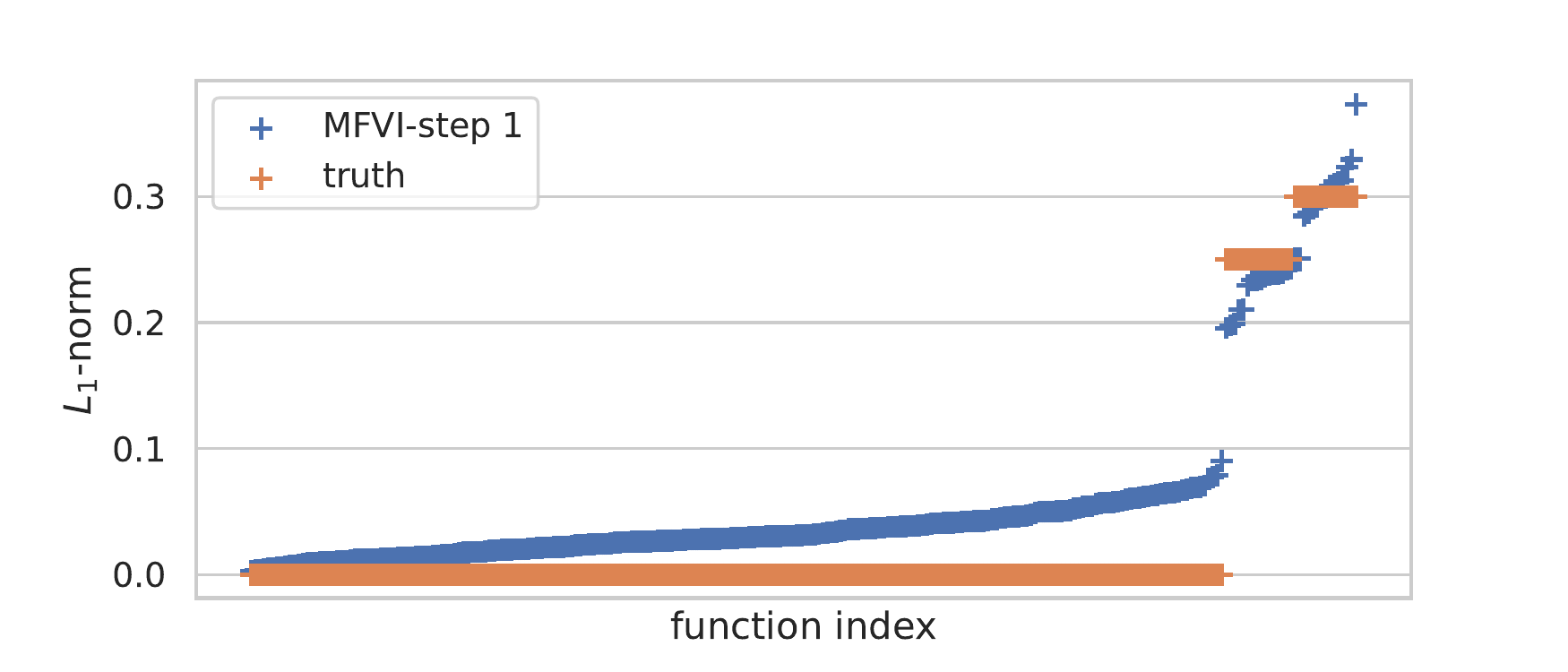}
    \caption{$K=16$}
    \end{subfigure}
            \begin{subfigure}[b]{0.49\textwidth}
    \includegraphics[width=\textwidth, trim=0.cm 0.cm 0cm  0.cm,clip]{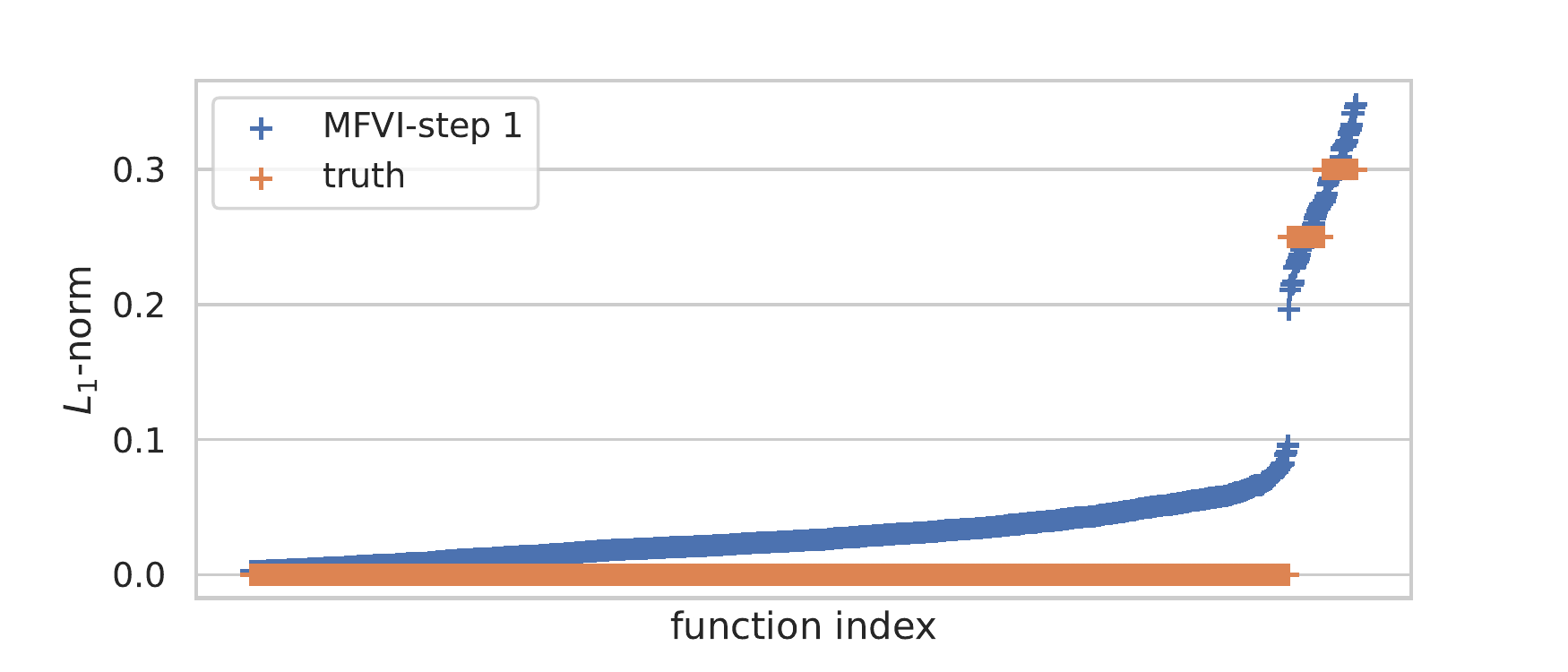}
    \caption{$K=32$}
    \end{subfigure}
                \begin{subfigure}[b]{0.49\textwidth}
    \includegraphics[width=\textwidth, trim=0.cm 0.cm 0cm  0.cm,clip]{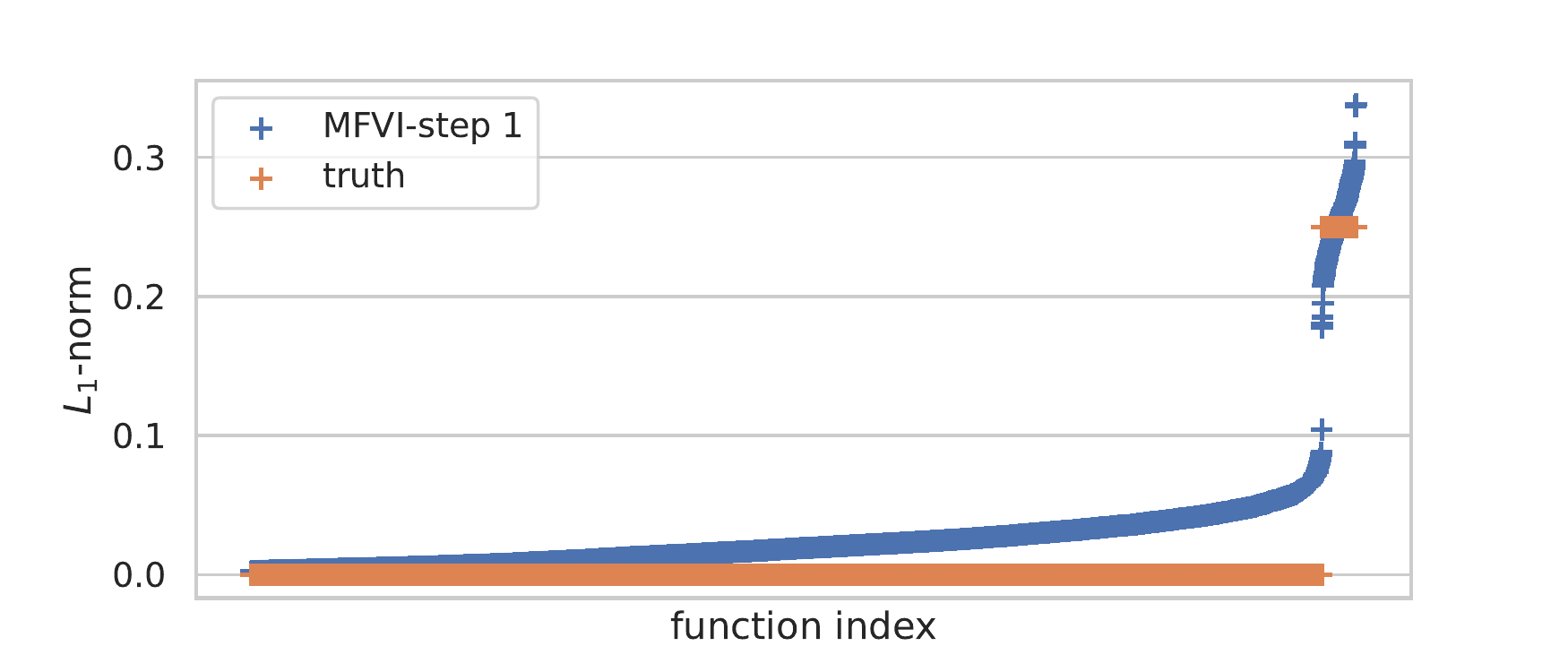}
    \caption{$K=64$}
    \end{subfigure}
\caption{Estimated $L_1$-norms after the first step of Algorithm \ref{alg:2step_adapt_cavi} (in blue), and ground-truth norms (in orange), plotted in increasing order, in the Inhibition scenario of Simulation 4, for the models with $K \in \{2,4,8,16, 32, 64\}$. } %The threshold in our algorithm $\eta_0 = 0.07$ is plotted in dotted red line.}
\label{fig:adaptive_VI_2step_norms_threshold_inhibition}
\end{figure}

\begin{figure}[hbt!]
\centering
\setlength{\tempwidth}{.2\linewidth}\centering
\settoheight{\tempheight}{\includegraphics[width=\tempwidth]{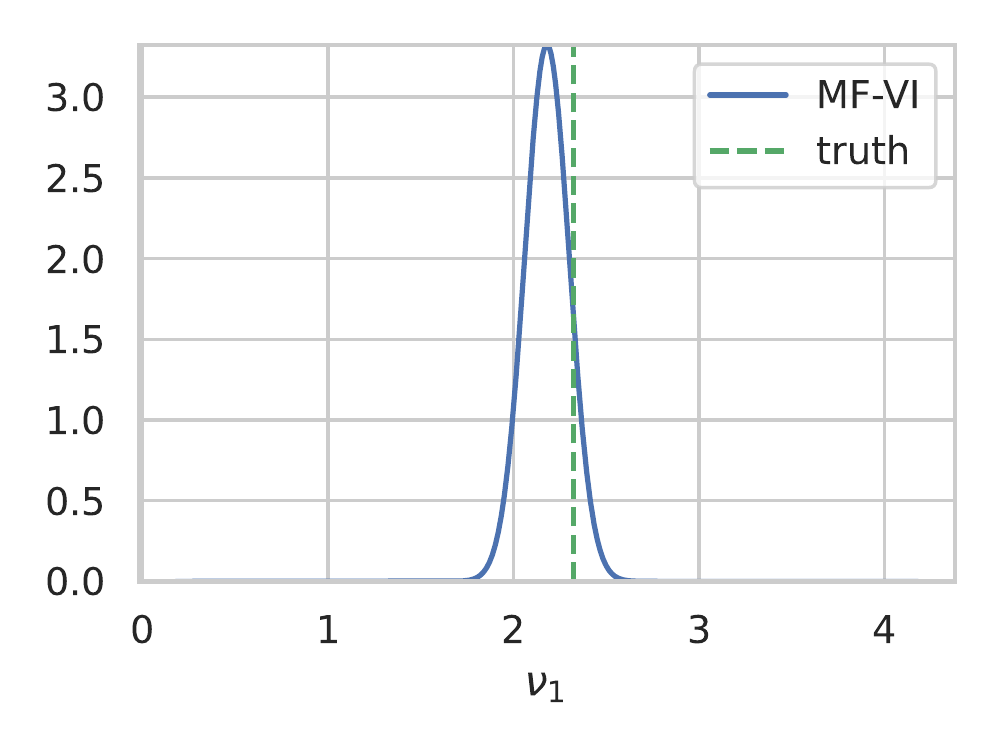}}
\hspace{-30mm} \fbox{\begin{minipage}  {\dimexpr 20mm} \begin{center} \itshape \large  \textbf{Inhibition} \end{center} \end{minipage}}
%\hspace{-30mm} 
\columnname{Background $\nu_1$}\hfil
\columnname{Interaction functions $h_{11}$ and $h_{21}$ }\\
\begin{minipage}{\dimexpr 20mm}  \vspace{-35mm} \flushright{\itshape \large  \textbf{$K=2$} } \end{minipage}
    \includegraphics[width=\tempwidth, trim=0.cm 0.cm 0cm 0cm,clip]{figs/simu_K2_estimated_nu_inh.pdf}\hfil
    \includegraphics[width=.6\linewidth, trim=0.cm 0.cm 0cm  1.cm,clip]{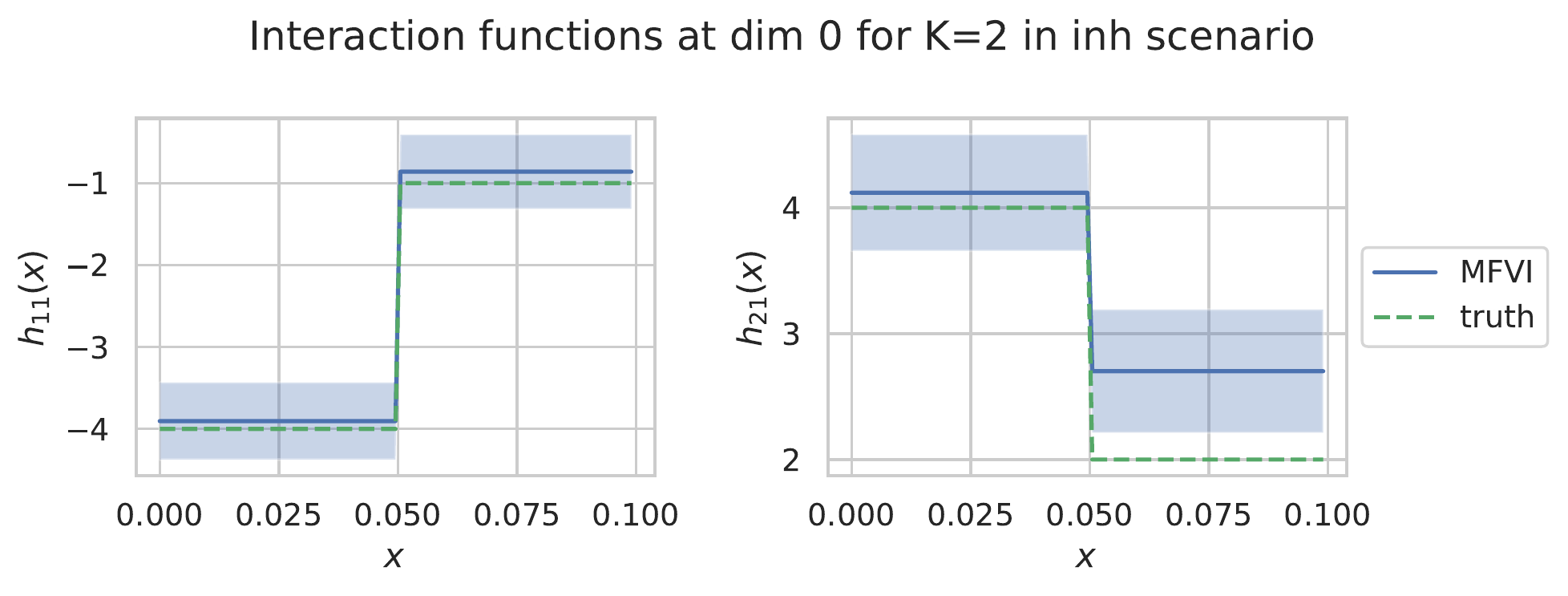}\\
\begin{minipage}{\dimexpr 20mm} \vspace{-35mm} \flushright{\itshape \large  \textbf{$K=4$} } \end{minipage}
    \includegraphics[width=\tempwidth, trim=0.cm 0.cm 0cm 0.cm,clip]{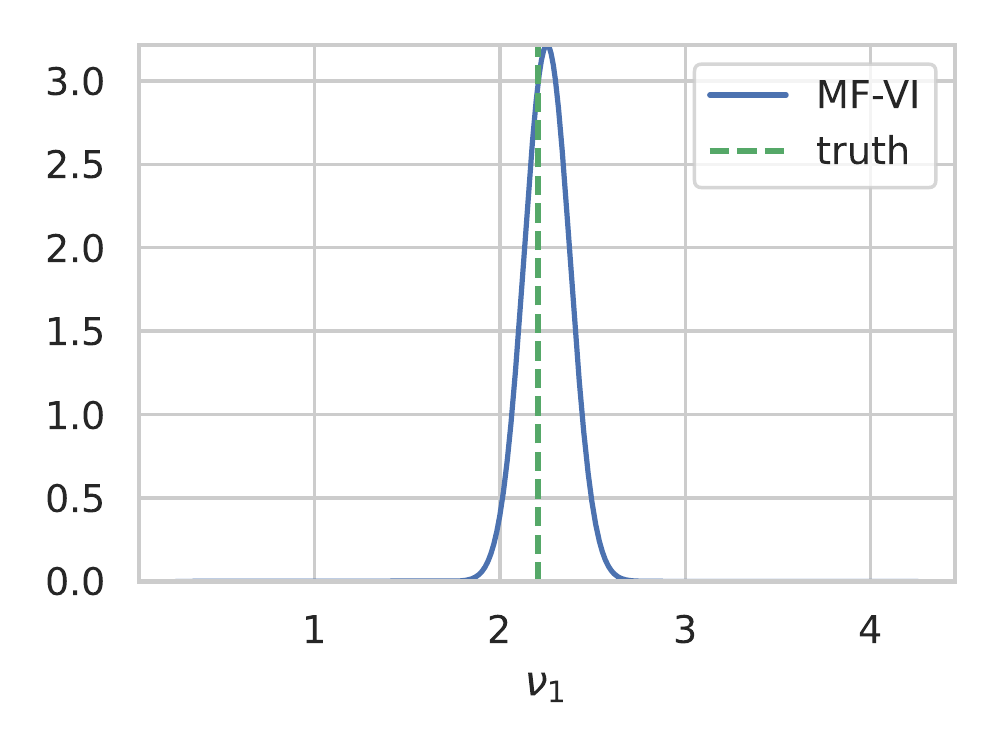}\hfil
    \includegraphics[width=.6\linewidth, trim=0.cm 0.cm 0cm  1.cm,clip]{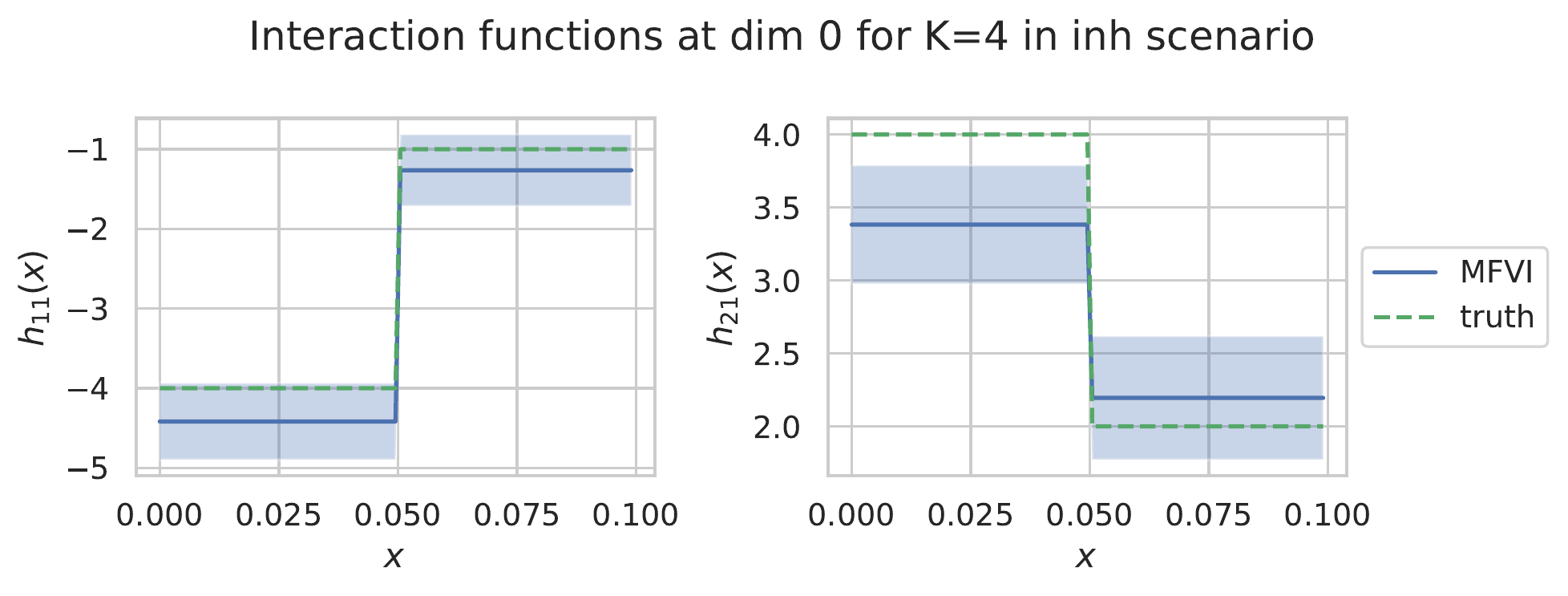}\\
\begin{minipage}{\dimexpr 20mm} \vspace{-35mm} \flushright{\itshape \large  \textbf{$K=8$} } \end{minipage}
    \includegraphics[width=\tempwidth, trim=0.cm 0.cm 0cm 0cm,clip]{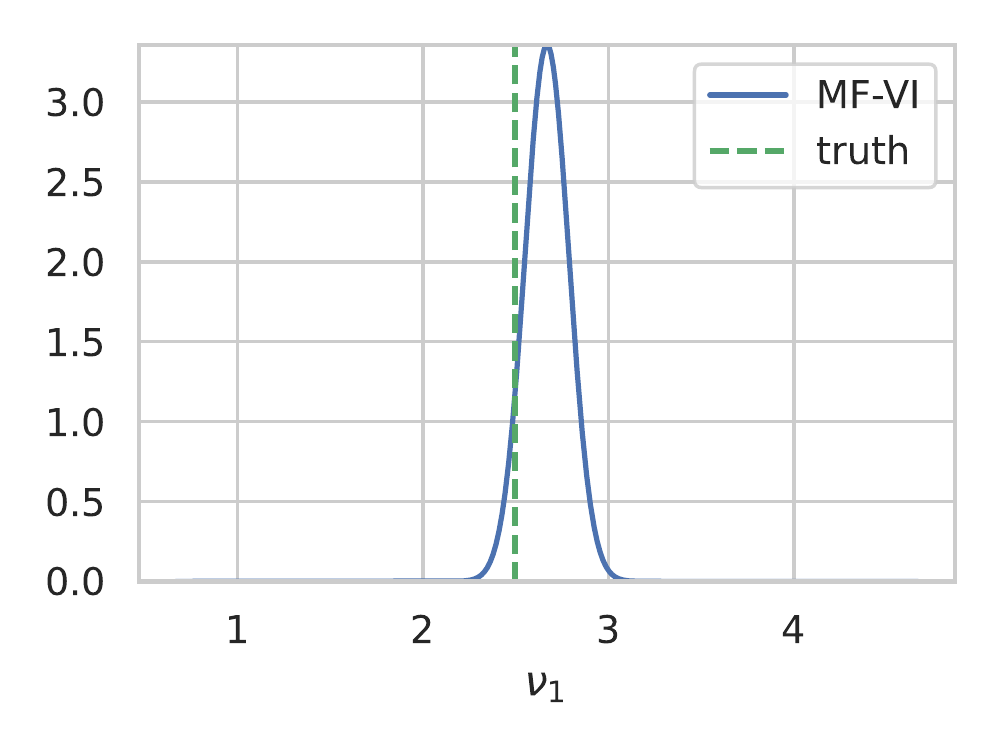}\hfil
    \includegraphics[width=.6\linewidth, trim=0.cm 0.cm 0cm  1.cm,clip]{figs/simu_K8_estimated_h_dim_0_exc.pdf}\\
\begin{minipage}{\dimexpr 20mm} \vspace{-35mm} \flushright{\itshape \large  \textbf{$K=16$} } \end{minipage}
    \includegraphics[width=\tempwidth, trim=0.cm 0.cm 0cm 0cm,clip]{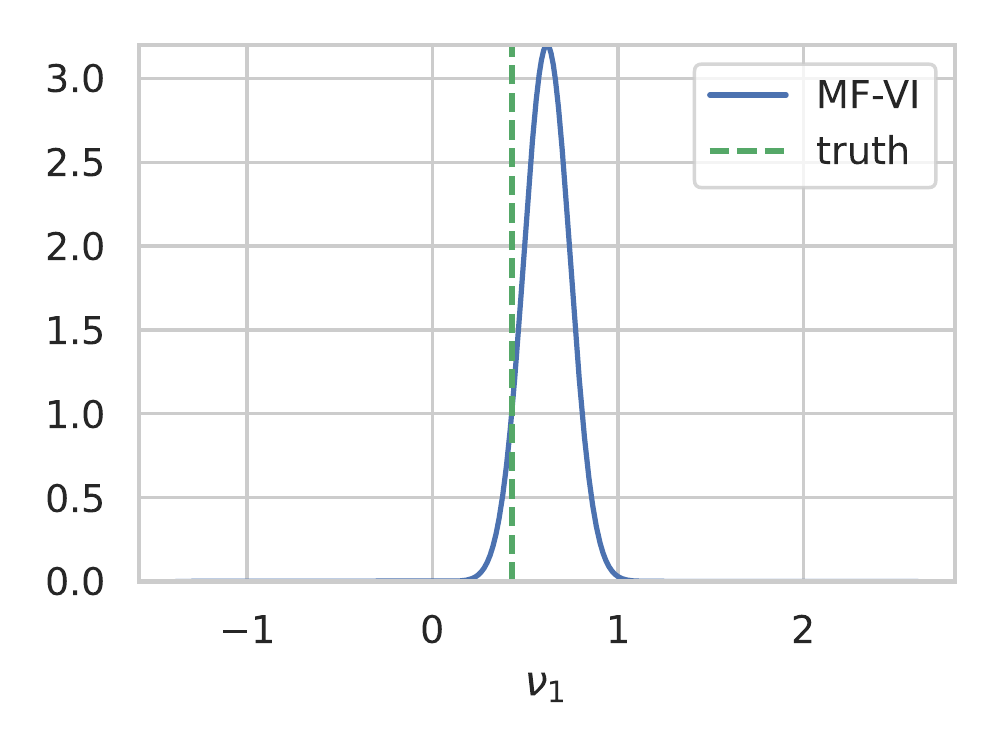}\hfil
    \includegraphics[width=.6\linewidth, trim=0.cm 0.cm 0cm  1.cm,clip]
    {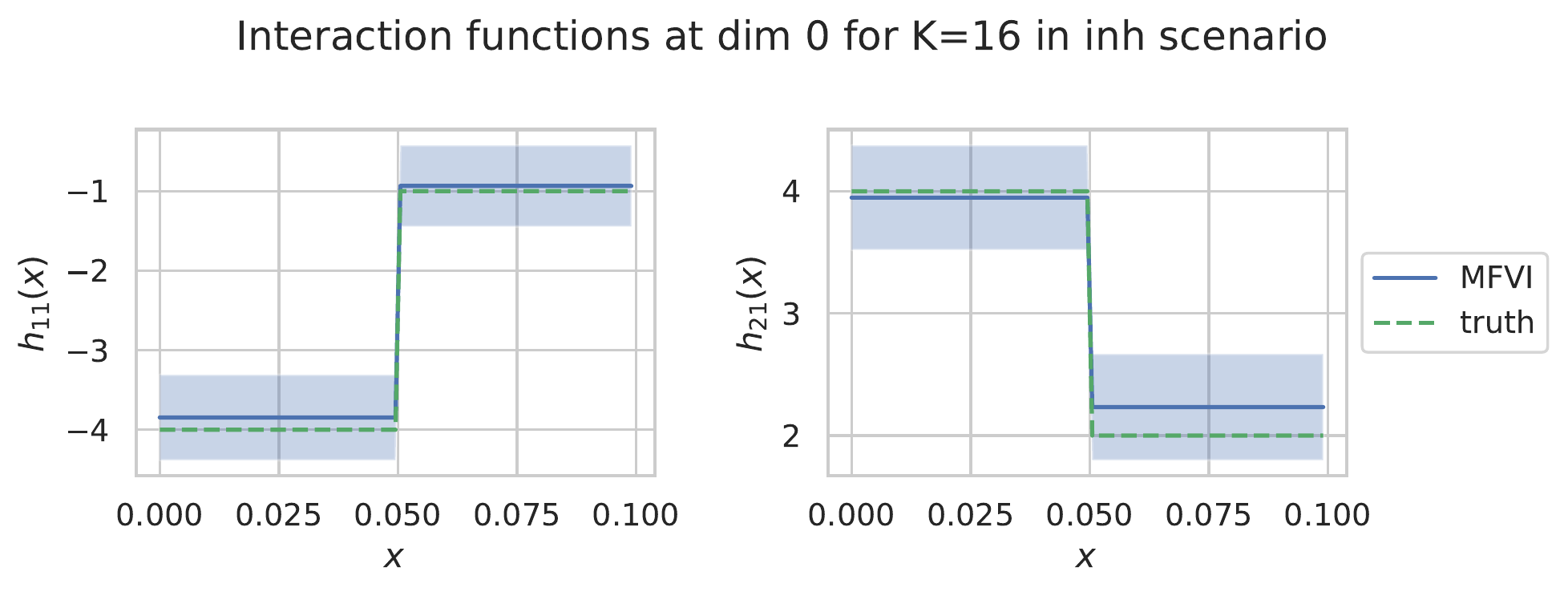}\\
\begin{minipage}{\dimexpr 20mm} \vspace{-35mm} \flushright{\itshape \large  \textbf{$K=32$} } \end{minipage}
    \includegraphics[width=\tempwidth, trim=0.cm 0.cm 0cm 0cm,clip]{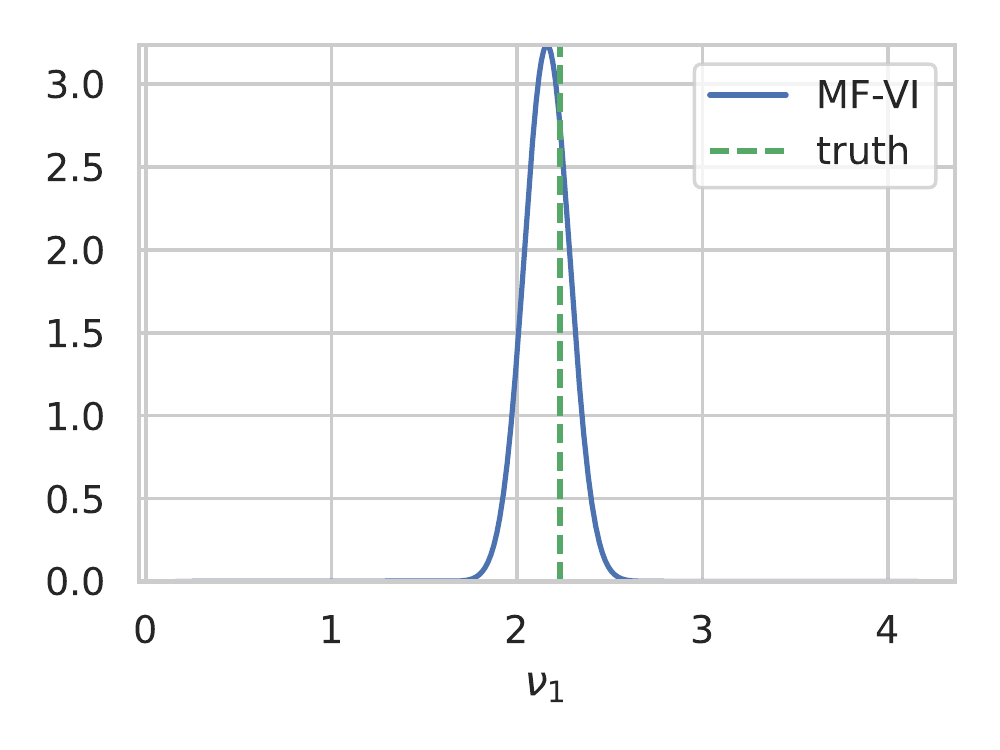}\hfil
    \includegraphics[width=.6\linewidth, trim=0.cm 0.cm 0cm  1.cm,clip]{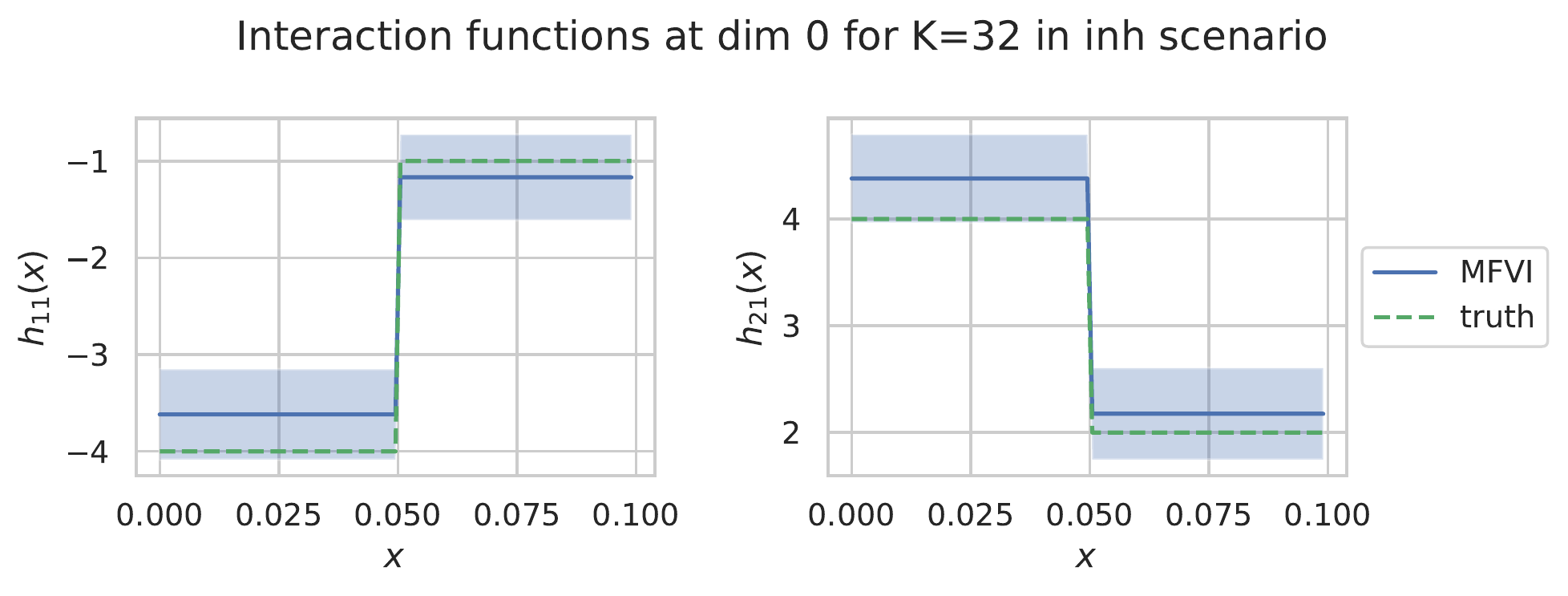}
\begin{minipage}{\dimexpr 20mm} \vspace{-35mm} \flushright{\itshape \large  \textbf{$K=64$} } \end{minipage}
    \includegraphics[width=\tempwidth, trim=0.cm 0.cm 0cm 0cm,clip]{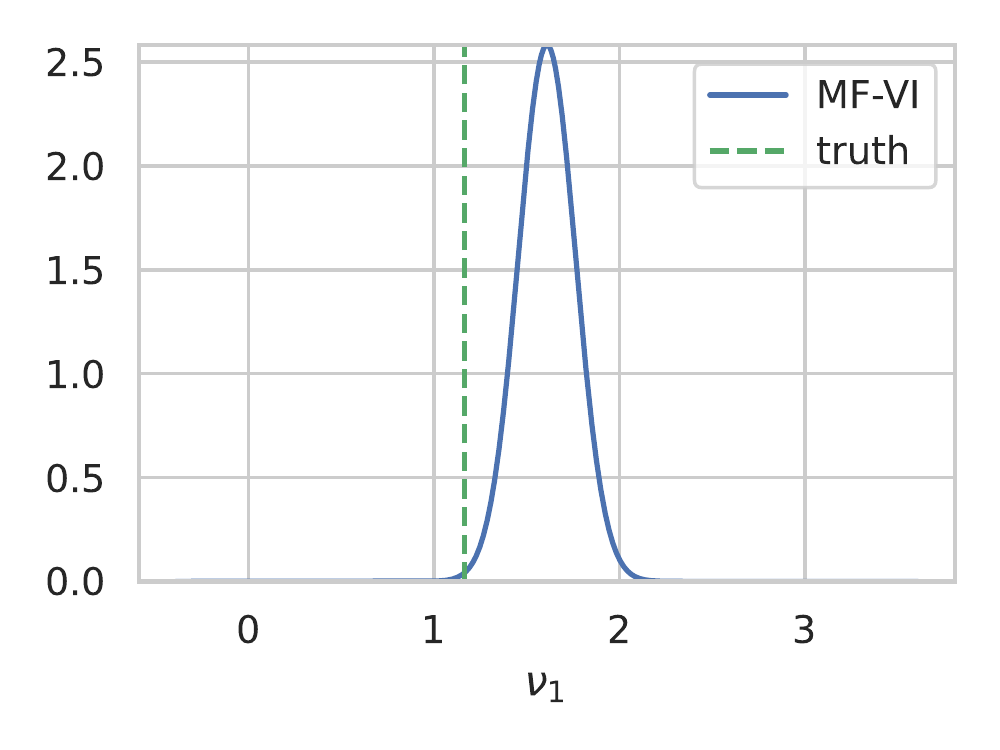}\hfil
    \includegraphics[width=.6\linewidth, trim=0.cm 0.cm 0cm  1.cm,clip]{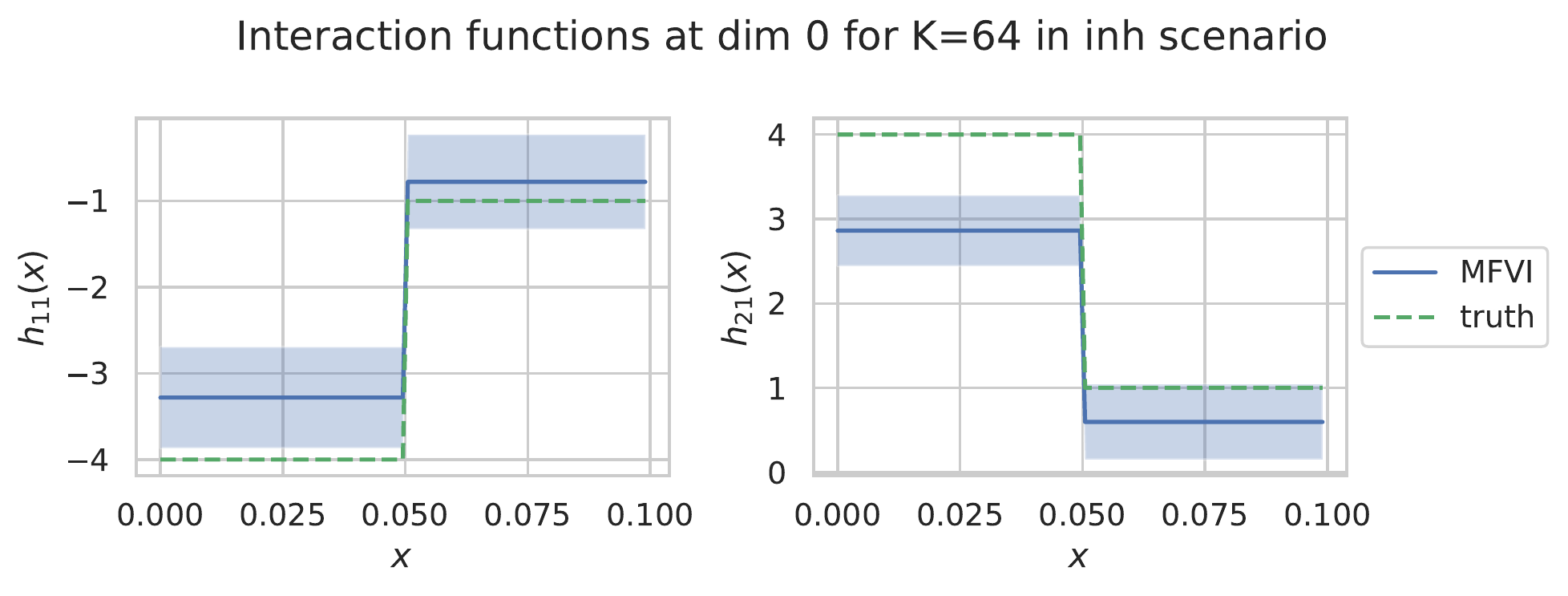}
\caption{Model-selection variational posterior distributions on $\nu_1$ (left column) and interaction functions $h_{11}$ and $ h_{21}$ (second and third columns) in the Inhibition scenario and multivariate sigmoid models of Simulation 4, computed with our two-step mean-field variational (MF-VI) algorithm (Algorithm \ref{alg:2step_adapt_cavi}). The different rows correspond to different multivariate settings $K=2,4,8,16,32, 64$.}
\label{fig:2step_adaptive_VI_inh_f}
\end{figure}

This section contains our results for the Inhibition setting of Simulation 4, i.e., the estimated graphs in (Figures \ref{fig:graphs_exc} and \ref{fig:graphs_inh}), the heatmaps of the risk on the interaction functions in Figure \ref{fig:adaptive_VI_2step_norms_inh}, the estimated $L_1$-norms after the first step of Algorithm \ref{alg:2step_adapt_cavi} in Figure \ref{fig:adaptive_VI_2step_norms_threshold_inhibition}, and the variational posterior distribution on the subset of the parameter in Figure \ref{fig:2step_adaptive_VI_inh_f}.

\FloatBarrier

\subsection{Simulation 5}\label{app:simu_5}

\begin{figure}
    \centering
    \begin{subfigure}[b]{\textwidth}
    \centering
    \includegraphics[width=0.49\textwidth, trim=0.cm 0.cm 0cm  0.cm,clip]{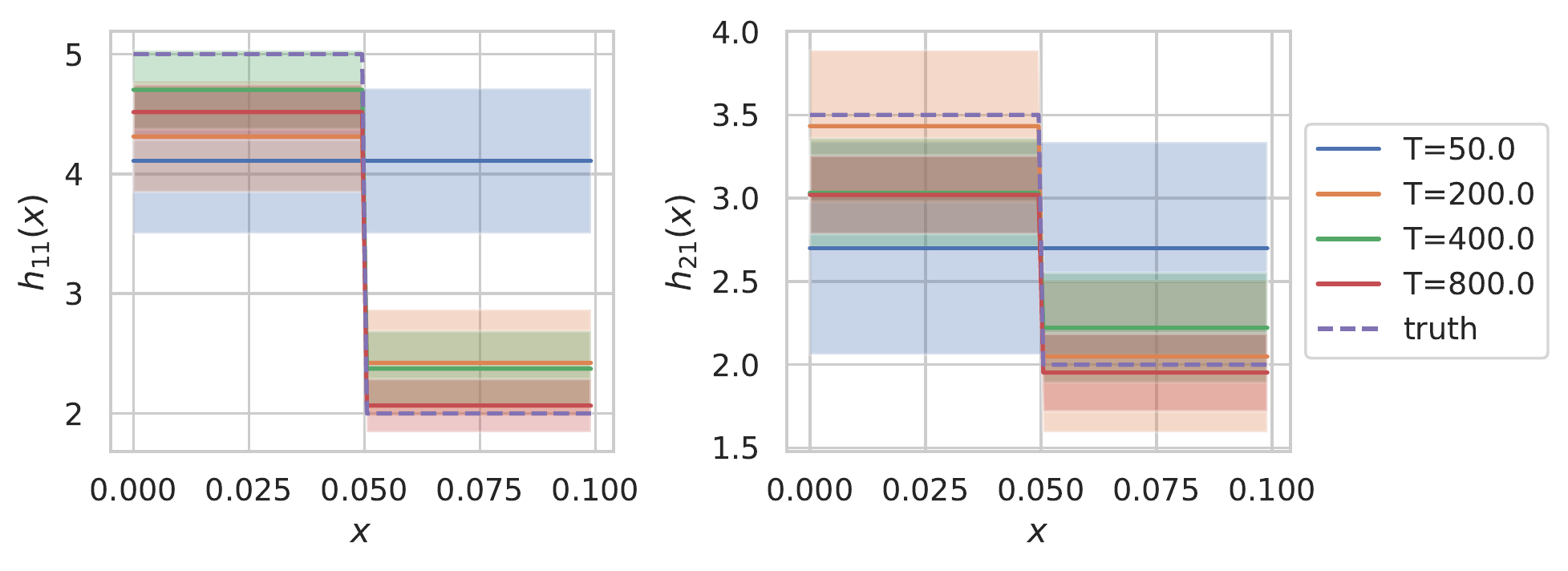}
    \includegraphics[width=0.49\textwidth, trim=0.cm 0.cm 0cm  0.cm,clip]{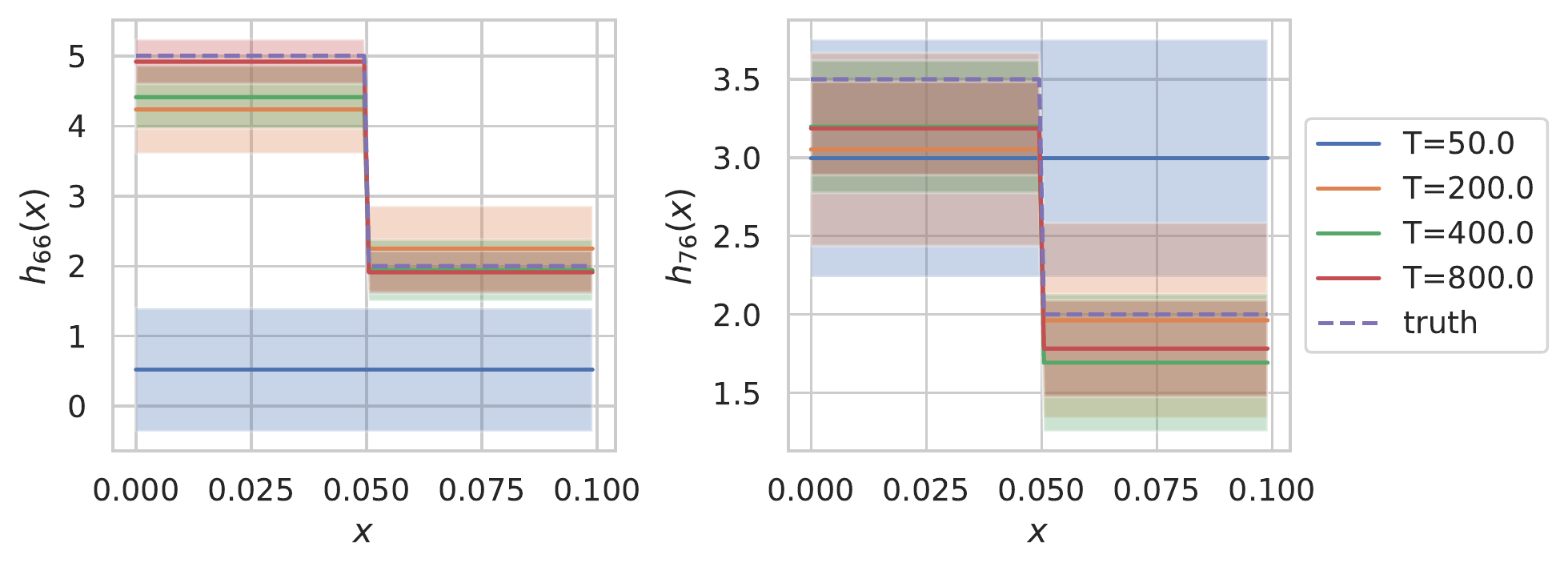}
    \caption{\emph{Excitation} scenario}
    \end{subfigure}
        \begin{subfigure}[b]{\textwidth}
    \centering
    \includegraphics[width=0.49\textwidth, trim=0.cm 0.cm 0cm  0.cm,clip]{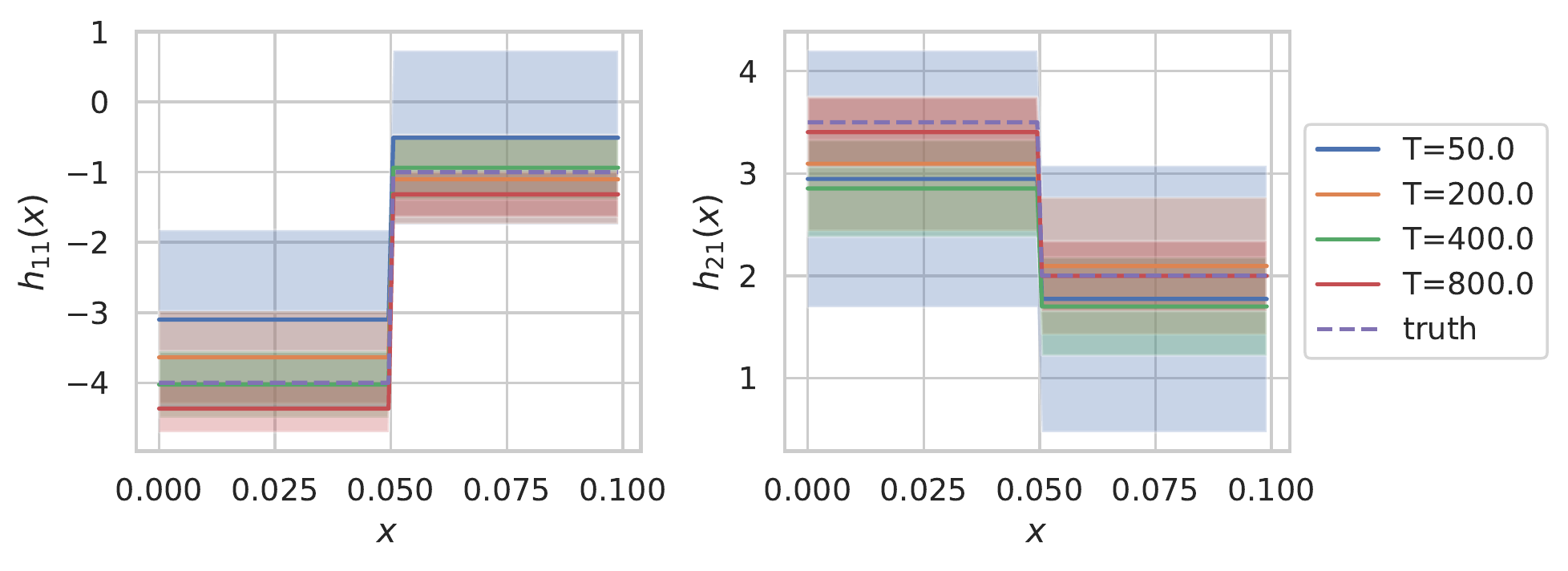}
        \includegraphics[width=0.49\textwidth, trim=0.cm 0.cm 0cm  0.cm,clip]{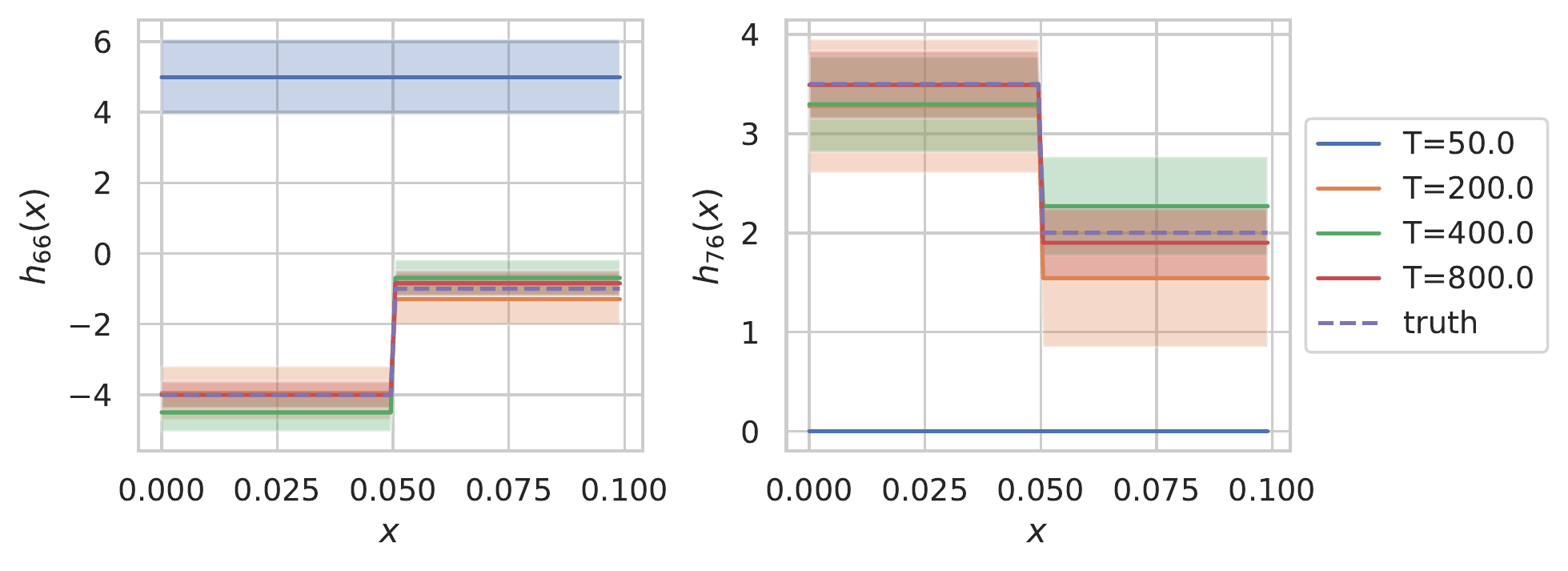}
    \caption{\emph{Inhibition} scenario}
    \end{subfigure}%
    \caption{Model-selection variational posterior on two interaction functions $h_{66}$ and $h_{76}$, for different observation lengths  $T \in \{50,200,400, 800\}$, in the  \emph{Excitation}  and \emph{Inhibition} scenarios in Simulation 5 with $K=10$. We note that in this simulation, the true number of basis functions is 2 and is well recovered  for all values of $T$. The estimation of these two interaction functions  is poor for the smallest $T$, however, it improves when $T$ increases.}
    \label{fig:est_h_T}
\end{figure}

In this  section, we report some characteristics of the simulated data in Simulation 5, in particular the number of points and excursions in each setting (see Table \ref{tab:simu_5_data}). Moreover, we report the plots of the posterior distribution in a subset of the parameter in Figure \ref{fig:est_h_T}.

\begin{table}[hbt!]
    \centering
\begin{tabular}{c|c|c|c|c}
\toprule
    Scenario & T & \# events & \# excursions & \# local excursions  \\
 \midrule
\multirow{4}{*}{Excitation}   & 50 & 2621 & 36 & 114  \\
  & 200 & 10,729 & 155 & 473  \\
    & 400 & 21,727 & 303 & 957  \\
      & 800 & 42,904 & 596 & 1921  \\
     \midrule
     \multirow{4}{*}{Inhibition}  & 50 & 1747 & 49 & 134  \\
  & 200 & 7019 & 222 & 529  \\
    & 400 & 13,819 & 466 & 1053  \\
      & 800 & 27,723 & 926 & 2118  \\
     \bottomrule
\end{tabular}
    \caption{Number of points and \emph{global} and average \emph{local}  excursions in the multidimensional data sets of Simulation 5 ($K=10$).}
    \label{tab:simu_5_data}
\end{table}

\FloatBarrier

\subsection{Simulation 6}\label{app:simu_6}

\begin{figure}
    \centering
    \begin{subfigure}[t]{0.8\textwidth}
    \includegraphics[width=\textwidth, trim=0.cm 0.cm 0cm  0.cm,clip]{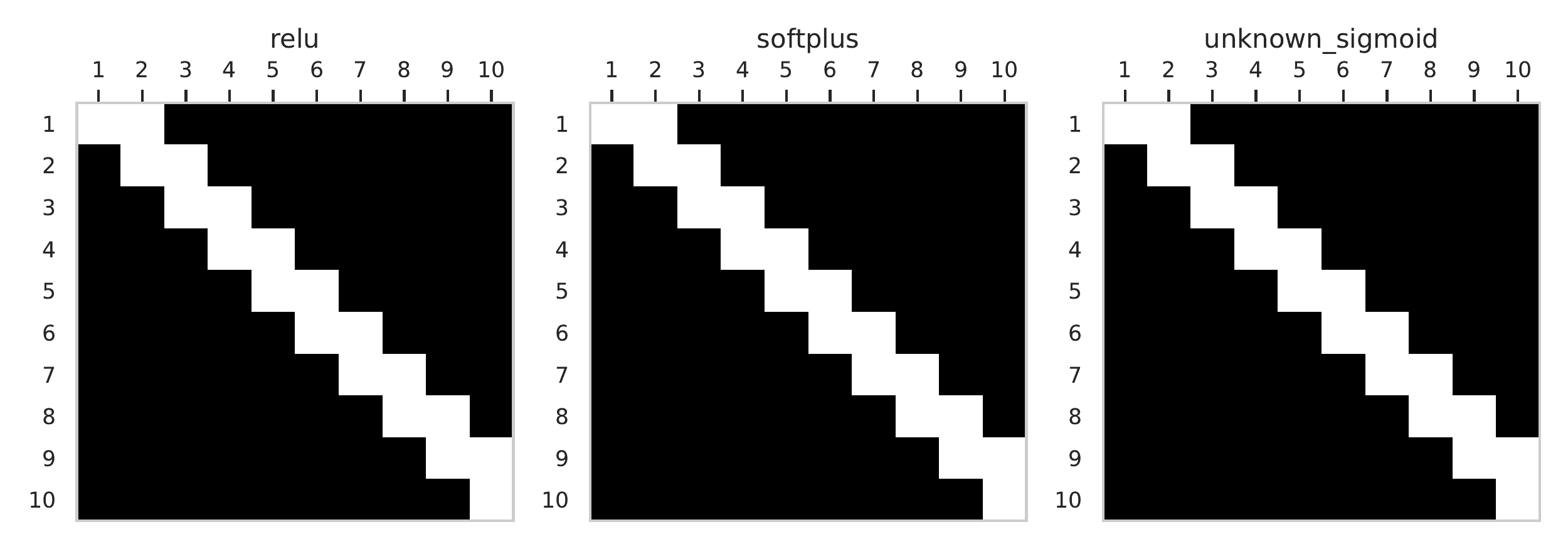}
    \caption{\emph{Excitation} scenario}
    \end{subfigure}
        \begin{subfigure}[b]{0.8\textwidth}
    \includegraphics[width=\textwidth, trim=0.cm 0.cm 0cm  0.cm,clip]{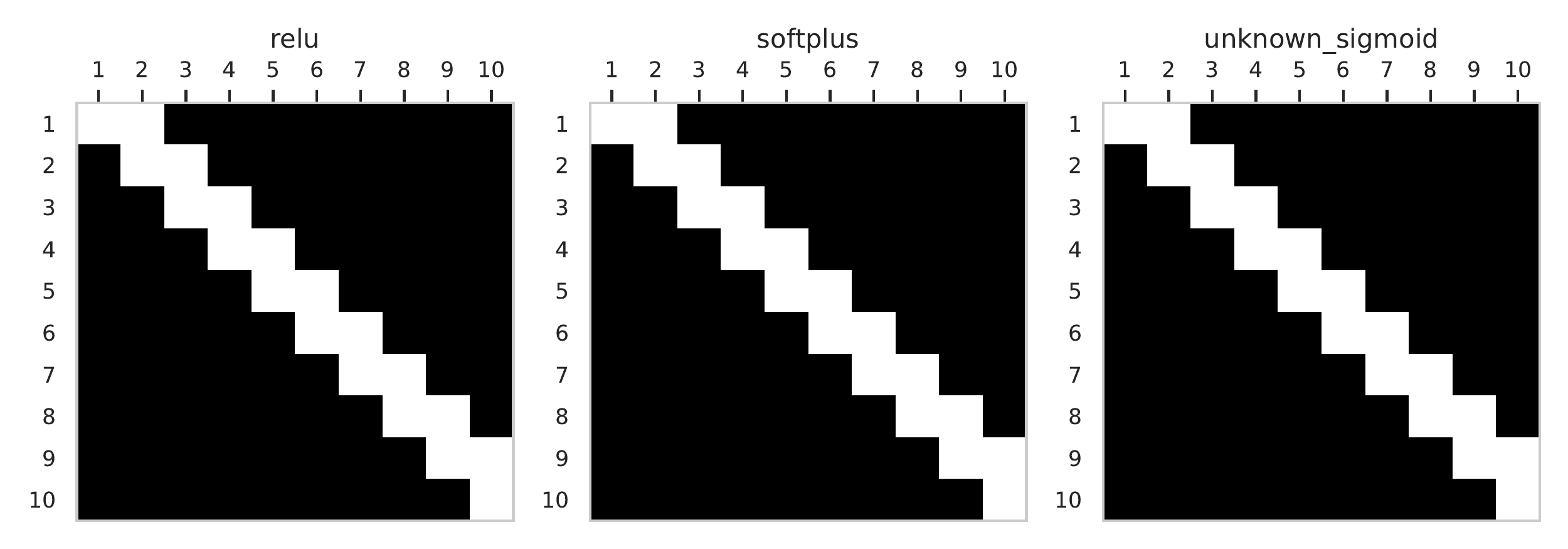}
    \caption{\emph{Inhibition} scenario}
    \end{subfigure}%
    \caption{Estimated graph after thresholding the $L_1$-norms using the ``gap" or ``slope change" heuristic, in the different settings of mis-specified link functions of Simulation 6, and in the \emph{Excitation} and \emph{Inhibition} scenarios. We observe that the true graph (with non-null principal and first off-diagonal) is correctly estimated for the ReLU mis-specification setting, while some errors happen in the two other link settings, in particular in the \emph{Inhibition} scenario.}
    \label{fig:graphs_mis}
\end{figure}

\begin{figure}
    \centering
    \begin{subfigure}[b]{\textwidth}
    \centering
    \includegraphics[width=0.8\textwidth, trim=0.cm 0.cm 0cm  0.cm,clip]{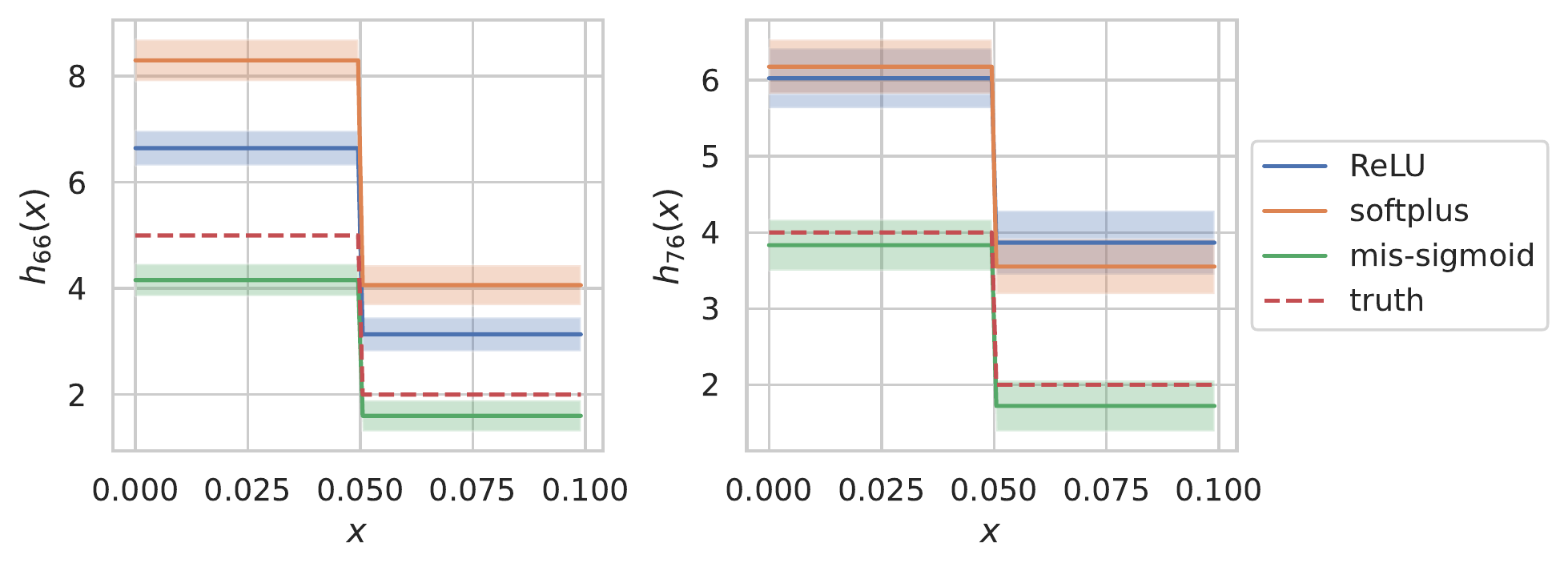}
    \caption{\emph{Excitation} scenario}
    \end{subfigure}
        \begin{subfigure}[b]{\textwidth}
        \centering
    \includegraphics[width=0.8\textwidth, trim=0.cm 0.cm 0cm  0.cm,clip]{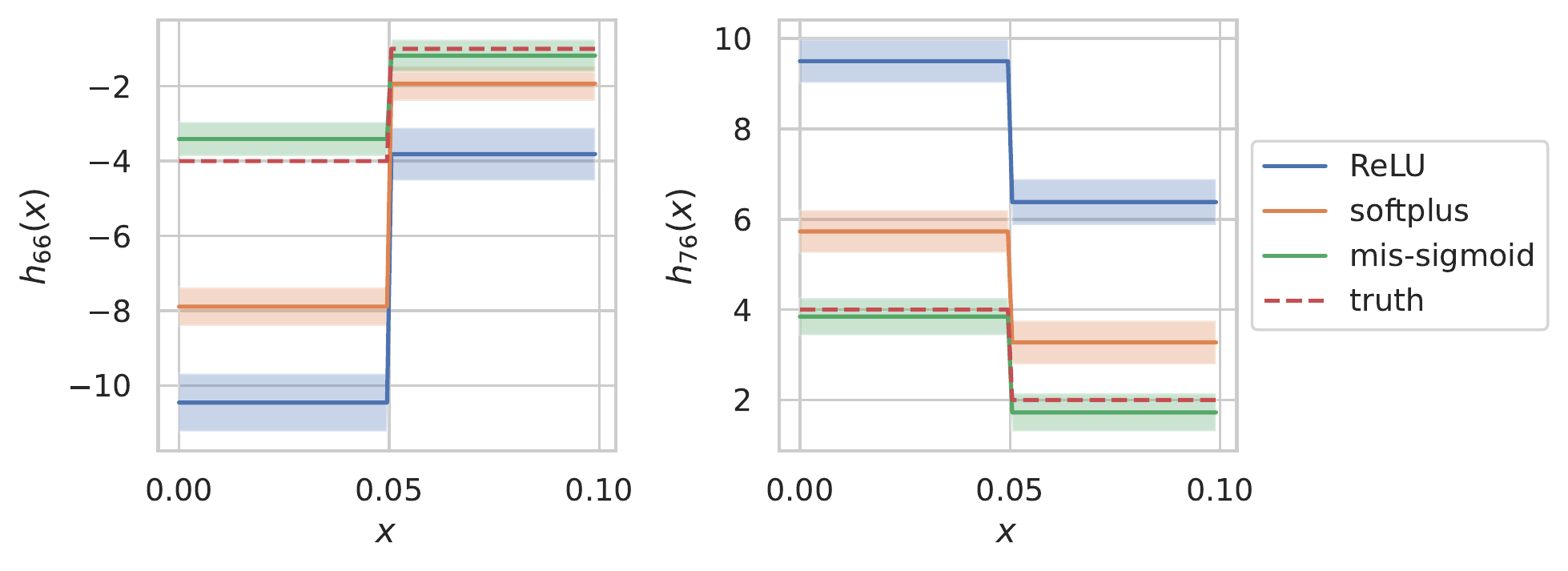}
    \caption{\emph{Inhibition} scenario}
    \end{subfigure}%
    \caption{Estimated interaction functions $h_{66}$ and $h_{76}$ in the mis-specified settings of Simulation 6, where the data is generated from a Hawkes model with ReLU, softplus, or a mis-specified link function, and in the \emph{Excitation} and \emph{Inhibition} scenarios. We note that the estimation of the interaction functions is deteriorated in these mis-specified cases, however the sign of the functions are still recovered. }
    \label{fig:est_h_mis_inh}
\end{figure}

\begin{figure}
    \centering
    \begin{subfigure}[b]{0.8\textwidth}
    \includegraphics[width=\textwidth, trim=0.cm 0.cm 0cm  0.cm,clip]{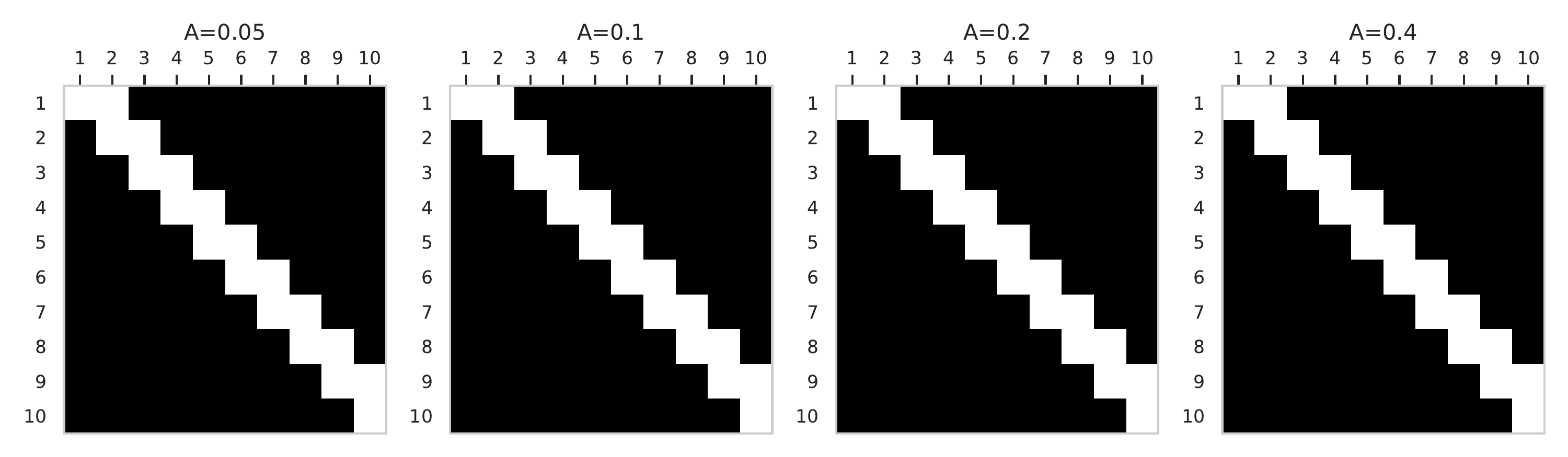}
    \caption{\emph{Excitation} scenario}
    \end{subfigure}
        \begin{subfigure}[b]{0.8\textwidth}
    \includegraphics[width=\textwidth, trim=0.cm 0.cm 0cm  0.cm,clip]{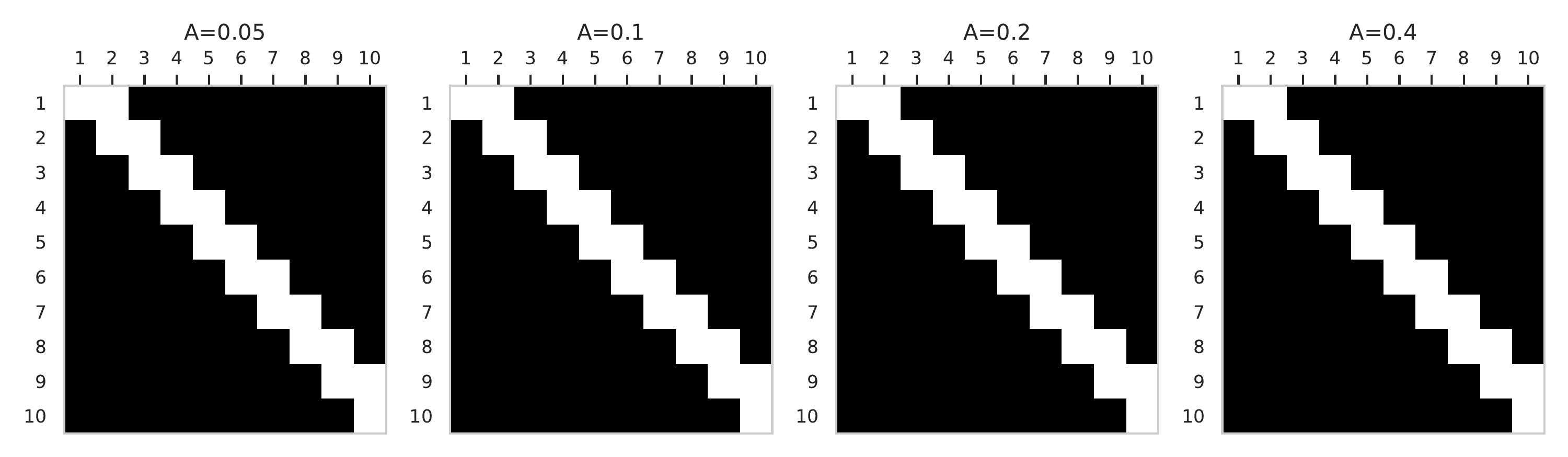}
    \caption{\emph{Inhibition} scenario}
    \end{subfigure}%
    \caption{Estimated graph after thresholding the $L_1$-norms, when using Algorithm \ref{alg:2step_adapt_cavi} with different support upper bounds $A'\in \{0.5, 0.1, 0.2, 0.4\}$, containing the true memory parameter $A=0.1$, in the settings of Simulation 7. We note that the true graph (with non-null principal and first off-diagonal) is correctly estimated in all cases, in the \emph{Excitation} scenario (first row) and in the \emph{Inhibition} scenario (second row).}
    \label{fig:graphs_A}
\end{figure}

\begin{figure}
    \centering
    \begin{subfigure}[b]{\textwidth}
    \includegraphics[width=\textwidth, trim=0.cm 0.cm 0cm  0.7cm,clip]{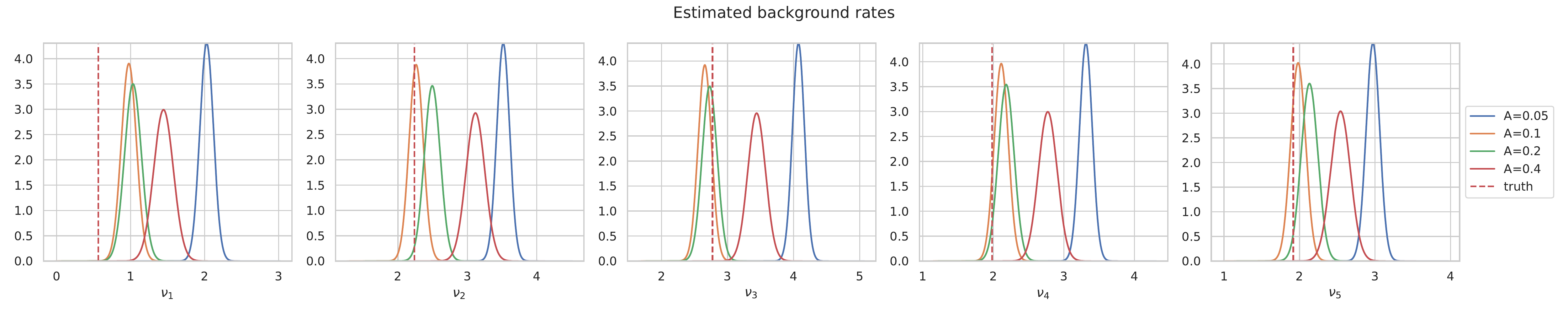}
    \caption{\emph{Excitation} scenario}
    \end{subfigure}
        \begin{subfigure}[b]{\textwidth}
    \includegraphics[width=\textwidth, trim=0.cm 0.cm 0cm  0.7cm,clip]{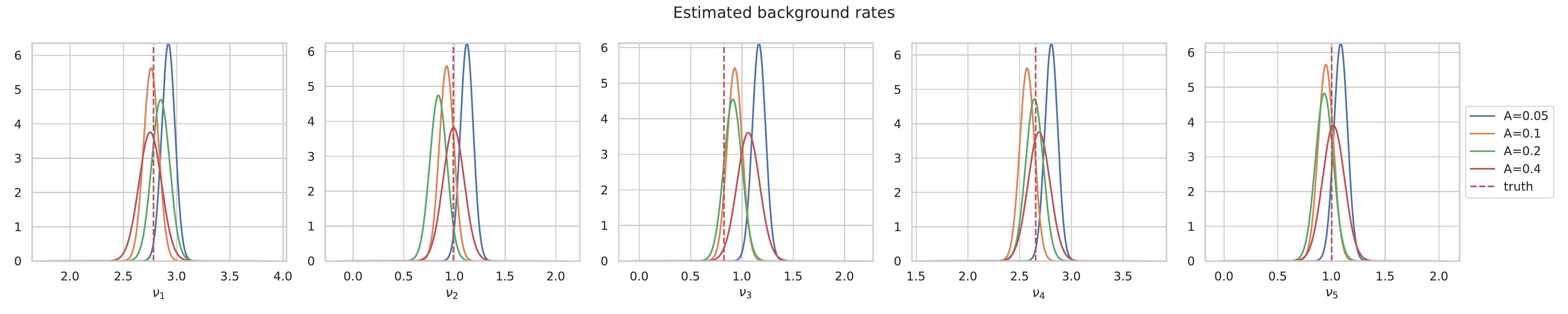}
    \caption{\emph{Inhibition} scenario}
    \end{subfigure}%
    \caption{Estimated background rates $\nu_k$ for $k=1,\dots, 5$ when using different values of the upper bound parameter $A \in \{0.05, 0.1, 0.2, 0.4\}$, in the two scenarios of Simulation 8. As expected, the background rates are better estimated in the well-specified setting $A=A_0=0.1$; nonetheless, when $A$ is not too far above $A_0$, the estimation does not deteriorate too much, in particular in the \emph{Inhibition} scenarios.}
    \label{fig:est_nu_A}
\end{figure}

This section contains the estimated graphs (Figures \ref{fig:graphs_mis} and \ref{fig:graphs_A}), the variational posterior distribution on a subset of the parameter (Figures \ref{fig:est_h_mis_inh} and \ref{fig:est_nu_A}), in the mis-specified settings of Simulation 6.

% \begin{figure}
%     \centering
%     \begin{subfigure}[b]{\textwidth}
%     \includegraphics[width=\textwidth, trim=0.cm 0.cm 0cm  0.7cm,clip]{figures/estimated_nu_D10_exc_mispecified.pdf}
%     \caption{\emph{Excitation} scenario}
%     \end{subfigure}
%         \begin{subfigure}[b]{\textwidth}
%     \includegraphics[width=\textwidth, trim=0.cm 0.cm 0cm  0.7cm,clip]{figures/estimated_nu_D10_inh_mispecified.pdf}
%     \caption{\emph{Inhibition} scenario}
%     \end{subfigure}%
%     \caption{Estimated background rates $\nu_k$ for $k=1,\dots, 5$ in the different settings of mis-specified link functions of Simulation 6, and in the \emph{Excitation} and \emph{Inhibition} scenarios. We observe that the background parameter is not well estimated in all mis-specified settings, and is worse in the softplus case.}
%     \label{fig:est_nu_mis}
% \end{figure}

\FloatBarrier

\end{document}